# Some series and integrals involving the Riemann zeta function, binomial coefficients and the harmonic numbers

## Volume V

Donal F. Connon

18 February 2008

**Abstract**

In this series of seven papers, predominantly by means of elementary analysis, we establish a number of identities related to the Riemann zeta function, including the following:

$$\varsigma(2n+1) = (-1)^{n+1} \frac{(2\pi)^{2n+1}}{(2n+1)!} \int_0^{1/2} B_{2n+1}(x) \cot(\pi x) dx$$

$$\frac{1}{2} \int_a^b p(x) dx = \sum_{n=0}^{\infty} \int_a^b p(x) \cos \alpha nx \, dx$$

$$\frac{1}{2} \int_a^b p(x) \cot(\alpha x/2) dx = \sum_{n=0}^{\infty} \int_a^b p(x) \sin \alpha nx \, dx$$

$$\int_0^{\pi/8} x \cot x \, dx = \frac{\pi}{16} \log\left[2-\sqrt{2}\right] + \frac{1}{8}\left[1-\sqrt{2}\right]G + \frac{1}{64}\left[\sqrt{2}\varsigma\left(2,\frac{1}{8}\right) - 2\left(\sqrt{2}+1\right)\pi^2\right]$$

$$\sum_{n=1}^{\infty} \frac{Ci(n\pi)}{n^4} = \varsigma(4)\left[\gamma + \log \pi - \frac{55}{32}\right] - \varsigma'(4)$$

$$\frac{\pi^2}{240} = \sum_{n=0}^{\infty} \frac{(n-1)\varsigma(2n)}{(n+1)(n+2)(n+3)2^n}$$

$$\sum_{n=1}^{\infty} \frac{H_n^{(1)} H_n^{(2)}}{2^n} = \frac{11}{6}\pi^2 \log 2 - \frac{179}{20}\varsigma(3) - \frac{1}{3}\log^3 2$$

$$\frac{\pi}{8}\log 2 + \frac{G}{2} = \frac{1}{\sqrt{2}} \sum_{n=0}^{\infty} \frac{1}{2^{3n}(2n+1)^2}\binom{2n}{n}$$

where $p(x)$ is a suitably behaved continuously differentiable function and $Ci(x)$ is the cosine integral.

Whilst the paper is mainly expository, some of the formulae reported in it are believed to be new, and the paper may also be of interest specifically due to the fact that most of the various identities have been derived by elementary methods.

## 5. AN APPLICATION OF THE BERNOULLI POLYNOMIALS

In (2.24) of Volume I we showed that for suitable functions

$$\sum_{n=1}^{\infty} \frac{1}{2^n} \sum_{k=1}^{n} \int_{a}^{b} \binom{n}{k} p(x) \sin 2\alpha kx \, dx = \int_{a}^{b} p(x) \cot \alpha x \, dx$$

It would be possible to replace $p(x)$ by $x^N$ in (2.24), for example, and carry out the detailed calculations as before, but considerable simplifications are achieved by employing the Bernoulli polynomials whose important properties are considered in Appendix A of Volume VI.

Therefore, in this part, we let $p(x) = B_{2N+1}(x)$, the selected range of integration is $[0, 1/2]$ and we use the generalised identity (2.24) with $\alpha = \pi$.

Since

(5.1) $\qquad B_{2N+1}(0) = B_{2N+1}(1/2) = 0$

the Riemann-Lebesgue lemma conditions will be satisfied at both end points of $[0, 1/2]$. We therefore have from (2.24)

(5.2) $\qquad \sum_{n=1}^{\infty} \frac{1}{2^n} \sum_{k=1}^{n} \int_{0}^{1/2} \binom{n}{k} B_{2N+1}(x) \sin 2\pi kx \, dx = \int_{0}^{1/2} B_{2N+1}(x) \cot \pi x \, dx$

and it will be noted that the right hand side appears in (1.13), an integral involving $\varsigma(2N+1)$.

Integration by parts gives

(5.3) $\qquad I_k = \int_{0}^{1/2} B_{2N+1}(x) \sin 2\pi kx \, dx$

(5.4) $\qquad = -B_{2N+1}(x) \frac{\cos 2\pi kx}{2\pi k} \bigg|_{0}^{1/2} + \frac{2N+1}{2\pi k} \int_{0}^{1/2} B_{2N}(x) \cos 2\pi kx \, dx$

The integrated part in (5.4) vanishes in view of the relations contained in (5.1). A further integration by parts gives

(5.5) $\qquad = \frac{(2N+1)}{2\pi k} \frac{B_{2N}(x) \sin 2\pi kx}{2\pi k} \bigg|_{0}^{1/2} - \frac{(2N+1)2N}{(2\pi k)^2} \int_{0}^{1/2} B_{2N-1}(x) \sin 2\pi kx \, dx$

$$= -\frac{(2N+1)2N}{(2\pi k)^2} \int_0^{1/2} B_{2N-1}(x) \sin 2\pi kx \, dx$$

Integrating by parts a total of $2N-1$ times, we obtain

(5.6) $$\int_0^{1/2} B_{2N}(x) \cos 2\pi kx \, dx = \frac{(-1)^{N-1}(2N)!}{2(2\pi k)^{2N}} + B_{2N} \int_0^{1/2} \cos 2\pi kx \, dx$$

and the last integral is zero. Therefore we get

(5.7) $$I_k = \frac{(-1)^{N-1}(2N+1)!}{2(2\pi k)^{2N+1}}$$

Substituting (5.7) in (5.2) gives

(5.8) $$\frac{(-1)^{N-1}(2N+1)!}{2(2\pi)^{2N+1}} \left\{ \sum_{n=1}^{\infty} \frac{1}{2^n} \sum_{k=1}^{n} \binom{n}{k} \frac{1}{k^{2N+1}} \right\} = \int_0^{1/2} B_{2N+1}(x) \cot \pi x \, dx$$

We have already shown in Lemma (3.4) that the series in parentheses is equal to $2\varsigma(2N+1)$ and therefore, upon rearranging (5.8) we have proved (1.13), namely

(5.9) $$\varsigma(2n+1) = (-1)^{n+1} \frac{(2\pi)^{2n+1}}{(2n+1)!} \int_0^{1/2} B_{2n+1}(x) \cot(\pi x) dx$$

Alternatively, the above could be used as a proof of Lemma 3.4 once we had determined a different proof for (1.13) (and an alternative proof of this is shown in (6.31)).

If we let $p(x) = B_{2N}(x) - B_{2N}$ (so that by definition $p(0) = 0$), and using the identity (2.23), it is easily shown that

(5.10) $$\frac{(-1)^{N-1}(2N)!}{2(2\pi)^{2N}} \left\{ \sum_{n=1}^{\infty} \frac{1}{2^n} \sum_{k=1}^{n} \binom{n}{k} \frac{1}{k^{2N}} \right\} = \int_0^{1/2} (B_{2N}(x) - B_{2N}) \, dx$$

and hence we have

(5.11) $$\sum_{n=1}^{\infty} \frac{1}{2^n} \sum_{k=1}^{n} \binom{n}{k} \frac{1}{k^{2N}} = 2\varsigma(2N) \qquad \square$$

# 6. TRIGONOMETRIC INTEGRAL IDENTITIES

As previously mentioned in Volume I, Berndt [19] gave an elementary evaluation of $\varsigma(2n)$ in his 1975 paper. Indeed, it was this paper which kindled my interest in this subject and, in particular, it was his use of a trigonometric function to show the relationship between $\varsigma(2n)$ and $B_{2n}$ that led me to consider using a function of the form $1/[1-\exp(ix)]$ (which, really, is only a mild variant of the generating function (1.8) used in defining the Bernoulli numbers).

We therefore consider the following identity (which is easily verified by multiplying the numerator and the denominator by the complex conjugate $(1-e^{-ix})$),

$$(6.1) \qquad \frac{1}{1-e^{ix}} = \frac{1}{2} - \frac{i}{2}\frac{\sin x}{1-\cos x} = \frac{1}{2} + \frac{i}{2}\cot(x/2)$$

and using (2.3) of Volume I we obtain

$$(6.1a) \qquad \frac{1}{1-e^{ix}} = \sum_{n=0}^{N} e^{inx} + R_N(x)$$

where $\quad R_N(x) = \dfrac{e^{i(N+1)x}}{1-e^{ix}} = \dfrac{1}{2}e^{i(N+1)x}\{1+i\cot(x/2)\} = \dfrac{ie^{i(N+\frac{1}{2})x}}{2\sin(x/2)}$

Separating the real and imaginary parts of (6.1a) produces the following two identities (the first of which is called Lagrange's trigonometric identity and contains the Dirichlet kernel $D_N(x)$ which is employed in the theory of Fourier series [145, p.49])

$$(6.2) \qquad \frac{1}{2} = \sum_{n=0}^{N}\cos nx - \frac{\sin(N+1/2)x}{2\sin(x/2)}$$

$$(6.3) \qquad \frac{1}{2}\cot(x/2) = \sum_{n=0}^{N}\sin nx + \frac{\cos(N+1/2)x}{2\sin(x/2)}$$

(and these equations may obviously be generalised by substituting $\alpha x$ in place of $x$ as in (2.23) and (2.24)).

Integration then gives us

$$\frac{1}{2}\int_a^b p(x)\bigl(1+i\cot(x/2)\bigr)dx = \sum_{n=0}^{N}\int_a^b p(x)e^{inx}dx + R_N$$

and, in what follows, $p(x)$ is assumed to be twice continuously differentiable on $[a,b]$ and we have

(6.3a)
$$R_N = \int_a^b p(x)R_N(x)dx = \frac{1}{2}\int_a^b p(x)\big(\cos(N+1)x + i\sin(N+1)x\big)\big(1 + i\cot(x/2)\big)dx$$

$$= \frac{1}{2}\int_a^b p(x)\left\{\cos(N+1)x - \frac{\cos(x/2)}{\sin(x/2)}\sin(N+1)x\right\}dx$$

$$+ \frac{i}{2}\int_a^b p(x)\left\{\sin(N+1)x + \frac{\cos(x/2)}{\sin(x/2)}\cos(N+1)x\right\}dx$$

Therefore, provided $\sin(x/2)$ has no zero in $[a,b]$, the Riemann-Lebesgue lemma (2.17) from Volume I tells us that

$$\lim_{N\to\infty} R_N = 0.$$

This will also be the case if $\sin(a/2) = 0$, provided $p(a) = 0$. From the above we can therefore derive the following trigonometric identities:

(6.4) $$\frac{1}{2}\int_a^b p(x)\,dx = \sum_{n=0}^{\infty}\int_a^b p(x)\cos nx\,dx$$

(6.4a) $$\frac{1}{2}\int_a^b p(x)\cot(x/2)\,dx = \sum_{n=0}^{\infty}\int_a^b p(x)\sin nx\,dx$$

and more generally we have

(6.5) $$\frac{1}{2}\int_a^b p(x)\,dx = \sum_{n=0}^{\infty}\int_a^b p(x)\cos\alpha nx\,dx$$

(6.5a) $$\frac{1}{2}\int_a^b p(x)\cot(\alpha x/2)\,dx = \sum_{n=0}^{\infty}\int_a^b p(x)\sin\alpha nx\,dx$$

Equations (6.5) and (6.5a) are valid provided (i) $\sin(\alpha x/2) \neq 0\ \forall\ x \in [a,b]$ or, alternatively, (ii) if $\sin(\alpha a/2) = 0$ then $p(a) = 0$ also.

Similarly, using the identity

(6.6) $$\frac{1}{1+e^{ix}} = \frac{1}{2} - \frac{i}{2}\frac{\sin x}{1+\cos x} = \frac{1}{2} - \frac{i}{2}\tan(x/2)$$

we may easily prove that

$$(6.7) \quad \frac{1}{2}\int_a^b p(x)dx = \sum_{n=0}^{\infty} \int_a^b p(x)(-1)^n \cos nx\, dx$$

$$(6.7a) \quad -\frac{1}{2}\int_a^b p(x)\tan(x/2)dx = \sum_{n=0}^{\infty} \int_a^b p(x)(-1)^n \sin nx\, dx$$

and more generally

$$(6.8) \quad \frac{1}{2}\int_a^b p(x)dx = \sum_{n=0}^{\infty} \int_a^b p(x)(-1)^n \cos\alpha nx\, dx$$

$$(6.8a) \quad -\frac{1}{2}\int_a^b p(x)\tan(\alpha x/2)dx = \sum_{n=0}^{\infty} \int_a^b p(x)(-1)^n \sin\alpha nx\, dx$$

From (6.6) we note that the denominator is $\cos(x/2)$ and hence (6.7) and (6.7a) are only valid provided either (i) $\cos(x/2)$ has no zero in $[a,b]$ or (ii) if $\cos(a/2)=0$, then $p(a)$ must also be zero.

Equations (6.8) and (6.8a) are valid provided (i) $\cos(\alpha x/2) \neq 0 \ \forall \ x \in [a,b]$ or, alternatively, (ii) if $\cos(\alpha a/2)=0$ then $p(a)=0$ also.

Equation (6.7) may be written as

$$(6.8b) \quad \frac{1}{2}\int_a^b p(x)dx + \sum_{n=1}^{\infty} \int_a^b p(x)(-1)^n \cos nx\, dx = 0$$

The following simple trigonometric identities are easily proved (for example see (A.18) in Appendix A of Volume VI)

$$(6.9) \quad \cot(x/2) + \tan(x/2) = \frac{2}{\sin x}$$

$$(6.10) \quad \cot(x/2) - \tan(x/2) = 2\cot x$$

Therefore, combining (6.5) and (6.7) produces the following identity

$$(6.11) \quad \int_a^b \frac{p(x)}{\sin x}dx = \sum_{n=0}^{\infty} \int_a^b p(x)\sin nx\, dx - \sum_{n=0}^{\infty} \int_a^b p(x)(-1)^n \sin nx\, dx$$

which simplifies to

$$(6.12) \quad \int_a^b \frac{p(x)}{\sin x}dx = 2\sum_{n=0}^{\infty} \int_a^b p(x)\sin(2n+1)x\, dx$$

Similarly, using (6.10) we obtain

(6.13) $$\int_a^b p(x)\cot x\,dx = \sum_{n=0}^{\infty}\int_a^b p(x)\sin nx\,dx + \sum_{n=0}^{\infty}\int_a^b p(x)(-1)^n \sin nx\,dx$$

which simplifies to

(6.14) $$\int_a^b p(x)\cot x\,dx = 2\sum_{n=1}^{\infty}\int_a^b p(x)\sin 2nx\,dx$$

It should be noted that in the above formulae we require either (i) both $\sin(x/2)$ and $\cos(x/2)$ have no zero in $[a,b]$ or (ii) if either $\sin(a/2)$ or $\cos(a/2)$ is equal to zero then $p(a)$ must also be zero. Condition (i) is equivalent to the requirement that $\sin x$ has no zero in $[a,b]$.
Note that (6.14) is equivalent to (6.5a) with $\alpha = 2$.

More generally we have

(6.14a) $$\int_a^b \frac{p(x)}{\sin \alpha x}\,dx = 2\sum_{n=0}^{\infty}\int_a^b p(x)\sin\left[(2n+1)\alpha x\right]dx$$

(6.14b) $$\int_a^b p(x)\cot \alpha x\,dx = 2\sum_{n=1}^{\infty}\int_a^b p(x)\sin 2\alpha nx\,dx$$

Equations (6.14a) and (6.14b) are valid provided (i) $\sin(\alpha x) \neq 0 \ \forall\ x \in [a,b]$ or, alternatively, (ii) if $\sin(\alpha a) = 0$ then $p(a) = 0$ also.

Duoandikoetxea [56a] has recently carried out some related work in this area.

Some examples of the application of these identities are shown below.

**Example 1:**

With (6.8b) and $p(x) = x$ we have

(6.15) $$\frac{1}{2}\int_0^{\pi/2} x\,dx + \sum_{n=1}^{\infty}\int_0^{\pi/2}(-1)^n x\cos nx\,dx = 0$$

A straightforward integration by parts gives

$$\int_0^{\pi/2} x\cos nx\,dx = \left.\frac{x\sin nx}{n} + \frac{\cos nx}{n^2}\right|_0^{\pi/2}$$

Accordingly we obtain

$$\frac{\pi^2}{16} + \sum_{n=1}^{\infty}\frac{(-1)^n}{(2n)^2} + \frac{\pi}{2}\sum_{n=1}^{\infty}\frac{(-1)^n}{2n-1} - \sum_{n=1}^{\infty}\frac{(-1)^n}{n^2} = 0$$

or equivalently using the Leibniz formula

$$\frac{\pi^2}{16} - \frac{1}{2^2}\varsigma_a(2) - \frac{\pi}{2}\frac{\pi}{4} + \varsigma_a(2) = 0$$

and, after a little algebra, we obtain the well-known result for the alternating zeta function

(6.16) $$\varsigma_a(2) = \frac{\pi^2}{12}$$

Whilst it must be acknowledged that similar results may be obtained from Fourier Theory, this treatment has the advantage of not requiring the full rigour and complexity of that theory. A very elementary proof of the identity $\sum_{n=1}^{\infty}\frac{1}{(2n+1)^2} = \frac{\pi^2}{8}$ is set out in Appendix D of Volume VI: this simple proof does not even require a direct knowledge of the Riemann-Lebesgue lemma and, as such, would be eminently suitable material for an elementary calculus class.

**Example 2:**

In fact, we do not even need to assume the veracity of the Leibniz formula in the above example because we may derive it directly as follows. With (6.8b) and $p(x) = 1$ we have

(6.17) $$\frac{1}{2}\int_0^{\pi/2} dx + \sum_{n=1}^{\infty}\int_0^{\pi/2}(-1)^n \cos nx\, dx = 0$$

and therefore, rather easily, we find that

(6.18) $$\frac{\pi}{4} = \sum_{n=0}^{\infty}\frac{(-1)^n}{2n+1}$$

**Example 3:**

Letting $p(x) = x^2$ in (6.8b) we obtain

$$\frac{1}{2}\int_0^{\pi/2} x^2 dx + \sum_{n=1}^{\infty}\left[\frac{(-1)^n 2x\cos nx}{n^2} - \frac{(-1)^n 2\sin nx}{n^3} + \frac{(-1)^n x^2 \sin nx}{n}\right]_0^{\pi/2} = 0$$

and, after some simplification, we obtain

(6.19) $$\sum_{n=0}^{\infty}\frac{(-1)^n}{(2n+1)^3}=\frac{\pi^3}{32}$$

This is a particular case of the general identity involving the Euler numbers which are defined in Appendix A (see (A.28)) of Volume VI.

**Example 4:**

As a further example, consider equation (6.14) with $p(x) = x^2$

$$\int_a^b p(x)\cot x\,dx = 2\sum_{n=1}^{\infty}\int_a^b p(x)\sin 2nx\,dx$$

We have

(6.20) $$\int_0^{\pi/2} x^2 \cot x\,dx = 2\sum_{n=1}^{\infty}\int_0^{\pi/2} x^2 \sin 2nx\,dx$$

$$= \frac{1}{2}\sum_{n=1}^{\infty}\frac{\cos 2nx}{n^3} - x^2\sum_{n=1}^{\infty}\frac{\cos 2nx}{n} + x\sum_{n=1}^{\infty}\frac{\sin 2nx}{n}\Bigg|_0^{\pi/2}$$

$$= \frac{1}{2}\sum_{n=1}^{\infty}\frac{(-1)^n}{n^3} - \frac{\pi^2}{4}\sum_{n=1}^{\infty}\frac{(-1)^n}{n} - \frac{1}{2}\varsigma(3)$$

Hence we see that

(6.20a) $$\int_0^{\pi/2} x^2 \cot x\,dx = -\frac{7}{8}\varsigma(3) + \frac{\pi^2}{4}\log 2$$

Therefore, reference to (3.9) shows that we have another elementary proof of the ubiquitous Euler integral identity (1.11) of 1772. Contrast the relative simplicity of this approach with the rather complex method originally employed by Euler (for example, see Ayoub's expository paper [15]).

**Example 5:**

We now employ (6.7a) with $p(x) = 1$ and $[a,b] = [0, \pi/4]$

$$-\frac{1}{2}\int_a^b p(x)\tan(x/2)dx = \sum_{n=0}^{\infty}\int_a^b p(x)(-1)^n \sin nx\,dx$$

(6.21) $$-\frac{1}{2}\int_0^{\pi/4} \tan(x/2)dx = \sum_{n=1}^{\infty}\int_0^{\pi/4} (-1)^n \sin nx\,dx$$

This results in

(6.22) $$\log\sqrt{\frac{1}{2}+\frac{\sqrt{2}}{4}} = \sum_{n=1}^{\infty}(-1)^{n+1}\frac{\cos(n\pi/4)}{n} - \log 2$$

since $\cos^2(\pi/8) = \frac{1}{2}[1+\cos(\pi/4)]$.

More generally, with $[a,b]=[0,t]$ we obtain the well-known Fourier series which, inter alia, is recorded in [130, p.148]

(6.22a) $$\sum_{n=1}^{\infty}(-1)^{n+1}\frac{\cos nt}{n} = \log[2\cos(t/2)]$$

**Example 6:**

With identity (6.4) and $p(x)=1$

$$\frac{1}{2}\int_a^b p(x)\,dx = \sum_{n=0}^{\infty}\int_a^b p(x)\cos nx\,dx$$

we have

(6.23) $$-\frac{1}{2}\int_{\pi/4}^{\pi/2} dx = \sum_{n=1}^{\infty}\int_{\pi/4}^{\pi/2}\cos nx\,dx$$

and we obtain

(6.24) $$\sum_{n=1}^{\infty}\frac{\sin(n\pi/4)}{n} = \frac{3\pi}{8}$$

More generally we have

$$-\frac{1}{2}\int_a^t dx = \sum_{n=1}^{\infty}\int_a^t \cos nx\,dx$$

and therefore we get

$$\frac{1}{2}(a-t) = \sum_{n=1}^{\infty}\frac{\sin nt - \sin na}{n}$$

Upon letting $a=\pi$ we obtain (7.5)

$$\frac{1}{2}(\pi - t) = \sum_{n=1}^{\infty} \frac{\sin nt}{n}$$

Note that, in the above example, since $p(x) = 1$ the interval of integration can not include 0. However, if we let $p(x) = x$ we can change the integral to include 0 as follows

**Example 7:**

$$-\frac{1}{2}\int_0^{\pi/4} x\,dx = \sum_{n=1}^{\infty} \int_0^{\pi/4} x \cos nx\,dx$$

Integration by parts gives us

(6.25) $$\sum_{n=1}^{\infty} \frac{\cos(n\pi/4)}{n^2} = \frac{11}{192}\pi^2$$

This is a particular case of the more general Fourier series [130, p.148]

(6.25a) $$\sum_{n=1}^{\infty} \frac{\cos nt}{n^2} = \frac{3t^2 - 6\pi + 2\pi^2}{12}$$

and this in turn is easily derived by using the identity

$$-\frac{1}{2}\int_0^t x\,dx = \sum_{n=1}^{\infty} \int_0^t x \cos nx\,dx$$

Using (6.5) gives us

$$-\frac{1}{2}\int_0^t B_3(x)\,dx = \sum_{n=1}^{\infty} \int_0^t B_3(x) \cos \alpha nx\,dx$$

where $B_3(x)$ is a Bernoulli polynomial and $B_3(0) = B_3 = 0$. We then have

$$\int_0^t B_3(x)\,dx = \frac{1}{4}[B_4(t) - B_4(0)]$$

With integration by parts we obtain

$$\int_0^t B_3(x) \cos \alpha nx\,dx = \left. \frac{B_3(x) \sin \alpha nx}{n\alpha} \right|_0^t - \frac{3}{n\alpha}\int_0^t B_2(x) \sin \alpha nx\,dx$$

$$= \frac{B_3(t)\sin \alpha nt}{n\alpha} - \frac{3}{n\alpha}\int_0^t B_2(x)\sin \alpha nx\, dx$$

$$\int_0^t B_2(x)\sin \alpha nx\, dx = -\frac{B_2(x)\cos \alpha nx}{n\alpha}\bigg|_0^t + \frac{2}{n\alpha}\int_0^t B_1(x)\cos \alpha nx\, dx$$

$$= -\frac{B_2(t)\cos \alpha nt}{n\alpha} + \frac{B_2}{n\alpha} + \frac{2}{n\alpha}\int_0^t B_1(x)\cos \alpha nx\, dx$$

$$\int_0^t B_1(x)\cos \alpha nx\, dx = \frac{B_1(x)\sin \alpha nx}{n\alpha}\bigg|_0^t - \frac{1}{n\alpha}\int_0^t B_0(x)\sin \alpha nx\, dx$$

$$= \frac{B_1(t)\sin \alpha nt}{n\alpha} - \frac{\cos \alpha nt - 1}{(n\alpha)^2}$$

This gives us

$$\int_0^t B_3(x)\cos \alpha nx\, dx =$$

$$\frac{B_3(t)\sin \alpha nt}{n\alpha} + 3\frac{B_2(t)\cos \alpha nt}{(n\alpha)^2} - 3\frac{B_2}{(n\alpha)^2} - 3!\frac{B_1(t)\sin \alpha nt}{(n\alpha)^3} + 3!\frac{\cos \alpha nt - 1}{(n\alpha)^4}$$

and hence we have

$$\frac{1}{8}[B_4 - B_4(t)] =$$

$$B_3(t)\sum_{n=1}^{\infty}\frac{\sin \alpha nt}{n\alpha} + 3B_2(t)\sum_{n=1}^{\infty}\frac{\cos \alpha nt}{(n\alpha)^2} - 3B_2\sum_{n=1}^{\infty}\frac{\cos \alpha nt}{(n\alpha)^2} - 3!B_1(t)\sum_{n=1}^{\infty}\frac{\sin \alpha nt}{(n\alpha)^3} + 3!\sum_{n=1}^{\infty}\frac{\cos \alpha nt - 1}{(n\alpha)^4}$$

With $\alpha = 2\pi$ we have

$$\frac{1}{8}[B_4 - B_4(t)] = B_3(t)\sum_{n=1}^{\infty}\frac{\sin 2n\pi t}{2n\pi} + 3[B_2(t) - B_2]\sum_{n=1}^{\infty}\frac{\cos 2n\pi t}{(2n\pi)^2}$$

$$-3!B_1(t)\sum_{n=1}^{\infty}\frac{\sin 2n\pi t}{(2n\pi)^3} + 3!\sum_{n=1}^{\infty}\frac{\cos 2n\pi t - 1}{(2n\pi)^4}$$

As noted in Apostol's book [13, p.338] we have

(6.25b) $$B_{2N}(t) = (-1)^{N+1} 2(2N)! \sum_{n=1}^{\infty} \frac{\cos 2n\pi t}{(2\pi n)^{2N}} \quad , (N = 1, 2, ...)$$

(6.25c) $$B_{2N+1}(t) = (-1)^{N+1} 2(2N+1)! \sum_{n=1}^{\infty} \frac{\sin 2n\pi t}{(2\pi n)^{2N+1}} \quad , (N = 0, 1, 2, ...)$$

and by substitution we then see that everything cancels out.

The above formulae for the Bernoulli polynomials are also direct consequences of the Hurwitz formula (4.4.228i) for the generalised zeta function.

**Example 8:**

Using (6.12) we have

(6.26) $$\int_0^{\pi/2} \frac{x}{\sin x} dx = 2 \sum_{n=0}^{\infty} \int_0^{\pi/2} x \sin(2n+1)x \, dx$$

(6.27) $$\int_0^{\pi/2} x \sin(2n+1)x \, dx = -\frac{x \cos(2n+1)x}{2n+1} + \frac{\sin(2n+1)x}{(2n+1)^2} \bigg|_0^{\pi/2}$$

$$= \frac{\cos n\pi}{(2n+1)^2} = \frac{(-1)^n}{(2n+1)^2}$$

Hence we obtain

(6.28) $$\int_0^{\pi/2} \frac{x}{\sin x} dx = 2 \sum_{n=0}^{\infty} \frac{(-1)^n}{(2n+1)^2} = 2G$$

where $G$ is Catalan's constant. This result is well-known and an alternative proof is contained, for example, in Bradley's website [33].

A connection with the gamma function is shown below.

$$\int_0^{\pi/2} \frac{x}{\sin x} dx = \pi \int_0^{1/2} \frac{\pi u}{\sin \pi u} du$$

and employing Euler's reflection formula (6.61) for the gamma function

$$\Gamma(u)\Gamma(1-u) = \frac{\pi}{\sin \pi u}$$

this becomes

$$= \pi \int_0^{1/2} u\Gamma(u)\Gamma(1-u)\,du = \pi \int_0^{1/2} \Gamma(1+u)\Gamma(1-u)\,du$$

Therefore we have

(6.28a) $$\int_0^{1/2} \Gamma(1+u)\Gamma(1-u)\,du = 2\frac{G}{\pi}$$

We have the well-known result [25, pp.195 & 254]

$$\int_0^{\pi/2} \frac{1}{\sqrt{\sin x}}\,dx = \frac{1}{4\sqrt{\pi}}\Gamma^2\left(\frac{1}{4}\right)$$

and therefore using (6.12) with $p(x) = \sqrt{\sin x}$ we get

$$\int_0^{\pi/2} \frac{1}{\sqrt{\sin x}}\,dx = \int_0^{\pi/2} \frac{\sqrt{\sin x}}{\sin x}\,dx = 2\sum_{n=0}^{\infty} \int_0^{\pi/2} \sqrt{\sin x}\,\sin(2n+1)x\,dx$$

Integration by parts gives us

$$\int_0^{\pi/2} \sqrt{\sin x}\,\sin(2n+1)x\,dx = -\sqrt{\sin x}\,\frac{\cos(2n+1)x}{2n+1}\bigg|_0^{\pi/2} + \frac{1}{2}\int_0^{\pi/2} \frac{\cos x}{\sqrt{\sin x}}\frac{\cos(2n+1)x}{2n+1}\,dx$$

$$= \frac{1}{2}\int_0^{\pi/2} \frac{\cos x}{\sqrt{\sin x}}\frac{\cos(2n+1)x}{2n+1}\,dx$$

Therefore we have

$$\sum_{n=0}^{\infty} \int_0^{\pi/2} \sqrt{\sin x}\,\sin(2n+1)x\,dx = \frac{1}{2}\int_0^{\pi/2} \frac{\cos x}{\sqrt{\sin x}} \sum_{n=0}^{\infty} \frac{\cos(2n+1)x}{2n+1}\,dx$$

From Tolstov [130, p.149] we have the Fourier series in the range $(0,\infty)$

$$\sum_{n=0}^{\infty} \frac{\cos(2n+1)x}{2n+1} = -\frac{1}{2}\log\tan\frac{x}{2}$$

and hence

$$\sum_{n=0}^{\infty} \int_0^{\pi/2} \sqrt{\sin x}\,\sin(2n+1)x\,dx = -\frac{1}{4}\int_0^{\pi/2} \frac{\cos x}{\sqrt{\sin x}}\log\tan\frac{x}{2}\,dx = \frac{1}{8\sqrt{\pi}}\Gamma^2\left(\frac{1}{4}\right)$$

The Wolfram Integrator expresses $\int \frac{\cos x}{\sqrt{\sin x}} \log \tan \frac{x}{2} dx$ as an elliptic integral of the first kind (and I have decided that such integrals are certainly beyond my remit!): this aspect may be worth exploring further to see if it gets us any nearer to obtaining a closed form expression for $\Gamma\left(\frac{1}{4}\right)$.

It should however be noted that a simple integration by parts gives us

$$\int_0^{\pi/2} \frac{\cos x}{\sqrt{\sin x}} \log \tan \frac{x}{2} dx = 2\sqrt{\sin x} \log \tan \frac{x}{2} \Big|_0^{\pi/2} - \int_0^{\pi/2} \sqrt{\sin x} \frac{\sec^2 \frac{x}{2}}{\tan \frac{x}{2}} dx$$

$$= -2 \int_0^{\pi/2} \frac{\sqrt{\sin x}}{\sin x} dx = -2 \int_0^{\pi/2} \frac{dx}{\sqrt{\sin x}}$$

and this completes the circle!

Using this method it is possible to evaluate more complex integrals, for example let $p(x) = \sinh x$ so that $p(0) = 0$. Using (6.12) we have

$$\int_0^t \frac{\sinh x}{\sin x} dx = 2 \sum_{n=0}^{\infty} \int_0^t \sinh x \sin(2n+1)x \, dx$$

$$\int_0^t \sinh x \sin(2n+1)x \, dx = \frac{\cosh x \sin(2n+1)x - (2n+1) \sinh x \cos(2n+1)x}{2n^2 + 2n + 1} \Big|_0^t$$

$$= \frac{\cosh t \sin(2n+1)t - (2n+1) \sinh t \cos(2n+1)t}{2n^2 + 2n + 1}$$

Therefore we obtain

$$\int_0^t \frac{\sinh x}{\sin x} dx = \sum_{n=0}^{\infty} \frac{\cosh t \sin(2n+1)t - (2n+1) \sinh t \cos(2n+1)t}{2n^2 + 2n + 1}$$

and for example

$$\int_0^{\pi/2} \frac{\sinh x}{\sin x} dx = \cosh(\pi/2) \sum_{n=0}^{\infty} \frac{(-1)^{n+1}}{2n^2 + 2n + 1}$$

**Example 9:**

Following on from the first part of Example 8 we have

$$\int_0^{\pi/2} \frac{x^2}{\sin x} dx = 2\sum_{n=0}^{\infty} \int_0^{\pi/2} x^2 \sin(2n+1)x \, dx$$

and carrying out the elementary integration we find that

(6.29) $$\int_0^{\pi/2} \frac{x^2}{\sin x} dx = 2\pi G - \frac{7}{2}\varsigma(3)$$

The evaluation of this integral is also contained in Bradley's website [33] and an alternative proof was provided by De Doelder [55] (see (6.30r) of this paper).

Combining (6.28) and (6.29) we have

(6.30) $$\int_0^{\pi/2} \frac{x(\pi - x)}{\sin x} dx = \frac{7}{2}\varsigma(3)$$

The following analysis was presented by Amigó [7] in 2001. He defined the constants $C(k)$ for $k \geq 2$ by

(6.30a) $$C(k) = \left(\frac{\pi}{2}\right)^{k-1} \sum_{j=0}^{\infty} \frac{\beta(2j+1)}{(2j+1)...(2j+k)}$$

where $\beta(j)$ is the Dirichlet Beta function defined by (A.27)

(6.30b) $$\beta(j) = \sum_{n=0}^{\infty} \frac{(-1)^n}{(2n+1)^j}$$

and for $j = 2$, $\beta(2)$ is known as Catalan's constant $G$. We also have the relationship

(6.30c) $$\beta(2j+1) = \frac{(-1)^j (\pi/2)^{2j+1} E_{2j}}{2(2j)!}$$

and from (A.29) we have

(6.30d) $$\sec(\pi x/2) = \frac{4}{\pi} \sum_{n=0}^{\infty} \beta(2n+1) x^{2n} \quad , \quad |x| < 1$$

We have using the Beta function

$$\frac{1}{(2j+1)\ldots(2j+k)} = \frac{(2j)!}{(2j+k)!} = \frac{B(k,2j+1)}{(k-1)!}$$

(6.30e)
$$= \frac{1}{(k-1)!}\int_0^1 t^{k-1}(1-t)^{2j}\,dt$$

and substituting this in (6.30a) we obtain

(6.30f)
$$C(k) = \left(\frac{\pi}{2}\right)^{k-1}\frac{1}{(k-1)!}\int_0^1 t^{k-1}\left[\sum_{j=0}^\infty \beta(2j+1)(1-t)^{2j}\right]dt$$

Using (6.30d) this becomes

$$= \frac{1}{2}\left(\frac{\pi}{2}\right)^k \frac{1}{(k-1)!}\int_0^1 t^{k-1}\sec\left[\frac{\pi(1-t)}{2}\right]dt$$

$$= \frac{1}{2}\frac{1}{(k-1)!}\int_0^1 \frac{\pi(\pi t/2)^{k-1}}{2\sin(\pi t/2)}\,dt$$

We therefore get

(6.30g)
$$C(k) = \frac{1}{2}\frac{1}{(k-1)!}\int_0^{\pi/2}\frac{t^{k-1}}{\sin t}\,dt$$

Therefore, with $k=2$ for example, we get by using (6.28)

$$C(2) = \frac{1}{2}\int_0^{\pi/2}\frac{t}{\sin t}\,dt = G$$

and hence we have

(6.30ga)
$$G = \frac{\pi}{2}\sum_{j=0}^\infty \frac{\beta(2j+1)}{(2j+1)(2j+2)}$$

(compare this with (6.54)).

With $k=3$ we get by using (6.30)

$$C(3) = \frac{1}{4}\int_0^{\pi/2}\frac{t^2}{\sin t}\,dt = \frac{\pi}{2}G - \frac{7}{8}\varsigma(3)$$

and therefore

(6.30h) $$\frac{\pi}{2}G - \frac{7}{8}\varsigma(3) = \left(\frac{\pi}{2}\right)^2 \sum_{j=0}^{\infty} \frac{\beta(2j+1)}{(2j+1)(2j+2)(2j+3)}$$

We have by partial fractions

$$\frac{1}{(u+1)(u+2)(u+3)} = \frac{1}{2}\left\{\frac{1}{u+1} - \frac{1}{u+2}\right\} - \frac{1}{2}\left\{\frac{1}{u+2} - \frac{1}{u+3}\right\}$$

$$= \frac{1}{2}\left\{\frac{1}{(u+1)(u+2)} - \frac{1}{(u+2)(u+3)}\right\}$$

Hence we have

$$\sum_{j=0}^{\infty} \frac{\beta(2j+1)}{(2j+1)(2j+2)(2j+3)} = \frac{1}{2}\sum_{j=0}^{\infty}\frac{\beta(2j+1)}{(2j+1)(2j+2)} - \frac{1}{2}\sum_{j=0}^{\infty}\frac{\beta(2j+1)}{(2j+2)(2j+3)}$$

and consequently we get

(6.30i) $$\left(\frac{\pi}{2}\right)^2 \sum_{j=0}^{\infty} \frac{\beta(2j+1)}{(2j+2)(2j+3)} = \frac{7}{8}\varsigma(3)$$

In a follow-up paper [7a] Amigó defined the constants $D(k)$ for $k \geq 2$ by

(6.30j) $$D(k) = 2\pi^{k-1} \sum_{j=0}^{\infty} \frac{\lambda(2j)}{2j\ldots(2j+k-1)}$$

where

(6.30k) $$\lambda(s) = \sum_{n=0}^{\infty} \frac{1}{(2n+1)^s} = (1-2^{-s})\varsigma(s)$$

Amigó then showed that

(6.30l) $$D(k) = \frac{1}{2}\frac{1}{(k-1)!}\int_0^{\pi} t^{k-1}\cot(t/2)\,dt = \frac{2^{k-1}}{(k-1)!}\int_0^{\pi/2} t^{k-1}\cot t\,dt$$

His proof follows:

Let $f(t) = \tan(\pi t/2)$. With repeated integration by parts we obtain

$$\int_0^1 t^{k-1} f(1-t)\,dt = \frac{f(0)}{k} + \frac{f'(0)}{k(k+1)} + \frac{f''(0)}{k(k+1)(k+2)} + \ldots$$

$$= \sum_{n=0}^{\infty} \frac{f^{(n)}(0)}{k(k+1)...(k+n)}$$

We have the Maclaurin expansion (see (A.19) in Appendix A)

(6.30m) $$\tan(\pi x/2) = \frac{4}{\pi} \sum_{n=0}^{\infty} \lambda(2n) x^{2n-1} \quad , \quad |x| < 1$$

and therefore

$$f^{(n)}(0) = \begin{cases} 0 & \text{if } n = 2j \\ 4(2n-1)\lambda(2n)/\pi & \text{if } n = 2j-1 \end{cases}$$

Hence we get

$$\frac{2^{k-1}}{(k-1)!} \int_0^{\pi/2} t^{k-1} \tan\left(\frac{\pi}{2} - t\right) dt = \frac{\pi^k}{2(k-1)!} \int_0^1 t^{k-1} \tan\frac{\pi}{2}(1-t) dt$$

$$= \frac{\pi^k}{2(k-1)!} \sum_{j=1}^{\infty} \frac{4(2j-1)\lambda(2j)/\pi}{k(k+1)...(k+2j-1)}$$

$$= 2\pi^{k-1} \sum_{j=1}^{\infty} \frac{\lambda(2j)}{2j...(2j+k-1)} = D(k)$$

It may also be noted that

$$\int_0^{\pi/2} t^{k-1} \tan\left(\frac{\pi}{2} - t\right) dt = \int_0^{\pi/2} t^{k-1} \cot t\, dt$$

and we therefore get

(6.30n) $$\frac{2^{k-1}}{(k-1)!} \int_0^{\pi/2} t^{k-1} \cot t\, dt = 2\pi^{k-1} \sum_{j=1}^{\infty} \frac{\lambda(2j)}{2j...(2j+k-1)}$$

These integrals are also evaluated in Muzaffar's recent paper [103ac].

In [7] and [7a] Amigó reports that the following integrals may be proved using a combination of integration by parts and induction

(6.30o)
$$\frac{1}{(2k-1)!} \int_0^{\pi/2} t^{2k-1} \sin(2n+1)t\, dt = \sum_{j=1}^{k} (-1)^{j+1} \frac{(\pi/2)^{2k-2j}}{(2k-2j)!} \frac{(-1)^n}{(2n+1)^{2j}}$$

$$\frac{1}{(2k)!}\int_0^{\pi/2} t^{2k}\sin(2n+1)t\,dt = \sum_{j=1}^{k}(-1)^{j+1}\frac{(\pi/2)^{2k+1-2j}}{(2k+1-2j)!}\frac{(-1)^n}{(2n+1)^{2j}} + (-1)^k\frac{1}{(2n+1)^{2k+1}}$$

$$\frac{2^{2k}}{(2k-1)!}\int_0^{\pi/2} t^{2k-1}\sin 2nt\,dt = \sum_{j=1}^{k}(-1)^{j}\frac{\pi^{2k-2j+1}}{(2k-2j+1)!}\frac{(-1)^n}{n^{2j-1}}$$

$$\frac{2^{2k+1}}{(2k)!}\int_0^{\pi/2} t^{2k}\sin 2nt\,dt = \sum_{j=1}^{k}(-1)^{j}\frac{\pi^{2k-2j+2}}{(2k-2j+2)!}\frac{(-1)^n}{n^{2j-1}} + (-1)^k\frac{1-(-1)^n}{n^{2k+1}}$$

and these integrals can obviously be used in conjunction with (6.5a) and (6.14a).

We have from Knopp [90, p.215]

(6.30p) $\qquad \left(\tan^{-1} x\right)^2 = \sum_{k=1}^{\infty}(-1)^{k-1}a_k\frac{x^{2k}}{k}$

where $\qquad a_k = 1 + \frac{1}{3} + \frac{1}{5} + \ldots \frac{1}{2k-1}$

Dividing (6.30p) by $x$ and integrating gives us

$$\int_0^t \frac{\left(\tan^{-1} x\right)^2}{x}\,dx = \frac{1}{2}\sum_{k=1}^{\infty}(-1)^{k-1}a_k\frac{t^{2k}}{k^2}$$

Using the substitution $u = \tan^{-1} x$ we have

$$\int_0^1 \frac{\left(\tan^{-1} x\right)^2}{x}\,dx = \frac{1}{4}\int_0^{\pi/2}\frac{x^2}{\sin x}\,dx$$

and using (6.30) we obtain

$$4\int_0^1 \frac{\left(\tan^{-1} x\right)^2}{x}\,dx = 2\pi G - \frac{7}{2}\varsigma(3)$$

We conclude that

(6.30q) $\qquad 2\sum_{k=1}^{\infty}\frac{(-1)^{k-1}}{k^2}\left(1 + \frac{1}{3} + \frac{1}{5} + \ldots \frac{1}{2k-1}\right) = 2\pi G - \frac{7}{2}\varsigma(3)$

This result was obtained by De Doelder [55] using contour integration and expressed in the (corrected) form

(6.30r) $$\sum_{k=1}^{\infty} \frac{(-1)^{k-1}}{k^2}\left[\psi\left(k+\frac{1}{2}\right)-\psi\left(\frac{1}{2}\right)\right]=2\pi G-\frac{7}{2}\varsigma(3)$$

and the equivalence arises because we have from [125, p.20]

$$\psi\left(k+\frac{1}{2}\right)=-\gamma-2\log 2+2\sum_{n=0}^{k-1}\frac{1}{2n+1}$$

$$\psi\left(\frac{1}{2}\right)=-\gamma-2\log 2$$

and hence $$\psi\left(k+\frac{1}{2}\right)-\psi\left(\frac{1}{2}\right)=2\sum_{n=0}^{k-1}\frac{1}{2n+1}$$

Knopp [90, p.216] also reports that

(6.30s) $$\frac{1}{2}\tan^{-1}x.\log(1+x^2)=\sum_{k=1}^{\infty}(-1)^{k-1}H_{2k}\frac{x^{2k+1}}{2k+1}$$

and this expansion merits further exploration. Letting $t=\tan^{-1}x$ we get

$$t\log\cos t=\sum_{k=1}^{\infty}(-1)^k H_{2k}\frac{\tan^{2k+1}t}{2k+1}$$

and upon integration we obtain

$$\int_0^u t\log\cos t\,dt=\sum_{k=1}^{\infty}(-1)^k\frac{H_{2k}}{2k+1}\int_0^u \tan^{2k+1}t\,dt$$

Using Wiener's formula (8.11j), which we consider in more detail in Section 8, results in

$$\int_0^u \tan^{2k+1}t\,dt=(-1)^k\log|\sec u|+(-1)^k\sum_{n=1}^{k}(-1)^n\frac{\tan^{2n}u}{2n}$$

However, term by term integration does not appear to result in a convergent series.

Using the substitution $y=\tan^{-1}t$ we have

$$\phi(x)=\int_0^x \frac{\tan^{-1}t}{t}\,dt=\sum_{n=0}^{\infty}(-1)^n\frac{x^{2n+1}}{(2n+1)^2}=\int_0^{2\tan^{-1}x}\frac{y}{\sin y}\,dy$$

and using (6.27) this becomes

(6.30t) $$\phi(x) = 2\sum_{n=0}^{\infty}\left[-2\tan^{-1}x\frac{\cos\left[(2n+1)2\tan^{-1}x\right]}{2n+1} + \frac{\sin\left[(2n+1)2\tan^{-1}x\right]}{(2n+1)^2}\right]$$

From [25a] we note that

(6.30u) $$\cos\left[(2n+1)2\tan^{-1}x\right] = \frac{1}{(1+x^2)^{2n+1}}\sum_{k=0}^{2n+1}(-1)^k\binom{4n+2}{2k}x^{2k}$$

(6.30v) $$\sin\left[(2n+1)2\tan^{-1}x\right] = \frac{1}{(1+x^2)^{2n+1}}\sum_{k=0}^{2n}(-1)^k\binom{4n+2}{2k+1}x^{2k+1}$$

and therefore we get

$$\phi(x) = \int_0^x \frac{\tan^{-1}t}{t}dt$$

$$= -4\tan^{-1}x\sum_{n=0}^{\infty}\left[\frac{1}{2n+1}\frac{1}{(1+x^2)^{2n+1}}\sum_{k=0}^{2n+1}(-1)^k\binom{4n+2}{2k}x^{2k}\right] + 2\sum_{n=0}^{\infty}\left[\frac{1}{(2n+1)^2}\frac{1}{(1+x^2)^{2n+1}}\sum_{k=0}^{2n}(-1)^k\binom{4n+2}{2k+1}x^{2k}\right]$$

With $x=1$ we obtain

(6.30w)

$$\phi(1) = \int_0^1 \frac{\tan^{-1}t}{t}dt = \int_0^{\pi/2}\frac{y}{\sin y}dy = 2G$$

$$= -\pi\sum_{n=0}^{\infty}\left[\frac{1}{(2n+1)2^{2n+1}}\sum_{k=0}^{2n+1}(-1)^k\binom{4n+2}{2k}\right] + 2\sum_{n=0}^{\infty}\left[\frac{1}{(2n+1)^2 2^{2n+1}}\sum_{k=0}^{2n}(-1)^k\binom{4n+2}{2k+1}\right]$$

From (6.30u) and (6.30v) with $x=1$ it is easily seen that

$$\cos\left[\frac{(2n+1)\pi}{2}\right] = 0 = \frac{1}{2^{2n+1}}\sum_{k=0}^{2n+1}(-1)^k\binom{4n+2}{2k}$$

$$\sin\left[\frac{(2n+1)\pi}{2}\right] = (-1)^n = \frac{1}{2^{2n+1}}\sum_{k=0}^{2n}(-1)^k\binom{4n+2}{2k+1}$$

and, disappointingly, we therefore obtain the rather uninspiring result that

$$G = \sum_{n=0}^{\infty}\frac{(-1)^n}{(2n+1)^2}.$$

It may be worthwhile exploring Euler's more rapidly convergent series (a short proof of which is contained in the most enjoyable paper by Castellanos [42] "The Ubiquitous $\pi$")

(6.30x) $$\tan^{-1} x = \sum_{n=0}^{\infty} \frac{2^{2n}(n!)^2}{(2n+1)!} \frac{x^{2n+1}}{(1+x)^{n+1}}$$

Division by $x$ and integration results in

$$\int_0^t \frac{\tan^{-1} x}{x} dx = \sum_{n=0}^{\infty} \frac{2^{2n}(n!)^2}{(2n+1)!} \int_0^t \frac{x^{2n}}{(1+x)^{n+1}} dx$$

We have

$$\int_0^t \frac{x^{2n}}{(1+x)^{n+1}} dx = \int_1^{1+t} \frac{(1-u)^{2n}}{u^{n+1}} du = \sum_{k=0}^{2n} (-1)^k \binom{2n}{k} \frac{(1+t)^{2k-n} - 1}{2k-n}$$

and therefore

$$\int_0^t \frac{\tan^{-1} x}{x} dx = \sum_{n=0}^{\infty} \frac{2^{2n}(n!)^2}{(2n+1)!} \sum_{k=0}^{2n} (-1)^k \binom{2n}{k} \frac{(1+t)^{2k-n} - 1}{2k-n}$$

Can this be simplified by using the binomial theorem on $(1+t)^{2k-n}$ and then miraculously finding some appropriate combinatorial relationship for the product of two binomial coefficients?

With $t = 1$ we obtain

$$2G = \int_0^1 \frac{\tan^{-1} x}{x} dx = \sum_{n=0}^{\infty} \frac{2^{2n}(n!)^2}{(2n+1)!} \sum_{k=0}^{2n} (-1)^k \binom{2n}{k} \frac{2^{2k-n} - 1}{2k-n}$$

One of the first papers published by Ramanujan [76, pp.40 & 337] considered the tangent integral [126, p.110] $\phi(x) = \int_0^x \frac{\tan^{-1} t}{t} dt$. In this paper, Ramanujan (as corrected by Berndt) showed that

(6.30y) $$\frac{1}{1^2} + \frac{1}{3^2} - \frac{1}{5^2} - \frac{1}{7^2} + \frac{1}{9^2} + \ldots = \sqrt{2}\phi(\sqrt{2}-1) + \frac{\pi}{4\sqrt{2}} \log(1+\sqrt{2})$$

See also (C.64) in connection with (6.30y).

In my formulation we have

$$\phi(\sqrt{2}-1) = \int_0^{\sqrt{2}-1} \frac{\tan^{-1} t}{t} dt = \int_0^{2\tan^{-1}(\sqrt{2}-1)} \frac{x}{\sin x} dx$$

and, since $\tan^{-1}(\sqrt{2}-1) = \frac{\pi}{8}$, this is equivalent to

$$\phi(\sqrt{2}-1) = \int_0^{\pi/4} \frac{x}{\sin x} dx$$

Using (6.26) we have

$$\int_0^{\pi/4} \frac{x}{\sin x} dx = 2 \sum_{n=0}^{\infty} \int_0^{\pi/4} x \sin(2n+1)x \, dx$$

$$\int_0^{\pi/4} x \sin(2n+1)x \, dx = -\frac{x \cos(2n+1)x}{2n+1} + \frac{\sin(2n+1)x}{(2n+1)^2} \Big|_0^{\pi/4}$$

$$= -\frac{\pi/4 \cos[(2n+1)\pi/4]}{2n+1} + \frac{\sin[(2n+1)\pi/4]}{(2n+1)^2}$$

Elementary trigonometry gives us

$$\cos[(2n+1)\pi/4] = \frac{1}{\sqrt{2}} \cos([n\pi/4] - \sin[n\pi/4])$$

$$\sin[(2n+1)\pi/4] = \frac{1}{\sqrt{2}} \sin([n\pi/4] + \cos[n\pi/4])$$

Therefore we have [130, p.149]

$$\sum_{n=0}^{\infty} \frac{\cos[(2n+1)\pi/4]}{2n+1} = -\frac{1}{2} \log \tan \frac{\pi}{8}$$

$$\sum_{n=0}^{\infty} \frac{\sin[(2n+1)\pi/4]}{(2n+1)^2} = \frac{1}{\sqrt{2}} \left\{ \frac{1}{1^2} + \frac{1}{3^2} - \frac{1}{5^2} - \frac{1}{7^2} + \frac{1}{9^2} + \ldots \right\}$$

and the Ramanujan/Berndt result (6.30y) easily follows.

In his letter to Hardy, Ramanujan [76, p.350] also reported that

$$\frac{2}{3} \int_0^1 \frac{\tan^{-1} x}{x} dx - \int_0^{2-\sqrt{3}} \frac{\tan^{-1} x}{x} dx = \frac{\pi}{12} \log(2+\sqrt{3})$$

See also the recent paper by Cho et al [43d] which considers the integrals $\int_0^t \frac{x^m}{\sin x} dx$ and in particular

$$\int_0^{\pi/3} \frac{x}{\sin x} dx = -\frac{\pi^2}{6}\log 3 - \frac{\sqrt{3}}{9}\pi^2 + \frac{\sqrt{3}}{18}\varsigma\left(2,\frac{1}{6}\right)$$

In their paper "The evaluation of character Euler double sums" [30b], Borwein, Zucker and Boersma make reference to the following integral

$$K_p = \int_0^{\pi/2} \frac{x^p}{\sin x} \log(2\sin x) dx$$

and we may obtain a series representation for this integral as follows. Using (6.12) we have

$$\int_0^{\pi/2} \frac{x^p}{\sin x} \log(2\sin x) dx = 2\sum_{n=0}^\infty \int_0^{\pi/2} x^p \log(2\sin x) \sin(2n+1)x \, dx$$

and employing (7.8) this becomes

$$= -2\sum_{n=0}^\infty \sum_{k=1}^\infty \int_0^{\pi/2} x^p \frac{\cos 2kx}{k} \sin(2n+1)x \, dx$$

We have the familiar trigonometric identity

$$\sin(2n+1)x \cos 2kx = \frac{1}{2}\left(\sin[2(n+k)+1]x + \sin[2(n-k)+1]x\right)$$

and hence we get

$$K_{2p-1} = -\sum_{n=0}^\infty \sum_{k=1}^\infty \frac{1}{k} \int_0^{\pi/2} x^{2p-1}\left(\sin[2(n+k)+1]x + \sin[2(n-k)+1]x\right) dx$$

We then use the formula (6.30o) quoted by Amigó

$$\frac{1}{(2p-1)!} \int_0^{\pi/2} x^{2p-1} \sin(2m+1)x \, dx = \sum_{j=1}^p (-1)^{j+1} \frac{(\pi/2)^{2p-2j}}{(2p-2j)!} \frac{(-1)^m}{(2m+1)^{2j}}$$

to obtain

(6.30z) $$K_{2p-1} = -(2p-1)!\sum_{n=0}^{\infty}\sum_{k=1}^{\infty}\frac{1}{k}\sum_{j=1}^{p}(-1)^{j+1}\frac{(\pi/2)^{2p-2j}}{(2p-2j)!}\frac{(-1)^{n+k}}{[2(n+k)+1]^{2j}}$$

$$-(2p-1)!\sum_{n=0}^{\infty}\sum_{k=1}^{\infty}\frac{1}{k}\sum_{j=1}^{p}(-1)^{j+1}\frac{(\pi/2)^{2p-2j}}{(2p-2j)!}\frac{(-1)^{n-k}}{[2(n-k)+1]^{2j}}$$

Can this expression be simplified?

**Example 10:**

With the Bernoulli and Euler polynomials (the relevant properties of which are set out in Appendix A), we can readily establish the following integral identities (the modus operandi is outlined in Section 5 together with the identities (6.5a), (6.8a) and (6.14a) respectively)

(6.31) $$\varsigma(2n+1) = (-1)^{n+1}\frac{(2\pi)^{2n+1}}{(2n+1)!}\int_{0}^{1/2}B_{2n+1}(x)\cot(\pi x)\,dx$$

(6.32) $$\varsigma(2n+1) = (-1)^{n+1}\frac{(2\pi)^{2n+1}}{(1-2^{-2n})(2n+1)!}\int_{0}^{1/2}B_{2n+1}(x)\tan(\pi x)\,dx$$

(6.33) $$\varsigma(2n+1) = (-1)^{n}\frac{\pi^{2n+1}}{2(1-2^{-(2n+1)})(2n)!}\int_{0}^{1/2}E_{2n}(x)\csc(\pi x)\,dx$$

The integral (6.31) appears in Abramowitz and Stegun [1, p.807] and the integrals (6.32) and (6.33) were obtained by Cvijović and Klinowski [50] in 2002 (and rediscovered by the author in 2004). Similar identities were also derived by Espinosa and Moll [59] in the form

(6.34) $$\int_{0}^{1}B_{2n}(x)\log\sin\pi x\,dx = (-1)^{n}\frac{(2n)!\varsigma(2n+1)}{(2\pi)^{2n}}$$

(6.35) $$\int_{0}^{1}B_{2n+1}(x)\log\sin\pi x\,dx = 0$$

and it is easily shown that equation (6.34) above is equivalent to (6.31) following a simple integration by parts, and making reference to (6.45). Alternative proofs are contained in (4.4.229v) in Volume IV.

We also have (see the Wolfram MathWorld website relating to the Riemann zeta function)

$$(6.36) \quad \varsigma(2n+1) = (-1)^{n+1} \frac{(2\pi)^{2n+1}}{2(2n+1)!} \int_0^1 B_{2n+1}(x) \cot(\pi x/2)\, dx$$

$$(6.37) \quad \varsigma(2n+1) = (-1)^{n} \frac{(2\pi)^{2n+1}}{2(2n+1)!} \int_0^1 B_{2n+1}(x) \tan(\pi x/2)\, dx$$

$$(6.38) \quad \varsigma(2n+1) = (-1)^{n} \frac{(2\pi)^{2n+1}}{4(2^{2n+1}-1)(2n)!} \int_0^1 E_{2n}(x) \cot(\pi x/2)\, dx$$

$$(6.39) \quad \varsigma(2n+1) = (-1)^{n} \frac{(2\pi)^{2n+1}}{4(2^{2n+1}-1)(2n)!} \int_0^1 E_{2n}(x) \tan(\pi x/2)\, dx$$

The reason why we have similar formulae for the integrals over two different intervals results from the well-known identity

$$(6.40) \quad \int_0^a f(x)\, dx = \int_0^{a/2} f(x)\, dx + \int_0^{a/2} f(a-x)\, dx$$

$$(6.41) \quad = 2\int_0^{a/2} f(x)\, dx \quad , \text{if } f(x) = f(a-x)$$

$$(6.42) \quad = 0 \quad , \text{if } f(x) = -f(a-x)$$

We know that $\cot x = -\cot(\pi - x)$ and hence $\cot \pi x = -\cot \pi(1-x)$ (using the substitution $x \to \pi x$). In addition from (A.14) we have

$$(6.43) \quad B_{2n+1}(x) = -B_{2n+1}(1-x)$$

Hence we obtain

$$(6.44) \quad B_{2n+1}(x) \cot \pi x = f(x) = f(1-x)$$

Therefore, for example, we have

$$(6.45) \quad \int_0^1 B_{2n+1}(x) \cot \pi x\, dx = 2\int_0^{1/2} B_{2n+1}(x) \cot \pi x\, dx$$

# A SURPRISING APPEARANCE BY THE BARNES DOUBLE GAMMA FUNCTION

**Example 11:**

The following well-known identity is proved in Appendix A of Volume VI

(6.46) $$x \cot x = \sum_{n=0}^{\infty} (-1)^n \frac{2^{2n} B_{2n}}{(2n)!} x^{2n} \quad , (|x| < \pi)$$

Combining this with Euler's formula (1.7)

(6.47) $$\varsigma(2n) = (-1)^{n+1} \frac{2^{2n-1} \pi^{2n} B_{2n}}{(2n)!} \quad , (n \geq 1)$$

and, letting $x \to \pi x$, we obtain

(6.48) $$\pi x \cot \pi x = -2 \sum_{n=0}^{\infty} \varsigma(2n) x^{2n} \quad , (|x| < 1)$$

Since the first term of the series (6.46) is equal to $B_0 = 1$, to be consistent with (6.48), we define $\varsigma(0) = -\frac{1}{2}$ (which, as mentioned in (3.11a), in fact also coincides with the value determined by the analytic continuation of $\varsigma(s)$ at $s = 0$). A different proof of (6.48) is shown in (6.139). With $x = 1/2$ we get

(6.48a) $$\sum_{n=1}^{\infty} \frac{\varsigma(2n)}{2^{n+1}} = 1$$

which may be compared with (E.43f) in Volume VI.

We now recall the basic identity (6.5a)

(6.49) $$\frac{1}{2} \int_a^b p(x) \cot(\alpha x / 2) \, dx = \sum_{n=0}^{\infty} \int_a^b p(x) \sin \alpha n x \, dx$$

which is valid provided $\sin(\alpha x / 2)$ has no zero in $[a, b]$ (or, alternatively, if $\sin(\alpha a / 2) = 0$ then $p(a) = 0$). We now let $\alpha = 2\pi$ in (6.49) and substitute (6.48) in (6.49) to obtain

(6.50) $$-\frac{1}{\pi} \int_a^b p(x) \sum_{n=0}^{\infty} \varsigma(2n) x^{2n-1} dx = \sum_{n=0}^{\infty} \int_a^b p(x) \sin 2\pi n x \, dx$$

The above formula can be used to generate an endless number of infinite series involving terms in $\varsigma(2n)$ by selecting appropriate functions for $p(x)$ and varying the interval of integration $[a,b]$. A number of such identities are recorded in the book by Srivastava and Choi, "Series Associated with the Zeta and Related Functions" [126] and a few examples are derived below.

Appendix A in Volume VI of this paper also lists series expansions for $\tan x$, $x/\sin x$ and $1/\cos x$ (and similar expansions for the corresponding hyperbolic functions can be easily derived): hence the methods employed below can be replicated to systematically produce a virtual encyclopaedia of identities involving $\varsigma(2n)$, $\lambda(2n)$ and $\beta(2n+1)$.

With $p(x) = x$ in (6.50) and $[a,b] = [0,1/2]$ we have

$$(6.51) \qquad -\frac{1}{\pi}\sum_{n=0}^{\infty}\frac{\varsigma(2n)}{(2n+1)2^{2n+1}} = \sum_{n=0}^{\infty}\int_0^{1/2} x\sin 2\pi nx\, dx$$

As in (3.3), integration by parts gives us

$$\int_0^{1/2} x\sin 2n\pi x\, dx = -\frac{x\cos 2n\pi x}{2n\pi} + \frac{\sin 2n\pi x}{(2n\pi)^2}\Big|_0^{1/2} = \frac{(-1)^{n+1}}{4n\pi}, (n\geq 1)$$

Therefore we have

$$\sum_{n=1}^{\infty}\int_0^{1/2} x\sin 2n\pi x\, dx = \frac{1}{4\pi}\sum_{n=1}^{\infty}\frac{(-1)^{n+1}}{n} = \frac{\log 2}{4\pi}$$

and hence we get

$$(6.52) \qquad \sum_{n=0}^{\infty}\frac{\varsigma(2n)}{(2n+1)2^{2n+1}} = -\frac{1}{4}\log 2$$

Or, alternatively, (with the summation starting at $n=1$) we have

$$(6.53) \qquad \sum_{n=1}^{\infty}\frac{\varsigma(2n)}{(2n+1)2^{2n+1}} = \frac{1}{4} - \frac{1}{4}\log 2$$

In his recent paper, Dalai [51] has derived a structurally similar identity to (6.52) involving the Dirichlet beta function defined in (A.27)

$$(6.54) \qquad \sum_{n=0}^{\infty}\frac{\beta(2n+1)}{(2n+1)2^{2n}} = \log(1+\sqrt{2})$$

The Dalai formula may also be obtained by integrating (6.30d)

$$\sec(\pi x/2) = \frac{4}{\pi}\sum_{n=0}^{\infty}\beta(2n+1)x^{2n} \quad , \ |x|<1$$

$$\int_0^t \sec(\pi x/2)dx = \frac{4}{\pi}\sum_{n=0}^{\infty}\frac{\beta(2n+1)}{2n+1}t^{2n+1}$$

$$\int_0^t \sec(\pi x/2)dx = \frac{\pi}{2}\int_0^{\pi t/2}\sec u\, du = \frac{\pi}{2}\log\left[\frac{1+\tan(u/2)}{1-\tan(u/2)}\right]_0^{\pi t/2}$$

$$= \frac{\pi}{2}\log\left[\frac{1+\tan(\pi t/4)}{1-\tan(\pi t/4)}\right]$$

and hence $t=1/2$ produces the required identity.

Perhaps the closed form formula for $\varsigma(3)$ contains a factor of $\frac{p}{q}\pi^2\log\left(1+\sqrt{2}\right)$!!

The following analysis corroborates formula (6.53). This part was mainly written in 2003 and, at that time, I knew next to nothing about the Barnes double gamma function; my analysis in Volume II(a) came about three years later.

In [126, p.212] we have

(6.55) $$\sum_{n=1}^{\infty}\frac{\varsigma(2n)z^{2n+1}}{(2n+1)} = \frac{1}{2}\left\{[1-\log(2\pi)]z + \log\frac{G(1+z)}{G(1-z)}\right\} \ , \ (|z|<1)$$

where $G(z)$ is the Barnes double gamma function [126, p.24]. In 1900 Barnes defined the G-function as

(6.55a) $$G(z+1) = (2\pi)^{z/2}\exp\left[-\frac{1}{2}\left(z^2+\gamma z^2+z\right)\right]\prod_{n=1}^{\infty}\left(1+\frac{z}{n}\right)^n e^{-z+z^2/2n}$$

With $z=1/2$ in (6.55) we have

(6.56) $$\sum_{n=1}^{\infty}\frac{\varsigma(2n)}{(2n+1)2^{2n+1}} = \frac{1}{2}\left\{\frac{1}{2}[1-\log(2\pi)] + \log\frac{G(3/2)}{G(1/2)}\right\}$$

We have from [126, p.25]

(6.57) $$G(z+1) = \Gamma(z)G(z)$$

and therefore with $z=1/2$ we obtain

(6.58) $$\frac{G(3/2)}{G(1/2)} = \Gamma(1/2)$$

From the definition of $\Gamma(z)$ in (4.3.1) we have

(6.59) $$\Gamma(1/2) = \int_0^\infty \frac{e^{-t}}{\sqrt{t}} dt$$

and, with the substitution $u = \sqrt{t}$, this becomes

(6.60) $$\Gamma(1/2) = 2\int_0^\infty e^{-u^2} du = \sqrt{\pi}$$

This is a very familiar integral and Sandham's elegant derivation is shown in Appendix B. Alternatively, it can be easily ascertained by substituting $z = 1/2$ in Euler's reflection formula for the gamma function (which is proved in Appendix C).

(6.61) $$\Gamma(z)\Gamma(1-z) = \frac{\pi}{\sin \pi z}$$

This formula actually indicates that $\Gamma(1/2) = \pm\sqrt{\pi}$: however, the negative answer is discarded because it is clear from the definition that $\Gamma(x) > 0 \ \forall x \in (0, \infty)$. Could the negative solution have any significance (recall Dirac's reason for postulating the existence of the anti-electron, the positron)?

For $z = 1/2$, the Srivastava and Choi result (6.56) is therefore in agreement with my simple analysis.

(6.62) $$\sum_{n=1}^\infty \frac{\varsigma(2n)}{(2n+1)2^{2n+1}} = \frac{1}{4} - \frac{1}{4}\log 2$$

With $p(x) = x$ and $[a,b] = [0, z]$ in (6.50) we have therefore shown that for $|z| < 1$

(6.63) $$\sum_{n=1}^\infty \frac{\varsigma(2n) z^{2n+1}}{(2n+1)} = \frac{1}{2}\left\{[1-\log(2\pi)]z + \log\frac{G(1+z)}{G(1-z)}\right\}$$

$$= \frac{z}{2} - \pi \sum_{n=1}^\infty \int_0^z x \sin 2n\pi x \, dx$$

(6.64) $$= \frac{z}{2} + \frac{z}{2}\sum_{n=1}^\infty \frac{\cos 2n\pi z}{n} - \frac{1}{4\pi}\sum_{n=1}^\infty \frac{\sin 2n\pi z}{n^2}$$

(6.65) $$= \frac{z}{2} - \frac{z}{2}\log(2\sin \pi z) - \frac{1}{4\pi}\text{Cl}_2(2\pi z)$$

where in (6.65) we have employed the familiar Fourier series, which will be easily derived in (7.8), and $\text{Cl}_2(t)$ is the Clausen function [100, p.101] defined by

$$(6.66) \qquad \text{Cl}_2(t) = -\int_0^t \log[2\sin(u/2)] \, du = \sum_{n=1}^\infty \frac{\sin nt}{n^2}$$

We therefore have

$$(6.67) \qquad \log\frac{G(1+z)}{G(1-z)} = z - z\log(2\sin\pi z) - z[1 - \log(2\pi)] - \frac{1}{2\pi}\text{Cl}_2(2\pi z)$$

$$(6.68) \qquad = z - z\log(2\sin\pi z) - z[1 - \log(2\pi)] + \int_0^z \log(2\sin\pi t) \, dt$$

Eliminating the factor of $z\log 2$ we obtain

$$(6.69) \qquad \log\frac{G(1+z)}{G(1-z)} = -z\log\left[\frac{\sin\pi z}{2\pi}\right] + \int_0^z \log(\sin\pi t) \, dt$$

which is in agreement with Adamchik's results in [5a] and [45]. More detailed information regarding the Barnes function may be found in Adamchik's papers [5a] and [6a].

We have therefore shown that

$$(6.69a) \qquad \log\frac{G(1+z)}{G(1-z)} = -z\log\left[\frac{\sin\pi z}{2\pi}\right] + \int_0^z \log(\sin\pi t) \, dt$$

and using integration by parts, we see that this is equivalent to the following integral formula originally due to Kinkelin (1860)

$$(6.69b) \qquad \log\frac{G(1+z)}{G(1-z)} = z\log(2\pi) - \int_0^z \pi t \cot\pi t \, dt$$

(this is recorded as an exercise in Whittaker and Watson [135, p.264]).

Referring to (6.48)

$$\pi z \cot\pi z = -2\sum_{n=0}^\infty \varsigma(2n) z^{2n} \qquad , \; (|x|<1)$$

and using [126, p.14]

$$\psi(1+z) - \psi(1-z) = \frac{1}{z} - \pi\cot\pi z$$

we see that

$$z\psi(1+z) - z\psi(1-z) = 2\sum_{n=1}^{\infty} \varsigma(2n) z^{2n}$$

since $\varsigma(0) = -1/2$.

Similarly (6.69b) may be expressed as

$$\log \frac{G(1+z)}{G(1-z)} = z[\log(2\pi) - 1] - \int_0^z [t\psi(1+t) - t\psi(1-t)] dt$$

Employing the definition of the Barnes function in (6.55a)

$$G(z+1) = (2\pi)^{z/2} \exp\left[-\frac{1}{2}(z^2 + \gamma z^2 + z)\right] \prod_{n=1}^{\infty} \left(1 + \frac{z}{n}\right)^n e^{-z+z^2/2n}$$

we get [126, p.32]

$$\Phi(z) = \frac{G(1+z)}{G(1-z)} = (2\pi)^z e^{-z} \prod_{n=1}^{\infty} \left(1 + \frac{z}{n}\right)^n \left(1 - \frac{z}{n}\right)^{-n} e^{-2z}$$

We then differentiate logarithmically and apply the expansion (C.42a)

$$\pi z \cot \pi z = 1 + 2\sum_{n=1}^{\infty} \frac{z^2}{z^2 - n^2}$$

to obtain

$$\frac{d}{dz} \log \Phi(z) = \log(2\pi) - \pi z \cot \pi z \quad (z \notin \mathbf{Z})$$

Integrating this, since $\Phi(0) = \frac{G(1)}{G(1)} = 1$, we obtain (6.69b). See a further proof in (4.3.87) in Volume II(a).

From (6.66) we see that

$$\mathrm{Cl}_2(2\pi t) = -2\pi \int_0^t \log[2\sin(\pi u)] \, du = \sum_{n=1}^{\infty} \frac{\sin 2\pi n t}{n^2}$$

and combining this with (6.69a) we get

(6.69c) $\quad \mathrm{Cl}_2(2\pi t) = \sum_{n=1}^{\infty} \frac{\sin 2\pi n t}{n^2} = -2\pi t \log\left[\frac{\sin \pi t}{\pi}\right] - 2\pi \log \frac{G(1+t)}{G(1-t)}$

Using the trigonometric series [48a] defined by

$$S_\nu(\alpha) = \sum_{n=0}^{\infty} \frac{\sin(2n+1)\alpha}{(2n+1)^\nu}$$

it is easily seen that

$$S_2(\alpha) = Cl_2(\alpha) - \frac{1}{4} Cl_2(2\alpha)$$

Hence we have

(6.69d) $\quad S_2(2\pi t) = \pi t \log[2\pi \cot(\pi t)] + \dfrac{\pi}{2} \log \dfrac{G(1+2t)}{G(1-2t)} - 2\pi \log \dfrac{G(1+t)}{G(1-t)}$

and this corrects the misprint in equation (4.10) of [45] (where $\pi$ is missing from the argument of the logarithm).

Letting $t = 1/4$ in (6.69c) we obtain the known result contained in [45]

(6.69e) $\quad \log \dfrac{G\left(\frac{5}{4}\right)}{G\left(\frac{3}{4}\right)} = \dfrac{1}{8}\log 2 + \dfrac{1}{4}\log \pi - \dfrac{G}{2\pi}$ where $G$ is Catalan's constant.

Letting $t = 1/4$ in (6.69d) we get

(6.69f) $\quad S_2(\pi/4) = \dfrac{\pi}{8} \log[2\pi \cot(\pi/8)] + \dfrac{\pi}{2} \log \dfrac{G\left(\frac{5}{4}\right)}{G\left(\frac{3}{4}\right)} - 2\pi \log \dfrac{G\left(\frac{9}{8}\right)}{G\left(\frac{7}{8}\right)}$

From (6.69b) we have

(6.69g) $\quad \log \dfrac{G\left(\frac{9}{8}\right)}{G\left(\frac{7}{8}\right)} = \dfrac{1}{8}\log(2\pi) - \int_0^{1/8} \pi t \cot \pi t \, dt$

Combining (6.69f) and (6.69g) we get

(6.69h)

$S_2(\pi/4) = \dfrac{\pi}{8} \log[2\pi \cot(\pi/8)] + \dfrac{\pi}{2} \log \dfrac{G\left(\frac{5}{4}\right)}{G\left(\frac{3}{4}\right)} - \dfrac{\pi}{4} \log(2\pi) + 2\pi \int_0^{1/8} \pi t \cot \pi t \, dt$

and this may be rewritten as

(6.69i) $\quad S_2(\pi/4) = \dfrac{\pi}{8}\log[2\pi\cot(\pi/8)] + \dfrac{\pi}{2}\log\dfrac{G\left(\dfrac{5}{4}\right)}{G\left(\dfrac{3}{4}\right)} - \dfrac{\pi}{4}\log(2\pi) + 2\int_0^{\pi/8} t\cot t\, dt$

We now employ (6.5a) with $\alpha = 2\pi$ to obtain

$$2\int_0^{\pi/8} t\cot t\, dt = 4\sum_{n=1}^{\infty}\int_0^{\pi/8} t\sin 2nt\, dt$$

Integration by parts gives

$$\int_0^{\pi/8} t\sin 2nt\, dt = -\dfrac{t\cos 2nt}{2n} + \dfrac{\sin 2nt}{4n^2}\bigg|_0^{\pi/8}$$

$$= -\dfrac{\pi}{16}\dfrac{\cos(n\pi/4)}{n} + \dfrac{1}{4}\dfrac{\sin(n\pi/4)}{n^2}$$

and therefore we get

$$2\int_0^{\pi/8} t\cot t\, dt = -\dfrac{\pi}{4}\sum_{n=1}^{\infty}\dfrac{\cos(n\pi/4)}{n} + \sum_{n=1}^{\infty}\dfrac{\sin(n\pi/4)}{n^2}$$

A quick glance forward to (7.8) shows that

$$\log[2\sin(x/2)] = -\sum_{n=1}^{\infty}\dfrac{\cos nx}{n}$$

and hence we have

(6.69j) $\quad 2\int_0^{\pi/8} t\cot t\, dt = \dfrac{\pi}{4}\log[2\sin(\pi/8)] + \sum_{n=1}^{\infty}\dfrac{\sin(n\pi/4)}{n^2}$

We have by separating the even and odd terms

$$\sum_{n=1}^{\infty}\dfrac{\sin(n\pi/4)}{n^2} = \sum_{n=0}^{\infty}\dfrac{\sin[(2n+1)\pi/4]}{(2n+1)^2} + \sum_{n=1}^{\infty}\dfrac{\sin[(2n)\pi/4]}{(2n)^2}$$

$$= S_2(\pi/4) + \dfrac{1}{4}G$$

and therefore we get

$$2\int_0^{\pi/8} t\cot t\,dt = \frac{\pi}{4}\log[2\sin(\pi/8)] + S_2(\pi/4) + \frac{1}{4}G$$

Substituting this in (6.69i) we obtain

(6.69k) $\quad \dfrac{\pi}{8}\log[2\pi\cot(\pi/8)] + \dfrac{\pi}{2}\log\dfrac{G\left(\frac{5}{4}\right)}{G\left(\frac{3}{4}\right)} - \dfrac{\pi}{4}\log(2\pi) + \dfrac{\pi}{4}\log[2\sin(\pi/8)] + \dfrac{1}{4}G = 0$

This is easily simplified by writing $\dfrac{1}{4}\log[\sin(\pi/8)] = \dfrac{1}{8}\log\left[\sin^2(\pi/8)\right]$ and, using (6.69g), rather disappointingly we do in fact find that everything cancels out! $G$ is proving as elusive as $\varsigma(3)$.

In (6.69j) we showed that

$$\int_0^{\pi/8} x\cot x\,dx = \frac{\pi}{8}\log\left[2\sin(\pi/8)\right] + \frac{1}{2}\sum_{n=1}^{\infty}\frac{\sin(n\pi/4)}{n^2}$$

Applying the half-angle trigonometric formula, this becomes

(6.69l) $\quad \displaystyle\int_0^{\pi/8} x\cot x\,dx = \frac{\pi}{16}\log\left[\sqrt{2}-2\right] + \frac{1}{2}\sum_{n=1}^{\infty}\frac{\sin(n\pi/4)}{n^2}$

We also have already shown that

$$\sum_{n=1}^{\infty}\frac{\sin(n\pi/4)}{n^2} = \sum_{n=0}^{\infty}\frac{\sin[(2n+1)\pi/4]}{(2n+1)^2} + \sum_{n=1}^{\infty}\frac{\sin[(2n)\pi/4]}{(2n)^2}$$

$$= S_2(\pi/4) + \frac{1}{4}G$$

In 2003, by a much more complex route, Choi, Srivastava and Adamchik [45] obtained the following representation

(6.69m)

$$\int_0^{\pi/8} x\cot x\,dx = \frac{\pi}{16}\log\left[(2-\sqrt{2})\pi\right] + \frac{1}{8}\left[1-\sqrt{2}\right]G + \frac{1}{64}\left[\sqrt{2}\varsigma\left(2,\tfrac{1}{8}\right) - 2\left(\sqrt{2}+1\right)\pi^2\right]$$

It therefore **seems** that we have

$$\sum_{n=1}^{\infty}\frac{\sin(n\pi/4)}{n^2}=\frac{\pi}{8}\log\pi+\frac{1}{4}\left[1-\sqrt{2}\right]G+\frac{1}{32}\left[\sqrt{2}\varsigma\left(2,\tfrac{1}{8}\right)-2\left(\sqrt{2}+1\right)\pi^2\right]$$

The presence of the term involving $\frac{\pi}{8}\log\pi$ certainly looked interesting!

Since the series $\sum_{n=1}^{\infty}\frac{\sin(n\pi/4)}{n^2}$ is absolutely convergent we may rearrange it as follows (using $\omega$ as a convenient abbreviation for $\sin(\pi/4)=\sqrt{2}/2$)

$$\sum_{n=1}^{\infty}\frac{\sin(n\pi/4)}{n^2}=\frac{\omega}{1^2}+\frac{1}{2^2}+\frac{\omega}{3^2}+\frac{0}{4^2}-\frac{\omega}{5^2}-\frac{1}{6^2}-\frac{\omega}{7^2}-\frac{0}{8^2}+\ldots$$

$$=\omega\left[\frac{1}{1^2}+\frac{1}{3^2}-\frac{1}{5^2}-\frac{1}{7^2}+\frac{1}{9^2}+\frac{1}{11^2}-\ldots\right]+\left[\frac{1}{2^2}-\frac{1}{6^2}+\frac{1}{10^2}-\frac{1}{14^2}+\frac{1}{18^2}-\frac{1}{22^2}+\ldots\right]$$

$$=\omega\left[\frac{1}{1^2}+\frac{1}{3^2}-\frac{1}{5^2}-\frac{1}{7^2}+\frac{1}{9^2}+\frac{1}{11^2}-\ldots\right]+\frac{1}{2^2}\left[\frac{1}{1^2}-\frac{1}{3^2}+\frac{1}{5^2}-\frac{1}{7^2}+\frac{1}{9^2}-\frac{1}{11^2}+\ldots\right]$$

$$=\omega\left[\frac{1}{1^2}+\frac{1}{3^2}-\frac{1}{5^2}-\frac{1}{7^2}+\frac{1}{9^2}+\frac{1}{11^2}-\ldots\right]+\frac{G}{4}$$

It is easily seen from the definition of the Hurwitz zeta function that

$$\frac{\sqrt{2}}{64}\varsigma\left(2,\frac{1}{8}\right)=\frac{\sqrt{2}}{64}\sum_{n=0}^{\infty}\frac{1}{\left(n+\frac{1}{8}\right)^2}=2\omega\sum_{n=0}^{\infty}\frac{1}{(8n+1)^2}$$

and this is certainly contained within the series in parentheses above. Similarly, $\pi^2$ may obviously be expressed as a multiple of $\varsigma(2)$. It seemed to me that someone more adept at series rearrangement may be able to deduce a relationship between $G$ and $\pi\log\pi$. Initially, it looked as if there may be gold in them there hills (but alas it turned out to be mere iron pyrites!).

Unfortunately, after I started looking at this I discovered that there was a misprint in equation (5.16) of [45]: the Choi, Srivastava and Adamchik equation (6.69m) should actually read as follows

(6.69n)
$$\int_0^{\pi/8}x\cot x\,dx=\frac{\pi}{16}\log\left[2-\sqrt{2}\right]+\frac{1}{8}\left[1-\sqrt{2}\right]G+\frac{1}{64}\left[\sqrt{2}\varsigma\left(2,\frac{1}{8}\right)-2\left(\sqrt{2}+1\right)\pi^2\right]$$

(i.e. the interesting factor of $\pi\log\pi$ is no longer there).

We now continue this exercise by giving an elementary proof of another integral which appears in [45]. Let us consider (6.5a) with $p(x) = x^2$ and $[a,b] = [0, \pi/4]$: we have

$$\frac{1}{2} \int_0^{\pi/4} x^2 \cot x \, dx = \sum_{n=1}^{\infty} \int_0^{\pi/4} x^2 \sin 2nx \, dx$$

Integration by parts gives

$$\int_0^{\pi/4} x^2 \sin 2nx \, dx = \left. \frac{\cos 2nx}{4n^3} - \frac{x^2 \cos 2nx}{2n} + \frac{x \sin 2nx}{2n^2} \right|_0^{\pi/4}$$

$$= \frac{1}{4n^3}[\cos(n\pi/2) - 1] - \frac{\pi^2}{32} \frac{\cos(n\pi/2)}{n} + \frac{\pi}{8} \frac{\sin(n\pi/2)}{n^2}$$

and therefore we get

$$\int_0^{\pi/4} x^2 \cot x \, dx = \frac{1}{2} \sum_{n=1}^{\infty} \frac{[\cos(n\pi/2) - 1]}{n^3} - \frac{\pi^2}{16} \sum_{n=1}^{\infty} \frac{\cos(n\pi/2)}{n} + \frac{\pi}{4} \sum_{n=1}^{\infty} \frac{\sin(n\pi/2)}{n^2}$$

We have (see [126, p.293])

$$\sum_{n=1}^{\infty} \frac{\cos(n\pi/2)}{n^s} = \frac{0}{1^s} - \frac{1}{2^s} + \frac{0}{3^s} + \frac{1}{4^s} - \ldots = -\frac{1}{2^s} \varsigma_a(s) = \frac{1}{2^s}(2^{1-s} - 1)\varsigma(s)$$

and hence the above equation easily simplifies to

(6.69o) $$\int_0^{\pi/4} x^2 \cot x \, dx = -\frac{35}{64} \varsigma(3) + \frac{\pi^2}{32} \log 2 + \frac{\pi G}{4}$$

which is equation (5.13) of [45]. As Jane Austin might say, it is a truth universally acknowledged that my method is substantially easier than the one employed by Choi, Srivastava and Adamchik in [45]!

Repeating the exercise again with $[0, \pi/6]$ as the interval of integration we obtain

$$\int_0^{\pi/6} x^2 \cot x \, dx = \frac{1}{2} \sum_{n=1}^{\infty} \frac{[\cos(n\pi/3) - 1]}{n^3} - \frac{\pi^2}{36} \sum_{n=1}^{\infty} \frac{\cos(n\pi/3)}{n} + \frac{\pi}{6} \sum_{n=1}^{\infty} \frac{\sin(n\pi/3)}{n^2}$$

From [126, p.293] we have

(6.69p) $$\sum_{n=1}^{\infty}\frac{\cos(n\pi/3)}{n^s} = \frac{1}{2}(6^{1-s}-3^{1-s}-2^{1-s}+1)\varsigma(s)$$

(6.69q) $$\sum_{n=1}^{\infty}\frac{\sin(n\pi/3)}{n^s} = \sqrt{3}\left[\frac{3^{-s}-1}{2}\varsigma(s)+6^{-s}\left\{\varsigma\left(s,\frac{1}{6}\right)+\varsigma\left(s,\frac{1}{3}\right)\right\}\right]$$

and hence

(6.69r) $$\sum_{n=1}^{\infty}\frac{\cos(n\pi/3)}{n^3} = \frac{1}{3}\varsigma(3)$$

(6.69s) $$\sum_{n=1}^{\infty}\frac{\sin(n\pi/3)}{n^2} = \sqrt{3}\left[-\frac{4}{9}\varsigma(2)+\frac{1}{36}\left\{\varsigma\left(2,\frac{1}{6}\right)+\varsigma\left(2,\frac{1}{3}\right)\right\}\right]$$

Accordingly we obtain

(6.69t)
$$\int_0^{\pi/6} x^2 \cot x\, dx = -\frac{1}{3}\varsigma(3)+\frac{\pi^2}{36}\log[2\sin(\pi/6)]+\frac{\sqrt{3}}{6}\pi\left[-\frac{4}{9}\varsigma(2)+\frac{1}{36}\left\{\varsigma\left(2,\frac{1}{6}\right)+\varsigma\left(2,\frac{1}{3}\right)\right\}\right]$$

Let us consider one further example (this stuff is addictive): we have

$$\int_0^{\pi/12} t\cot t\, dt = 2\sum_{n=1}^{\infty}\int_0^{\pi/12} t\sin 2nt\, dt$$

Integration by parts gives

$$\int_0^{\pi/12} t\sin 2nt\, dt = -\frac{t\cos 2nt}{2n}+\frac{\sin 2nt}{4n^2}\bigg|_0^{\pi/12}$$

$$= -\frac{\pi}{24}\frac{\cos(n\pi/6)}{n}+\frac{1}{4}\frac{\sin(n\pi/6)}{n^2}$$

and therefore we get

$$\int_0^{\pi/12} t\cot t\, dt = -\frac{\pi}{12}\sum_{n=1}^{\infty}\frac{\cos(n\pi/6)}{n}+\frac{1}{2}\sum_{n=1}^{\infty}\frac{\sin(n\pi/6)}{n^2}$$

Reference again to (7.8) results in

$$\int_0^{\pi/12} t \cot t \, dt = -\frac{\pi}{12} \log[2\sin(\pi/12)] + \frac{1}{2} \sum_{n=1}^{\infty} \frac{\sin(n\pi/6)}{n^2}$$

Using the half-angle trigonometric formula

$$\sin^2(\pi/12) = \frac{1}{2}[1 - \cos(\pi/6)] = \frac{1}{4}\left[2 - \sqrt{3}\right]$$

this becomes

(6.69u) $$\int_0^{\pi/12} t \cot t \, dt = -\frac{\pi}{24} \log\left[2 - \sqrt{3}\right] + \frac{1}{2} \sum_{n=1}^{\infty} \frac{\sin(n\pi/6)}{n^2}$$

This should be compared with equation (5.22) in [45]: however, I am not convinced that the latter is correct because, inter alia, it contains a factor of $\pi \log \pi$.

We saw in (4.3.158) in Volume II(a) that

$$\int_0^x \pi u \cot \pi u \, du = \varsigma'(-1, x) - \varsigma'(-1, 1-x) + x \log(2 \sin \pi x)$$

and therefore we have

$$\int_0^x \pi u \cot \pi u \, du = \frac{1}{\pi} \int_0^{\pi x} t \cot t \, dt$$

We have

$$\int_0^{\pi x} t \cot t \, dt = 2 \sum_{n=1}^{\infty} \int_0^{\pi x} t \sin 2nt \, dt = -\sum_{n=1}^{\infty} \frac{\pi x \cos(2n\pi x)}{n} + \frac{1}{2} \sum_{n=1}^{\infty} \frac{\sin(2n\pi x)}{n^2}$$

and we therefore see that

$$\varsigma'(-1, x) - \varsigma'(-1, 1-x) + x \log(2 \sin \pi x) = -\sum_{n=1}^{\infty} \frac{x \cos(2n\pi x)}{n} + \frac{1}{2\pi} \sum_{n=1}^{\infty} \frac{\sin(2n\pi x)}{n^2}$$

Then, having regard to (7.8), we see that

(6.69v) $$\varsigma'(-1, x) - \varsigma'(-1, 1-x) = \frac{1}{2\pi} \sum_{n=1}^{\infty} \frac{\sin(2n\pi x)}{n^2}$$

Let us now integrate (6.64) to obtain

(6.70) $$\sum_{n=1}^{\infty}\frac{\varsigma(2n)}{(2n+1)(2n+2)2^{2n+2}}=\frac{1}{8}+\frac{1}{2}\int_0^{1/2}\sum_{n=1}^{\infty}\frac{z\cos 2n\pi z}{n}\,dz-\frac{1}{4\pi}\int_0^{1/2}\sum_{n=1}^{\infty}\frac{\sin 2n\pi z}{n^2}\,dz$$

Completing the integration, we have yet another rediscovery of the now well-known result (see the discussion in Section 8 for more details of this zeta expansion)

(6.71) $$\varsigma(3)=-\frac{4\pi^2}{7}\sum_{n=0}^{\infty}\frac{\varsigma(2n)}{(2n+1)(2n+2)2^{2n}}$$

As a further example, we now let $p(x)=x^2$ and $[a,b]=[0,1/2]$ in (6.50) and obtain

(6.72) $$-\frac{1}{\pi}\sum_{n=0}^{\infty}\frac{\varsigma(2n)}{(2n+2)2^{2n+2}}=\sum_{n=1}^{\infty}\int_0^{1/2}x^2\sin 2\pi nx\,dx$$

Integration by parts shows that

(6.73) $$\int_0^{1/2}x^2\sin 2\pi nx\,dx=\left.\frac{\cos 2\pi nx}{4(\pi n)^3}-\frac{x^2\cos 2\pi nx}{2\pi n}+\frac{x\sin 2\pi nx}{2(\pi n)^2}\right|_0^{1/2}$$

$$=\frac{(-1)^n}{4(\pi n)^3}-\frac{1}{4(\pi n)^3}-\frac{(-1)^n}{8\pi n}\quad,n\geq 1$$

Therefore we have

(6.74) $$\sum_{n=1}^{\infty}\int_0^{1/2}x^2\sin 2\pi nx\,dx=\sum_{n=1}^{\infty}\left\{\frac{(-1)^n}{4(\pi n)^3}-\frac{1}{4(\pi n)^3}+\frac{(-1)^{n+1}}{8\pi n}\right\}$$

(6.75) $$=-\frac{1}{4\pi^3}\varsigma_a(3)-\frac{1}{4\pi^3}\varsigma(3)+\frac{1}{8\pi}\log 2$$

(6.76) $$=-\frac{7}{16\pi^3}\varsigma(3)+\frac{1}{8\pi}\log 2$$

Hence for (6.72) we have

(6.77) $$\frac{1}{8}\sum_{n=1}^{\infty}\frac{\varsigma(2n)}{(n+1)2^{2n}}=\frac{1}{16}-\pi\sum_{n=1}^{\infty}\int_0^{1/2}x^2\sin 2\pi nx\,dx$$

(6.78) $$=\frac{1}{16}+\frac{7}{16\pi^2}\varsigma(3)-\frac{1}{8}\log 2$$

Therefore we obtain

(6.79) $$\sum_{n=1}^{\infty}\frac{\varsigma(2n)}{(n+1)2^{2n}}=\frac{1}{2}+\frac{7}{2\pi^2}\varsigma(3)-\log 2$$

An expression for the left hand side of (6.79) is also contained in the book by Srivastava and Choi [126, p.217]: in this book [126, p.215] the authors give the following formula in the case where $|z|<|a|$

(6.80) $$\sum_{k=1}^{\infty}\frac{\varsigma(2k,a)z^{2k+2}}{k+1}=2(a-1)\log G(a)-2(a-1)^2\log\Gamma(a)+[1-\log(2\pi)]\frac{z^2}{2}$$

$$+(a-1)^2\log\Gamma(a+z)\Gamma(a-z)+(z-a+1)\log G(a+z)$$

$$-(z+a-1)\log G(a-z)-\int_0^z\log G(t+a)dt-\int_0^{-z}\log G(t+a)dt$$

Taking $a=1$ they report the following result which is valid for $|z|<1$

(6.81)

$$\sum_{k=1}^{\infty}\frac{\varsigma(2k)z^{2k+2}}{k+1}=[1-\log(2\pi)]\frac{z^2}{2}+z\log\frac{G(1+z)}{G(1-z)}-\int_0^z\log G(1+t)dt-\int_0^{-z}\log G(1+t)dt$$

and taking $z=1/2$ they obtain [126, p.217]

(6.82) $$\sum_{n=1}^{\infty}\frac{\varsigma(2n)}{(n+1)2^{2n}}=\frac{1}{2}+\log\left(2^{-1}B^{14}\right)$$

where $B$ is defined in [126, p.36] as

(6.83) $$B=\lim_{n\to\infty}\left[\sum_{k=1}^n k^2\log k-\left(\frac{n^3}{3}+\frac{n^2}{2}+\frac{n}{6}\right)\log n+\frac{n^3}{9}-\frac{n}{12}\right]$$

$B$ is one of the generalised Glaisher-Kinkelin constants which are referred to in (4.4.211). Mercifully, we are told in [126, p.100] that (see also Appendix F of Volume VI)

(6.84) $$\log B=-\varsigma'(-2)=\frac{\varsigma(3)}{4\pi^2}$$

A proof of this is also given in Volume III (and a further proof is presented in Volume VI).

Using (6.84) we may write (6.82) as

(6.85) $$\sum_{n=1}^{\infty}\frac{\varsigma(2n)}{(n+1)2^{2n}} = \frac{1}{2} + \frac{7\varsigma(3)}{2\pi^2} - \log 2$$

which is in agreement with my considerably more elementary analysis (6.79).

## A CONNECTION WITH SINE AND COSINE INTEGRALS

We may write (6.50) as

(6.86) $$\sum_{n=1}^{\infty}\int_a^b \varsigma(2n)p(x)x^{2n-1}dx = \frac{1}{2}\int_a^b \frac{p(x)}{x}dx - \pi\sum_{n=1}^{\infty}\int_a^b p(x)\sin 2n\pi x\,dx$$

Let $p(x) = x$ and integrate over the interval $[0,t]$ to obtain

(6.87) $$\sum_{n=1}^{\infty}\frac{\varsigma(2n)t^{2n+1}}{2n+1} = \frac{1}{2}t - \pi\sum_{n=1}^{\infty}\left\{-\frac{t\cos 2n\pi t}{2n\pi} + \frac{\sin 2n\pi t}{(2n\pi)^2}\right\}$$

Now divide (6.87) by $t$ and integrate the result over the interval $[0,x]$ (it may be shown that term by term integration is valid).

(6.88) $$\sum_{n=1}^{\infty}\frac{\varsigma(2n)x^{2n+1}}{(2n+1)^2} = \frac{1}{2}x - \sum_{n=1}^{\infty}\int_0^x\left\{-\frac{\cos 2n\pi t}{2n} + \frac{\sin 2n\pi t}{4\pi n^2 t}\right\}dt$$

(6.89) $$= \frac{1}{2}x + \frac{1}{4\pi}\sum_{n=1}^{\infty}\frac{\sin 2n\pi x}{n^2} - \frac{1}{4\pi}\sum_{n=1}^{\infty}\frac{1}{n^2}\int_0^x \frac{\sin 2n\pi t}{t}dt$$

Let us now consider the integral

(6.90) $$\int_0^x \log t \cdot \cos nt\,dt = \log t\,\frac{\sin nt}{n}\bigg|_0^x - \int_0^x \frac{\sin nt}{nt}dt$$

We have

$$\lim_{t\to 0}\log t \sin nt = \lim_{t\to 0}\frac{\log t}{1/\sin nt}$$

and using L'Hôpital's rule

$$= -\lim_{t\to 0}\frac{\sin^2 nt}{nt\cos nt} = -\lim_{t\to 0}\frac{\sin nt}{nt}\tan nt = 0$$

Therefore we get

$$\int_0^x \log t \cdot \cos nt \, dt = \frac{\sin nx \log x}{n} - \int_0^x \frac{\sin nt}{nt} dt$$

$$= \frac{\sin nx \log x}{n} - \frac{1}{n} \int_0^{nx} \frac{\sin u}{u} du$$

and hence we get

(6.90a) $$\int_0^x \log t \cdot \cos nt \, dt = \frac{\sin nx \log x}{n} - \frac{Si(nx)}{n}$$

(6.90ai) $$\int_\alpha^x \log t \cdot \cos nt \, dt = \frac{\sin nx \log x}{n} - \frac{\sin n\alpha \log \alpha}{n} - \frac{Si(nx)}{n} + \frac{Si(n\alpha)}{n}$$

where $Si(x)$ is the sine integral function defined by G&R [74, p.878] and [1, p.231] as

(6.90b) $$Si(x) = \int_0^x \frac{\sin t}{t} dt \quad , \quad Si(0) = 0$$

For reference we may note that

(6.90bi) $$Si(ax) = \int_0^x \frac{\sin at}{t} dt$$

We therefore have

(6.90c) $$\int_0^\pi \log t \cdot \cos nt \, dt = -\frac{Si(n\pi)}{n}$$

We have the well-known integral from Fourier series analysis

(6.90d) $$\frac{\pi}{2} = \int_0^\infty \frac{\sin t}{t} dt$$

and therefore defining

(6.90di) $$si(x) = Si(x) - \frac{\pi}{2}$$

we have

(6.90e) $$si(x) = \int_0^x \frac{\sin t}{t} dt - \int_0^\infty \frac{\sin t}{t} dt = -\int_x^\infty \frac{\sin t}{t} dt$$

With $x = 1/2$ in (6.89) we obtain

$$\sum_{n=1}^{\infty} \frac{\varsigma(2n)}{(2n+1)^2 2^{2n}} = \frac{1}{2}\left[1 - \frac{1}{\pi}\sum_{n=1}^{\infty}\frac{1}{n^2}\int_0^{1/2}\frac{\sin 2n\pi t}{t}dt\right]$$

$$= \frac{1}{2}\left[1 - \frac{1}{\pi}\sum_{n=1}^{\infty}\frac{1}{n^2}\int_0^{n\pi}\frac{\sin u}{u}du\right]$$

and hence we have

(6.90f) $$\sum_{n=1}^{\infty}\frac{\varsigma(2n)}{(2n+1)^2 2^{2n}} = \frac{1}{2}\left[1 - \frac{1}{\pi}\sum_{n=1}^{\infty}\frac{Si(n\pi)}{n^2}\right]$$

We have from (6.89)

(6.90fi) $$\sum_{n=1}^{\infty}\frac{\varsigma(2n)x^{2n+1}}{(2n+1)^2} = \frac{1}{2}x + \frac{1}{4\pi}\sum_{n=1}^{\infty}\frac{\sin 2n\pi x}{n^2} - \frac{1}{4\pi}\sum_{n=1}^{\infty}\frac{Si(2n\pi x)}{n^2}$$

We note that we could multiply (6.90fi) by $x^p$ and then integrate to obtain an expression for $\sum_{n=1}^{\infty}\frac{\varsigma(2n)x^{2n+2+p}}{(2n+2+p)(2n+1)^2}$.

From (6.90c) we have

$$\sum_{n=1}^{\infty}\int_0^{\pi}\log t\frac{\cos nt}{n}dt = -\sum_{n=1}^{\infty}\frac{Si(n\pi)}{n^2}$$

and using (7.8) we get

$$\sum_{n=1}^{\infty}\int_0^{\pi}\log t\frac{\cos nt}{n}dt = -\int_0^{\pi}\log t\log[2\sin(t/2)]dt$$

Thus we have

(6.90g) $$\sum_{n=1}^{\infty}\frac{Si(n\pi)}{n^2} = \int_0^{\pi}\log t\log[2\sin(t/2)]dt$$

We have

$$\int_0^{\pi}\log t\log[2\sin(t/2)]dt = \log 2\int_0^{\pi}\log t\, dt + \int_0^{\pi}\log t\log\sin(t/2)\, dt$$

and integration by parts gives us

$$\int_0^\pi \log t \log \sin(t/2)\, dt = (t\log t - t)\log \sin(t/2)\Big|_0^\pi - \frac{1}{2}\int_0^\pi (t\log t - t)\cot(t/2)\, dt$$

$$= -\frac{1}{2}\int_0^\pi (t\log t - t)\cot(t/2)\, dt$$

We now employ the basic identity (6.5a) to give us

(6.90h) $\quad \int_0^\pi (t\log t - t)\cot(t/2)\, dt = 2\sum_{n=1}^\infty \int_0^\pi (t\log t - t)\sin nt\, dt$

It is left as an exercise for the reader to show why (6.90h) satisfies the conditions for (6.5a).

We have

$$\int_0^\pi (t\log t - t)\sin nt\, dt = -(t\log t - t)\frac{\cos nt}{n}\Big|_0^\pi + \int_0^\pi \log t \frac{\cos nt}{n}\, dt$$

$$= \pi(\log \pi - 1)\frac{(-1)^{n+1}}{n} - \frac{Si(n\pi)}{n^2}$$

where we have used (6.90c). Completing the summation we obtain

$$\sum_{n=1}^\infty \int_0^\pi (t\log t - t)\sin nt\, dt = \pi(\log \pi - 1)\sum_{n=1}^\infty \frac{(-1)^{n+1}}{n} - \sum_{n=1}^\infty \frac{Si(n\pi)}{n^2}$$

and hence we get the same result as before

$$\sum_{n=1}^\infty \int_0^\pi \log t \frac{\cos nt}{n}\, dt = -\pi(\log \pi - 1)\log 2 + \pi(\log \pi - 1)\sum_{n=1}^\infty \frac{(-1)^{n+1}}{n} - \sum_{n=1}^\infty \frac{Si(n\pi)}{n^2} = -\sum_{n=1}^\infty \frac{Si(n\pi)}{n^2}$$

It has been shown in [45ae] that (an alternative proof is given later in this paper):

(6.91) $\quad \sum_{n=1}^\infty \frac{Si(n\pi)}{n^3} = \frac{5}{72}\pi^3$

and this is tantalisingly similar in structure to the series under consideration in (6.90f).

We have from Bromwich's book [36b, p.312] and also [1, p.232]

$$si(x) = -\int_0^{\pi/2} e^{-x\cos t} \cos(x\sin t)\,dt$$

$$si(n\pi) = -\int_0^{\pi/2} e^{-n\pi\cos t} \cos(n\pi \sin t)\,dt = -\frac{1}{2}\int_0^{\pi/2} e^{-n\pi\cos t}\left[e^{in\pi\sin t} + e^{-in\pi\sin t}\right]dt$$

We then see that

$$\sum_{n=1}^{\infty} \frac{si(n\pi)}{n^p} = -\frac{1}{2}\sum_{n=1}^{\infty} \frac{1}{n^p}\int_0^{\pi/2} \left(e^{n\pi[-\cos t+i\sin t]} + e^{n\pi[-\cos t-i\sin t]}\right)dt$$

$$= -\frac{1}{2}\sum_{n=1}^{\infty} \frac{1}{n^p}\int_0^{\pi/2} \left(e^{n\pi[-\cos t+i\sin t]} + e^{n\pi[-\cos t-i\sin t]}\right)dt$$

and hence we get

$$\sum_{n=1}^{\infty} \frac{si(n\pi)}{n^p} = -\frac{1}{2}\int_0^{\pi/2}\left(Li_p\left[e^{\pi[-\cos t+i\sin t]}\right] + Li_p\left[e^{\pi[-\cos t-i\sin t]}\right]\right)dt$$

A formula corresponding to (6.90f) does not appear in the mammoth list compiled by Srivastava and Choi in [126]: they do report the result [126, p.223]

(6.92) $$\sum_{n=1}^{\infty} \frac{\varsigma(2n)}{n(2n+1)2^{2n}} = -1 + \log \pi$$

(a proof of this identity is given in (6.107g)). Later in this paper in (6.123) we show that a related series is given by

(6.92a)
$$\int_0^1 \log\Gamma(x+1)\log[2\sin(\pi x)]\,dx = \frac{1}{2\pi}\sum_{n=1}^{\infty} \frac{si(2n\pi)}{n^2} = \frac{1}{2\pi}\sum_{n=1}^{\infty} \frac{Si(2n\pi)}{n^2} - \frac{1}{4}\varsigma(2)$$

and this gives some indication of the complexity involved in this area.

If (6.89) is valid for $x=1$, we may have

(6.92b) $$\sum_{n=1}^{\infty} \frac{\varsigma(2n)}{(2n+1)^2} = \frac{1}{2} - \frac{1}{4\pi}\sum_{n=1}^{\infty} \frac{Si(2n\pi)}{n^2} \quad (?)$$

and, with our unproven assumption in mind, we may then have

(6.92c) $$\int_0^1 \log \Gamma(x+1) \log[2\sin(\pi x)] dx = 1 - \frac{1}{4}\varsigma(2) - 2\sum_{n=1}^{\infty} \frac{\varsigma(2n)}{(2n+1)^2} \quad (?)$$

Reference to (6.48) shows that

$$\pi \int_0^1 \frac{dt}{t} \int_0^t x \cot \pi x \, dx = -2\sum_{n=1}^{\infty} \frac{\varsigma(2n)}{(2n+1)^2}$$

See also (6.117j) where, some three years later, I showed that

$$\frac{1}{2\pi^2} \sum_{n=1}^{\infty} \frac{Si(2n\pi)}{n^2} = \log A - \frac{1}{4}$$

A series similar to (6.92) was posed as a question in The American Mathematical Monthly in 1965 and, Danese, one of the solvers [75a], noted that a more general identity, known as Burnside's formula, is reported in Higher Transcendental Functions by Erdélyi et al.

(6.93) $$\sum_{n=1}^{\infty} \frac{\varsigma(2n,z)}{n(2n+1)2^{2n}} = 2\left(z-\frac{1}{2}\right)\log\left(z-\frac{1}{2}\right) - 2\left(z-\frac{1}{2}\right) + \log 2\pi - 2\log \Gamma(z)$$

where $\varsigma(s,z)$ is the Hurwitz zeta function (the relationship becomes clear by letting $z = 1$ in (6.93)). The popularity of this series continued through into 1987 where it reappeared as a problem in the same journal [130a], and was solved in an entirely different manner.

Andersen [45a] proved (6.91) by using Parseval's generalised identity

(6.94) $$\frac{1}{\pi} \int_{-\pi}^{\pi} f(x)g(x) dx = \frac{1}{2}a_0 \alpha_0 + \sum_{n=1}^{\infty}(a_n \alpha_n + b_n \beta_n)$$

I first came across this identity in 1999 when I purchased a 1930 edition of Carslaw's book [41] on Fourier series in a car-boot sale for the princely sum of 10 pence. Glancing through this book, almost immediately (on 2 December 1999, to be exact) I realised for the first time that one could potentially evaluate $\varsigma(3)$ by identifying suitable functions $f(x)$ and $g(x)$, with Fourier coefficients $(a_n, b_n)$ and $(\alpha_n, \beta_n)$ respectively, which satisfied the following conditions

$$a_n = \frac{1}{n} \qquad b_n = 0$$

$$\alpha_n = \frac{1}{n^2} \qquad \beta_n = 0$$

That was the very start of my zeta quest! Of course I never did find any such functions; somewhat later, I realised that a simple modification would suffice, viz

$$a_n = \frac{(-1)^n}{n} \qquad b_n = 0$$

$$\alpha_n = \frac{(-1)^n}{n^2} \qquad \beta_n = 0$$

(at that stage of my knowledge, I didn't even know that the alternating zeta function was an appropriate object of study). Some further perusal of Carslaw's inspiring book then revealed that suitable functions were $f(x) = x$ and $g(x) = \log \cos x$ over an appropriate interval: all that remained was to evaluate an integral of the form

$\int x \log \cos x \, dx$. Little did I know way back in 1999 how difficult was the task that then confronted me, and indeed had confronted countless mathematicians for the past three centuries!

In this connection, it may be worthwhile noting from (E.46) of Volume V that

$$\int_0^\pi \log \Gamma(x/\pi) \sin 2nx \, dx = \frac{1}{2n}(\gamma + \log 2n\pi)$$

Let us now move on: from (6.90c) we might deduce that

$$\sum_{n=1}^\infty \int_0^\pi \log t . (-1)^n \cos 2nt \, dt = -\sum_{n=1}^\infty (-1)^n \frac{Si(2n\pi)}{2n}$$

and reference to (6.8) would imply that

$$\sum_{n=1}^\infty \int_0^\pi \log t . (-1)^n \cos 2nt \, dt = -\int_0^\pi \log t \, dt = \pi(1 - \log \pi)$$

Therefore we might have

$$\sum_{n=1}^\infty (-1)^n \frac{Si(2n\pi)}{n} = 2\pi(\log \pi - 1) \quad (?)$$

However we would be wrong because the necessary conditions for the Riemann-Lebesgue lemma are not satisfied by $\log t$ at $t = 0$. Nielsen instead gave the following expression in his 1906 book, "Theorie des Integrallogarithmus und verwanter tranzendenten" [104a, p.79]

(6.94a) $\qquad \sum_{n=1}^\infty (-1)^n \frac{si(2n\pi)}{n\pi} = \log \sqrt{2} - 1$

where $si(x) = Si(x) - \frac{\pi}{2}$. This may be written as

$$\sum_{n=1}^{\infty}(-1)^n \frac{Si(2n\pi)}{n\pi} - \frac{1}{2}\sum_{n=1}^{\infty}\frac{(-1)^n}{n} = \log\sqrt{2} - 1$$

and hence

(6.94ai) $$\sum_{n=1}^{\infty}(-1)^n \frac{Si(2n\pi)}{n} = -\pi$$

Alternatively, using (6.90ai) in the interval $[a,b]$ we get for $0 < a < b < \pi/2$

$$\int_a^b \log t \cdot \cos 2nt \, dt = \frac{\sin 2nb \log b}{2n} - \frac{\sin 2na \log a}{2n} - \frac{Si(2nb)}{2n} + \frac{Si(2na)}{2n}$$

The necessary conditions for the Riemann-Lebesgue lemma to apply to (6.8) are now satisfied in that interval and in conjunction with (6.94a) we may obtain an expression for $\sum_{n=1}^{\infty}(-1)^n \frac{Si(2n\alpha)}{n}$.

From (6.8) we have

$$-\frac{1}{2}\int_a^b \log t \, dt = \sum_{n=1}^{\infty}\int_a^b \log t (-1)^n \cos 2nt \, dt$$

and hence we get

$$-\frac{1}{2}[b\log b - b - a\log a + a] =$$

$$\frac{1}{2}\log b \sum_{n=1}^{\infty}(-1)^n \frac{\sin 2nb}{n} - \frac{1}{2}\log a \sum_{n=1}^{\infty}(-1)^n \frac{\sin 2na}{n} - \frac{1}{2}\sum_{n=1}^{\infty}(-1)^n \frac{Si(2nb)}{n} + \frac{1}{2}\sum_{n=1}^{\infty}(-1)^n \frac{Si(2na)}{n}$$

We have from Tolstov [130, p.148] for $-\pi/2 < \alpha < \pi/2$

$$\sum_{n=1}^{\infty}(-1)^n \frac{\sin 2n\alpha}{n} = -\alpha$$

and using (6.94ai) we may write the above equation as

$$\frac{1}{2}b - \frac{1}{2}a = -\frac{1}{2}\sum_{n=1}^{\infty}(-1)^n \frac{Si(2nb)}{n} + \frac{1}{2}\sum_{n=1}^{\infty}(-1)^n \frac{Si(2na)}{n}$$

Hence by symmetry we have

$$b = -\sum_{n=1}^{\infty}(-1)^n \frac{Si(2nb)}{n} + c$$

and since $Si(0) = 0$ we have $c = 0$ which then gives us

(6.94a) $$\frac{x}{2} = \frac{\pi}{2}\log 2 - \sum_{n=1}^{\infty}(-1)^n \frac{si(nx)}{n}$$

as noted in [104a, p.83] for $|x| < \pi$.

In the same book (which may be viewed on the internet), Nielsen [104a, p.78] states that for $0 < x < 1$

(6.94b) $$\log x + \log(1-x) - \log(2\sin \pi x) + 2 = -\frac{2}{\pi}\sum_{n=1}^{\infty}\frac{si(2n\pi)\cos 2n\pi x}{n}$$

(equation (6.94b) corrects the misprint in Nielsen's book).

It is easily seen that (6.94a) is derived from (6.94b) by letting $x = 1/2$.

The proof of (6.94b), which is quite straightforward, is given below. We have the Fourier series for $\log x$

(6.94c) $$\log x = \frac{1}{2}a_0 + \sum_{n=1}^{\infty}(a_n \cos nx + b_n \sin nx)$$

where, as proved below, we have

(6.94d) $$a_n = \frac{1}{\pi}\int_0^{2\pi} \log x . \cos nx\, dx = -\frac{1}{n\pi}\left[si(2n\pi) + \frac{\pi}{2}\right]$$

(6.94e) $$b_n = \frac{1}{\pi}\int_0^{2\pi} \log x . \sin nx\, dx = \frac{1}{n\pi}\left[Ci(2n\pi) - \gamma - \log(2n\pi)\right]$$

With integration by parts we get

$$\int_\alpha^x \log t . \sin at\, dt = -\log t \frac{\cos at}{a}\bigg|_\alpha^x + \int_\alpha^x \frac{\cos at}{at}\, dt$$

$$= -\log t \frac{\cos at}{a}\bigg|_\alpha^x + \int_\alpha^x \frac{\cos at - 1}{at} dt + \int_\alpha^x \frac{1}{at} dt$$

$$= -\log x \frac{\cos ax}{a} + \int_\alpha^{ax} \frac{\cos u - 1}{u} du + \frac{\log \alpha}{a}(\cos a\alpha - 1) + \frac{\log x}{a}$$

Therefore, in the limit as $\alpha \to 0$, we get

$$\int_0^x \log t \cdot \sin at \, dt = -\log x \frac{\cos ax}{a} + \int_0^{ax} \frac{\cos u - 1}{u} du + \frac{\log x}{a}$$

and, as reported in Nielsen's book [104a, p.12], we have using (6.94g)

(6.94f) $\qquad \int_0^x \log t \cdot \sin at \, dt = \frac{1}{a}\left[Ci(ax) - \gamma - \cos(ax)\log x - \log a\right]$

where the cosine integral $Ci(x)$ is defined in G&R [74, p.878] and [1, p.231] as

(6.94g) $\qquad Ci(x) = \gamma + \log x + \int_0^x \frac{\cos t - 1}{t} dt = \gamma + \log x + \sum_{n=1}^\infty \frac{(-1)^n x^{2n}}{2n(2n)!}$

(where, in the final part, we have simply substituted the Maclaurin series for the integrand). We also have

(6.94ga) $\qquad Ci(x) = -\int_x^\infty \frac{\cos t}{t} dt$

and this more clearly shows the connection with $si(x)$ in (6.90e). $Ci(x)$ is frequently designated as $ci(x)$ in other works such as G&R [74]. We also have the following representations recently given by Harris [76c] in terms of spherical Bessel functions

$$Si(x) = x \sum_{n=0}^\infty \left[j_n\left(\frac{x}{2}\right)\right]^2$$

$$Ci(x) = \gamma + \log x + \sum_{n=0}^\infty a_n \left[j_n\left(\frac{x}{2}\right)\right]^2$$

where $a_n = -(2n+1)\left(1 - (-1)^n + H_n\right)$ and $a_0 = 0$.

From Nielsen's book [104a, p.6] we have

$$Ci(x) = \frac{\mathrm{li}(e^{ix}) + \mathrm{li}(e^{-ix})}{2}$$

$$si(x) = \frac{\mathrm{li}(e^{ix}) - \mathrm{li}(e^{-ix})}{2i}$$

in terms of the logarithmic integral $\mathrm{li}(x) = \int_0^x \frac{du}{\log u}$.

After that digression we therefore obtain

(6.94h) $\quad b_n = \dfrac{1}{\pi}\displaystyle\int_0^{2\pi} \log x \cdot \sin nx\, dx = \dfrac{1}{n\pi}\big[Ci(2n\pi) - \gamma - \log(2n\pi)\big]$

The coefficients $a_n$ can be obtained in the same manner as (6.90c). We also have

(6.94i) $\quad a_0 = \dfrac{1}{\pi}\displaystyle\int_0^{2\pi} \log x\, dx = \log(2\pi) - 1$

and making use of (7.5) and (7.8)

$$\frac{1}{2}(\pi - x) = \sum_{n=1}^{\infty} \frac{\sin nx}{n} \quad (0 < x < 2\pi)$$

$$\log\big[2\sin(x/2)\big] = -\sum_{n=1}^{\infty} \frac{\cos nx}{n} \quad (0 < x < 2\pi)$$

we therefore obtain the Fourier series (where we have made the substitution $x \to 2\pi x$)

(6.94j)
$$1 + \log x - \frac{1}{2}\log(2\sin \pi x) + (\gamma + \log 2\pi)\left(\frac{1}{2} - x\right) =$$
$$\sum_{n=1}^{\infty} \frac{[Ci(2n\pi) - \log n]\sin 2n\pi x - si(2n\pi)\cos 2n\pi x}{n\pi}$$

In (6.94j) we let $x \to 1 - x$ and combine the resulting identity with (6.94j) to obtain (6.94b)

(6.94ji) $\quad \log x + \log(1 - x) - \log(2\sin \pi x) + 2 = -\dfrac{2}{\pi}\displaystyle\sum_{n=1}^{\infty} \dfrac{si(2n\pi)\cos 2n\pi x}{n}$

and subtraction results in

(6.94jii)  $$\log x - \log(1-x) + (\gamma + \log 2\pi)(1-2x) = 2\sum_{n=1}^{\infty} \frac{[Ci(2n\pi) - \log n]\sin 2n\pi x}{n\pi}$$

If we multiply (6.94ji) by $x$ and integrate we obtain

$$\int_0^{1/2} x[\log x + \log(1-x) - \log(2\sin \pi x) + 2]dx = -\frac{2}{\pi}\sum_{n=1}^{\infty} \frac{si(2n\pi)}{n}\int_0^{1/2} x\cos 2n\pi x\, dx$$

Therefore we get

$$\left\{-\frac{1}{4}x^2 - \frac{1}{2}x^2 \log x\right\} + \left\{-\frac{1}{2}x - \frac{1}{4}x^2 + \frac{1}{2}x^2 \log(1-x) - \frac{1}{2}\log(1-x)\right\} + x^2 \Bigg|_0^{1/2}$$

$$-\frac{1}{8}\log 2 - \int_0^{1/2} x\log\sin\pi x\, dx = -\frac{2}{\pi}\sum_{n=1}^{\infty} \frac{si(2n\pi)}{n}\left[\frac{(-1)^n - 1}{4\pi^2 n^2}\right]$$

Hence, using our old favourite equation (1.11) from Volume I, we have

$$-\frac{1}{8} + \frac{1}{8}\log 2 - \frac{1}{\pi^2}\left[\frac{7}{16}\varsigma(3) - \frac{\pi^2}{8}\log 2\right] = -\frac{1}{2\pi^3}\sum_{n=1}^{\infty} \frac{si(2n\pi)[(-1)^n - 1]}{n^3}$$

and this becomes

(6.94jiv)  $$\frac{1}{4}\pi^3 + \frac{7}{8}\pi\varsigma(3) = \sum_{n=1}^{\infty} \frac{si(2n\pi)[(-1)^n - 1]}{n^3}$$

Since $si(x) = Si(x) - \frac{\pi}{2}$ we have

$$\frac{1}{4}\pi^3 + \frac{7}{8}\pi\varsigma(3) = \sum_{n=1}^{\infty} \frac{Si(2n\pi)[(-1)^n - 1]}{n^3} - \frac{\pi}{2}\sum_{n=1}^{\infty} \frac{[(-1)^n - 1]}{n^3}$$

$$= -2\sum_{n=0}^{\infty} \frac{Si[2(2n+1)\pi]}{(2n+1)^3} + \pi\sum_{n=0}^{\infty} \frac{1}{(2n+1)^3}$$

$$= -2\sum_{n=0}^{\infty} \frac{Si[2(2n+1)\pi]}{(2n+1)^3} + \frac{7}{8}\pi\varsigma(3)$$

Therefore we have

(6.94k)  $$\sum_{n=0}^{\infty} \frac{Si[2(2n+1)\pi]}{(2n+1)^3} = -\frac{1}{8}\pi^3$$

This agrees with [45ae] where, using Parseval's formula, it was shown that

(6.94ki) $$\sum_{n=1}^{\infty}\frac{Si(n\pi)}{n^3} = \left(\frac{1}{8}-\frac{1}{18}\right)\pi^3 \text{ and } \sum_{n=1}^{\infty}(-1)^n\frac{Si(n\pi)}{n^3} = -\frac{1}{18}\pi^3$$

and hence combining the above two series we obtain (6.94k).

It is a pity, but not altogether surprising, that the interesting factor involving $\varsigma(3)$ cancels out! Other interesting identities may be obtained by multiplying (6.94b) by different powers of $x$, or perhaps by a different function altogether.

We have using integration by parts

(6.94l) $$\int (x^3 \log x)\cos ax\, dx = \frac{(a^2x^2-11)\cos ax}{a^4} - \frac{5x\sin ax}{a^3} + 6\frac{Ci(ax)}{a^4}$$

$$+ \frac{-6\cos ax + 3a^2x^2\cos ax - 6ax\sin ax + a^3x^3\sin ax}{a^4}\log x$$

Therefore we have

$$\int_0^1 (x^3 \log x)\cos n\pi x\, dx = 11\left[\frac{1}{n^4\pi^4} - \frac{(-1)^n}{n^4\pi^4}\right] + \frac{(-1)^n}{n^2\pi^2} + 6\frac{Ci(n\pi)}{n^4\pi^4} - \lim_{\varepsilon\to 0}\frac{6}{n^4\pi^4}\left[Ci(n\pi\varepsilon) - \cos(n\pi\varepsilon)\log\varepsilon\right]$$

Using (6.94f) it is clear that

$$\lim_{x\to 0}\int_0^x \log t.\sin at\, dt = \frac{1}{a}\lim_{x\to 0}\left[Ci(ax)-\gamma-\cos(ax)\log x - \log a\right] = 0$$

This in turn implies that

(6.94m) $$\lim_{x\to 0}\left[Ci(ax) - \cos(ax)\log x\right] = \gamma + \log a$$

This limit may also be derived directly from the definition of $Ci(x)$ in (6.94g) as follows

$$Ci(ax) - \cos(ax)\log x = \gamma + \log ax - \cos(ax)\log x + \int_0^{ax}\frac{\cos t - 1}{t}dt$$

$$= \gamma + \log a + \log x[1-\cos(ax)] + \int_0^{ax}\frac{\cos t - 1}{t}dt$$

and then take the limit as $x\to 0$.

Therefore we have

(6.94n)
$$\int_0^1 (x^3 \log x) \cos n\pi x \, dx = 11 \left[ \frac{1}{n^4 \pi^4} - \frac{(-1)^n}{n^4 \pi^4} \right] + \frac{(-1)^n}{n^2 \pi^2} + 6 \frac{Ci(n\pi)}{n^4 \pi^4} - \frac{6}{n^4 \pi^4} [\gamma + \log n\pi]$$

Since $p(x) = x^m \log x$ is twice continuously differentiable on $[0,1]$ for $m \geq 3$, it is clear that $p(x) = x^3 \log x$ meets the conditions of the Riemann-Lebesgue lemma on that interval (also noting that $p(0) = 0$ as required). We will show in (6.94q) that (6.5) and (6.5a) may in fact be employed in the case where $m = 2$.

We may therefore employ equation (6.5) which is noted below for ease of reference

$$\sum_{n=1}^{\infty} \int_a^b p(x) \cos \alpha nx \, dx = -\frac{1}{2} \int_a^b p(x) dx$$

Hence we get

$$\sum_{n=1}^{\infty} \int_0^1 (x^3 \log x) \cos n\pi x \, dx = -\frac{1}{2} \int_0^1 x^3 \log x \, dx = \frac{x^4}{32} - \frac{1}{8} x^3 \log x \Big|_0^1 = \frac{1}{32}$$

Alternatively, completing the summation of (6.94n) we have

$$\sum_{n=1}^{\infty} \int_0^1 (x^3 \log x) \cos n\pi x \, dx =$$

$$\frac{11}{\pi^4} [\varsigma(4) + \varsigma_a(4)] - \frac{1}{\pi^2} \varsigma_a(2) + \frac{6}{\pi^4} \sum_{n=1}^{\infty} \frac{Ci(n\pi)}{n^4} - \frac{6}{\pi^4} \varsigma(4) [\gamma + \log \pi] + \frac{6}{\pi^4} \varsigma'(4)$$

In conclusion, we have

$$\frac{11}{\pi^4} [\varsigma(4) + \varsigma_a(4)] - \frac{1}{\pi^2} \varsigma_a(2) + \frac{6}{\pi^4} \sum_{n=1}^{\infty} \frac{Ci(n\pi)}{n^4} - \frac{6}{\pi^4} \varsigma(4) [\gamma + \log \pi] + \frac{6}{\pi^4} \varsigma'(4) = \frac{1}{32}$$

which implies that

(6.94o) $\quad \frac{11}{576} \pi^4 + \sum_{n=1}^{\infty} \frac{Ci(n\pi)}{n^4} = \varsigma(4)[\gamma + \log \pi] - \varsigma'(4)$

Using the definition of $Ci(x)$ in (6.94g) this becomes

$$\frac{11}{576} \pi^4 + \sum_{n=1}^{\infty} \frac{1}{n^4} \int_0^{n\pi} \frac{\cos t - 1}{t} dt = 0$$

or equivalently

$$\frac{11}{576}\pi^4 + \sum_{n=1}^{\infty}\frac{1}{n^4}\int_0^{\pi}\frac{\cos nu - 1}{u}du = 0$$

With integration by parts we see that

$$\int_0^{\pi}\frac{\cos nu - 1}{u}du = \left[(-1)^n - 1\right]\log\pi + n\int_0^{\pi}\sin nu \log u\, du$$

and we therefore obtain

$$\frac{11}{576}\pi^4 - \varsigma_a(4) - \varsigma(4) + \sum_{n=1}^{\infty}\frac{1}{n^3}\int_0^{\pi}\sin nu \log u\, du = 0$$

Assuming that we may interchange the order of summation and integration we have

$$\sum_{n=1}^{\infty}\frac{1}{n^3}\int_0^{\pi}\sin nu \log u\, du = \int_0^{\pi}\sum_{n=1}^{\infty}\frac{\sin nu}{n^3}\log u$$

We have the Fourier series from [130, p.148]

$$\sum_{n=1}^{\infty}\frac{\sin nu}{n^3} = \frac{1}{12}(u^3 - 3\pi u^2 + 2\pi^2 u)$$

We have seen before that

$$\int u^n \log u\, du = \frac{u^{n+1}}{(n+1)^2}\left[(n+1)\log u - 1\right]$$

and therefore

$$\frac{1}{12}\int_0^{\pi}(u^3 - 3\pi u^2 + 2\pi^2 u)\log u\, du = \frac{1}{48}\pi^4 \log\pi - \frac{11}{576}\pi^4$$

This then gives us another derivation of (6.94o).

In a similar way we have

$$\int (x^2 \log x)\sin ax\, dx = \frac{3\cos ax}{a^3} + \frac{x\sin ax}{a^2} + \frac{2\cos ax + 2ax\sin ax - a^2 x^2 \cos ax}{a^3}\log x - \frac{2Ci(ax)}{a^3}$$

Therefore we get using (6.94m)

$$\int_0^1 (x^2 \log x)\sin n\pi x\, dx = 3\frac{(-1)^n}{n^3\pi^3} - 3\frac{1}{n^3\pi^3} - \frac{2Ci(n\pi)}{n^3\pi^3} + \lim_{\varepsilon\to 0}\frac{2}{n^3\pi^3}\left[Ci(n\pi\varepsilon) - \cos(n\pi\varepsilon)\log\varepsilon\right]$$

$$= 3\frac{(-1)^n}{n^3\pi^3} - 3\frac{1}{n^3\pi^3} - \frac{2Ci(n\pi)}{n^3\pi^3} + \frac{2}{n^3\pi^3}[\gamma + \log n\pi]$$

Completing the summation of the above integral we get

$$\sum_{n=1}^{\infty}\int_0^1 (x^2 \log x)\frac{\sin n\pi x}{n}dx =$$

$$-\frac{3}{\pi^3}[\varsigma(4) + \varsigma_a(4)] - \frac{2}{\pi^3}\sum_{n=1}^{\infty}\frac{Ci(n\pi)}{n^4} + \frac{2}{\pi^3}\varsigma(4)[\gamma + \log \pi] - \frac{2}{\pi^3}\varsigma'(4)$$

From (7.5) we have for $x \in (0,2)$

$$\frac{\pi}{2}(1-x) = \sum_{n=1}^{\infty}\frac{\sin n\pi x}{n}$$

and therefore interchanging the order of summation and integration we get (assuming that this operation is valid)

$$\sum_{n=1}^{\infty}\int_0^1 (x^2 \log x)\frac{\sin n\pi x}{n}dx = \frac{\pi}{2}\int_0^1 (x^2 \log x)(1-x)dx$$

(since $x^2 \log x$ is zero at $x = 0$, the Fourier series (7.5) may also be employed at $x = 0$).

Using (3.207) we have

$$\int_0^1 x^2 \log x\, dx = -\frac{1}{9} \qquad \int_0^1 x^3 \log x\, dx = -\frac{1}{16}$$

and therefore

$$\sum_{n=1}^{\infty}\int_0^1 (x^2 \log x)\frac{\sin n\pi x}{n}dx = -\frac{7\pi}{288}$$

$$= -\frac{3}{\pi^3}[\varsigma(4) + \varsigma_a(4)] - \frac{2}{\pi^3}\sum_{n=1}^{\infty}\frac{Ci(n\pi)}{n^4} + \frac{2}{\pi^3}\varsigma(4)[\gamma + \log \pi] - \frac{2}{\pi^3}\varsigma'(4)$$

Accordingly, we have another proof of the equation (6.94o) involving $\sum_{n=1}^{\infty}\frac{Ci(n\pi)}{n^4}$ and $\varsigma'(4)$.

As a further example we have

$$\int (x^3 \log x) \sin ax \, dx = \frac{5x \cos ax}{a^3} + \frac{-11 + a^2 x^2 \sin ax}{a^4}$$

$$- \frac{-6ax \cos ax + a^3 x^3 \cos ax + 6 \sin ax - 3a^2 x^2 \sin ax}{a^4} \log x + \frac{6 Si(ax)}{a^4}$$

Therefore we obtain

(6.94p) $$\int_0^1 (x^3 \log x) \sin n\pi x \, dx = 5 \frac{(-1)^n}{n^3 \pi^3} + 6 \frac{Si(n\pi)}{n^4 \pi^4}$$

Using (6.5a) we get

$$\int_0^1 x^3 \log x \cot(\pi x/2) \, dx = \sum_{n=1}^{\infty} \int_0^1 (x^3 \log x) \sin n\pi x \, dx$$

and hence we have

(6.94q) $$\int_0^1 x^3 \log x \cot(\pi x/2) \, dx = -\frac{5}{\pi^3} \varsigma_a(3) + \frac{6}{\pi^4} \sum_{n=1}^{\infty} \frac{Si(n\pi)}{n^4}$$

Since from (6.90c) $Si(n\pi) = -n \int_0^{\pi} \log t . \cos nt \, dt$ we have

$$\sum_{n=1}^{\infty} \frac{Si(n\pi)}{n^4} = \int_0^{\pi} \log t . \sum_{n=1}^{\infty} \frac{\cos nt}{n^3} dt$$

In [130, p.148] we find

$$\sum_{n=1}^{\infty} \frac{\cos nt}{n^3} = \int_0^t du \int_0^u \log[2 \sin(v/2)] dv$$

but further simplification does not seem possible.

Clearly, the method may be extended for other powers of $x$.

Let us now reconsider the remainder term in (6.3a) with $p(x) = x^2 \log x$

$$R_N = \frac{1}{2} \int_a^b p(x) \left( \cos(N+1)x + i \sin(N+1)x \right) \left( 1 + i \cot(x/2) \right) dx$$

We have with $[a, b] = [0, t]$

$$2\operatorname{Im}(R_N) = \int_0^t x^2 \log x \{\sin(N+1)x + \cot(x/2)\cos(N+1)x\} dx$$

and integration by parts gives us

$$\int_0^t x^2 \log x \sin(N+1)x\, dx = -\frac{\cos(N+1)x}{N+1} x^2 \log x \Big|_0^t + \int_0^t \frac{\cos(N+1)x}{N+1}(2x\log x + x)dx$$

It is clear that the above terms remain finite for a fixed value of $N$ and that the value of each of the integrated term and the integral approaches zero as $N \to \infty$.

Another integration by parts gives us

$$\int_0^t x^2 \log x \cot(x/2)\cos(N+1)x\, dx =$$

$$\frac{\sin(N+1)x}{N+1} x^2 \log x \cot(x/2)\Big|_0^t - \int_0^t \frac{\sin(N+1)x}{N+1}\left[\cot(x/2)(x + 2x\log x) - \frac{x^2 \log x}{2\sin^2(x/2)}\right]dx$$

With regard to the integrated part we see that

$$I_N = \frac{\sin(N+1)x}{N+1} x^2 \log x \cot(x/2) = \frac{\sin(N+1)x}{(N+1)} \frac{x}{\sin(x/2)} \cos(x/2) x \log x$$

and hence it is easily seen that $I_N$ is finite at $x = 0$ and that $\lim_{x \to 0} I_N = 0$. Furthermore, with regard to the above integral, we have

$$\frac{\sin(N+1)x}{N+1}\left[\cot(x/2)(x + 2x\log x) - \frac{x^2 \log x}{2\sin^2(x/2)}\right] =$$

$$\frac{\sin(N+1)x}{(N+1)} \frac{x}{\sin(x/2)}\cos(x/2) + \frac{\sin(N+1)x}{(N+1)\sin(x/2)} x \log x \left[2\cos(x/2) - \frac{(x/2)}{\sin(x/2)}\right]$$

and it is clear that this expression remains finite as $x \to 0$.

Therefore we have $\lim_{N \to \infty}[\operatorname{Im}(R_N)] = 0$ and similarly we can show that $\lim_{N \to \infty}[\operatorname{Re}(R_N)] = 0$. Accordingly, we have shown that (6.5) and (6.5a) are valid for $p(x) = x^2 \log x$ notwithstanding that $x^2 \log x$ is not twice continuously differentiable on $[0,t]$ (and therefore does not meet the specific conditions for our proof of the Riemann-Lebesgue lemma). We therefore have

$$\frac{1}{2}\int_0^t x^2 \log x\, dx = \sum_{n=0}^{\infty} \int_0^t x^2 \log x \cos \alpha nx\, dx$$

$$\frac{1}{2}\int_0^t x^2 \log x \cot(\alpha x/2) x\, dx = \sum_{n=1}^{\infty} \int_0^t x^2 \log x \sin \alpha nx\, dx$$

As shown above we have

$$\int_0^1 (x^2 \log x) \sin n\pi x\, dx = 3\frac{(-1)^n}{n^3\pi^3} - 3\frac{1}{n^3\pi^3} - \frac{2Ci(n\pi)}{n^3\pi^3} + \frac{2}{n^3\pi^3}[\gamma + \log n\pi]$$

and we then obtain

(6.94r) $$\frac{1}{2}\int_0^1 (x^2 \log x)\cot(\pi x/2)dx =$$

$$-\frac{3}{\pi^3}[\varsigma_a(3) + \varsigma(3)] - \frac{2}{\pi^3}\sum_{n=1}^{\infty}\frac{Ci(n\pi)}{n^3} + \frac{2}{\pi^3}\varsigma(3)[\gamma + \log \pi] - \frac{2}{\pi^3}\varsigma'(3)$$

Mathematica was not able to evaluate the above integral.

An alternative proof of (6.91) is given below. Using integration by parts gives us

$$\int (x^2 \log x)\cos ax\, dx =$$

$$\frac{x\cos ax}{a^2} - \frac{3\sin ax}{a^3} + \frac{2ax\cos ax - 2\sin ax + a^2 x^2 \sin ax}{a^3}\log x + \frac{2Si(ax)}{a^3}$$

and therefore

(6.94s) $$\int_0^1 (x^2 \log x)\cos n\pi x\, dx = \frac{(-1)^n}{n^2\pi^2} + \frac{2Si(n\pi)}{n^3\pi^3}$$

Therefore, using the basic identity (6.5) with $\alpha = \pi$ we obtain with the assistance of (3.207)

$$\sum_{n=0}^{\infty}\int_0^1 (x^2 \log x)\cos n\pi x\, dx = \frac{1}{2}\int_0^1 x^2 \log x\, dx = -\frac{x^3}{18} + \frac{1}{6}x^3 \log x\Big|_0^1 = -\frac{1}{18}$$

Accordingly we have

$$\sum_{n=1}^{\infty}\left[\frac{(-1)^n}{n^2\pi^2} + \frac{2Si(n\pi)}{n^3\pi^3}\right] - \frac{1}{18} = -\frac{1}{18}$$

and, since $\sum_{n=1}^{\infty} \frac{(-1)^n}{n^2 \pi^2} = -\frac{1}{12}$, this easily simplifies to prove (6.91).

$$\sum_{n=1}^{\infty} \frac{Si(n\pi)}{n^3} = \frac{5}{72}\pi^3$$

Completing the summation of (6.94s) we get

$$\sum_{n=1}^{\infty} \frac{1}{n} \int_0^1 (x^2 \log x) \cos n\pi x \, dx = -\frac{1}{12} + \frac{2}{\pi^3} \sum_{n=1}^{\infty} \frac{Si(n\pi)}{n^4}$$

Using the Fourier series [130, p.148] we have

$$\sum_{n=1}^{\infty} \frac{1}{n} \int_0^1 (x^2 \log x) \cos n\pi x \, dx = -\int_0^1 (x^2 \log x) \log[2\sin(\pi x/2)] \, dx$$

and therefore we obtain

(6.94t) $$\int_0^1 (x^2 \log x) \log[2\sin(\pi x/2)] \, dx = \frac{1}{12} - \frac{2}{\pi^3} \sum_{n=1}^{\infty} \frac{Si(n\pi)}{n^4}$$

Similarly, we have

$$\sum_{n=1}^{\infty} \frac{1}{n^2} \int_0^1 (x^2 \log x) \cos n\pi x \, dx = -\frac{1}{\pi^2} \varsigma_a(4) + \frac{2}{\pi^3} \sum_{n=1}^{\infty} \frac{Si(n\pi)}{n^5}$$

and [130, p.148] gives us

$$\sum_{n=1}^{\infty} \frac{1}{n^2} \int_0^1 x^2 \log x \cos n\pi x \, dx = \frac{\pi^2}{12} \int_0^1 x^2 \log x \left[3x^2 - 6x + 2\right] dx$$

Thus we may easily derive a closed ended formula for $\sum_{n=1}^{\infty} \frac{Si(n\pi)}{n^5}$.

We also have for example

$$\int (x \log x) \cos ax \, dx = \frac{\cos ax}{a^2} + \frac{\cos ax + ax \sin ax}{a^2} \log x - \frac{Ci(ax)}{a^2}$$

$$\int (x \log x) \sin ax \, dx = \frac{\sin ax}{a^2} - \frac{ax \cos ax - \sin ax}{a^2} \log x - \frac{Si(ax)}{a^2}$$

In particular we have

$$\int_0^1 (x\log x)\cos n\pi x\,dx = \frac{(-1)^n}{\pi^2 n^2} - \frac{1}{\pi^2 n^2} - \frac{1}{\pi^2}\frac{Ci(n\pi)}{n^2} + \frac{1}{\pi^2 n^2}(\gamma + \log n\pi)$$

Therefore, completing the summation gives us

$$\int_0^1 (x\log x)\sum_{n=1}^{\infty}\frac{\cos n\pi x}{n}\,dx = -\frac{1}{\pi^2}\varsigma_a(2) - \frac{1}{\pi^2}\varsigma(2) - \frac{1}{\pi^2}\sum_{n=1}^{\infty}\frac{Ci(n\pi)}{n^2} + \frac{1}{\pi^2}\sum_{n=1}^{\infty}\frac{(\gamma + \log n\pi)}{n^2}$$

and we therefore get

$$\int_0^1 (x\log x)\log[2\sin(\pi x/2)]\,dx =$$

$$\frac{1}{\pi^2}\varsigma_a(2) + \frac{1}{\pi^2}\varsigma(2) + \frac{1}{\pi^2}\sum_{n=1}^{\infty}\frac{Ci(n\pi)}{n^2} - \frac{1}{\pi^2}\sum_{n=1}^{\infty}\frac{(\gamma + \log n\pi)}{n^2}$$

$$= \frac{1}{\pi^2}\varsigma_a(2) + \frac{1}{\pi^2}\varsigma(2) + \frac{1}{\pi^2}\sum_{n=1}^{\infty}\int_0^{\pi}\frac{\cos nu - 1}{un^2}\,du + \frac{1}{\pi^2}\varsigma'(2)$$

Using

$$\int_0^{\pi}\frac{\cos nu - 1}{u}\,du = [(-1)^n - 1]\log\pi + n\int_0^{\pi}\sin nu\,\log u\,du$$

we then obtain

$$\sum_{n=1}^{\infty}\int_0^{\pi}\frac{\cos nu - 1}{un^2}\,du = -[\varsigma_a(2) + \varsigma(2)]\log\pi + \int_0^{\pi}\sum_{n=1}^{\infty}\frac{\sin nu}{n}\log u\,du$$

$$= -[\varsigma_a(2) + \varsigma(2)]\log\pi + \int_0^{\pi}\frac{\pi - u}{2}\log u\,du$$

$$= -[\varsigma_a(2) + \varsigma(2)]\log\pi - \frac{3}{8}\pi^2 + \frac{1}{4}\pi^2\log\pi$$

Therefore we get

(6.94u)
$$\int_0^1 (x\log x)\log[2\sin(\pi x/2)]\,dx = [\varsigma_a(2) + \varsigma(2)]\left[1 - \frac{1}{\pi^2}\log\pi\right] - \frac{3}{8} + \frac{1}{4}\log\pi + \frac{1}{\pi^2}\varsigma'(2)$$

$$= [\varsigma_a(2) + \varsigma(2)] - \frac{3}{8} + \frac{1}{\pi^2}\varsigma'(2)$$

We also have

$$\int_0^1 (x\log x)\sin n\pi x\,dx = -\frac{1}{\pi^2}\frac{Si(n\pi)}{n^2}$$

Some similar relations with the $\log\Gamma(x)$ function are considered in (6.108) et seq.

Let us now divide (6.87) by $t^2$ and integrate the result over the interval $[a, x]$ to obtain

$$\sum_{n=1}^{\infty}\frac{\varsigma(2n)(x^{2n}-a^{2n})}{2n(2n+1)} = \frac{1}{2}(\log x - \log a) - \pi\int_a^x\sum_{n=1}^{\infty}\left\{-\frac{\cos 2n\pi t}{2n\pi t}+\frac{\sin 2n\pi t}{(2n\pi)^2 t^2}\right\}dt$$

$$= \frac{1}{2}(\log x - \log a) + \pi\int_a^x\sum_{n=1}^{\infty}\frac{d}{dt}\left\{\frac{\sin 2n\pi t}{(2n\pi)^2 t}\right\}dt$$

$$= \frac{1}{2}(\log x - \log a) + \frac{\pi}{x}\sum_{n=1}^{\infty}\frac{\sin 2n\pi x}{(2n\pi)^2} - \frac{\pi}{a}\sum_{n=1}^{\infty}\frac{\sin 2n\pi a}{(2n\pi)^2}$$

Letting $a = 1/2$ we get

$$\frac{1}{2}\sum_{n=1}^{\infty}\frac{\varsigma(2n)x^{2n}}{n(2n+1)} - \frac{1}{2}\sum_{n=1}^{\infty}\frac{\varsigma(2n)}{n(2n+1)2^{2n}} = \frac{1}{2}\log 2x + \frac{\pi}{x}\sum_{n=1}^{\infty}\frac{\sin 2n\pi x}{(2n\pi)^2}$$

From (6.92) we have

$$\sum_{n=1}^{\infty}\frac{\varsigma(2n)}{n(2n+1)2^{2n}} = -1 + \log\pi$$

and we see that

$$\sum_{n=1}^{\infty}\frac{\varsigma(2n)x^{2n}}{n(2n+1)} = -1 + \log 2\pi x + \frac{2\pi}{x}\sum_{n=1}^{\infty}\frac{\sin 2n\pi x}{(2n\pi)^2}$$

From [130, p.148] we have

$$\sum_{n=1}^{\infty}\frac{\sin nx}{n^2} = -\int_0^x\log\left[2\sin\frac{t}{2}\right]dt$$

and accordingly we obtain

(6.94v) $$\sum_{n=1}^{\infty}\frac{\varsigma(2n)x^{2n+1}}{n(2n+1)} = -x + x\log 2\pi x - \frac{1}{2\pi}\int_{0}^{2\pi x}\log\left[2\sin\frac{t}{2}\right]dt$$

Srivastava and Choi [126, p.223] report that for $|x| < 1$

$$\sum_{n=1}^{\infty}\frac{\varsigma(2n)x^{2n+1}}{n(2n+1)} = [1-\log(2\pi)]x + x\log\Gamma(1-x)\Gamma(1+x) + \log\frac{G(1+x)}{G(1-x)}$$

The Wolfram Integrator gives us

$$\int\log\left[2\sin\frac{t}{2}\right]dt = -t\log(1-e^{it}) + t\log\left[2\sin\frac{t}{2}\right] + i\left[\frac{t^2}{4} + Li_2(e^{it})\right]$$

and thus

$$\int_{0}^{2\pi x}\log\left[2\sin\frac{t}{2}\right]dt = -2\pi x\log(1-e^{i2\pi x}) + 2\pi x\log[2\sin\pi x] + i\left[\pi^2 x^2 + Li_2(e^{i2\pi x})\right]$$

**Example 12:**

We now let $p(x) = 1$ in (6.50) and use the interval $[a,t]$ where $a > 0$ so that the conditions for the Riemann-Lebesgue lemma are satisfied at $x = a$: we have

(6.95) $$\sum_{n=1}^{\infty}\int_{a}^{t}\varsigma(2n)x^{2n-1}dx = \frac{1}{2}\int_{a}^{t}\frac{1}{x}dx - \pi\sum_{n=1}^{\infty}\int_{a}^{t}\sin 2n\pi x\,dx$$

Hence we get

(6.96)

$$\sum_{n=1}^{\infty}\frac{\varsigma(2n)}{2n}t^{2n} - \sum_{n=1}^{\infty}\frac{\varsigma(2n)}{2n}a^{2n} = \frac{1}{2}(\log t - \log a) + \pi\sum_{n=1}^{\infty}\left\{\frac{\cos 2n\pi t}{2n\pi} - \frac{\cos 2n\pi a}{2n\pi}\right\}$$

Therefore, using the Fourier series which will be derived in (7.8), we obtain

(6.97) $$\sum_{n=1}^{\infty}\frac{\varsigma(2n)}{n}t^{2n} = \log t - \log a - \log\sin\pi t + \log\sin\pi a + \sum_{n=1}^{\infty}\frac{\varsigma(2n)}{n}a^{2n}$$

Since the logarithm is a continuous function, in the limit as $a \to 0$, we have

(6.98) $$\lim_{a\to 0}\log\frac{\sin\pi a}{a} = \log\lim_{a\to 0}\left(\frac{\pi\sin\pi a}{\pi a}\right) = \log\pi$$

and hence we have

(6.99) $$\sum_{n=1}^{\infty}\frac{\varsigma(2n)}{n}t^{2n}=\log\pi t-\log\sin\pi t=\log\frac{\pi t}{\sin\pi t}$$

This identity is also recorded in [126, p.160] as

(6.99a) $$\sum_{n=1}^{\infty}\frac{\varsigma(2n)}{n}t^{2n}=\log\Gamma(1+t)+\log\Gamma(1-t)$$

and upon using the Euler reflection formula (6.61) it is easily seen that (6.99) and (6.99a) are equivalent. Alternatively, the above example could be presented as a further proof of the Euler reflection formula.

Let us now multiply (6.99) by $t$ and integrate to obtain

(6.100) $$\sum_{n=1}^{\infty}\frac{\varsigma(2n)x^{2n+2}}{n(2n+2)}=\int_{0}^{x}t\log\pi t\,dt-\int_{0}^{x}t\log\sin\pi t\,dt$$

(6.100a) $$x^2\sum_{n=1}^{\infty}\frac{\varsigma(2n)x^{2n}}{n(n+1)}=x^2\log\pi x-\frac{1}{2}x^2-2\int_{0}^{x}t\log\sin\pi t\,dt$$

or equivalently,

(6.100b) $$x^2\sum_{n=1}^{\infty}\frac{\varsigma(2n)x^{2n}}{n(n+1)}=x^2\log\pi x-\frac{1}{2}x^2-\frac{2}{\pi^2}\int_{0}^{\pi x}u\log\sin u\,dt$$

Using the Riemann-Lebesgue lemma, we give a simple derivation of the following integral in Section 8

(6.101)
$$2\int_{0}^{y}t\log\sin t\,dt=\sum_{n=1}^{\infty}\left(\frac{y^2\cos 2ny}{n}-\frac{\cos 2ny}{2n^3}-\frac{y\sin 2ny}{n^2}\right)+y^2\log\sin y+\frac{1}{2}\varsigma(3)$$

(6.101a) $$2\int_{0}^{\pi x}t\log\sin t\,dt=\sum_{n=1}^{\infty}\left(\frac{\pi^2 x^2\cos 2n\pi x}{n}-\frac{\cos 2n\pi x}{2n^3}-\frac{\pi x\sin 2n\pi x}{n^2}\right)$$

$$+\pi^2 x^2\log\sin\pi x+\frac{1}{2}\varsigma(3)$$

With $x=1/2$ (6.100b) becomes

(6.102) $$\sum_{n=1}^{\infty}\frac{\varsigma(2n)}{n(n+1)2^{2n}}=\log\frac{\pi}{2}-\frac{1}{2}-8\int_{0}^{1/2}t\log\sin\pi t\,dt$$

(6.102a) $$= \log\frac{\pi}{2} - \frac{1}{2} - \frac{8}{\pi^2}\int_0^{\pi/2} t\log\sin t\, dt$$

Therefore using the Euler integral (1.11) we have

(6.102b) $$\sum_{n=1}^{\infty}\frac{\varsigma(2n)}{n(n+1)2^{2n}} = \log\pi - \frac{1}{2} - \frac{7}{2\pi^2}\varsigma(3)$$

Employing partial fractions we obtain

(6.102c) $$\sum_{n=1}^{\infty}\frac{\varsigma(2n)}{n(n+1)2^{2n}} = \sum_{n=1}^{\infty}\frac{\varsigma(2n)}{n2^{2n}} - \sum_{n=1}^{\infty}\frac{\varsigma(2n)}{(n+1)2^{2n}}$$

Using (6.99) we have for $x = 1/2$

(6.102d) $$\sum_{n=1}^{\infty}\frac{\varsigma(2n)}{n2^{2n}} = \log\frac{\pi}{2}$$

and this is reported in [25, p.131]. Using (6.102b) we have

(6.102e) $$\sum_{n=1}^{\infty}\frac{\varsigma(2n)}{(n+1)2^{2n}} = \frac{1}{2} - \log 2 + \frac{7\varsigma(3)}{2\pi^2}$$

and we have previously seen this in (6.85). When I was doing these calculations, I initially thought that I could combine the above equations so as to obtain a closed form expression for $\varsigma(3)$ but, alas, the interesting part cancels out!

Quite clearly equation (6.100) can be generalised to produce similar infinite series for $\varsigma(2m+1)$ by multiplying (6.99) by $t^m$ or, perhaps more simply, by the corresponding Bernoulli polynomial. Some generalised formulae of a similar nature (the Ramanujan-Yoshimoto formula) are contained in [76b]: this formula is far reaching, and apparently the formulae contained in some 180 pages of the book by Srivastava and Choi [126] can readily be deduced from it.

**Example 13:**

Let us now revisit (6.1a).

$$\frac{1}{1-e^{ix}} = \sum_{n=0}^{N} e^{inx} + R_N(x) \text{ where } R_N(x) = \frac{ie^{i(N+\frac{1}{2})x}}{2\sin(x/2)}$$

Upon differentiation we obtain

$$\frac{ie^{ix}}{(1-e^{ix})^2} = i\sum_{n=0}^{N} ne^{inx} + R'_N(x)$$

and hence we get

$$-\frac{1}{4\sin^2(x/2)} = \frac{e^{ix}}{(1-e^{ix})^2} = \sum_{n=0}^{N} n e^{inx} - iR'_N(x)$$

Multiplying by $p(x)$ and integrating we have

$$-\frac{1}{4}\int_a^b \frac{p(x)}{\sin^2(x/2)} dx = \int_a^b \sum_{n=0}^{N} np(x) e^{inx} dx - i\int_a^b p(x) R'_N(x) dx$$

Using integration by parts we have (assuming that the relevant functions are Riemann integrable)

$$i\int_a^b p(x) R'_N(x) dx = ip(x) R_N(x) \Big|_a^b - i\int_a^b p'(x) R_N(x) dx$$

Therefore, if (i) $p(a) = p(b) = 0$ and (ii) $\sin(x/2) \neq 0$ for all $x \in [a,b]$ and $p'(x)$ is bounded on $[a,b]$ we have using the Riemann-Lebesgue lemma

$$i\int_a^b p(x) R'_N(x) dx = 0$$

Hence, assuming that the above conditions are met, as $N \to \infty$ we have

(6.103a) $$\sum_{n=0}^{\infty} \int_a^b p(x) n \cos nx\, dx = -\frac{1}{4} \int_a^b \frac{p(x)}{\sin^2(x/2)} dx$$

(6.103b) $$\sum_{n=0}^{\infty} \int_a^b p(x) n \sin nx\, dx = 0$$

Equation (6.103a) reminded me of Stark's 1974 paper [124] where, from quite the opposite direction and using the Fejér kernel $F_n(x)$, he showed that

(6.103c) $$\varsigma_a(3) = \frac{\pi^2}{6} \log 2 - \frac{1}{3\pi} \int_0^{\pi/2} \frac{t^4}{\sin^2 t} dt$$

(6.103d) $$\varsigma_a(5) = \frac{7\pi^3}{360} \int_0^{\pi/2} \frac{t^2}{\sin^2 t} dt - \frac{\pi}{18} \int_0^{\pi/2} \frac{t^4}{\sin^2 t} dt + \frac{2}{45\pi} \int_0^{\pi/2} \frac{t^6}{\sin^2 t} dt$$

where the Fejér kernel is defined by

(6.103e) $$F_n(x) = \frac{1}{2(n+1)} \left[ \frac{\sin\left(\frac{n+1}{2}\right)x}{\sin\left(\frac{x}{2}\right)} \right]^2$$

These equations may also be derived directly from (6.103a) by letting $p(x) = x(x-\pi/2)^r$ and $[a,b] = [0, \pi/2]$. Whilst Stark's paper focuses on the interval $[0, \pi/2]$, equation (6.103a) may be employed to derive countless identities with the more general interval $[a,b]$ and indeed many more suitable functions $p(x)$. See also [125aa]. The integrals in (6.103c) and (6.103d) are also evaluated in Muzaffar's recent paper [103ac].

More generally we have

(6.103f) $$\sum_{n=0}^{\infty} \int_a^b p(x) n \cos \alpha n x \, dx = -\frac{1}{4} \int_a^b \frac{p(x)}{\sin^2(\alpha x/2)} dx$$

With integration by parts it is easily seen that

$$\int_a^b \frac{p(x)}{\sin^2(\alpha x/2)} dx = -\frac{2p(x)\cot(\alpha x/2)}{\alpha}\bigg|_a^b + \frac{2}{\alpha}\int_a^b p'(x)\cot(\alpha x/2) dx$$

$$= \frac{2}{\alpha} \int_a^b p'(x) \cot(\alpha x/2) dx$$

since we have postulated that $p(a) = p(b) = 0$. In what follows we let $\alpha = 2$ and hence we have

$$\sum_{n=0}^{\infty} \int_a^b p(x) n \cos 2nx \, dx = -\frac{1}{4}\int_a^b \frac{p(x)}{\sin^2 x} dx = -\frac{1}{4}\int_a^b p'(x) \cot x \, dx$$

Therefore, using (6.46) we have

$$\sum_{n=0}^{\infty} \int_a^b p(x) n \cos 2nx \, dx = -\frac{1}{4}\int_a^b \frac{p(x)}{\sin^2 x} dx = \frac{1}{2}\sum_{n=0}^{\infty} \int_a^b p'(x) \frac{\varsigma(2n)}{x}\left(\frac{x}{\pi}\right)^n dx$$

We now let $[a,b] = [0, \pi/2]$ and $p(x) = x^2(x-\pi/2)^2$ so that $p(0) = p(\pi/2) = 0$

$$\sum_{n=0}^{\infty} \int_0^{\pi/2} x^2(x-\pi/2)^2 n \cos 2nx \, dx = -\frac{1}{4}\int_0^{\pi/2} \frac{x^2(x-\pi/2)^2}{\sin^2 x} dx$$

$$= \sum_{n=0}^{\infty} \int_0^{\pi/2} (x-\pi/2)(2x-\pi/2)\varsigma(2n)\left(\frac{x}{\pi}\right)^n dx$$

The integrals are easily evaluated and we get

$$\varsigma(4) - \varsigma_a(4) = \frac{\pi^2}{3} \sum_{n=0}^{\infty} \frac{(n-1)\varsigma(2n)}{(n+1)(n+2)(n+3)2^n}$$

which may be expressed as

(6.103g) $$\frac{\pi^2}{240} = \sum_{n=0}^{\infty} \frac{(n-1)\varsigma(2n)}{(n+1)(n+2)(n+3)2^n}$$

where, as before, we have used $\varsigma(0) = -1/2$.

**Example 14:**

It is elementary to prove the identity

(6.104) $$\frac{1}{1-ie^{ix}} = \frac{1}{2} + \frac{i}{2}\frac{\cos x}{1+\sin x}$$

and using (2.3) we have

$$\frac{1}{1-ie^{ix}} = 1 + \sum_{n=1}^{2N} \left(ie^{ix}\right)^n + R_{2N+1}(x)$$

where

$$R_{2N+1}(x) = \frac{(ie^{ix})^{2N+1}}{1-ie^{ix}}$$

Now it is easily seen that

$$\sum_{n=1}^{2N}\left(ie^{ix}\right)^n = \sum_{n=1}^{N} i^{2n}e^{2inx} + \sum_{n=1}^{N} i^{2n-1}e^{(2n-1)ix}$$

$$= \sum_{n=1}^{N}(-1)^n\{\cos 2nx + i\sin 2nx\} - i\sum_{n=1}^{N}(-1)^n\{\cos(2n-1)x + i\sin(2n-1)x\}$$

$$= \sum_{n=1}^{N}(-1)^n\{\cos 2nx + \sin(2n-1)x\} + i\sum_{n=1}^{N}(-1)^n\{\sin 2nx - \cos(2n-1)x\}$$

Employing the same integration technique as before, we have

$$\int_a^b p(x)\left\{\frac{1}{2}+\frac{i}{2}\frac{\cos x}{1+\sin x}\right\}dx = \int_a^b p(x)dx + \int_a^b p(x)\sum_{n=1}^N(-1)^n\{\cos 2nx + \sin(2n-1)x\}dx$$

$$+i\int_a^b p(x)\sum_{n=1}^N(-1)^n\{\sin 2nx - \cos(2n-1)x\}dx + \int_a^b p(x)R_{2N+1}(x)dx$$

Now we have

$$R_{2N+1} = \int_a^b p(x)R_{2N+1}(x)dx$$

$$= \int_a^b p(x)i(-1)^N\left(\cos(2N+1)x + i\sin(2N+1)x\right)\left(\frac{1}{2}+\frac{i}{2}\frac{\cos x}{1+\sin x}\right)dx$$

and from the analysis carried out in proving (2.16b) it is clear that $R_{2N+1} \to 0$ as $N \to \infty$ provided that $p'(x)$ is bounded on $[a,b]$ and $\sin x \neq -1$ on $[a,b]$.

We therefore obtain the following identities

(6.105) $$-\frac{1}{2}\int_a^b p(x)dx = \sum_{n=1}^\infty (-1)^n \int_a^b p(x)\{\cos 2nx + \sin(2n-1)x\}dx$$

(6.105a) $$\frac{1}{2}\int_a^b p(x)\left\{\frac{\cos x}{1+\sin x}\right\}dx = \sum_{n=1}^\infty (-1)^n \int_a^b p(x)\{\sin 2nx - \cos(2n-1)x\}dx$$

For example, let $p(x) = x$ in (6.105a) over the range $[0,\pi]$. Then we have

$$\frac{1}{2}\int_0^\pi x\left\{\frac{\cos x}{1+\sin x}\right\}dx$$

$$= \sum_{n=1}^\infty (-1)^n \left\{-\frac{\cos(2n-1)x}{(2n-1)^2} - \frac{x\sin(2n-1)x}{2n-1} - \frac{x\cos 2nx}{2n} + \frac{\sin 2nx}{2n}\right\}\Big|_0^\pi$$

$$= \sum_{n=1}^\infty (-1)^n \left\{-\frac{(-1)^{2n-1}}{(2n-1)^2} - \frac{1}{(2n-1)^2} - \frac{\pi}{2n}\right\}$$

$$= -2\sum_{n=1}^\infty \frac{(-1)^n}{(2n-1)^2} + \frac{\pi}{2}\sum_{n=1}^\infty \frac{(-1)^n}{n}$$

Therefore we get

(6.106) $$\frac{1}{2}\int_0^\pi x\left\{\frac{\cos x}{1+\sin x}\right\}dx = -2G + \frac{\pi}{2}\log 2$$

where $G$ is Catalan's constant. This is in agreement with the formula in G&R [74, p.442, 3.791.2]. The main attraction of the current method is that the given integral may be evaluated systematically over many different intervals of integration.

An alternative proof is given below. Employing integration by parts, we obtain

$$I = \int_0^\pi x \left\{ \frac{\cos x}{1+\sin x} \right\} dx = x \log(1+\sin x) \Big|_0^\pi - \int_0^\pi \log(1+\sin x) dx$$

$$= -\int_0^\pi \log(1+\sin x) dx$$

We have

$$\log(1+\sin x) = \log(\sin x/2 + \cos x/2)^2$$

$$= 2\log \sin\left[(x/2) + \cos(x/2)\right]$$

$$= 2\log\left\{\sqrt{2} \sin\left(\frac{x}{2} + \frac{\pi}{4}\right)\right\}$$

$$= \log 2 + 2\log \sin\left(\frac{x}{2} + \frac{\pi}{4}\right)$$

Therefore we have

$$I = -\pi \log 2 - 2\int_0^\pi \log \sin\left(\frac{x}{2} + \frac{\pi}{4}\right) dx$$

and, using the substitution $u + \frac{\pi}{2} = \frac{x}{2} + \frac{\pi}{4}$, we get

$$I = -\pi \log 2 - 8\int_0^{\pi/4} \log \cos u \, du$$

We now use (6.8a) with $\alpha = 2$ and $p(x) = x$ to obtain

$$-\frac{1}{2}\int_0^t x \tan x \, dx = \sum_{n=0}^\infty \int_0^t (-1)^n x \sin 2nx \, dx$$

With integration by parts we have

$$\int_0^t x \tan x \, dx = -t \log \cos t + \int_0^t \log \cos x \, dx$$

$$\int_0^t x \sin 2nx \, dx = -\frac{t \cos 2nt}{2n} + \frac{\sin 2nt}{4n^2} \quad , n \geq 1$$

Hence we get

$$\frac{1}{2} t \log \cos t - \frac{1}{2} \int_0^t \log \cos x \, dx = \sum_{n=1}^{\infty} (-1)^n \left[ -\frac{t \cos 2nt}{2n} + \frac{\sin 2nt}{4n^2} \right]$$

$$= \frac{t}{2} \sum_{n=1}^{\infty} (-1)^{n+1} \frac{\cos 2nt}{n} - \frac{1}{4} \sum_{n=1}^{\infty} (-1)^{n+1} \frac{\sin 2nt}{n^2}$$

We have already seen in (6.22a) that

$$\sum_{n=1}^{\infty} (-1)^{n+1} \frac{\cos nt}{n} = \log[2 \cos(t/2)]$$

and therefore we get

(6.107) $$\int_0^t \log \cos x \, dx = \frac{1}{2} \sum_{n=1}^{\infty} (-1)^{n+1} \frac{\sin 2nt}{n^2} - t \log 2$$

This may be written in its more familiar form [130, p.148]

(6.107a) $$\int_0^t \log\left[ 2 \cos \frac{x}{2} \right] dx = \sum_{n=1}^{\infty} (-1)^{n+1} \frac{\sin nt}{n^2}$$

With $t = \pi/4$ we obtain

$$\int_0^{\pi/4} \log \cos x \, dx = \frac{1}{2} \sum_{n=1}^{\infty} (-1)^{n+1} \frac{\sin(n\pi/2)}{n^2} - \frac{\pi}{4} \log 2$$

$$= \frac{1}{2} \sum_{n=0}^{\infty} \frac{(-1)^{n+1}}{(2n+1)^2} - \frac{\pi}{4} \log 2$$

(6.107b) $$= \frac{G}{2} - \frac{\pi}{4} \log 2$$

This integral is well-known and is contained, inter alia, in G&R [74, p.526].

Lobachevsky's function $L(x)$ is defined in G&R [74, p.883] as

$$L(x) = -\int_0^x \log \cos t \, dt$$

and using (6.107) it is easily seen that

$$L(\pi/2) = -\int_0^{\pi/2} \log \cos t \, dt = \frac{\pi}{2} \log 2$$

Using (6.5a) it can similarly be shown that the Clausen function $\text{Cl}_2(t)$ referred to in (6.66) is

(6.107c) $\qquad \text{Cl}_2(t) = -\int_0^t \log\left[2 \sin \frac{x}{2}\right] dx = \sum_{n=1}^{\infty} \frac{\sin nt}{n^2} \quad , 0 \le x \le 2\pi$

and hence we have

(6.107d) $\qquad \int_0^{\pi/4} \log \sin x \, dx = -\frac{G}{2} - \frac{\pi}{4} \log 2$

We have another Lobachevsky function defined by

(6.107e) $\qquad \Lambda(x) = -\int_0^x \log|2 \sin t| \, dt = \frac{1}{2} \sum_{n=1}^{\infty} \frac{\sin 2nx}{n^2}$

This function is odd with period $\pi$ and its graph vaguely resembles a slightly distorted sine function. Now recall (6.46)

$$\cot x = \frac{1}{x} \sum_{n=0}^{\infty} (-1)^n \frac{2^{2n} B_{2n}}{(2n)!} x^{2n} \quad , (|x| < \pi)$$

Integrating this we obtain

$$\int_a^t \left(\cot x - \frac{1}{x}\right) dx = \sum_{n=1}^{\infty} (-1)^n \frac{2^{2n} B_{2n}}{2n(2n)!} (t^{2n} - a^{2n})$$

Using L'Hôpital's rule it is easily shown that $\lim_{x \to 0}\left[\cot x - \frac{1}{x}\right] = 0$. We have

$$\lim_{x \to 0}\left[\cot x - \frac{1}{x}\right] = \lim_{x \to 0}\left[\left(\frac{x}{\sin x} \cos x - 1\right)/x\right] = \frac{\lim_{x \to 0}\left[\frac{x}{\sin x} \cos x - 1\right]}{\lim_{x \to 0} x} \approx \frac{0}{0}$$

and applying L'Hôpital's rule again we have

$$\frac{\lim_{x\to 0}\left[\dfrac{x}{\sin x}\cos x - 1\right]}{\lim_{x\to 0} x} = \lim_{x\to 0}\left[\frac{\cos x}{\sin x} - \frac{x}{\sin^2 x}\right] = \lim_{x\to 0}\left[\frac{\cos x - x/\sin x}{\sin x}\right] \approx \frac{0}{0}$$

Yet a further application of L'Hôpital's rule gives us

$$\lim_{x\to 0}\left[\frac{\cos x - x/\sin x}{\sin x}\right] = \frac{\lim_{x\to 0}\left[-\sin x - \left(\dfrac{\sin x - x\cos x}{\sin^2 x}\right)\right]}{\cos x}$$

We now need to consider

$$\lim_{x\to 0}\left[\frac{\sin x - x\cos x}{\sin^2 x}\right] \approx \frac{0}{0} = \frac{\lim_{x\to 0}(x\sin x)}{\lim_{x\to 0}(2\sin x\cos x)} = \frac{\lim_{x\to 0}(x\cos x + \sin x)}{\lim_{x\to 0} 2(\cos^2 x - \sin^2 x)} = 0$$

and the required limit then follows.

We have

$$\int_a^t \left(\cot x - \frac{1}{x}\right)dx = \log\sin t - \log t - \log\frac{\sin a}{a}$$

and hence we get as $a \to 0$

(6.107ei) $\qquad \displaystyle\int_0^t \left(\cot x - \frac{1}{x}\right)dx = \log\sin t - \log t = \sum_{n=1}^{\infty}(-1)^n \frac{2^{2n} B_{2n}}{2n(2n)!} t^{2n}$

Integrating once more we get

$$\int_0^x \left[\log\sin t - \log t\right]dt = \sum_{n=1}^{\infty}(-1)^n \frac{2^{2n} B_{2n}}{2n(2n+1)!} x^{2n+1}$$

and therefore we have

(6.107f) $\qquad \Lambda(x) = x - x\log 2x + x\sum_{n=1}^{\infty}(-1)^{n+1}\frac{B_{2n}}{2n(2n+1)!}(2x)^{2n}$

I came across this identity whilst searching for the Lobachevsky function on the internet: the source was a homework paper on hyperbolic geometry by D. Calegari. Letting $x = \pi/2$, we have $\Lambda(\pi/2) = 0$ and deduce that

(6.107g) $\qquad -1 + \log\pi = \displaystyle\sum_{n=1}^{\infty}(-1)^{n+1}\frac{B_{2n}\pi^{2n}}{2n(2n+1)!}$

and, using (1.7) from Volume I,

$$\varsigma(2n) = (-1)^{n+1} \frac{2^{2n-1} \pi^{2n} B_{2n}}{(2n)!}$$

it is easily seen that this is equivalent to (6.92).

We also have from (6.107e)

$$\frac{1}{2}\sum_{n=1}^{\infty} \frac{\sin 2nx}{n^2} = x - x\log 2x + x\sum_{n=1}^{\infty} (-1)^{n+1} \frac{B_{2n}}{2n(2n+1)!}(2x)^{2n}$$

$$= x - x\log 2x + x\sum_{n=1}^{\infty} \frac{\varsigma(2n)}{n(2n+1)}\left(\frac{x}{\pi}\right)^{2n}$$

Integrating one more time we obtain

$$\frac{1}{2}\int_0^t \sum_{n=1}^{\infty} \frac{\sin 2nx}{n^2} dx = \int_0^t \left[ x - x\log 2x + x\sum_{n=1}^{\infty} \frac{\varsigma(2n)}{n(2n+1)}\left(\frac{x}{\pi}\right)^{2n} \right] dx$$

Therefore we obtain

(6.107h) $$\frac{1}{4}\varsigma(3) - \frac{1}{4}\sum_{n=1}^{\infty} \frac{\cos 2nt}{n^3} = \frac{3}{4}t^2 - \frac{1}{2}t^2\log(2t) + \frac{1}{2}\sum_{n=1}^{\infty} \frac{\varsigma(2n)t^{2n+2}}{n(2n+1)(n+1)\pi^{2n}}$$

With $t = \pi/2$ we get

$$\frac{1}{4}\left[\varsigma(3) + \varsigma_a(3)\right] = \frac{3}{16}\pi^2 - \frac{1}{8}\pi^2\log\pi + \frac{\pi^2}{8}\sum_{n=1}^{\infty} \frac{\varsigma(2n)}{n(2n+1)(n+1)2^{2n}}$$

and, with a little algebra, we obtain

(6.107i) $$\sum_{n=1}^{\infty} \frac{\varsigma(2n)}{n(2n+1)(n+1)2^{2n}} = \frac{7}{2}\frac{\varsigma(3)}{\pi^2} - \frac{3}{2} + \log\pi$$

In [126, p.229] we find that

(6.107j) $$\sum_{n=1}^{\infty} \frac{\varsigma(2n)}{n(2n+1)(n+1)2^{2n}} = -\frac{3}{2} + \log\left(\pi.B^{14}\right)$$

where $B$ is one of the generalised Glaisher-Kinkelin constants which we referred to previously in (6.83). We have from (6.84)

$$\log B = -\varsigma'(-2) = \frac{\varsigma(3)}{4\pi^2}$$

and it is easily seen that (6.107i) and (6.107j) are equivalent.

Integrating (6.107h) one more time we see that

(6.107k)

$$\frac{1}{4}x\varsigma(3) - \frac{1}{8}\sum_{n=1}^{\infty}\frac{\sin 2nx}{n^4} = \frac{1}{4}x^3 - \frac{1}{18}x^3\left[3\log(2x) - 1\right] + \frac{1}{8}\sum_{n=1}^{\infty}\frac{\varsigma(2n)x^{2n+3}}{2n(2n+1)(2n+2)(2n+3)\pi^{2n}}$$

and with $x = \pi/2$ we get

(6.107l) $$\sum_{n=1}^{\infty}\frac{\varsigma(2n)}{2n(2n+1)(2n+2)(2n+3)2^{2n}} = \frac{\varsigma(3)}{2\pi^2} + \frac{1}{12}\log\pi - \frac{11}{72}$$

This was originally derived by Wilton in 1922 (see [126, p.148]).

Integrating (6.107ei) results in

$$\int_0^x t\log\sin t \, dt - \int_0^x t\log t \, dt = \sum_{n=1}^{\infty}(-1)^n \frac{2^{2n}B_{2n}}{2n(2n+2)(2n)!}x^{2n+2}$$

and with $x = \pi/2$ and using Euler's equation (1.11) from Volume I

$$\int_0^{\pi/2} t\log\sin t \, dt = \frac{7}{16}\varsigma(3) - \frac{\pi^2}{8}\log 2$$

we obtain

$$\frac{7}{16}\varsigma(3) - \frac{\pi^2}{8}\log 2 - \left(\frac{\pi^2}{8}\log\pi - \frac{\pi^2}{8}\log 2 - \frac{\pi^2}{16}\right) = \frac{\pi^2}{4}\sum_{n=1}^{\infty}(-1)^n \frac{2^{2n-1}B_{2n}\pi^{2n}}{n(2n+2)(2n)!2^{2n}}$$

We then have

$$\frac{7}{16}\varsigma(3) - \frac{\pi^2}{8}\log\pi + \frac{\pi^2}{16} = -\frac{\pi^2}{4}\sum_{n=1}^{\infty}\frac{\varsigma(2n)}{n(2n+2)2^{2n}}$$

$$= \frac{\pi^2}{8}\sum_{n=1}^{\infty}\frac{\varsigma(2n)}{(n+1)2^{2n}} - \frac{\pi^2}{8}\sum_{n=1}^{\infty}\frac{\varsigma(2n)}{n2^{2n}}$$

Using (6.79)

$$\frac{1}{2} - \log 2 + \frac{7\varsigma(3)}{2\pi^2} = \sum_{n=1}^{\infty}\frac{\varsigma(2n)}{(n+1)2^{2n}}$$

we obtain (as reported in [126, p.162])

(6.107m) $$\log\frac{\pi}{2} = \sum_{n=1}^{\infty}\frac{\varsigma(2n)}{n2^{2n}}$$

## SOME INTEGRALS INVOLVING THE CLAUSEN FUNCTION

Further integral forms of the Clausen function are derived below. Define $A(\theta)$ by

$$A(\theta) = \int_0^1 \tan^{-1}\left(\frac{x\sin\theta}{1-x\cos\theta}\right)\frac{dx}{x}$$

Then, differentiating under the integral sign, we easily obtain

$$A'(\theta) = \int_0^1 \frac{\cos\theta - x}{1-2x\cos\theta + x^2}dx = -\frac{1}{2}\int_0^1 \frac{-2\cos\theta + 2x}{1-2x\cos\theta + x^2}dx$$

$$= -\frac{1}{2}\log(1-2x\cos\theta + x^2)\Big|_0^1 = -\frac{1}{2}\log 2(1-\cos\theta)$$

$$= -\log[2\sin(\theta/2)]$$

Therefore, using (7.8) we get

$$A'(\theta) = -\log[2\sin(\theta/2)] = \sum_{n=1}^{\infty}\frac{\cos n\theta}{n}$$

and, since $A(0) = 0$ we obtain

$$A(\theta) = \sum_{n=1}^{\infty}\frac{\sin n\theta}{n^2}$$

Accordingly, we have

(6.107n) $$Cl_2(\theta) = -\int_0^{\theta}\log\left[2\sin\frac{x}{2}\right]dx = \sum_{n=1}^{\infty}\frac{\sin n\theta}{n^2} = \int_0^1 \tan^{-1}\left(\frac{x\sin\theta}{1-x\cos\theta}\right)\frac{dx}{x}$$

Employing integration by parts we get

$$\int_0^1 \tan^{-1}\left(\frac{x\sin\theta}{1-x\cos\theta}\right)\frac{dx}{x} = \tan^{-1}\left(\frac{x\sin\theta}{1-x\cos\theta}\right)\log x\Big|_0^1 - \sin\theta\int_0^1 \frac{\log x}{1-2x\cos\theta + x^2}dx$$

Hence we also have

(6.107o) $$\operatorname{Cl}_2(\theta) = -\sin\theta \int_0^1 \frac{\log x}{1 - 2x\cos\theta + x^2} dx$$

and therefore with $\theta = \pi/2$ we get

(6.107p) $$\operatorname{Cl}_2(\pi/2) = -\int_0^1 \frac{\log x}{1 + x^2} dx = G$$

Using the Poisson kernel

$$\frac{1}{2} + \sum_{n=1}^\infty x^n \cos n\theta = \frac{1}{2} \frac{1 - x^2}{1 - 2x\cos\theta + x^2}$$

we have

(6.107pi) $$\sum_{n=1}^\infty x^{n-1} \cos n\theta = \frac{\cos\theta - x}{1 - 2x\cos\theta + x^2}$$

We may the write (6.107o) as

$$-\frac{\operatorname{Cl}_2(\theta)}{\sin\theta} = \int_0^1 \sum_{n=1}^\infty \frac{x^{n-1} \log x}{\cos\theta - x} \cos n\theta \, dx$$

The Wolfram Integrator can evaluate $\int_0^1 \frac{x^{n-1} \log x}{\cos\theta - x} dx$ but I have not explored the consequences of this in any detail.

$$\int_0^1 \frac{x^{n-1} \log x}{a - x} dx = -\operatorname{Li}_2\left(\frac{1}{a}\right) a^{n-1} + a^{n-1}\left[\frac{1}{a1^2} + \frac{1}{a^2 2^2} + \ldots + \frac{1}{a^{n-1}(n-1)^2}\right]$$

Integrating (6.107o) gives us

$$\operatorname{Cl}_3(u) - \varsigma(3) = \int_0^u \operatorname{Cl}_2(\theta) d\theta = \frac{1}{2} \int_0^1 \frac{\log x}{x} dx \int_0^u \frac{2x\sin\theta}{1 - 2x\cos\theta + x^2} d\theta$$

We see that

$$\int_0^u \frac{2x\sin\theta}{1 - 2x\cos\theta + x^2} d\theta = \log[1 - 2x\cos u + x^2] - 2\log(1 - x)$$

and hence we get

$$Cl_3(u) - \varsigma(3) = \frac{1}{2}\int_0^1 \frac{\log x \log[1-2x\cos u + x^2]}{x}dx - \int_0^1 \frac{\log x \log(1-x)}{x}dx$$

We have already seen that

$$\int \frac{\log x \log(1-x)}{x}dx = Li_3(x) - Li_2(x)\log x$$

and we therefore obtain

$$Cl_3(u) = \frac{1}{2}\int_0^1 \frac{\log x \log[1-2x\cos u + x^2]}{x}dx$$

In particular we have with $u = \pi/2$

$$Cl_3(\pi/2) = \frac{1}{2}\int_0^1 \frac{\log x \log[1+x^2]}{x}dx = -\frac{3}{32}\varsigma(3)$$

where we have used the formula for $Cl_3(\pi/2)$ in Volume IV. This may also be derived using the formula in Bromwich's book [36b, p.187]

$$\log[1-2x\cos u + x^2] = -2\sum_{n=1}^{\infty}\frac{x^n \cos nu}{n}$$

which is valid for $0 \leq x \leq 1$ and $0 < u < 2\pi$ (and which is easily obtained by integrating (6.107pi)).

The Wolfram Integrator kindly evaluates the following integral

$$2\int \frac{\log x \log[1-ax+x^2]}{x}dx = -\log\left[1+\frac{2x}{\sqrt{a^2-4}}\right]\log^2 x - \log\left[1-\frac{2x}{a+\sqrt{a^2-4}}\right]\log^2 x$$

$$+ \log[1-ax+x^2]\log^2 x - 2Li_2\left[\frac{2x}{a-\sqrt{a^2-4}}\right]\log x - 2Li_2\left[\frac{2x}{a+\sqrt{a^2-4}}\right]\log x$$

$$+ 2Li_3\left[\frac{2x}{a-\sqrt{a^2-4}}\right] + 2Li_3\left[\frac{2x}{a+\sqrt{a^2-4}}\right]$$

and we obtain the definite integral

$$2\int_0^t \frac{\log x \log[1-ax+x^2]}{x}dx = -\log\left[1+\frac{2t}{\sqrt{a^2-4}}\right]\log^2 t - \log\left[1-\frac{2t}{a+\sqrt{a^2-4}}\right]\log^2 t$$

$$+\log[1-at+t^2]\log^2 t - 2Li_2\left[\frac{2t}{a-\sqrt{a^2-4}}\right]\log t - 2Li_2\left[\frac{2t}{a+\sqrt{a^2-4}}\right]\log t$$

$$+2Li_3\left[\frac{2t}{a-\sqrt{a^2-4}}\right] + 2Li_3\left[\frac{2t}{a+\sqrt{a^2-4}}\right]$$

With $t=1$ we obtain a more compact result

$$\int_0^1 \frac{\log x \log[1-ax+x^2]}{x}dx = Li_3\left[\frac{2}{a-\sqrt{a^2-4}}\right] + Li_3\left[\frac{2}{a+\sqrt{a^2-4}}\right]$$

and hence we have with $a = 2\cos u$

$$Cl_3(u) = \frac{1}{2}\int_0^1 \frac{\log x \log[1-2x\cos u + x^2]}{x}dx = \frac{1}{2}Li_3\left[\frac{1}{\cos u - i\sin u}\right] + \frac{1}{2}Li_3\left[\frac{1}{\cos u + i\sin u}\right]$$

This gives us

$$Cl_3(u) = \frac{1}{2}Li_3\left[e^{iu}\right] + \frac{1}{2}Li_3\left[e^{-iu}\right]$$

and, upon taking the real part, we simply come right back to the definition of the Clausen function, namely

$$Cl_3(u) = \sum_{n=1}^\infty \frac{\cos nu}{n^3}$$

More generally [51b] we have

(6.107q) $\quad Cl_{2n}(\theta) = -\frac{\sin \theta}{(2n-1)!}\int_0^1 \frac{\log^{2n-1} x}{1-2x\cos\theta + x^2}dx$

(6.107r) $\quad Cl_{2n+1}(\theta) = -\frac{1}{(2n)!}\int_0^1 \frac{(x-\cos\theta)\log^{2n} x}{1-2x\cos\theta + x^2}dx$

A novel way of deriving the above two identities has recently been given by Efthimiou [58aa]. His method is shown below.

Since from (4.4.28) $\frac{1}{n^s} = \frac{1}{\Gamma(s)}\int_0^\infty e^{-nt}t^{s-1}dt$ we may write

$$\sum_{n=1}^\infty \frac{\cos nx}{n^s} = \frac{1}{\Gamma(s)}\sum_{n=1}^\infty \cos nx \int_0^\infty e^{-nt}t^{s-1}dt$$

$$= \frac{1}{\Gamma(s)} \int_0^\infty \left( \sum_{n=1}^\infty \cos nx \, e^{-nt} \right) t^{s-1} dt$$

$$= \frac{1}{\Gamma(s)} \int_0^\infty \frac{e^{-t}(\cos x - e^{-t})}{1 - 2\cos x \, e^{-t} + e^{-2t}} t^{s-1} dt$$

Letting $u = e^{-t}$ we obtain

(6.107ri) $\quad \displaystyle\sum_{n=1}^\infty \frac{\cos nx}{n^s} = \frac{(-1)^s}{\Gamma(s)} \int_0^1 \frac{(u - \cos x) \log^{s-1} u}{1 - 2u \cos x + u^2} du$

Similarly we obtain

(6.107rii) $\quad \displaystyle\sum_{n=1}^\infty \frac{\sin nx}{n^s} = \frac{(-1)^{s-1} \sin x}{\Gamma(s)} \int_0^1 \frac{\log^{s-1} u}{1 - 2u \cos x + u^2} du$

We may also write (6.107rii) as

(6.107riii) $\quad \displaystyle\sum_{n=1}^\infty \frac{\sin nx}{n^s} = \frac{\sin x}{\Gamma(s)} \int_0^1 \frac{\log^{s-1}(1/u)}{1 - 2u \cos x + u^2} du$

and then differentiate with respect to $s$ to obtain

(6.107rv) $\quad -\displaystyle\sum_{n=1}^\infty \frac{\log n}{n^s} \sin nx = \frac{\sin x}{\Gamma(s)} \int_0^1 \frac{\log^{s-1}(1/u) \log \log(1/u)}{1 - 2u \cos x + u^2} du$

$$-\frac{\psi(s) \sin x}{\Gamma(s)} \int_0^1 \frac{\log^{s-1}(1/u)}{1 - 2u \cos x + u^2} du$$

or alternatively using (6.107riii)

(6.107rv) $\quad -\displaystyle\sum_{n=1}^\infty \frac{\log n}{n^s} \sin nx = \frac{\sin x}{\Gamma(s)} \int_0^1 \frac{\log^{s-1}(1/u) \log \log(1/u)}{1 - 2u \cos x + u^2} du - \psi(s) \sum_{n=1}^\infty \frac{\sin nx}{n^s}$

With $s = 1$ we have

$$-\sum_{n=1}^\infty \frac{\log n}{n} \sin nx = \sin x \int_0^1 \frac{\log \log(1/u)}{1 - 2u \cos x + u^2} du + \gamma \sum_{n=1}^\infty \frac{\sin nx}{n}$$

and with $x \to 2\pi x$ this becomes

$$-\sum_{n=1}^{\infty}\frac{\log n}{n}\sin 2\pi nx=\sin 2\pi x\int_{0}^{1}\frac{\log\log(1/u)}{1-2u\cos 2\pi x+u^{2}}du+\gamma\sum_{n=1}^{\infty}\frac{\sin 2\pi nx}{n}$$

Then referring to Kummer's Fourier series (E.44a)

$$\log\Gamma(x)=\frac{1}{2}\log\frac{\pi}{\sin\pi x}+\frac{1}{2}(1-2x)[\gamma+\log(2\pi)]+\frac{1}{\pi}\sum_{n=1}^{\infty}\frac{\log n}{n}\sin 2\pi nx$$

we see that

$$\sin 2\pi x\int_{0}^{1}\frac{\log\log(1/u)}{1-2u\cos 2\pi x+u^{2}}du=$$

$$-\frac{1}{2}\gamma\pi(1-2x)+\frac{1}{2}\pi\log\frac{\pi}{\sin\pi x}+\frac{1}{2}\pi(1-2x)[\gamma+\log(2\pi)]-\pi\log\Gamma(x)$$

We therefore obtain for $0<x<1$

(6.107rvi)

$$\sin 2\pi x\int_{0}^{1}\frac{\log\log(1/u)}{1-2u\cos 2\pi x+u^{2}}du=\frac{1}{2}\pi\log\frac{\pi}{\sin\pi x}+\frac{1}{2}\pi(1-2x)\log(2\pi)-\pi\log\Gamma(x)$$

It is easily seen that both sides of the above equation vanish at $x=1/2$. With $x=1/4$ we have

(6.107rvii) $$\int_{0}^{1}\frac{\log\log(1/u)}{1+u^{2}}du=\frac{3}{4}\pi\log\pi+\frac{1}{2}\pi\log 2-\pi\log\Gamma\left(\frac{1}{4}\right)$$

One is then immediately reminded of Adamchik's result (C.57) in Volume VI

$$\int_{0}^{1}\frac{u^{p-1}}{1+u^{n}}\log\log\left(\frac{1}{u}\right)du=\frac{\gamma+\log(2n)}{2n}\left[\psi\left(\frac{p}{2n}\right)-\psi\left(\frac{n+p}{2n}\right)\right]+\frac{1}{2n}\left[\varsigma'\left(1,\frac{p}{2n}\right)-\varsigma'\left(1,\frac{n+p}{2n}\right)\right]$$

where with $p=1$ and $n=2$ we get

$$\int_{0}^{1}\frac{\log\log(1/u)}{1+u^{2}}du=\frac{1}{4}[\gamma+2\log 2]\left[\psi\left(\frac{1}{4}\right)-\psi\left(\frac{3}{4}\right)\right]+\frac{1}{4}\left[\varsigma'\left(1,\frac{1}{4}\right)-\varsigma'\left(1,\frac{3}{4}\right)\right]$$

We have previously seen in (4.3.161c) in Volume II(a) that

$$\varsigma'\left(1,\frac{1}{4}\right)-\varsigma'\left(1,\frac{3}{4}\right)=\pi\left[\gamma+4\log 2+3\log\pi-4\log\Gamma\left(\frac{1}{4}\right)\right]$$

and it is well known that [126, p.20]

$$\psi\left(\frac{1}{4}\right)-\psi\left(\frac{3}{4}\right)=-\pi$$

and hence we obtain

$$\int_0^1 \frac{\log\log(1/u)}{1+u^2}du=-\frac{\pi}{4}[\gamma+2\log 2]+\frac{1}{4}\pi\left[\gamma+4\log 2+3\log\pi-4\log\Gamma\left(\frac{1}{4}\right)\right]$$

which is the same as (6.107rvii).

With $x=1/6$ in (6.107rvi) we get

$$\frac{\sqrt{3}}{2}\int_0^1 \frac{\log\log(1/u)}{1-u+u^2}du=\frac{1}{2}\pi\log(2\pi)+\frac{\pi}{3}\log(2\pi)-\pi\log\Gamma\left(\frac{1}{6}\right)$$

$$=\frac{5}{6}\pi\log(2\pi)-\pi\log\Gamma\left(\frac{1}{6}\right)$$

which is in agreement with Adamchik's paper [2a].

We may also write (6.107rii) as

$$\sum_{n=1}^{\infty}\frac{\cos nx}{n^s}=-\frac{1}{\Gamma(s)}\int_0^1 \frac{(u-\cos x)\log^{s-1}(1/u)}{1-2u\cos x+u^2}du$$

and then differentiate with respect to $s$ to obtain

$$-\sum_{n=1}^{\infty}\frac{\log n}{n^s}\cos nx=-\frac{1}{\Gamma(s)}\int_0^1 \frac{(u-\cos x)\log^{s-1}(1/u)\log\log(1/u)}{1-2u\cos x+u^2}du$$

$$+\frac{\psi(s)}{\Gamma(s)}\int_0^1 \frac{(u-\cos x)\log^{s-1}(1/u)}{1-2u\cos x+u^2}du$$

and this may be written as

(6.107rviii)

$$\frac{1}{\Gamma(s)}\int_0^1 \frac{(u-\cos x)\log^{s-1}(1/u)\log\log(1/u)}{1-2u\cos x+u^2}du=\sum_{n=1}^{\infty}\frac{\log n}{n^s}\cos nx-\psi(s)\sum_{n=1}^{\infty}\frac{\cos nx}{n^s}$$

For example, we may note that $\sum_{n=1}^{\infty} \frac{\log n}{n^2} \cos(n\pi/2)$ is known from (4.4.229j) in Volume IV.

Reference should also be made to the 2002 paper by Koyama and Kurokawa, "Kummer's formula for the multiple gamma functions" [93a] where they show that

(6.107rvix)
$$\log \Gamma_2(x) = -\frac{1}{2\pi^2} \sum_{n=1}^{\infty} \frac{\log n}{n^2} \cos 2\pi nx - \frac{\log(2\pi) + \gamma - 1}{2\pi^2} \sum_{n=1}^{\infty} \frac{\cos 2\pi nx}{n^2}$$

$$+ \frac{1}{4\pi} \sum_{n=1}^{\infty} \frac{\sin 2\pi nx}{n^2} + (1-x) \log \Gamma_1(x)$$

(6.107rvx)
$$\log \Gamma_3(x) = -\frac{1}{4\pi^3} \sum_{n=1}^{\infty} \frac{\log n}{n^3} \sin 2\pi nx - \frac{2\log(2\pi) + 2\gamma - 3}{8\pi^3} \sum_{n=1}^{\infty} \frac{\sin 2\pi nx}{n^3}$$

$$+ \frac{1}{8\pi^2} \sum_{n=1}^{\infty} \frac{\cos 2\pi nx}{n^3} + \left(\frac{3}{2} - x\right) \log \Gamma_2(x) - \frac{1}{2}(1-x)^2 \log \Gamma_1(x)$$

where $\Gamma_1(x) = \frac{\Gamma(x)}{\sqrt{2\pi}}$. These equations will then enable us to evaluate the integral in (6.107rviii) when $s = 2$ and the integral in (6.107rv) when $s = 3$.

We could also differentiate (6.107rvi) with respect to $x$ to gives us the integral

$$-4\pi \sin^2 2\pi x \int_0^1 \frac{u \log \log(1/u)}{[1 - 2u \cos 2\pi x + u^2]^2} du + 2\pi \cos 2\pi x \int_0^1 \frac{\log \log(1/u)}{1 - 2u \cos 2\pi x + u^2} du =$$

$$\frac{1}{2} \pi^2 \cot \pi x - \pi \log(2\pi) - \pi \psi(x)$$

and we note that Adamchik [2a] has also evaluated integrals such as

$$\int_0^1 \frac{u \log \log(1/u)}{[1 + u^2]^2} du$$

Differentiating (6.107riii) with respect to $x$ to gives us

$$\sum_{n=1}^{\infty} \frac{\cos nx}{n^{s-1}} = -\frac{2\sin^2 x}{\Gamma(s)} \int_0^1 \frac{u \log^{s-1}(1/u)}{[1 - 2u \cos x + u^2]^2} du + \frac{\cos x}{\Gamma(s)} \int_0^1 \frac{\log^{s-1}(1/u)}{1 - 2u \cos x + u^2} du$$

In addition, we could also integrate (6.107rv) with respect to $x$.

$\square$

Chen and Khalili [43d] have recently shown that

(6.107s) $\quad B_{2n+2}(\theta) = (-1)^{n+1} \dfrac{(2n+1)(2n+2)}{(2\pi)^{2n+2}} \int_0^1 \dfrac{\log[1 - 2x\cos(2\pi\theta) + x^2] \log^{2n} x}{x} dx$

and differentiation results in

(6.107t) $\quad B_{2n+1}(\theta) = (-1)^{n+1} \dfrac{(2n+1)\sin(2\pi\theta)}{(2\pi)^{2n+1}} \int_0^1 \dfrac{\log^{2n} x}{1 - 2x\cos(2\pi\theta) + x^2} dx$

We note from (6.107q) that

$$\operatorname{Cl}_{2n}(2\pi\theta) = -\dfrac{\sin(2\pi\theta)}{(2n-1)!} \int_0^1 \dfrac{\log^{2n-1} x}{1 - 2x\cos(2\pi\theta) + x^2} dx$$

Cvijović [49a] has recently considered related integrals for the Bernoulli and Euler polynomials. See also another recent paper by Cvijović [49b].

We have from (6.36)

$$\varsigma(2n+1) = (-1)^{n+1} \dfrac{(2\pi)^{2n+1}}{2(2n+1)!} \int_0^1 B_{2n+1}(x) \cot(\pi x/2) dx$$

and therefore

$$\varsigma(3) = \dfrac{(2\pi)^3}{12} \int_0^1 B_3(x) \cot(\pi x/2) dx$$

Using (6.107t) we have

$$B_3(x) = \dfrac{3\sin(2\pi x)}{(2\pi)^3} \int_0^1 \dfrac{\log t}{1 - 2t\cos(2\pi x) + t^2} dt$$

and therefore we obtain

$$\varsigma(3) = \dfrac{1}{4} \int_0^1 \sin(2\pi x) \int_0^1 \dfrac{\log t}{1 - 2t\cos(2\pi x) + t^2} \cot(\pi x/2) dx dt$$

$$= \int_0^1 \int_0^1 \dfrac{\sin(\pi x/2)\cos^2(\pi x/2)\cos(\pi x)\log t}{1 - 2t\cos(2\pi x) + t^2} dx dt$$

Mathematica may be able to evaluate the integral in $x$. Alternatively we have

$$\int \frac{\log t}{1-2t\cos(2\pi x)+t^2} dt = \frac{\log t}{2i\sin(2\pi x)} \left( \log\left[1-\frac{t}{\exp(2i\pi x)}\right] - \log\left[\frac{t}{-\exp(2i\pi x)}-1\right] \right)$$

$$-\frac{1}{2i\sin(2\pi x)} \left( Li_2\left[1-\frac{t}{\exp(2i\pi x)}\right] - Li_2\left[\frac{t}{-\exp(2i\pi x)}-1\right] \right)$$

$$= \frac{\log t}{2i\sin(2\pi x)} \log\frac{[t\exp(-2i\pi x)-1]}{[t\exp(-2i\pi x)+1]}$$

$$-\frac{1}{2i\sin(2\pi x)} \left( Li_2\left[1-\frac{t}{\exp(2i\pi x)}\right] - Li_2\left[\frac{t}{-\exp(2i\pi x)}-1\right] \right)$$

As mentioned in Cvijović's paper [49b], the Chebyshev polynomials defined by the following generating functions may also be usefully employed

$$T_n(x) = \cos(n\cos^{-1} x) \qquad T_n(\cos\theta) = \cos n\theta$$

$$U_n(\cos\theta) = \frac{\sin(n+1)\theta}{\sin\theta}$$

$$\frac{1-t^2}{1-2xt+t^2} = 1 + 2\sum_{n=1}^{\infty} T_n(x)t^n$$

$$\frac{1}{1-2xt+t^2} = \sum_{n=0}^{\infty} U_n(x)t^n$$

$$\frac{1-t^2}{1-2t\cos(2\pi x)+t^2} = 1 + \frac{2t(\cos(2\pi x)-1)}{1-2t\cos(2\pi x)+t^2}$$

**Example 15:**

The well-known Maclaurin expansion for the log gamma function is derived in (E.22n) of Volume VI

(6.108) $$\log\Gamma(1+x) = -\gamma x + \sum_{n=2}^{\infty} (-1)^n \frac{\varsigma(n)}{n} x^n, \quad -1 < x \leq 1$$

With reference to (6.14b) we have

(6.109) $$\int_a^b p(x)\cot\alpha x\, dx = 2\sum_{n=1}^{\infty} \int_a^b p(x)\sin 2n\alpha x$$

As mentioned previously, it should be noted that in the above formulae we require either (i) both $\sin(\alpha x/2)$ and $\cos(\alpha x/2)$ have no zero in $[a,b]$ or (ii) if either $\sin(\alpha a/2)$ or $\cos(\alpha a/2)$ is equal to zero then $p(a)$ must also be zero. Condition (i) is equivalent to the requirement that $\sin \alpha x$ has no zero in $[a,b]$.

Since $\log \Gamma(1) = \log \Gamma(2) = 0$, it is clear that $p(x) = \log \Gamma(x+1)$ satisfies the necessary conditions for the Riemann-Lebesgue lemma on the interval $[0,1]$ with $\alpha = \pi$. Accordingly, substituting (6.108) in (6.109) we have

$$(6.110) \quad \int_0^1 \left[ -\gamma x + \sum_{n=2}^{\infty} (-1)^n \frac{\varsigma(n)}{n} x^n \right] \cot \pi x \, dx = 2 \sum_{n=1}^{\infty} \int_0^1 \log \Gamma(x+1) \sin 2n\pi x \, dx$$

and

$$\int_0^1 \log \Gamma(x+1) \cot \pi x \, dx = 2 \sum_{n=1}^{\infty} \int_0^1 \log \Gamma(x+1) \sin 2n\pi x \, dx$$

The corresponding result for (6.5) is shown in (6.120) below.

The following integral is recorded in G&R [74, p.650] for $a > 0$

$$\int_0^1 \log \Gamma(x+a) \sin 2n\pi x \, dx = -\frac{1}{2\pi n} \left[ \log a + \cos(2n\pi a) Ci(2n\pi a) - \sin(2n\pi a) si(2n\pi a) \right]$$

but I believe that the correct version should be

(6.111)
$$\int_0^1 \log \Gamma(x+a) \sin 2n\pi x \, dx = -\frac{1}{2\pi n} \left[ \log a - \cos(2n\pi a) Ci(2n\pi a) + \sin(2n\pi a) si(2n\pi a) \right]$$

and with $a = 1$ this implies that

$$(6.111a) \quad \int_0^1 \log \Gamma(x+1) \sin 2n\pi x \, dx = \frac{Ci(2n\pi)}{2\pi n}$$

where $Ci(x)$ is the cosine integral defined for $x > 0$ in G&R [74, p.878] by

$$Ci(x) = -\int_x^{\infty} \frac{\cos t}{t} dt = \gamma + \log x + \int_0^x \frac{\cos t - 1}{t} dt$$

We therefore obtain from (6.109)

$$(6.111b) \quad \int_0^1 \log \Gamma(x+1) \cot \pi x \, dx = 2 \sum_{n=1}^{\infty} \int_0^1 \log \Gamma(x+1) \sin 2n\pi x \, dx = \frac{1}{\pi} \sum_{n=1}^{\infty} \frac{Ci(2n\pi)}{n}$$

An alternative proof is shown below. We have

$$\int_0^1 \log \Gamma(x+1)\sin 2n\pi x\, dx = \int_0^1 \log x \sin 2n\pi x\, dx + \int_0^1 \log \Gamma(x)\sin 2n\pi x\, dx$$

and in (E.46) in Volume VI we show that (G&R [74, p.650])

(6.111bi) $\quad \int_0^1 \log \Gamma(x)\sin 2n\pi x\, dx = \dfrac{\gamma + \log 2n\pi}{2\pi n}$

Since from (6.94f)

$$\int_0^1 \log x \sin 2n\pi x\, dx = \dfrac{Ci(2n\pi)}{2\pi n} - \dfrac{\gamma + \log 2n\pi}{2\pi n}$$

we have

(6.111c) $\quad \int_0^1 \log \Gamma(x+1)\sin 2n\pi x\, dx = \dfrac{Ci(2n\pi)}{2\pi n}$

in agreement with (6.111a).

It therefore appears that two of the signs in (6.111) are recorded incorrectly in G&R [74, p.650] (and also in "Integrals and Series", Volume 2, p.60 by Prudnikov et al). In this regard, I note that both Havil [78, p.126] and Elizalde [58c] define $Ci(x)$ as the negative of (6.94g).

Equation (6.111) also applies in the limit $a \to 0$ because as we shall see in (6.117ei)

$$\lim_{y \to 0}[\cos y\, Ci(y) - \log y] = \gamma$$

and we therefore obtain another proof of (6.111bi).

In passing, we note that integration by parts gives us

$$\int_0^1 \log \Gamma(x+1)\cot \pi x\, dx = \dfrac{1}{\pi}\log \Gamma(x+1)\log \sin \pi x \Big|_0^1 - \dfrac{1}{\pi}\int_0^1 \psi(1+x)\log \sin \pi x\, dx$$

We note from L'Hôpital's rule that

$$\lim_{x \to 1}[\log \Gamma(x+1)\log \sin \pi x] = \lim_{x \to 1}\left[\dfrac{\log \sin \pi x}{1/\log \Gamma(x+1)}\right]$$

$$= \lim_{x \to 1}\left[-\frac{\pi \cot \pi x \log^2 \Gamma(x+1)}{\psi(x+1)}\right]$$

$$= \lim_{x \to 1}\left[-\frac{\pi \log^2 \Gamma(x+1)}{\psi(2)\sin \pi x}\right]$$

We also have

$$\lim_{x \to 1}\left[\frac{\log \Gamma(x+1)}{\sin \pi x}\log \Gamma(x+1)\right] = \lim_{x \to 1}\left[\frac{\psi(x+1)}{\pi \cos \pi x}\right]\lim_{x \to 1}\log \Gamma(x+1) = 0$$

Similarly we have

$$\lim_{x \to 0}\left[\frac{\log \Gamma(x+1)}{\sin \pi x}\log \Gamma(x+1)\right] = \lim_{x \to 0}\left[\frac{\pi x}{\sin \pi x}\frac{\log \Gamma(x+1)}{\pi x}\log \Gamma(x+1)\right]$$

and since $\lim_{x \to 0}\left[\frac{\log \Gamma(x+1)}{\pi x}\right] = \lim_{x \to 0}\left[\frac{\psi(x+1)}{\pi}\right]$ we obtain $\lim_{x \to 0}[\log \Gamma(x+1)\log \sin \pi x] = 0$

Hence we have

$$\int_0^1 \log \Gamma(x+1)\cot \pi x\, dx = -\frac{1}{\pi}\int_0^1 \psi(1+x)\log \sin \pi x\, dx$$

Completing the summation of (6.111c) we obtain

$$\sum_{n=1}^{\infty}\frac{1}{n}\int_0^1 \log \Gamma(x+1)\sin 2n\pi x\, dx = \frac{1}{2\pi}\sum_{n=1}^{\infty}\frac{Ci(2n\pi)}{n^2}$$

Using (7.5) we get

$$\sum_{n=1}^{\infty}\frac{1}{n}\int_0^1 \log \Gamma(x+1)\sin 2n\pi x\, dx = \int_0^1 \log \Gamma(x+1)\sum_{n=1}^{\infty}\frac{\sin 2n\pi x}{n}\, dx$$

$$= \frac{\pi}{2}\int_0^1 (1-2x)\log \Gamma(x+1)\, dx$$

$$= \frac{\pi}{2}\int_0^1 \log \Gamma(x+1)\, dx - \pi\int_0^1 x\log \Gamma(x+1)\, dx$$

In the above we have tacitly assumed that interchanging the order of integration and summation is valid and that the Fourier series expansion (7.5) may be validly used at both end points because the multiplication factor $\log \Gamma(x+1)$ is zero at both end points.

We have from [126, p.32]

(6.112) $\displaystyle\int_0^z \log \Gamma(1+t)\,dt = \frac{1}{2}[\log(2\pi)-1]z - \frac{z^2}{2} + z\log\Gamma(1+z) - \log G(1+z)$

and integration by parts readily gives us

$$\int_0^x t\log\Gamma(1+t)\,dt = t\left[\frac{1}{2}[\log(2\pi)-1]t - \frac{t^2}{2} + t\log\Gamma(1+t) - \log G(1+t)\right]\Big|_0^x$$

$$-\int_0^x \left[\frac{1}{2}[\log(2\pi)-1]t - \frac{t^2}{2} + t\log\Gamma(1+t) - \log G(1+t)\right]dt$$

Hence we obtain

(6.112a) $\displaystyle 2\int_0^x t\log\Gamma(1+t)\,dt = \frac{1}{4}[\log(2\pi)-1]x^2 - \frac{1}{3}x^3 + x^2\log\Gamma(1+x)$

$$-x\log G(1+x) + \int_0^x \log G(1+t)\,dt$$

Therefore we get

(6.113) $\displaystyle 2\int_0^1 t\log\Gamma(1+t)\,dt = \frac{1}{4}[\log(2\pi)-1] - \frac{1}{3} - \log G(2) + \int_0^1 \log G(1+t)\,dt$

From (6.57) we have $G(2) = \Gamma(1)G(1) = 1$ and Srivastava and Choi [126, p.217] report that

(6.114)
$$\int_0^1 \log G(a+t)\,dt = \frac{a(1-a)}{2} + \frac{a}{2}\log 2\pi + a\log\Gamma(a) - \log G(a+1) + \log\left[2^{-\frac{1}{4}}\pi^{-\frac{1}{4}}e^{\frac{1}{12}}A^{-2}\right]$$

and hence, as originally discovered by Barnes [126, p.37], we have

(6.115)
$$\int_0^1 \log G(1+t)\,dt = \frac{1}{2}\log 2\pi + \log\left[2^{-\frac{1}{4}}\pi^{-\frac{1}{4}}e^{\frac{1}{12}}A^{-2}\right] = \frac{1}{12} + \frac{1}{4}\log 2\pi - 2\log A$$

This results in

(6.116) $$2\int_0^1 t\log\Gamma(1+t)\,dt = \frac{1}{2}[\log(2\pi)-1] - 2\log A$$

This was also obtained by Espinosa and Moll [59].

I subsequently discovered that (6.116) is a particular case of the more general formula (6.126) for $\int_0^z t\log\Gamma(a+t)\,dt$ given by Choi and Srivastava [45ab] in 2000 in terms of the multiple gamma functions.

(6.116a) $$2\int_0^z t\log\Gamma(1+t)\,dt = \left(\frac{1}{4} - 2\log A\right)z + \left(\frac{1}{2}\log(2\pi) - \frac{1}{4}\right)z^2$$
$$-\frac{1}{2}z^2 + z^2\log\Gamma(1+z) - \log G(1+z) - 2\log\Gamma_3(1+z)$$

Using the formula (6.112) we also obtain (see also (C.43b) in Volume VI)

$$\int_0^1 \log\Gamma(1+t)\,dt = \int_0^1 \log[t\Gamma(t)]\,dt = \int_0^1 \log t\,dt + \int_0^1 \log\Gamma(t)\,dt = -1 + \frac{1}{2}\log(2\pi)$$

and we finally deduce that

(6.117) $$\sum_{n=1}^\infty \frac{Ci(2n\pi)}{n^2} = 2\pi^2\left(\log A - \frac{1}{4}\right) = -2\pi^2\left(\varsigma'(-1) + \frac{1}{6}\right)$$

where in the last part we have used (4.4.225). The set of series $\sum_{n=1}^\infty \frac{Ci(2n\pi)}{n^{2N}}$ may also be evaluated in a similar way.

We have from (6.94g)

$$Ci(2n\pi) = \gamma + \log 2n\pi + \int_0^{2\pi} \frac{\cos nu - 1}{u}\,du = \gamma + \log 2n\pi + n\int_0^{2\pi} \sin nu \log u\,du$$

and therefore

$$\sum_{n=1}^\infty \frac{Ci(2n\pi)}{n^2} = [\gamma + \log 2\pi]\varsigma(2) - \varsigma'(2) + \int_0^{2\pi} \sum_{n=1}^\infty \frac{\sin nu}{n}\log u\,du$$

Unfortunately, the Fourier series $\sum_{n=1}^\infty \frac{\sin nu}{n} = \frac{\pi - u}{2}$ is only valid in the interval $(0, 2\pi)$ and we cannot proceed any further in this direction.

□

Elizalde [58c] reported in 1985 that for $a > 0$

(6.117a) $\varsigma'(-1,a) =$

$$-\varsigma(-1,a)\log a - \frac{1}{4}a^2 + \frac{1}{12} - \frac{1}{2\pi^2}\sum_{n=1}^{\infty}\frac{1}{n^2}[\cos(2n\pi a)Ci(2n\pi a) + \sin(2n\pi a)si(2n\pi a)]$$

where $\varsigma(s,a)$ is the Hurwitz zeta function. Since $si(x) = Si(x) - \frac{\pi}{2}$ this may be written as

(6.117ai) $\varsigma'(-1,a) = -\varsigma(-1,a)\log a - \frac{1}{4}a^2 + \frac{1}{12} + \frac{1}{4\pi}\sum_{n=1}^{\infty}\frac{\sin(2n\pi a)}{n^2}$

$$-\frac{1}{2\pi^2}\sum_{n=1}^{\infty}\frac{1}{n^2}[\cos(2n\pi a)Ci(2n\pi a) + \sin(2n\pi a)Si(2n\pi a)]$$

and we therefore have with $a = 1$

$$\varsigma'(-1,1) = \varsigma'(-1) = -\frac{1}{6} - \frac{1}{2\pi^2}\sum_{n=1}^{\infty}\frac{Ci(2n\pi)}{n^2}$$

which concurs with (6.117). I was quite pleased to come across this reference in 2007 since I had not previously confronted any similar analysis during my limited research.

Readers should note that both of the functions $Si(x)$ and $si(x)$ are employed in the various formulae set out in the following section.

With $a = 1/2$ we get

$$\varsigma'\left(-1,\frac{1}{2}\right) = \frac{1}{24}\log 2 + \frac{1}{48} - \frac{1}{2\pi^2}\sum_{n=1}^{\infty}(-1)^n\frac{Ci(n\pi)}{n^2}$$

and from (4.3.140) we have

$$\varsigma'\left(-1,\frac{1}{2}\right) = -\frac{1}{24}\log 2 - \frac{1}{2}\varsigma'(-1)$$

This then gives us (see also (6.130))

(6.117b) $\frac{1}{2\pi^2}\sum_{n=1}^{\infty}(-1)^n\frac{Ci(n\pi)}{n^2} = \frac{1}{12}\log 2 + \frac{1}{48} + \frac{1}{2}\varsigma'(-1)$

Elizalde [58c] gave little indication of the source of the identity but differentiation of (6.117a) sheds more light on the subject: we have

$$\frac{\partial}{\partial a}\varsigma'(-1,a) = \left(a-\frac{1}{2}\right)\log a + \frac{1}{2}\left(a^2 - a + \frac{1}{6}\right)\frac{1}{a} - \frac{1}{2}a - \frac{1}{2\pi^2}\frac{\varsigma(2)}{a}$$

$$-\frac{1}{\pi}\sum_{n=1}^{\infty}\frac{1}{n}[-\sin(2n\pi a)Ci(2n\pi a) + \cos(2n\pi a)si(2n\pi a)]$$

since $\frac{d}{dx}Ci(x) = \frac{\cos x}{x}$ and $\frac{d}{dx}si(x) = \frac{\sin x}{x}$. This simplifies to

$$\frac{\partial}{\partial a}\varsigma'(-1,a) = \left(a-\frac{1}{2}\right)\log a - \frac{1}{2} + \frac{1}{\pi}\sum_{n=1}^{\infty}\frac{1}{n}[\sin(2n\pi a)Ci(2n\pi a) - \cos(2n\pi a)si(2n\pi a)]$$

We have previously seen in the analysis following (4.3.117) that

$$\frac{\partial}{\partial a}\frac{\partial}{\partial s}\varsigma(s,a) = -\varsigma(s+1,a) - s\varsigma'(s+1,a)$$

and hence

$$\frac{\partial}{\partial a}\varsigma'(-1,a) = -\varsigma(0,a) + \varsigma'(0,a)$$

Then using Lerch's identity (4.3.116) in Volume II(a)

$$\varsigma'(0,a) = \log\Gamma(a) - \frac{1}{2}\log(2\pi)$$

this becomes

$$\frac{\partial}{\partial a}\varsigma'(-1,a) = -\varsigma(0,a) + \log\Gamma(a) - \frac{1}{2}\log(2\pi)$$

We then obtain

$$-\varsigma(0,a) + \log\Gamma(a) - \frac{1}{2}\log(2\pi)$$

$$= \left(a-\frac{1}{2}\right)\log a - \frac{1}{2} + \frac{1}{\pi}\sum_{n=1}^{\infty}\frac{1}{n}[\sin(2n\pi a)Ci(2n\pi a) - \cos(2n\pi a)si(2n\pi a)]$$

which, for $0 < a < 1$, simplifies to

(6.117c) $\quad \log\Gamma(a) =$

$$\frac{1}{2}\log(2\pi) + \left(a-\frac{1}{2}\right)\log a - a + \frac{1}{\pi}\sum_{n=1}^{\infty}\frac{1}{n}[\sin(2n\pi a)Ci(2n\pi a) - \cos(2n\pi a)si(2n\pi a)]$$

$$= \frac{1}{2}\log(2\pi) + \left(a - \frac{1}{2}\right)\log a - a + \frac{1}{2}\sum_{n=1}^{\infty}\frac{\cos(2n\pi a)}{n}$$

$$+ \frac{1}{\pi}\sum_{n=1}^{\infty}\frac{1}{n}[\sin(2n\pi a)Ci(2n\pi a) - \cos(2n\pi a)Si(2n\pi a)]$$

and using (7.8a) this becomes

$$\log \Gamma(a) = \frac{1}{2}\log(2\pi) + \left(a - \frac{1}{2}\right)\log a - a - \frac{1}{2}\log[2\sin(\pi a)]$$

$$+ \frac{1}{\pi}\sum_{n=1}^{\infty}\frac{1}{n}[\sin(2n\pi a)Ci(2n\pi a) - \cos(2n\pi a)Si(2n\pi a)]$$

The formula (6.117c) was given by Nörlund in [105, p.114]. When $a = 1/2$ we obtain

$$\frac{\pi}{2}(1 - \log 2) = \sum_{n=1}^{\infty}\frac{(-1)^n}{n}si(n\pi)$$

which is a particular case of (6.94a). With $a = 1/4$ we get

$$\log \Gamma\left(\frac{1}{4}\right) = \frac{1}{2}\log(2\pi) + \frac{1}{2}\log 2 - \frac{1}{4} - \frac{1}{4}\log 2$$

$$+ \frac{1}{\pi}\sum_{n=1}^{\infty}\frac{1}{n}[\sin(n\pi/2)Ci(n\pi/2) - \cos(n\pi/2)Si(n\pi/2)]$$

The following is posed as a question in Whittaker & Watson [135, p.261]: Prove that for all values of $a$ except negative real values

(6.117ci)  $$\log \Gamma(a) = \frac{1}{2}\log(2\pi) + \left(a - \frac{1}{2}\right)\log a - a + \frac{1}{\pi}\sum_{n=1}^{\infty}\int_0^{\infty}\frac{\sin(2n\pi x)}{n(x+a)}dx$$

and this was attributed by Stieltjes to Bourguet. Equation (6.117ci) may also be derived using the Euler-Maclaurin summation formula (see in particular Knopp's book [90, p.530]).

By differentiation we can easily see that

$$\frac{d}{dx}\left(\cos(2n\pi a)Si[2n\pi(x+a)] - \sin(2n\pi a)Ci[2n\pi(x+a)]\right) = \frac{\sin(2n\pi x)}{x+a}$$

and we therefore have

$$\int_0^M \frac{\sin(2n\pi x)}{x+a} dx = \left(\cos(2n\pi a) Si[2n\pi(x+a)] - \sin(2n\pi a) Ci[2n\pi(x+a)]\right)\Big|_0^M$$

$$= \cos(2n\pi a)\{Si[2n\pi(M+a)] - Si[2n\pi a]\} - \sin(2n\pi a)\{Ci[2n\pi(M+a)] - Ci[2n\pi a]\}$$

From (6.90b) and (6.90d) we see that

$$\lim_{M\to\infty} Si[2n\pi(M+a)] = \frac{\pi}{2}$$

and from (6.9ga) we have

$$\lim_{M\to\infty} Ci[2n\pi(M+a)] = 0$$

Hence we obtain as $M \to \infty$

$$\int_0^\infty \frac{\sin(2n\pi x)}{x+a} dx = \cos(2n\pi a)\left\{\frac{\pi}{2} - Si(2n\pi a)\right\} + \sin(2n\pi a) Ci(2n\pi a)$$

and reference to (6.90e) shows that this is equal to

$$= -\cos(2n\pi a) si(2n\pi a) + \sin(2n\pi a) Ci(2n\pi a)$$

Therefore from Bourguet's formula we have

$$\log \Gamma(a) =$$

$$\frac{1}{2}\log(2\pi) + \left(a - \frac{1}{2}\right)\log a - a + \frac{1}{\pi}\sum_{n=1}^\infty \frac{1}{n}[\sin(2n\pi a) Ci(2n\pi a) - \cos(2n\pi a) si(2n\pi a)]$$

which we have already seen in (6.117c).

This may be expressed as

$$\log \Gamma(a) = \frac{1}{2}\log(2\pi) + \left(a - \frac{1}{2}\right)\log a - a + \sum_{n=1}^\infty C_n$$

and we have with $a \to 2a$

$$\log \Gamma(2a) = \frac{1}{2}\log(2\pi) + \left(2a - \frac{1}{2}\right)\log(2a) - 2a + 2\sum_{n=1}^\infty C_{2n}$$

Since $2\sum_{n=1}^\infty C_{2n} = \sum_{n=1}^\infty C_n + \sum_{n=1}^\infty (-1)^n C_n$ we have

$$\log \Gamma(2a) = \left(2a - \frac{1}{2}\right)\log(2a) - a + \log \Gamma(a) - \left(a - \frac{1}{2}\right)\log a + \sum_{n=1}^{\infty}(-1)^n C_n$$

and using Legendre's duplication formula for the gamma function [126, p.7]

$$\sqrt{\pi}\,\Gamma(2a) = 2^{2a-1}\Gamma(a)\Gamma\left(a + \frac{1}{2}\right)$$

this results in

(6.117ca)    $\log \Gamma(a + 1/2) =$

$$\frac{1}{2}\log(2\pi) + a\log a - a + \frac{1}{\pi}\sum_{n=1}^{\infty}\frac{(-1)^n}{n}[\sin(2n\pi a)Ci(2n\pi a) - \cos(2n\pi a)si(2n\pi a)]$$

which is also reported by Nörlund [105, p.114]. Since $si(0) = -\pi/2$ this identity may be easily verified for $a = 0$. With $a = 1/2$ we obtain

$$\sum_{n=1}^{\infty}\frac{si(n\pi)}{n} = \frac{\pi}{2}\log \pi - \frac{\pi}{2}$$

which is contained in [104a, p,82]. Letting $a = 1$ we have

$$\sum_{n=1}^{\infty}\frac{si(2n\pi)}{n} = \frac{3}{2}\pi \log 2 - \pi$$

We also have using (6.117c)

(6.117cai)

$$\log \Gamma(2a) = \left(2a - \frac{1}{2}\right)\log 2 + \log \Gamma(a) + a\log\left(a + \frac{1}{2}\right) - a - \frac{1}{2}$$

$$+ \frac{1}{\pi}\sum_{n=1}^{\infty}\frac{(-1)^n}{n}[\sin(2n\pi a)Ci(2n\pi a + n\pi) - \cos(2n\pi a)si(2n\pi a + n\pi)]$$

and with $a = 1/2$ we obtain (6.120).

Using (7.8) for $0 < a < 1$ (6.117c) may be written as

(6.117ci)    $\log \Gamma(a) =$

$$= \frac{1}{2}\log(2\pi) + \left(a - \frac{1}{2}\right)\log a - a - \frac{1}{2}\log[2\sin(\pi a)]$$

$$+\frac{1}{\pi}\sum_{n=1}^{\infty}\frac{1}{n}[\sin(2n\pi a)Ci(2n\pi a)-\cos(2n\pi a)Si(2n\pi a)]$$

This may be written more compactly as

$$\log\left[\frac{\Gamma^2(a)\sin(\pi a)}{\pi a^{2a-1}}\right]=\frac{2}{\pi}\sum_{n=1}^{\infty}\frac{1}{n}[\sin(2n\pi a)Ci(2n\pi a)-\cos(2n\pi a)Si(2n\pi a)]$$

With $a=1$ we get

(6.117cii) $$\sum_{n=1}^{\infty}\frac{si(2n\pi)}{n}=\frac{\pi}{2}\log(2\pi)-\pi$$

which is given by Nielsen [104a, p.79], and with $a=1/2$ we get (6.94a)

(6.117ciii) $$\sum_{n=1}^{\infty}(-1)^n\frac{si(n\pi)}{n}=\frac{\pi}{2}\log 2-\frac{\pi}{2}$$

which is a particular case of Nielsen's formula [104a, p.83] for $x=\pi$

(6.117civ) $$\sum_{n=1}^{\infty}(-1)^n\frac{si(nx)}{n}=\frac{\pi}{2}\log 2-\frac{x}{2}$$

$$\sum_{n=1}^{\infty}(-1)^n\frac{Si(nx)}{n}=\pi\log 2-\frac{x}{2}$$

Note that Nielsen indicated that (6.117civ) was only valid for $x\in(-\pi,\pi)$.

Letting

$$h(a)=\sin(2n\pi a)Ci(2n\pi a)-\cos(2n\pi a)Si(2n\pi a)$$

we see that

$$h'(a)=\sin(2n\pi a)\frac{Ci(2n\pi a)}{a}+2n\pi\cos(2n\pi a)Ci(2n\pi a)$$

$$-\cos(2n\pi a)\frac{Si(2n\pi a)}{a}+2n\pi\sin(2n\pi a)Si(2n\pi a)$$

We therefore have

$$h'(1)=2n\pi\,Ci(2n\pi)-si(2n\pi)$$

$$h'(1/2)=2n\pi(-1)^n Ci(n\pi)-2(-1)^n si(2n\pi)$$

and by differentiating (6.117ca) we see that

$$\psi\left(\frac{3}{2}\right) = \frac{1}{\pi}\sum_{n=1}^{\infty}\frac{(-1)^n}{n}[2n\pi\, Ci(2n\pi) - si(2n\pi)]$$

$$= 2\sum_{n=1}^{\infty}(-1)^n Ci(2n\pi) - \frac{1}{\pi}\sum_{n=1}^{\infty}\frac{(-1)^n si(2n\pi)}{n}$$

Therefore, using (6.121)

$$\sum_{n=1}^{\infty}(-1)^n \frac{si(2n\pi)}{n} = \frac{\pi}{2}\log 2 - \pi$$

and $\psi\left(\dfrac{3}{2}\right) = 2 - \gamma - 2\log 2$, we obtain

$$\sum_{n=1}^{\infty}(-1)^n Ci(2n\pi) = \frac{1}{2}2 - \frac{1}{2}\gamma - \frac{3}{4}\log 2$$

With $a = 1/2$ we get

$$\psi(1) = -\gamma = -\log 2 + 2\sum_{n=1}^{\infty}Ci(n\pi) - \frac{2}{\pi}\sum_{n=1}^{\infty}\frac{si(n\pi)}{n}$$

Then referring to (6.117ca) we obtain (see also (6.117s))

$$2\sum_{n=1}^{\infty}Ci(n\pi) = \log(2\pi) - \gamma - 1$$

□

Integrating (6.117civ) results in

(6.117cv) $\qquad x\sum_{n=1}^{\infty}(-1)^n \dfrac{si(nx)}{n} + \sum_{n=1}^{\infty}(-1)^n \dfrac{\cos(nx)-1}{n^2} = \dfrac{\pi x}{2}\log 2 - \dfrac{1}{4}x^2$

or equivalently

$$x\sum_{n=1}^{\infty}(-1)^n \frac{Si(nx)}{n} + \sum_{n=1}^{\infty}(-1)^n \frac{\cos(nx)-1}{n^2} = \pi x\log 2 - \frac{1}{4}x^2$$

where we have used

$$\int Si(ax)dx = xSi(ax) + \frac{1}{a}\cos x$$

Multiplying (6.117civ) by $x$ and subtracting (6.117cv) gives us

$$\sum_{n=1}^{\infty}(-1)^n \frac{\cos(nx)-1}{n^2} = -\frac{1}{2}x^2$$

and this may be written as

$$\sum_{n=1}^{\infty}(-1)^{n+1} \frac{\cos(nx)}{n^2} = \frac{\pi^2 - 3x^2}{12}$$

which is contained in [130, p.148]. Integrating (6.117cv) results in

$$\frac{1}{2}x^2 \sum_{n=1}^{\infty}(-1)^n \frac{Si(nx)}{n} + \frac{1}{2}x\sum_{n=1}^{\infty}(-1)^n \frac{\cos(nx)}{n^2} - \frac{1}{2}\sum_{n=1}^{\infty}(-1)^n \frac{\sin(nx)}{n^3} + \sum_{n=1}^{\infty}(-1)^n \frac{\sin(nx)}{n^3} - x\varsigma_a(2)$$

$$= \frac{\pi x^2}{2}\log 2 - \frac{1}{12}x^3$$

where we have used

$$\int x\, Si(ax)\, dx = \frac{1}{2}x^2 Si(ax) + \frac{x\cos ax}{2a} - \frac{\sin ax}{2a^2}$$

$$\frac{1}{2}x^2 \sum_{n=1}^{\infty}(-1)^n \frac{Si(nx)}{n} + \frac{1}{2}x\sum_{n=1}^{\infty}(-1)^n \frac{\cos(nx)}{n^2} + \frac{1}{2}\sum_{n=1}^{\infty}(-1)^n \frac{\sin(nx)}{n^3} - x\varsigma_a(2)$$

$$= \frac{\pi x^2}{2}\log 2 - \frac{1}{12}x^3$$

We then have

$$\frac{1}{2}x^2\left[\pi \log 2 - \frac{x}{2}\right] + \frac{1}{2}x\sum_{n=1}^{\infty}(-1)^n \frac{\cos(nx)}{n^2} + \frac{1}{2}\sum_{n=1}^{\infty}(-1)^n \frac{\sin(nx)}{n^3} - x\varsigma_a(2)$$

$$= \frac{\pi x^2}{2}\log 2 - \frac{1}{12}x^3$$

This becomes

$$\frac{1}{2}x\sum_{n=1}^{\infty}(-1)^n \frac{\cos(nx)}{n^2} + \frac{1}{2}\sum_{n=1}^{\infty}(-1)^n \frac{\sin(nx)}{n^3} - x\varsigma_a(2) = \frac{1}{6}x^3$$

and this is a simple consequence of the Fourier series contained in [130, p.148].

More generally we have the Fourier sine and cosine transforms [58d] for $|\arg a| < \pi$ and $y > 0$

$$\int_0^\infty \frac{\sin(xy)}{x+a} dx = -\cos(ay)si(ay) + \sin(ay)Ci(ay)$$

$$\int_0^\infty \frac{\cos(xy)}{x+a} dx = -\sin(ay)si(ay) - \cos(ay)Ci(ay)$$

and we may make the following summation

$$\sum_{a=1}^\infty \frac{1}{a} \int_0^\infty \frac{\sin(xy)}{x+a} dx = \sum_{a=1}^\infty \frac{1}{a} [-\cos(ay)si(ay) + \sin(ay)Ci(ay)]$$

We have from (E.14) in Volume VI

$$\sum_{a=1}^\infty \frac{1}{a(x+a)} = \frac{\psi(x)+\gamma}{x} + \frac{1}{x^2}$$

and therefore we obtain

$$\int_0^\infty \left[\frac{\psi(x)+\gamma}{x} + \frac{1}{x^2}\right] \sin(xy) dx = \sum_{a=1}^\infty \frac{1}{a} [-\cos(ay)si(ay) + \sin(ay)Ci(ay)]$$

We have from (6.90d)

$$\int_0^\infty \frac{\sin(xy)dx}{x} = \frac{\pi}{2}$$

and hence we have

$$\int_0^\infty \left[\frac{\psi(x)}{x} + \frac{1}{x^2}\right] \sin(xy) dx = -\frac{\gamma\pi}{2} + \sum_{a=1}^\infty \frac{1}{a}[-\cos(ay)si(ay) + \sin(ay)Ci(ay)]$$

Since $\sum_{a=1}^\infty \frac{1}{a(x+a)} = \frac{\psi(x)+\gamma}{x} + \frac{1}{x^2}$ we note that $\lim_{x\to 0}\left[\frac{\psi(x)+\gamma}{x} + \frac{1}{x^2}\right] = \varsigma(2)$ and this

may also be seen from the Maclaurin expansion

$$\psi(1+x) = -\gamma + \sum_{n=2}^\infty (-1)^n \varsigma(n) x^{n-1}$$

Abramowitz and Stegun [1, p.232] define auxiliary functions

$$f(x) = -\cos x\, si(x) + \sin x\, Ci(x) = \int_0^\infty \frac{\sin y}{y+x}\, dy$$

$$g(x) = -\cos x\, Ci(x) - \sin x\, si(x) = \int_0^\infty \frac{\cos y}{y+x}\, dy$$

and report that for $\text{Re}(x) > 0$

$$f(x) = \int_0^\infty \frac{e^{-xu}}{1+u^2}\, du$$

$$g(x) = \int_0^\infty \frac{u e^{-xu}}{1+u^2}\, du$$

Letting $t = xy$ we see that

$$\int_0^\infty \frac{\sin(xy)}{x+a}\, dx = \int_0^\infty \frac{\sin t}{t+ay}\, dt$$

We also see that

$$\sum_{n=1}^\infty \frac{1}{n^p} \int_0^\infty \frac{e^{-nxu}}{1+u^2}\, du = \sum_{n=1}^\infty \frac{1}{n^p}[-\cos(nx)si(nx) + \sin(nx)Ci(nx)]$$

and therefore we obtain

$$\int_0^\infty \frac{Li_p\left[e^{-xu}\right]}{1+u^2}\, du = \sum_{n=1}^\infty \frac{1}{n^p}[-\cos(nx)si(nx) + \sin(nx)Ci(nx)]$$

Similarly we have

$$\int_0^\infty \frac{u Li_p\left[e^{-xu}\right]}{1+u^2}\, du = -\sum_{n=1}^\infty \frac{1}{n^p}[\sin(nx)si(nx) + \cos(nx)Ci(nx)]$$

$\square$

Integrating Elizalde's identity (6.117ai) gives us

$$\int_0^x \varsigma'(-1,a)\, da =$$

$$\frac{1}{72}x[-4x^2 + 9x - 6 + 6(x-1)(2x-1)\log x] - \frac{1}{12}x^3 + \frac{1}{12}x - \frac{1}{8\pi^2}\sum_{n=1}^\infty \frac{\cos(2n\pi x)}{n^3}$$

$$+\frac{1}{8\pi^2}\varsigma(3)-\frac{1}{4\pi^3}\sum_{n=1}^{\infty}\frac{1}{n^3}[\sin(2n\pi x)Ci(2n\pi x)-\cos(2n\pi x)Si(2n\pi x)]$$

where we have used the integral (easily derived using integration by parts)

$$\int Ci(x)dx = xCi(x)-\sin x \quad \text{and} \quad \int Si(x)dx = xSi(x)+\cos x$$

In evaluating the integral at $a=0$, we have used the fact that $Si(0)=0$ and from (6.94g) we have

$$\sin x\, Ci(x) = \gamma \sin x + \sin x \log x + \sin x \int_0^x \frac{\cos t - 1}{t}dt$$

$$= \gamma \sin x + \frac{\sin x}{x} x \log x + \sin x \int_0^x \frac{\cos t - 1}{t}dt$$

We therefore see that

$$\lim_{x \to 0} \sin x\, Ci(x) = 0$$

From (4.3.131) in Volume II(a) we have

$$n\int_0^x \varsigma'(1-n,u)\,du = \frac{B_{n+1}-B_{n+1}(x)}{n(n+1)}+\varsigma'(-n,x)-\varsigma'(-n)$$

which gives us for $n=2$

$$\int_0^x \varsigma'(-1,a)\,da = -\frac{1}{12}B_3(x)+\frac{1}{2}\varsigma'(-2,x)-\frac{1}{2}\varsigma'(-2)$$

Therefore we obtain

$$\frac{1}{72}x[-4x^2+9x-6+6(x-1)(2x-1)\log x]-\frac{1}{12}x^3+\frac{1}{12}x-\frac{1}{8\pi^2}\sum_{n=1}^{\infty}\frac{\cos(2n\pi x)}{n^3}$$

$$+\frac{1}{8\pi^2}\varsigma(3)-\frac{1}{4\pi^3}\sum_{n=1}^{\infty}\frac{1}{n^3}[\sin(2n\pi x)Ci(2n\pi x)-\cos(2n\pi x)Si(2n\pi x)]$$

$$=-\frac{1}{12}B_3(x)+\frac{1}{2}\varsigma'(-2,x)-\frac{1}{2}\varsigma'(-2)$$

This is easily simplified to

(6.117d)

$$\frac{1}{12}x(x-1)(2x-1)\log x - \frac{5}{36}x^3 + \frac{1}{8}x^2 - \frac{1}{8\pi^2}\sum_{n=1}^{\infty}\frac{\cos(2n\pi x)}{n^3} + \frac{1}{8\pi^2}\varsigma(3)$$

$$-\frac{1}{4\pi^3}\sum_{n=1}^{\infty}\frac{1}{n^3}[\sin(2n\pi x)Ci(2n\pi x) - \cos(2n\pi x)Si(2n\pi x)]$$

$$= -\frac{1}{12}B_3(x) + \frac{1}{2}\varsigma'(-2,x) - \frac{1}{2}\varsigma'(-2)$$

Equation (6.117d) could be integrated to produce an identity involving $\varsigma'(-3,x)$ and so on. This could also be written in several other formats using for example (6.25c)

$$B_3(x) = \frac{3}{2\pi^3}\sum_{n=1}^{\infty}\frac{\sin 2n\pi x}{n^3}$$

and (A.15)

$$B_3(x) = \frac{1}{4}x(x-1)(2x-1)$$

Alternatively, we could employ the generalised Clausen functions $Cl_N(t)$ defined by [126, p.115] as

$$Cl_{2N}(t) = \sum_{n=1}^{\infty}\frac{\sin nt}{n^{2N}} \qquad Cl_{2N+1}(t) = \sum_{n=1}^{\infty}\frac{\cos nt}{n^{2N+1}}$$

With $x = 1$ in (6.117d) we see that

(6.117e) $$\frac{1}{18}\pi^3 = \sum_{n=1}^{\infty}\frac{Si(2n\pi)}{n^3}$$

and this concurs with the result previously obtained in (6.94k).

With $x = 1/2$ in (6.117d) we see that

$$\frac{1}{72} + \frac{1}{8\pi^2}[\varsigma_a(3) + \varsigma(3)] + \frac{1}{4\pi^3}\sum_{n=1}^{\infty}(-1)^n\frac{Si(n\pi)}{n^3} = -\frac{1}{12}B_3\left(\frac{1}{2}\right) + \frac{1}{2}\varsigma'\left(-2,\frac{1}{2}\right) - \frac{1}{2}\varsigma'(-2)$$

We know from (6.94ki) that $\sum_{n=1}^{\infty}(-1)^n\frac{Si(n\pi)}{n^3} = -\frac{1}{18}\pi^3$ and from (F.8b) we have $\varsigma'(-2) = -\frac{\varsigma(3)}{4\pi^2}$. Therefore we obtain

$$\varsigma'\left(-2,\frac{1}{2}\right) = \frac{3\varsigma(3)}{16\pi^2}$$

and we have previously seen this in (4.3.168d) in Volume II(a).

## YET ANOTHER DERIVATION OF GOSPER'S INTEGRAL

Integration of (6.117c) results in

$$\int_\varepsilon^x \log\Gamma(a) =$$

$$\frac{1}{2}(x-\varepsilon)\log(2\pi) + \frac{1}{4}x[2-x+2(x-1)\log x] - \frac{1}{4}\varepsilon[2-\varepsilon+2(\varepsilon-1)\log\varepsilon] - \frac{1}{2}x^2 + \frac{1}{2}\varepsilon^2$$

$$+ \frac{1}{4\pi}\sum_{n=1}^\infty \frac{\sin(2n\pi x)}{n^2} - \frac{1}{4\pi}\sum_{n=1}^\infty \frac{\sin(2n\pi\varepsilon)}{n^2}$$

$$- \frac{1}{2\pi^2}\sum_{n=1}^\infty \frac{1}{n^2}\left[\cos(2n\pi x)Ci(2n\pi x) + \sin(2n\pi x)Si(2n\pi x) - \log(2n\pi x)\right]$$

$$+ \frac{1}{2\pi^2}\sum_{n=1}^\infty \frac{1}{n^2}\left[\cos(2n\pi\varepsilon)Ci(2n\pi\varepsilon) + \sin(2n\pi\varepsilon)Si(2n\pi\varepsilon) - \log(2n\pi\varepsilon)\right]$$

We note from (6.94g) that upon multiplying by $\cos y$

$$\cos y\, Ci(y) = \gamma\cos y + \cos y\log y + \cos y\int_0^y \frac{\cos t - 1}{t}dt$$

and we therefore see that

(6.117ei)

$$\lim_{y\to 0}[\cos y\, Ci(y) - \log y] = \lim_{y\to 0}\left[\gamma\cos y + \log y[\cos y - 1] + \cos y\int_0^y \frac{\cos t - 1}{t}dt\right] = \gamma$$

Therefore as $\varepsilon \to 0$ we have

(6.117f) $$\int_0^x \log\Gamma(a)\,da =$$

$$\frac{1}{2}x\log(2\pi) + \frac{1}{4}x[2-x+2(x-1)\log x] - \frac{1}{2}x^2 + \frac{1}{4\pi}\sum_{n=1}^\infty \frac{\sin(2n\pi x)}{n^2}$$

$$-\frac{1}{2\pi^2}\sum_{n=1}^{\infty}\frac{1}{n^2}\left[\cos(2n\pi x)Ci(2n\pi x)+\sin(2n\pi x)Si(2n\pi x)-\log(2n\pi x)\right]+\frac{\gamma\varsigma(2)}{2\pi^2}$$

$$=\frac{1}{2}x\log(2\pi)+\frac{1}{4}x\left[2-x+2(x-1)\log x\right]-\frac{1}{2}x^2+\frac{1}{4\pi}\sum_{n=1}^{\infty}\frac{\sin(2n\pi x)}{n^2}+\frac{1}{12}\log x$$

$$-\frac{1}{2\pi^2}\sum_{n=1}^{\infty}\frac{1}{n^2}\left[\cos(2n\pi x)Ci(2n\pi x)+\sin(2n\pi x)Si(2n\pi x)\right]+\frac{\varsigma(2)\log(2\pi)}{2\pi^2}-\frac{\varsigma'(2)}{2\pi^2}+\frac{\gamma\varsigma(2)}{2\pi^2}$$

and with a little algebra and using (F.7) in Volume VI

$$\varsigma'(-1)=\frac{1}{12}(1-\gamma-\log 2\pi)+\frac{1}{2\pi^2}\varsigma'(2)$$

we obtain

(6.117fi) $$\int_0^x \log\Gamma(a)\,da=$$

$$=\frac{1}{2}x\log(2\pi)+\frac{1}{4}x\left[2-x+2(x-1)\log x\right]-\frac{1}{2}x^2+\frac{1}{4\pi}\sum_{n=1}^{\infty}\frac{\sin(2n\pi x)}{n^2}+\frac{1}{12}\log x$$

$$-\frac{1}{2\pi^2}\sum_{n=1}^{\infty}\frac{1}{n^2}\left[\cos(2n\pi x)Ci(2n\pi x)+\sin(2n\pi x)Si(2n\pi x)\right]+\frac{1}{12}-\varsigma'(-1)$$

Letting $x=1$ in (6.117fi) gives us back (6.117), and with $x=1/2$ we get

$$\int_0^{1/2}\log\Gamma(a)\,da=\frac{1}{4}\log(2\pi)+\frac{1}{16}+\frac{1}{24}\log 2-\frac{1}{2\pi^2}\sum_{n=1}^{\infty}\frac{(-1)^n}{n^2}Ci(n\pi)+\frac{1}{12}-\varsigma'(-1)$$

and using (6.130b) we see that

(6.117fii) $$\int_0^{1/2}\log\Gamma(a)\,da=\frac{5}{24}\log 2+\frac{1}{4}\log\pi+\frac{3}{2}\log A$$

which is reported in [126, p.35]. With $x=1/4$ we obtain

$$\int_0^{1/4}\log\Gamma(a)\,da=$$

$$=\frac{1}{8}\log(2\pi)+\frac{5}{64}+\frac{G}{4\pi}+\frac{1}{48}\log 2$$

$$-\frac{1}{2\pi^2}\sum_{n=1}^{\infty}\frac{1}{n^2}\left[\cos(n\pi/2)Ci(n\pi/2)+\sin(n\pi/2)Si(n\pi/2)\right]+\frac{1}{12}-\varsigma'(-1)$$

From [126, p.35] we have

$$\int_0^{1/4}\log\Gamma(a)\,da=\frac{1}{8}\log 2+\frac{1}{8}\log\pi+\frac{9}{8}\log A+\frac{G}{4\pi}$$

and we therefore obtain

(6.117fiii)

$$\frac{1}{2\pi^2}\sum_{n=1}^{\infty}\frac{1}{n^2}\left[\cos(n\pi/2)Ci(n\pi/2)+\sin(n\pi/2)Si(n\pi/2)\right]=\frac{5}{64}+\frac{1}{48}\log 2-\frac{1}{8}\log A$$

Using (6.117ai) we may write (6.117fi) as Gosper's integral (4.3.129a) in Volume II(a)

(6.117g) $$\int_0^x \log\Gamma(a)da=\frac{1}{2}x(1-x)+\frac{1}{2}x\log(2\pi)+\varsigma'(-1,x)-\varsigma'(-1)$$

In fact, using (6.117f) with $x=1$ and Raabe's integral $\int_0^1\log\Gamma(a)=\frac{1}{2}\log(2\pi)$ we can easily derive (F.7), which was there obtained using Riemann's functional equation for the zeta function.

Referring to (6.117c), we obtain upon multiplication by $a$ and then integrating

$$\int_0^x a\log\Gamma(a)=\frac{1}{4}x^2\log(2\pi)+\frac{1}{72}x^2\left(9-8x+6(4x-3)\log x\right)-\frac{1}{3}x^3$$

$$+\frac{1}{8\pi^2}\sum_{n=1}^{\infty}\frac{\cos(2n\pi x)-1}{n^3}+\frac{1}{4\pi}x\sum_{n=1}^{\infty}\frac{\sin(2n\pi x)}{n^2}$$

$$+\frac{1}{4\pi^3}\sum_{n=1}^{\infty}\frac{1}{n^3}[\sin(2n\pi x)Ci(2n\pi x)-\cos(2n\pi x)Si(2n\pi x)]$$

$$-\frac{x}{2\pi^2}\sum_{n=1}^{\infty}\frac{1}{n^2}[\sin(2n\pi x)Ci(2n\pi x)+\cos(2n\pi x)Si(2n\pi x)]+\frac{x}{2\pi^2}\sum_{n=1}^{\infty}\frac{1}{n^2}$$

and using (6.117d) we may write this as

(6.117h)
$$\int_0^x a\log\Gamma(a) = \frac{1}{4}x^2\log(2\pi) + \frac{1}{72}x^2\left(9-8x+6(4x-3)\log x\right) - \frac{1}{3}x^3$$

$$+ \frac{1}{8\pi^2}\sum_{n=1}^{\infty}\frac{\cos(2n\pi x)-1}{n^3} + \frac{1}{4\pi}x\sum_{n=1}^{\infty}\frac{\sin(2n\pi x)}{n^2}$$

$$+ \frac{1}{12}x(x-1)(2x-1)\log x - \frac{5}{36}x^3 + \frac{1}{8}x^2 - \frac{1}{8\pi^2}\sum_{n=1}^{\infty}\frac{\cos(2n\pi x)}{n^3} + \frac{1}{8\pi^2}\varsigma(3)$$

$$+ \frac{1}{12}B_3(x) - \frac{1}{2}\varsigma'(-2,x) + \frac{1}{2}\varsigma'(-2)$$

$$- \frac{x}{2\pi^2}\sum_{n=1}^{\infty}\frac{1}{n^2}[\sin(2n\pi x)Ci(2n\pi x) + \cos(2n\pi x)Si(2n\pi x)] + \frac{x}{2\pi^2}\sum_{n=1}^{\infty}\frac{1}{n^2}$$

and with some cancellations we have

$$\int_0^x a\log\Gamma(a) = \frac{1}{4}x^2\log(2\pi) + \frac{1}{8}x + \frac{1}{8}x^2 - \frac{1}{2}x^3 + \frac{1}{12}x(6x^2 - 6x + 1)\log x$$

$$+ \frac{1}{4\pi}x\sum_{n=1}^{\infty}\frac{\sin(2n\pi x)}{n^2} - \frac{1}{2}\varsigma'(-2,x) + \frac{1}{2}\varsigma'(-2)$$

$$- \frac{x}{2\pi^2}\sum_{n=1}^{\infty}\frac{1}{n^2}[\sin(2n\pi x)Ci(2n\pi x) + \cos(2n\pi x)Si(2n\pi x)]$$

We also have from (6.116a)

$$\int_0^x a\log\Gamma(a)\,da = \frac{1}{2}\left(\frac{1}{4} - 2\log A\right)x + \frac{1}{2}\left(\frac{1}{2}\log(2\pi) - \frac{1}{4}\right)x^2$$

$$+ \frac{1}{2}x^2\log\Gamma(x) - \frac{1}{2}\log G(1+x) - \log\Gamma_3(1+x)$$

and equating the above two integrals gives us

$$\frac{1}{8}x + \frac{1}{4}x^2 - \frac{1}{2}x^3 + \frac{1}{12}x(6x^2 - 6x + 1)\log x$$

$$+ \frac{1}{4\pi}x\sum_{n=1}^{\infty}\frac{\sin(2n\pi x)}{n^2} - \frac{1}{2}\varsigma'(-2, x) + \frac{1}{2}\varsigma'(-2)$$

(6.117i)

$$-\frac{x}{2\pi^2}\sum_{n=1}^{\infty}\frac{1}{n^2}[\sin(2n\pi x)Ci(2n\pi x) + \cos(2n\pi x)Si(2n\pi x)] =$$

$$\frac{1}{2}\left(\frac{1}{4} - 2\log A\right)x + \frac{1}{2}x^2\log\Gamma(x) - \frac{1}{2}\log G(1+x) - \log\Gamma_3(1+x)$$

Letting $x = 1$ results in

(6.117j) $\qquad \dfrac{1}{2\pi^2}\sum_{n=1}^{\infty}\dfrac{Si(2n\pi)}{n^2} = \log A - \dfrac{1}{4}$

and hence we may have a closed form expression for the series in (6.92b). Letting $x = 1/2$ will result in the following expression for $\sum_{n=1}^{\infty}(-1)^n\dfrac{Si(n\pi)}{n^2}$ involving $\log\Gamma_3(3/2)$.

$$\frac{1}{4\pi^2}\sum_{n=1}^{\infty}\frac{(-1)^n Si(n\pi)}{n^2} =$$

$$\frac{1}{48}\log 2 - \frac{1}{2}\varsigma'\left(-2, \frac{1}{2}\right) + \frac{1}{2}\varsigma'(-2) - \frac{1}{2}\log A + \frac{1}{16}\log\pi - \frac{1}{2}\log G(3/2) - \log\Gamma_3(3/2)$$

which may be written as

$$\frac{1}{4\pi^2}\sum_{n=1}^{\infty}\frac{(-1)^n Si(n\pi)}{n^2} =$$

(6.117ji)

$$\frac{1}{48}\log 2 - \frac{7\varsigma(3)}{32\pi^2} - \frac{1}{2}\log A + \frac{1}{16}\log\pi - \frac{1}{2}\log G(3/2) - \log\Gamma_3(3/2)$$

Using (6.117o)

$$\log\Gamma_3(3/2) = -\frac{1}{16}\log\pi + \frac{7}{32}\frac{\varsigma(3)}{\pi^2}$$

and $\log G(3/2) = \log G(1/2) + \log\Gamma(1/2)$ and (6.117ki)

$$\log G(1/2) = \frac{1}{24}\log 2 - \frac{1}{4}\log \pi + \frac{3}{2}\varsigma'(-1)$$

we obtain

$$\frac{1}{4\pi^2}\sum_{n=1}^{\infty}\frac{(-1)^n Si(n\pi)}{n^2} =$$

$$\frac{1}{48}\log 2 - \frac{7\varsigma(3)}{16\pi^2} - \frac{1}{2}\log A + \frac{1}{8}\log \pi - \frac{1}{2}\log G(3/2)$$

and hence we have

(6.117jii) $$\frac{1}{4\pi^2}\sum_{n=1}^{\infty}\frac{(-1)^n Si(n\pi)}{n^2} = -\frac{7\varsigma(3)}{16\pi^2} + \frac{1}{4}\log A - \frac{3}{48}$$

□

Using (4.3.85)

$$\int_0^x \log \Gamma(a)\,da = \frac{1}{2}x(1-x) + \frac{x}{2}\log(2\pi) - \log G(1+x) + x\log\Gamma(x)$$

we could alternatively write (6.117f) as

$$\frac{1}{2}x(1-x) + \frac{x}{2}\log(2\pi) - \log G(1+x) + x\log\Gamma(x)$$

$$= \frac{1}{2}x\log(2\pi) + \frac{1}{4}x\left[2 - x + 2(x-1)\log x\right] - \frac{1}{2}x^2 + \frac{1}{4\pi}\sum_{n=1}^{\infty}\frac{\sin(2n\pi x)}{n^2} + \frac{1}{12}\log x$$

$$-\frac{1}{2\pi^2}\sum_{n=1}^{\infty}\frac{1}{n^2}\left[\cos(2n\pi x)Ci(2n\pi x) + \sin(2n\pi x)Si(2n\pi x)\right] + \frac{\varsigma(2)\log(2\pi)}{2\pi^2} - \frac{\varsigma'(2)}{2\pi^2} + \frac{\gamma\varsigma(2)}{2\pi^2}$$

This easily simplifies to

(6.117k) $$-\log G(1+x) + x\log\Gamma(x) = \frac{1}{4}x\left[-x + 2(x-1)\log x\right] + \frac{Cl_2(2\pi x)}{4\pi} + \frac{1}{12}\log x$$

$$-\frac{1}{2\pi^2}\sum_{n=1}^{\infty}\frac{1}{n^2}\left[\cos(2n\pi x)Ci(2n\pi x) + \sin(2n\pi x)Si(2n\pi x)\right] + \frac{1}{12} - \varsigma'(-1)$$

where we have used (F.7).

Letting $x = 1$ results in

$$\frac{1}{2\pi^2}\sum_{n=1}^{\infty}\frac{Ci(2n\pi)}{n^2}=-\frac{1}{6}-\varsigma'(-1)$$

which we saw previously in (6.117).

With $x=1/2$ we have

$$-\log G(3/2)+\frac{1}{4}\log\pi=-\frac{1}{16}+\frac{1}{8}\log 2-\frac{1}{2\pi^2}\sum_{n=1}^{\infty}(-1)^n\frac{Ci(n\pi)}{n^2}+\frac{1}{12}-\varsigma'(-1)-\frac{1}{12}\log 2$$

and this becomes

$$-\log G(1/2)-\frac{1}{2}\log\pi+\frac{1}{4}\log\pi=-\frac{1}{16}+\frac{1}{24}\log 2-\frac{1}{2\pi^2}\sum_{n=1}^{\infty}(-1)^n\frac{Ci(n\pi)}{n^2}+\frac{1}{12}-\varsigma'(-1)$$

Employing (6.117b) then shows that

(6.117ki) $\qquad \log G(1/2)=\frac{1}{24}\log 2-\frac{1}{4}\log\pi+\frac{3}{2}\varsigma'(-1)$

With $x=1/4$ in (6.117k) we get

$$-\log G(5/4)+\frac{1}{4}\log\Gamma(1/4)=\frac{1}{16}\left[-\frac{1}{4}+3\log 2\right]+\frac{Cl_2(\pi/2)}{4\pi}-\frac{1}{6}\log 2$$

$$-\frac{1}{2\pi^2}\sum_{n=1}^{\infty}\frac{1}{n^2}\left[\cos(n\pi/2)Ci(n\pi/2)+\sin(n\pi/2)Si(n\pi/2)\right]+\frac{1}{12}-\varsigma'(-1)$$

and using (6.117fiii)

$$\frac{1}{2\pi^2}\sum_{n=1}^{\infty}\frac{1}{n^2}\left[\cos(n\pi/2)Ci(n\pi/2)+\sin(n\pi/2)Si(n\pi/2)\right]=\frac{5}{64}+\frac{1}{48}\log 2-\frac{1}{8}\log A$$

we may write this as

$$-\log G(5/4)+\frac{1}{4}\log\Gamma(1/4)=\frac{1}{16}\left[-\frac{1}{4}+3\log 2\right]+\frac{Cl_2(\pi/2)}{4\pi}-\frac{1}{6}\log 2$$

$$-\frac{5}{64}-\frac{1}{48}\log 2+\frac{1}{8}\log A+\frac{1}{12}-\varsigma'(-1)$$

From (6.107p) we know that $Cl_2(\pi/2)=G$ and we have

$$\log G(5/4)=\log G(1/4)+\log\Gamma(1/4)$$

and hence we have

(6.117kii) $$\log G(1/4) = \frac{3}{32} - \frac{G}{4\pi} - \frac{9}{8}\log A - \frac{3}{4}\log \Gamma(1/4)$$

as is shown by Srivastava and Choi in [126, p.30].

Combining (6.117k) and (6.117ai) results in the Gosper/Vardi identity (4.3.126)

$$\log G(u+1) - u \log \Gamma(u) = \varsigma'(-1) - \varsigma'(-1,u)$$

whereupon letting $u = 1/4$ results in (6.117kii).

As shown below, we may perform a relatively painless integration of (6.117k) to evaluate $\int_0^t \log G(1+x) du$.

We have in an instant from the Wolfram Integrator

(6.117l) $$\int_0^t [\cos(2n\pi x) Ci(2n\pi x) + \sin(2n\pi x) Si(2n\pi x)] dx =$$

$$\frac{1}{2n\pi}[\sin(2n\pi t) Ci(2n\pi t) - \cos(2n\pi t) Si(2n\pi t)]$$

and we obtain

$$-\int_0^t \log G(1+x) du + \int_0^t x \log \Gamma(x) = -\frac{1}{12}t^3 + \frac{1}{72}t^2(9 - 4t + 6(2t-3)\log t)$$

$$+ \frac{1}{12}(t \log t - t) - \frac{1}{8\pi^2}\sum_{n=1}^{\infty} \frac{\cos(2n\pi t) - 1}{n^3}$$

$$- \frac{1}{4\pi^3}\sum_{n=1}^{\infty}\frac{1}{n^3}[\sin(2n\pi t) Ci(2n\pi t) - \cos(2n\pi t) Si(2n\pi t)] + \frac{1}{12}t - \varsigma'(-1)t$$

Using (4.3.85) and (4.3.87c) from Volume II(a)

$$\int_0^t \log G(1+x) dx =$$

$$\left(\frac{1}{4} - 2\log A\right)t + \frac{1}{4}\log(2\pi)t^2 - \frac{1}{6}t^3 + (t-1)\log G(1+t) - 2\log \Gamma_3(1+t)$$

where the triple gamma function $\Gamma_3(z)$ is the unique meromorphic function which satisfies each of the following properties

(i) $\Gamma_3(1) = 1$

(ii) $\Gamma_3(z+1) = G(z)\Gamma_3(z)$

(iii) For $x \geq 1$, $\Gamma_3(x)$ is infinitely differentiable and $\dfrac{d^4}{dx^4}\log\Gamma_3(x) \geq 0$

we then determine that

$$-\left(\frac{1}{4} - 2\log A\right)t - \frac{1}{4}\log(2\pi)t^2 + \frac{1}{6}t^3 - (t-1)\log G(1+t) + 2\log\Gamma_3(1+t)$$

$$+\frac{1}{2}\left(\frac{1}{4} - 2\log A\right)t + \frac{1}{2}\left(\frac{1}{2}\log(2\pi) - \frac{1}{4}\right)t^2$$

$$+\frac{1}{2}t^2\log\Gamma(t) - \frac{1}{2}\log G(1+t) - \log\Gamma_3(1+t)$$

$$= -\frac{1}{12}t^3 + \frac{1}{72}t^2(9 - 4t + 6(2t-3)\log t) + \frac{1}{12}t\log t - \frac{1}{8\pi^2}\sum_{n=1}^{\infty}\frac{\cos(2n\pi t) - 1}{n^3}$$

$$-\frac{1}{4\pi^3}\sum_{n=1}^{\infty}\frac{1}{n^3}\left[\sin(2n\pi t)Ci(2n\pi t) - \cos(2n\pi t)Si(2n\pi t)\right] - \varsigma'(-1)t$$

which simplifies to

$$t\log A - \frac{1}{8}t - \frac{1}{4}t^2 + \frac{11}{36}t^3 - (t-1)\log G(1+t) + \log\Gamma_3(1+t)$$

$$+\frac{1}{2}t^2\log\Gamma(t) - \frac{1}{2}\log G(1+t)$$

$$= \frac{1}{12}t(2t^2 - 3t + 1)\log t - \frac{1}{8\pi^2}\sum_{n=1}^{\infty}\frac{\cos(2n\pi t) - 1}{n^3}$$

$$-\frac{1}{4\pi^3}\sum_{n=1}^{\infty}\frac{1}{n^3}\left[\sin(2n\pi t)Ci(2n\pi t) - \cos(2n\pi t)Si(2n\pi t)\right] - \varsigma'(-1)t$$

Using (4.4.225)

$$\log A = \frac{1}{12} - \varsigma'(-1)$$

this becomes

(6.117m)

$$\frac{11}{36}t^3 - \frac{1}{24}t - \frac{1}{4}t^2 - \left(t - \frac{1}{2}\right)\log G(1+t) + \log \Gamma_3(1+t) + \frac{1}{2}t^2 \log \Gamma(t)$$

$$= \frac{1}{12}t(t-1)(2t-1)\log t - \frac{1}{8\pi^2}\sum_{n=1}^{\infty}\frac{\cos(2n\pi t)-1}{n^3}$$

$$-\frac{1}{4\pi^3}\sum_{n=1}^{\infty}\frac{1}{n^3}[\sin(2n\pi t)Ci(2n\pi t) - \cos(2n\pi t)Si(2n\pi t)]$$

Letting $t = 1$ gives us

(6.117n) $\quad \sum_{n=1}^{\infty}\frac{Si(2n\pi)}{n^3} = \frac{\pi^3}{18}$

in agreement with (6.117e).

With $t = 1/2$ we have

$$\frac{1}{4\pi^3}\sum_{n=1}^{\infty}(-1)^n \frac{Si(n\pi)}{n^3} = -\frac{1}{72} + \frac{1}{16}\log \pi + \log \Gamma_3(3/2) - \frac{7}{32}\frac{\varsigma(3)}{\pi^2}$$

and using (6.94ki)

(6.117o) $\quad \sum_{n=1}^{\infty}(-1)^n \frac{Si(n\pi)}{n^3} = -\frac{1}{18}\pi^3$

we see that

(6.117p) $\quad \log \Gamma_3(3/2) = -\frac{1}{16}\log \pi + \frac{7}{32}\frac{\varsigma(3)}{\pi^2}$

as reported in [45ab] and [126, p.44].

In (6.117d) we saw that

$$\frac{1}{12}t(t-1)(2t-1)\log t - \frac{5}{36}t^3 + \frac{1}{8}t^2 - \frac{1}{8\pi^2}\sum_{n=1}^{\infty}\frac{\cos(2n\pi x)-1}{n^3}$$

$$-\frac{1}{4\pi^3}\sum_{n=1}^{\infty}\frac{1}{n^3}[\sin(2n\pi t)Ci(2n\pi t) - \cos(2n\pi t)Si(2n\pi t)]$$

$$= -\frac{1}{12}B_3(t) + \frac{1}{2}\varsigma'(-2,t) - \frac{1}{2}\varsigma'(-2)$$

and comparing this with (6.117m) gives us

(6.117pi)

$$\varsigma'(-2,t) - \varsigma'(-2) = \frac{5}{6}t^3 - t^2 + \frac{1}{6}t - (2t-1)\log G(1+t) + 2\log \Gamma_3(1+t) + t^2 \log \Gamma(t)$$

Differentiating (6.117c) (and boldly assuming that the procedure is valid) easily results in

(6.117q) $$\psi(a) = \log a - \frac{1}{2a} + 2\sum_{n=1}^{\infty}[\cos(2n\pi a)Ci(2n\pi a) + \sin(2n\pi a)si(2n\pi a)]$$

which appears in Nörlund's book [105, p.108]. Letting $a=1$ results in

(6.117r) $$\frac{1}{2} - \gamma = 2\sum_{n=1}^{\infty} Ci(2n\pi)$$

and this corrects the corresponding formula given by Nielsen [104a, p.80]. This formula was also used in Volume II(b). With $a = 1/2$ we get

$$\psi\left(\frac{1}{2}\right) = -\log 2 - 1 + 2\sum_{n=1}^{\infty}(-1)^n Ci(n\pi)$$

and hence we have

(6.117s) $$2\sum_{n=1}^{\infty}(-1)^n Ci(n\pi) = 1 - \gamma - \log 2$$

Differentiating (6.117q) (and again assuming that the operation is valid) results in

$$\psi'(a) = \frac{1}{a} + \frac{1}{2a^2} + 2\sum_{n=1}^{\infty}\left[\frac{1}{a} - 2n\pi \sin(2n\pi a)Ci(2n\pi a) + 2n\pi \cos(2n\pi a)si(2n\pi a)\right]$$

but this series does not appear to be convergent.

We have from G&R [74, p.650] for $a > 0$

(6.118)

$$\int_0^1 \log \Gamma(x+a) \cos 2n\pi x\, dx = -\frac{1}{2\pi n}\left[\sin(2n\pi a)Ci(2n\pi a) + \cos(2n\pi a)si(2n\pi a)\right]$$

and with $a=1$ we get for $n \geq 1$

(6.119) $$\int_0^1 \log \Gamma(x+1) \cos 2n\pi x \, dx = -\frac{si(2n\pi)}{2\pi n}$$

This may be easily confirmed as follows. We have

$$\int_0^1 \log \Gamma(x+1) \cos 2n\pi x \, dx = \int_0^1 \log x \cos 2n\pi x \, dx + \int_0^1 \log \Gamma(x) \cos 2n\pi x \, dx$$

and in (E.44b) in Volume VI we have proved that

(6.119i) $$\int_0^1 \log \Gamma(x) \cos 2n\pi x \, dx = \frac{1}{4n}$$

It may be seen from (6.90a) that

$$\int_0^1 \log x \cos 2n\pi x \, dx = -\frac{Si(2n\pi)}{2n\pi} = -\frac{si(2n\pi)}{2n\pi} - \frac{1}{4n}$$

and (6.119) follows directly.

Equation (6.118) also applies in the limit $a \to 0$ because we have $\lim_{x \to 0} \sin x \, Ci(x) = 0$ and $si(0) = -\frac{\pi}{2}$ and we therefore obtain (6.119i).

Furthermore, we may employ (6.5) to give us

$$-\frac{1}{2}\int_0^1 \log \Gamma(x+1) \, dx = \sum_{n=1}^\infty \int_0^1 \log \Gamma(x+1) \cos 2n\pi x \, dx = -\frac{1}{2\pi} \sum_{n=1}^\infty \frac{si(2n\pi)}{n}$$

and we therefore get

(6.120) $$\sum_{n=1}^\infty \frac{si(2n\pi)}{n} = \frac{\pi}{2} \log(2\pi) - \pi$$

and this concurs with Nielsen's result [104a, p.79]. Nielsen also records that

(6.121) $$\sum_{n=1}^\infty (-1)^n \frac{si(2n\pi)}{n} = \frac{\pi}{2} \log 2 - \pi$$

but, in view of the singularity in the interval of integration, we are not able to employ (6.8) to derive this result (see (6.94a) instead).

We recall (6.111) and (6.118)

$$\int_0^1 \log \Gamma(x+a)\sin 2n\pi x\, dx = -\frac{1}{2\pi n}\left[\log a - \cos(2n\pi a)Ci(2n\pi a) + \sin(2n\pi a)si(2n\pi a)\right]$$

$$\int_0^1 \log \Gamma(x+a)\cos 2n\pi x\, dx = -\frac{1}{2\pi n}\left[\sin(2n\pi a)Ci(2n\pi a) + \cos(2n\pi a)si(2n\pi a)\right]$$

and multiplying each by $\sin(2n\pi a)$ and $\cos(2n\pi a)$ respectively we get

$$\sin(2n\pi a)\int_0^1 \log \Gamma(x+a)\sin 2n\pi x\, dx =$$

$$-\frac{1}{2\pi n}\left[\log a \sin(2n\pi a) - \cos(2n\pi a)\sin(2n\pi a)Ci(2n\pi a) + \sin^2(2n\pi a)si(2n\pi a)\right]$$

$$\cos(2n\pi a)\int_0^1 \log \Gamma(x+a)\cos 2n\pi x\, dx =$$

$$-\frac{1}{2\pi n}\left[\cos(2n\pi a)\sin(2n\pi a)Ci(2n\pi a) + \cos^2(2n\pi a)si(2n\pi a)\right]$$

Adding these two equations results in

$$\int_0^1 \log \Gamma(x+a)[\sin(2n\pi a)\sin 2n\pi x + \cos(2n\pi a)\cos 2n\pi x]dx = -\frac{1}{2\pi n}\left[\log a \sin(2n\pi a) + si(2n\pi a)\right]$$

or equivalently

(6.121a) $$\int_0^1 \log \Gamma(x+a)\cos[2n\pi(x-a)]dx = -\frac{1}{2\pi n}\left[\log a \sin(2n\pi a) + si(2n\pi a)\right]$$

With $a=0$ we immediately obtain (6.119i)

$$\int_0^1 \log \Gamma(x)\cos 2n\pi x\, dx = \frac{1}{4n}$$

and $a = 1/2$ gives us

(6.121b) $$\int_0^1 \log \Gamma(x+1/2)\cos 2n\pi x\, dx = (-1)^{n+1}\frac{si(n\pi)}{2\pi n}$$

in accordance with (6.118).

Similarly subtraction gives us

(6.121c) $$\int_0^1 \log \Gamma(x+a) \cos[2n\pi(x+a)] dx =$$

$$\frac{1}{4\pi n}\left[2\log a \sin(2n\pi a) - \sin(4n\pi a)Ci(2n\pi a) - \cos(4n\pi a)si(2n\pi a)\right]$$

and with $a=1$ and $a=1/2$ we see that

$$\int_0^1 \log \Gamma(x+1) \cos 2n\pi x \, dx = \frac{si(2n\pi)}{4\pi n}$$

$$\int_0^1 \log \Gamma(x+1/2) \cos 2n\pi x \, dx = (-1)^n \frac{si(n\pi)}{4\pi n}$$

Completing the summation of (6.119) we obtain

(6.122) $$\sum_{n=1}^\infty \frac{1}{n} \int_0^1 \log \Gamma(x+1) \cos 2n\pi x \, dx = -\frac{1}{2\pi} \sum_{n=1}^\infty \frac{si(2n\pi)}{n^2}$$

Using (7.8) we get (assuming that it is valid to interchange the order of summation and integration)

$$\sum_{n=1}^\infty \frac{1}{n} \int_0^1 \log \Gamma(x+1) \cos 2n\pi x \, dx = \int_0^1 \log \Gamma(x+1) \sum_{n=1}^\infty \frac{\cos 2n\pi x}{n} dx$$

$$= -\int_0^1 \log \Gamma(x+1) \log[2\sin(\pi x)] dx$$

and hence we get

(6.123) $$\int_0^1 \log \Gamma(x+1) \log[2\sin(\pi x)] dx = \frac{1}{2\pi} \sum_{n=1}^\infty \frac{si(2n\pi)}{n^2}$$

In the above we have tacitly assumed that interchanging the order of integration and summation is valid and that it is appropriate to use the Fourier series expansion at both end points because the multiplication factor $\log \Gamma(x+1)$ is zero at both of those end points.

Similarly, we may also deduce that

$$\sum_{n=1}^\infty \frac{1}{n} \int_0^1 \log \Gamma(x+a) \cos 2n\pi x \, dx = \int_0^1 \log \Gamma(x+a) \sum_{n=1}^\infty \frac{\cos 2n\pi x}{n} dx$$

$$= -\int_0^1 \log \Gamma(x+a) \log[2\sin(\pi x)] dx$$

and using (6.118) we have

$$\sum_{n=1}^{\infty} \frac{1}{n} \int_0^1 \log \Gamma(x+a) \cos 2n\pi x \, dx = -\frac{1}{2\pi} \sum_{n=1}^{\infty} \frac{1}{n^2} [\sin(2n\pi a) Ci(2n\pi a) + \cos(2n\pi a) si(2n\pi a)]$$

Therefore we see that

(6.124)

$$\frac{1}{2\pi} \sum_{n=1}^{\infty} \frac{1}{n^2} [\sin(2n\pi a) Ci(2n\pi a) + \cos(2n\pi a) si(2n\pi a)] = \int_0^1 \log \Gamma(x+a) \log[2\sin(\pi x)] dx$$

We now use (6.111) again

$$\int_0^1 \log \Gamma(x+a) \sin 2n\pi x \, dx = -\frac{1}{2\pi n} [\log a - \cos(2n\pi a) Ci(2n\pi a) + \sin(2n\pi a) si(2n\pi a)]$$

and with $a = 1/2$ we get

(6.125) $$\int_0^1 \log \Gamma(x+1/2) \sin 2n\pi x \, dx = \frac{1}{2\pi n} [\log 2 + (-1)^n Ci(n\pi)]$$

Completing the summation we obtain

$$\sum_{n=1}^{\infty} \frac{1}{n} \int_0^1 \log \Gamma(x+1/2) \sin 2n\pi x \, dx = \frac{\log 2}{2\pi} \varsigma(2) + \frac{1}{2\pi} \sum_{n=1}^{\infty} (-1)^n \frac{Ci(n\pi)}{n^2}$$

Using (7.5) we get as before

$$\sum_{n=1}^{\infty} \frac{1}{n} \int_0^1 \log \Gamma(x+1/2) \sin 2n\pi x \, dx = \int_0^1 \log \Gamma(x+1/2) \sum_{n=1}^{\infty} \frac{\sin 2n\pi x}{n} dx$$

$$= \frac{\pi}{2} \int_0^1 (1-2x) \log \Gamma(x+1/2) dx$$

$$= \frac{\pi}{2} \left( \int_0^1 \log \Gamma(x+1/2) dx - 2\int_0^1 x \log \Gamma(x+1/2) dx \right)$$

We have from [45ab] and [125, p.209]

$$(6.126) \quad 2\int_0^z t \log \Gamma(a+t)\,dt = \left(\frac{1}{4} - \frac{1}{2}a + \frac{1}{2}a^2 - 2\log A\right)z + \left(\frac{1}{2}\log(2\pi) - \frac{1}{2}a + \frac{1}{4}\right)z^2$$

$$-\frac{1}{2}z^2 + z^2 \log \Gamma(z+a) - (a-1)^2 [\log \Gamma(z+a) - \log \Gamma(a)] + (2a-3)[\log G(z+a) - \log G(a)]$$

$$-2[\log \Gamma_3(z+a) - \log \Gamma_3(a)]$$

With $z = 1$ and $a = 1/2$ we get

$$2\int_0^1 t \log \Gamma(t+1/2)\,dt =$$

$$-\frac{3}{8} - 2\log A + \frac{1}{2}\log(2\pi) + \log \Gamma(3/2) - \frac{1}{4}\log \frac{\Gamma(3/2)}{\Gamma(1/2)} - 2\log \frac{G(3/2)}{G(1/2)} - 2\log \frac{\Gamma_3(3/2)}{\Gamma_3(1/2)}$$

Since $\Gamma_3(z+1) = G(z)\Gamma_3(z)$ we have $\Gamma_3(3/2) = G(1/2)\Gamma_3(1/2)$ and therefore

$$\log \frac{\Gamma_3(3/2)}{\Gamma_3(1/2)} = \log G(1/2)$$

As mentioned previously in Volume III the value of $G(1/2)$ was originally determined by Barnes [17aa] in 1899 as (see also [126, p.26])

$$(6.127) \quad G(1/2) = A^{-\frac{3}{2}} \pi^{-\frac{1}{4}} e^{\frac{1}{8}} 2^{\frac{1}{24}}$$

where $A$ is the Glaisher-Kinkelin constant. This result is derived above in (6.117ki) and also in (4.4.228ti) of Volume IV.

Using (6.58) we have $\log \frac{G(3/2)}{G(1/2)} = \log \Gamma(1/2)$ and we know that $\Gamma(3/2) = \frac{1}{2}\Gamma(1/2)$. Therefore we obtain

$$(6.128) \quad 2\int_0^1 t \log \Gamma(t+1/2)\,dt = -\frac{3}{8} - 2\log A - \frac{1}{4}\log 2 - 2\log G(1/2)$$

Using the following formula in [126, p.207] which is due to Barnes

$$(6.128a) \quad \int_0^z \log \Gamma(t+a)\,dt =$$

$$\frac{1}{2}[\log(2\pi) + 1 - 2a]z - \frac{z^2}{2} + z\log \Gamma(z+a) + (a-1)\log[\Gamma(z+a) - \log \Gamma(a)]$$

$$-\left[\log G(z+a) - \log G(a)\right]$$

we get

$$\int_0^1 \log \Gamma(t+1/2)\,dt = \frac{1}{2}\log(2\pi) - \frac{1}{2} + \log\Gamma(3/2) - \frac{1}{2}\log\frac{\Gamma(3/2)}{\Gamma(1/2)} - \log\frac{G(3/2)}{G(1/2)}$$

Therefore we have

(6.129) $$\int_0^1 \log\Gamma(t+1/2)\,dt = \frac{1}{2}\log\pi - \frac{1}{2}$$

Alternatively we could have used (C.43a).

$$\int_0^1 \log\Gamma(t+a)\,dt = \int_a^{1+a} \log\Gamma(u)\,du = \frac{1}{2}\log(2\pi) + a\log a - a$$

Hence we have

$$\frac{\pi}{2}\left(\int_0^1 \log\Gamma(t+1/2)\,dt - 2\int_0^1 t\log\Gamma(t+1/2)\,dt\right) =$$

$$\frac{\pi}{2}\left(\frac{1}{2}\log\pi + \frac{1}{4}\log 2 - \frac{1}{8} + 2\log A + 2\log G(1/2)\right)$$

and therefore we obtain

$$\frac{\log 2}{2\pi}\varsigma(2) + \frac{1}{2\pi}\sum_{n=1}^\infty (-1)^n \frac{Ci(n\pi)}{n^2} =$$

$$\frac{\pi}{2}\left(\frac{1}{2}\log\pi + \frac{1}{4}\log 2 - \frac{1}{8} + 2\log A + 2\log G(1/2)\right)$$

$$= \frac{\pi}{2}\left(\frac{1}{3}\log 2 + \frac{1}{8} - \log A\right)$$

This then gives us

(6.130) $$\sum_{n=1}^\infty (-1)^n \frac{Ci(n\pi)}{n^2} = \pi^2\left[\frac{1}{6}\log 2 + \frac{1}{8} - \log A\right]$$

We saw in (6.117b) that

$$\sum_{n=1}^{\infty}(-1)^n \frac{Ci(n\pi)}{n^2} = \pi^2\left[\frac{1}{24} + \frac{1}{6}\log 2 + \varsigma'(-1)\right]$$

and this is equivalent to (6.130).

□

We now let $a = 1/2$ in (6.118) to obtain

$$\int_0^1 \log \Gamma(x+1/2)\cos 2n\pi x\, dx = (-1)^{n+1}\frac{si(n\pi)}{2\pi n}$$

Summation gives us

$$\sum_{n=1}^{\infty}\frac{1}{n}\int_0^1 \log \Gamma(x+1/2)\cos 2n\pi x\, dx = \frac{1}{2\pi}\sum_{n=1}^{\infty}(-1)^{n+1}\frac{si(n\pi)}{n^2}$$

and using [130, p.148] (see also (7.8a)) for $0 < x < 1$

$$\sum_{n=1}^{\infty}\frac{\cos nx}{n} = -\log[2\sin(x/2)]$$

we obtain

(6.130a) $$\int_0^1 \log \Gamma(x+1/2)\log[2\sin(\pi x)]\, dx = -\frac{1}{2\pi}\sum_{n=1}^{\infty}(-1)^{n+1}\frac{si(n\pi)}{n^2}$$

We may then use (6.117ji) to obtain a closed form result.

Similarly we obtain

(6.130b) $$\int_0^1 \log \Gamma(x+1/2)(2x^2 - 2x + 1)dx = \frac{3}{\pi^3}\sum_{n=1}^{\infty}(-1)^{n+1}\frac{si(n\pi)}{n^3}$$

and from (6.94ki) we have $\sum_{n=1}^{\infty}(-1)^n \frac{Si(n\pi)}{n^3} = -\frac{1}{18}\pi^3$. This then gives us

(6.130bi) $$\int_0^1 \log \Gamma(x+1/2)(2x^2 - 2x + 1)dx = \frac{\pi}{6} - \frac{3}{2\pi^2}\varsigma_a(3)$$

It should however be noted that the validity of the above formulae is dubious because of the restriction to the interval (0,1).

In [126, p.210] we have a formula for $\int_0^z x^2 \log \Gamma(x+a)\,dx$ which involves the integral of the triple gamma function $\int_0^z \log \Gamma_3(x+a)\,dx$ and therefore a little bit of algebraic manipulation will give us $\sum_{n=1}^{\infty} (-1)^{n+1} \frac{si(n\pi)}{n^3}$ in terms of $\int_0^1 \log \Gamma_3(x+1/2)\,dx$ (and having regard to (6.130bi) this will also involve the term $\frac{3}{2\pi^2} \varsigma_a(3)$).

We have from (6.108)

$$\int_0^1 \log \Gamma(1+x)\,dx = \int_0^1 \left[ -\gamma x + \sum_{n=2}^{\infty} (-1)^n \frac{\varsigma(n)}{n} x^n \right] \cot \pi x\,dx$$

$$= -\gamma \int_0^1 x \cot x\,dx + \sum_{n=2}^{\infty} \frac{\varsigma(n)}{n} \int_0^1 x^n \cot \pi x\,dx$$

From [59] we have the inversion formula for the Bernoulli polynomials

(6.130c) $$x^n = \frac{1}{n+1} \sum_{j=0}^n \binom{n+1}{j} B_j(x)$$

and we therefore have

(6.130d) $$\int_0^1 x^n \cot \pi x\,dx = \frac{1}{n+1} \sum_{j=0}^n \binom{n+1}{j} \int_0^1 B_j(x) \cot \pi x\,dx$$

The above integral may be evaluated using the method indicated in (6.31) et seq.

In passing, we mention that G&R [74, p.428] contains the integral

(6.130e) $$\int_0^{\pi/2} x^n \cot x\,dx = \left(\frac{\pi}{2}\right)^n \left[ \frac{1}{n} - 2 \sum_{k=1}^{\infty} \frac{\varsigma(2k)}{4^k(n+2k)} \right]$$

and in [126, p.35] we have

(6.130f) $$\int_0^{1/2} \log \Gamma(1+x)\,dx = -\frac{1}{2} - \frac{7}{24} \log 2 + \frac{1}{4} \log \pi + \frac{3}{2} \log A$$

(6.130g) $$\int_0^{1/4} \log \Gamma(1+x)\,dx = -\frac{1}{4} - \frac{3}{8} \log 2 + \frac{1}{8} \log \pi + \frac{9}{8} \log A + \frac{G}{4\pi}$$

Employing (6.108) we obtain

$$(6.130\mathrm{h}) \qquad \int_0^{1/2} \log \Gamma(1+x)\,dx = -\frac{1}{8}\gamma + \sum_{n=2}^{\infty}(-1)^n \frac{\varsigma(n)}{n(n+1)2^{n+1}}$$

$$= -\frac{1}{2} - \frac{7}{24}\log 2 + \frac{1}{4}\log \pi + \frac{3}{2}\log A$$

in accordance with [126, p.223]. We also have

$$(6.130\mathrm{i}) \qquad \int_0^{1/4} \log \Gamma(1+x)\,dx = -\frac{1}{32}\gamma + \sum_{n=2}^{\infty}(-1)^n \frac{\varsigma(n)}{n(n+1)4^{n+1}}$$

$$= -\frac{1}{4} - \frac{3}{8}\log 2 + \frac{1}{8}\log \pi + \frac{9}{8}\log A + \frac{G}{4\pi}$$

Differentiating the "corrected" version of (6.111) with respect to $a$

$$\int_0^1 \log \Gamma(x+a) \sin 2n\pi x\,dx = -\frac{1}{2\pi n}\left[\log a - \cos(2n\pi a)Ci(2n\pi a) + \sin(2n\pi a)si(2n\pi a)\right]$$

we obtain

(6.130j)

$$\int_0^1 \psi(x+a)\sin 2n\pi x\,dx = -\frac{1}{2\pi n}\begin{bmatrix}\dfrac{1}{a} - \cos(2n\pi a)\dfrac{\cos(2n\pi a)}{a} + 2n\pi \sin(2n\pi a)Ci(2n\pi a) \\ + \sin(2n\pi a)\dfrac{\sin(2n\pi a)}{a} + 2n\pi \cos(2n\pi a)si(2n\pi a)\end{bmatrix}$$

With $a=1$ we get

$$(6.130\mathrm{k}) \qquad \int_0^1 \psi(x+1)\sin 2n\pi x\,dx = -si(2n\pi)$$

Similarly, differentiating (6.118) we obtain

$$(6.130\mathrm{l}) \qquad \int_0^1 \psi(x+1)\cos 2n\pi x\,dx = -Ci(2n\pi)$$

Consideration may be given to a summation of the above two integrals.

Using integration by parts we see that for $a \geq 0$

$$\int_0^1 \psi(x+a)\sin 2n\pi x\, dx = \sin 2n\pi x \log\Gamma(x+a)\Big|_0^1 - 2n\pi \int_0^1 \Gamma(x+a)\cos 2n\pi x\, dx$$

We have for $a \geq 0$

$$\lim_{x\to 0}\sin 2n\pi x \log\Gamma(x+a) = \lim_{x\to 0}\frac{\sin 2n\pi x}{x} x\log\Gamma(x+a) = 0$$

and hence we have

(6.130m) $\quad \int_0^1 \psi(x+a)\sin 2n\pi x\, dx = -2n\pi \int_0^1 \Gamma(x+a)\cos 2n\pi x\, dx$

Similarly we have for $a > 0$

(6.130m) $\quad \int_0^1 \psi(x+a)\cos 2n\pi x\, dx = \log\Gamma(1+a) - \log\Gamma(a) + 2n\pi \int_0^1 \Gamma(x+a)\sin 2n\pi x\, dx$

Nielsen [104a, p.80] reports these integrals in the case $a = 1$ as

$$\int_0^1 \psi(x+1)\sin 2n\pi x\, dx = -2n\pi \int_0^1 \log\Gamma(x+1)\cos 2n\pi x\, dx$$

$$\int_0^1 \psi(x+1)\cos 2n\pi x\, dx = 2n\pi \int_0^1 \log\Gamma(x+1)\sin 2n\pi x\, dx$$

$\square$

We now revisit (6.94ji)

$$\log x + \log(1-x) - \log(2\sin\pi x) + 2 = -\frac{2}{\pi}\sum_{n=1}^{\infty}\frac{si(2n\pi)\cos 2n\pi x}{n}$$

and integrate to obtain

$$u\log u + (u-1)\log(1-u) - \int_0^u \log(2\sin\pi x)dx = -\frac{1}{\pi^2}\sum_{n=1}^{\infty}\frac{si(2n\pi)\sin 2n\pi u}{n^2}$$

From this we immediately see that

$$\int_0^{1/2}\log(\sin\pi x)dx = -\frac{1}{2}\log 2$$

$$\int_0^1 \log(\sin \pi x)dx = -\log 2$$

Using (4.3.158a) in Volume II(a) we get

(6.130o)
$$u\log u + (u-1)\log(1-u) + \varsigma'(-1,u) - \varsigma'(-1,1-u) = -\frac{1}{\pi^2}\sum_{n=1}^{\infty}\frac{si(2n\pi)\sin 2n\pi u}{n^2}$$

For $0 \leq u \leq 1$ this may also be expressed as

(6.130p)
$$u\log u + (u-1)\log(1-u) + \frac{1}{2\pi}\sum_{n=1}^{\infty}\frac{\sin 2n\pi u}{n^2} = -\frac{1}{\pi^2}\sum_{n=1}^{\infty}\frac{si(2n\pi)\sin 2n\pi u}{n^2}$$

and a further integration results in

(6.130q)

$$\frac{1}{2}u^2[\log u + \log(1-u)] + \left[\frac{1}{2}-u\right]\log(1-u) - \frac{1}{2}u(u-1) - \frac{1}{4\pi^2}\sum_{n=1}^{\infty}\frac{\cos 2n\pi u - 1}{n^3} =$$

$$\frac{1}{2\pi^3}\sum_{n=1}^{\infty}\frac{si(2n\pi)[\cos 2n\pi u - 1]}{n^3}$$

It is easily seen that the above equation is valid at $u=0$ and $u=1$. With $u=1/2$ we obtain

$$\frac{1}{8} + \frac{\varsigma_a(3) + \varsigma(3)}{4\pi^2} = \frac{1}{2\pi^3}\sum_{n=1}^{\infty}\frac{si(2n\pi)[\cos 2n\pi u - 1]}{n^3} = -\frac{1}{\pi^3}\sum_{n=1}^{\infty}\frac{si[2(2n+1)\pi]}{(2n+1)^3}$$

which we have seen before in (6.94jiv).

We now consider (6.94jii) for $0 < x < 1$

$$\log x - \log(1-x) + (\gamma + \log 2\pi)(1-2x) = 2\sum_{n=1}^{\infty}\frac{[Ci(2n\pi) - \log n]\sin 2n\pi x}{n\pi}$$

where integration results in

(6.130r)
$$u\log u + (1-u)\log(1-u) + (\gamma + \log 2\pi)(u - u^2) = -\frac{1}{\pi^2}\sum_{n=1}^{\infty}\frac{[Ci(2n\pi) - \log n][\cos 2n\pi u - 1]}{n^2}$$

With $u = 1/2$ we obtain

$$\frac{1}{4}(\gamma + \log 2\pi) - \log 2 = -\frac{1}{\pi^2}\sum_{n=1}^{\infty}\frac{[Ci(2n\pi) - \log n][(-1)^n - 1]}{n^2}$$

and we therefore have

(6.130s) $$\frac{1}{4}(\gamma + \log 2\pi) - \log 2 = \frac{2}{\pi^2}\sum_{n=0}^{\infty}\frac{[Ci[2(2n+1)\pi] - \log(2n+1)]}{(2n+1)^2}$$

We have

$$(1-2^{-s})\varsigma(s) = \sum_{n=0}^{\infty}\frac{1}{(2n+1)^s}$$

and therefore

$$(1-2^{-s})\varsigma'(s) + 2^{-s}\varsigma(s)\log 2 = -\sum_{n=0}^{\infty}\frac{\log(2n+1)}{(2n+1)^s}$$

We then obtain

$$\frac{1}{4}(\gamma + \log 2\pi) - \log 2 = \frac{2}{\pi^2}\sum_{n=0}^{\infty}\frac{Ci[2(2n+1)\pi]}{(2n+1)^2} + \frac{2}{\pi^2}\left[\frac{3}{4}\varsigma'(2) + \frac{1}{4}\varsigma(2)\log 2\right]$$

$$= \frac{2}{\pi^2}\sum_{n=0}^{\infty}\frac{Ci[2(2n+1)\pi]}{(2n+1)^2} + \frac{3\varsigma'(2)}{2\pi^2} + \frac{1}{12}\log 2$$

and using (F.7) this becomes

(6.130si) $$\frac{2}{\pi^2}\sum_{n=0}^{\infty}\frac{Ci[2(2n+1)\pi]}{(2n+1)^2} = \frac{1}{4} - 3\varsigma'(-1) - \frac{13}{12}\log 2$$

□

We recall (6.117ca)

$$\log \Gamma(a+1/2) =$$

$$\frac{1}{2}\log(2\pi) + a\log a - a + \frac{1}{\pi}\sum_{n=1}^{\infty}\frac{(-1)^n}{n}[\sin(2n\pi a)Ci(2n\pi a) - \cos(2n\pi a)si(2n\pi a)]$$

$$= \frac{1}{2}\log(2\pi) + a\log a - a + \frac{1}{\pi}\sum_{n=1}^{\infty}\frac{(-1)^n}{n}[\sin(2n\pi a)Ci(2n\pi a) - \cos(2n\pi a)Si(2n\pi a)]$$

$$+ \frac{1}{2}\sum_{n=1}^{\infty}\frac{(-1)^n \cos(2n\pi a)}{n}$$

and we may note from [130, p.148] that

$$-\log[2\cos\pi a] = \sum_{n=1}^{\infty} \frac{(-1)^n \cos(2n\pi a)}{n} \qquad \text{for } -\frac{1}{2} < a < \frac{1}{2}$$

We have the integral

$$\int [\sin(2n\pi a)Ci(2n\pi a) - \cos(2n\pi a)Si(2n\pi a)]da =$$

$$-\frac{1}{2n\pi}[\cos(2n\pi a)Ci(2n\pi a) + \sin(2n\pi a)Si(2n\pi a)] + \frac{1}{2n\pi}\log(2n\pi a)$$

Therefore using (6.117ei) we see that

$$\int_0^z [\sin(2n\pi a)Ci(2n\pi a) - \cos(2n\pi a)Si(2n\pi a)]da =$$

$$-\frac{1}{2n\pi}[\cos(2n\pi z)Ci(2n\pi z) + \sin(2n\pi z)Si(2n\pi z)] + \frac{\gamma + \log(2n\pi z)}{2n\pi}$$

and hence we get

(6.130t) $\int_0^z \log\Gamma(a+1/2)\,da =$

$$\frac{1}{2}\log(2\pi)z + \frac{1}{2}z^2\log z - \frac{3}{4}z^2 - \frac{1}{2\pi^2}\sum_{n=1}^{\infty}\frac{(-1)^n}{n^2}[\cos(2n\pi z)Ci(2n\pi z) + \sin(2n\pi z)Si(2n\pi z)]$$

$$+\frac{1}{2\pi^2}\sum_{n=1}^{\infty}\frac{(-1)^n[\gamma+\log(2n\pi z)]}{n^2} + \frac{1}{4\pi}\sum_{n=1}^{\infty}\frac{(-1)^n\sin(2n\pi z)}{n^2}$$

$$= \frac{1}{2}\log(2\pi)z + \frac{1}{2}z^2\log z - \frac{3}{4}z^2 - \frac{1}{2\pi^2}\sum_{n=1}^{\infty}\frac{(-1)^n}{n^2}[\cos(2n\pi z)Ci(2n\pi z) + \sin(2n\pi z)Si(2n\pi z)]$$

$$-\frac{[\gamma+\log(2\pi z)]}{2\pi^2}\varsigma_a(2) + \frac{1}{2\pi^2}\varsigma_a'(2) + \frac{1}{4\pi}\sum_{n=1}^{\infty}\frac{(-1)^n\sin(2n\pi z)}{n^2}$$

With $z = 1$ we get

$$\int_0^1 \log\Gamma(a+1/2)\,da =$$

$$= \frac{1}{2}\log(2\pi) - \frac{3}{4} - \frac{1}{2\pi^2} \sum_{n=1}^{\infty} \frac{(-1)^n Ci(2n\pi)}{n^2} - \frac{[\gamma + \log(2\pi)]}{2\pi^2} \varsigma_a(2) + \frac{1}{2\pi^2} \varsigma_a'(2)$$

Since

$$\int_0^z \log \Gamma(a+1/2)\, da = \int_{1/2}^{z+1/2} \log \Gamma(a)\, da$$

$$= \int_0^{z+1/2} \log \Gamma(a)\, da - \int_0^{1/2} \log \Gamma(a)\, da$$

we may use Alexeiewsky's theorem (4.3.85) in Volume II(a)

$$\int_0^u \log \Gamma(x)\, dx = \frac{1}{2} u(1-u) + \frac{u}{2} \log(2\pi) - \log G(u+1) + u \log \Gamma(u)$$

to obtain

$$\int_0^z \log \Gamma(a+1/2)\, da = \frac{1}{2}\left(z+\frac{1}{2}\right)\left(\frac{1}{2}-z\right) + \frac{1}{2}\left(z+\frac{1}{2}\right)\log(2\pi) - \log G\left(z+\frac{3}{2}\right)$$

$$+ \left(z+\frac{1}{2}\right)\log \Gamma\left(z+\frac{1}{2}\right) - \frac{1}{8} - \frac{1}{4}\log(2\pi) + \log G\left(\frac{3}{2}\right) - \frac{1}{2}\log \Gamma\left(\frac{1}{2}\right)$$

$$= \frac{1}{2}\left(z+\frac{1}{2}\right)\left(\frac{1}{2}-z\right) + \frac{1}{2}\left(z+\frac{1}{2}\right)\log(2\pi) - \log G\left(z+\frac{3}{2}\right)$$

$$+ \left(z+\frac{1}{2}\right)\log \Gamma\left(z+\frac{1}{2}\right) - \frac{1}{8} - \frac{1}{4}\log \pi - \frac{5}{24}\log 2 + \frac{3}{2}\varsigma'(-1)$$

$$= \frac{1}{2}\left(z+\frac{1}{2}\right)\left(\frac{1}{2}-z\right) + \frac{1}{2} z \log(2\pi) - \log G\left(z+\frac{3}{2}\right)$$

$$+ \left(z+\frac{1}{2}\right)\log \Gamma\left(z+\frac{1}{2}\right) - \frac{1}{8} - \frac{1}{4}\log \pi + \frac{1}{24}\log 2 + \frac{3}{2}\varsigma'(-1)$$

Therefore we have

$$\frac{1}{2}\left(z+\frac{1}{2}\right)\left(\frac{1}{2}-z\right) + \frac{1}{2} z \log(2\pi) - \log G\left(z+\frac{3}{2}\right) + \left(z+\frac{1}{2}\right)\log \Gamma\left(z+\frac{1}{2}\right)$$

$$-\frac{1}{8} + \frac{1}{24}\log 2 + \frac{3}{2}\varsigma'(-1) =$$

$$\frac{1}{2}z\log(2\pi) + \frac{1}{2}z^2\log z - \frac{3}{4}z^2 - \frac{1}{2\pi^2}\sum_{n=1}^{\infty}\frac{(-1)^n}{n^2}[\cos(2n\pi z)Ci(2n\pi z) + \sin(2n\pi z)Si(2n\pi z)]$$

$$-\frac{[\gamma + \log(2\pi z)]}{2\pi^2}\varsigma_a(2) + \frac{1}{2\pi^2}\varsigma_a'(2) + \frac{1}{4\pi}\sum_{n=1}^{\infty}\frac{(-1)^n \sin(2n\pi z)}{n^2}$$

which may be slightly simplified to

(6.130u)

$$\frac{1}{2}\left(z+\frac{1}{2}\right)\left(\frac{1}{2}-z\right) - \log G\left(z+\frac{3}{2}\right) + \left(z+\frac{1}{2}\right)\log \Gamma\left(z+\frac{1}{2}\right)$$

$$-\frac{1}{8} + \frac{1}{24}\log 2 + \frac{3}{2}\varsigma'(-1) =$$

$$\frac{1}{2}z^2\log z - \frac{3}{4}z^2 - \frac{1}{2\pi^2}\sum_{n=1}^{\infty}\frac{(-1)^n}{n^2}[\cos(2n\pi z)Ci(2n\pi z) + \sin(2n\pi z)Si(2n\pi z)]$$

$$-\frac{[\gamma + \log(2\pi z)]}{2\pi^2}\varsigma_a(2) + \frac{1}{2\pi^2}\varsigma_a'(2) + \frac{1}{4\pi}\sum_{n=1}^{\infty}\frac{(-1)^n \sin(2n\pi z)}{n^2}$$

Letting $z = 1/2$ we get

$$-\frac{1}{8} + \frac{1}{24}\log 2 + \frac{3}{2}\varsigma'(-1) = -\frac{1}{8}\log 2 - \frac{3}{16} - \frac{1}{2\pi^2}\sum_{n=1}^{\infty}\frac{Ci(n\pi)}{n^2} - \frac{[\gamma + \log \pi]}{2\pi^2}\varsigma_a(2) + \frac{1}{2\pi^2}\varsigma_a'(2)$$

$$\frac{1}{2\pi^2}\sum_{n=1}^{\infty}\frac{Ci(n\pi)}{n^2} = -\frac{1}{6}\log 2 - \frac{1}{16} - \frac{3}{2}\varsigma'(-1) - \frac{[\gamma + \log \pi]}{2\pi^2}\varsigma_a(2) + \frac{1}{2\pi^2}\varsigma_a'(2)$$

$$\sum_{n=1}^{\infty}\frac{Ci(n\pi)}{n^2} = -\frac{1}{3}\pi^2\log 2 - \frac{1}{8}\pi^2 - 3\pi^2\varsigma'(-1) - [\gamma + \log \pi]\varsigma_a(2) + \varsigma_a'(2)$$

From (F.8i) in Volume VI we have

$$\varsigma_a'(2) = \varsigma'(-1)\pi^2 - \frac{\pi^2}{12} + \frac{\pi^2}{12}[\gamma + \log \pi] + \frac{\pi^2}{6}\log 2$$

and we then obtain

$$\sum_{n=1}^{\infty}\frac{Ci(n\pi)}{n^2} = -\frac{1}{6}\pi^2\log 2 - \frac{17}{24}\pi^2 - 2\pi^2\varsigma'(-1)$$

This is consistent with the following two formulae

$$\sum_{n=1}^{\infty}(-1)^n \frac{Ci(n\pi)}{n^2} = \frac{1}{6}\pi^2 \log 2 + \frac{\pi^2}{24} + \pi^2 \varsigma'(-1)$$

$$\sum_{n=1}^{\infty} \frac{Ci(2n\pi)}{n^2} = 2\pi^2\left(\log A - \frac{1}{4}\right) = -2\pi^2\left[\frac{1}{6} + \varsigma'(-1)\right] = -\frac{1}{3}\pi^2 - 2\pi^2 \varsigma'(-1)$$

after employing the identity

$$\sum_{n=1}^{\infty}(-1)^n \frac{Ci(n\pi)}{n^2} + \sum_{n=1}^{\infty} \frac{Ci(n\pi)}{n^2} = \frac{1}{2}\sum_{n=1}^{\infty} \frac{Ci(2n\pi)}{n^2}$$

□

We recall Elizalde's formula (6.117a)

$$\varsigma'(-1,t) =$$

$$-\varsigma(-1,t)\log t - \frac{1}{4}t^2 + \frac{1}{12} - \frac{1}{2\pi^2}\sum_{n=1}^{\infty}\frac{1}{n^2}[\cos(2n\pi t)Ci(2n\pi t) + \sin(2n\pi t)si(2n\pi t)]$$

and in (4.4.229i) in Volume IV we showed that

$$\varsigma'(-1,t) = -2\sum_{n=1}^{\infty}\frac{\log(2\pi n) + \gamma - 1}{(2\pi n)^2}\cos 2n\pi t + \sum_{n=1}^{\infty}\frac{\sin 2n\pi t}{4\pi n^2}$$

In (4.4.229m) we then showed that

$$-\varsigma(-1,t)\log t - \frac{1}{4}t^2 + \frac{1}{12} - \frac{1}{2\pi^2}\sum_{n=1}^{\infty}\frac{1}{n^2}[\cos(2n\pi t)Ci(2n\pi t) + \sin(2n\pi t)Si(2n\pi t)] =$$

$$-\frac{1}{2\pi^2}\sum_{n=1}^{\infty}\frac{\log(2\pi) + \log n + \gamma - 1}{n^2}\cos 2n\pi t$$

□

Some of the identities involving the sine and cosine integrals are set out below (and it seems that this is a subject ripe for more systematic research).

(6.120) $$\sum_{n=1}^{\infty}\frac{si(2n\pi)}{n} = \frac{\pi}{2}\log(2\pi) - \pi$$

(6.121) $$\sum_{n=1}^{\infty}(-1)^n \frac{si(2n\pi)}{n} = \frac{\pi}{2}\log 2 - \pi$$

(6.94ai) $$\sum_{n=1}^{\infty}(-1)^n \frac{Si(2n\pi)}{n} = -\pi$$

(6.94a) $$\sum_{n=1}^{\infty}(-1)^n \frac{si(nx)}{n} = \frac{\pi}{2}\log 2 - \frac{x}{2}$$

(6.94ji) $$\log x + \log(1-x) - \log(2\sin\pi x) + 2 = -\frac{2}{\pi}\sum_{n=1}^{\infty}\frac{si(2n\pi)\cos 2n\pi x}{n}$$

(6.94jii) $$\log x - \log(1-x) + (\gamma + \log 2\pi)(1-2x) = 2\sum_{n=1}^{\infty}\frac{[Ci(2n\pi) - \log n]\sin 2n\pi x}{n\pi}$$

(6.117c) $$\log \Gamma(a) =$$

$$\frac{1}{2}\log(2\pi) + \left(a - \frac{1}{2}\right)\log a - a + \frac{1}{\pi}\sum_{n=1}^{\infty}\frac{1}{n}[\sin(2n\pi a)Ci(2n\pi a) - \cos(2n\pi a)si(2n\pi a)]$$

$$= \frac{1}{2}\log(2\pi) + \left(a - \frac{1}{2}\right)\log a - a + \frac{1}{2}\sum_{n=1}^{\infty}\frac{\cos(2n\pi a)}{n}$$

$$+ \frac{1}{\pi}\sum_{n=1}^{\infty}\frac{1}{n}[\sin(2n\pi a)Ci(2n\pi a) - \cos(2n\pi a)Si(2n\pi a)]$$

(6.117ca) $$\log \Gamma(a+1/2) =$$

$$\frac{1}{2}\log(2\pi) + a\log a - a + \frac{1}{\pi}\sum_{n=1}^{\infty}\frac{(-1)^n}{n}[\sin(2n\pi a)Ci(2n\pi a) - \cos(2n\pi a)si(2n\pi a)]$$

$$\sum_{n=1}^{\infty}\frac{si(n\pi)}{n} = \frac{\pi}{2}\log\pi - \frac{\pi}{2}$$

$$\sum_{n=1}^{\infty}(-1)^n \frac{si(2n\pi)}{n} = \frac{3}{2}\pi\log 2 - \pi$$

(6.117cii) $$\sum_{n=1}^{\infty}\frac{si(2n\pi)}{n} = \frac{\pi}{2}\log(2\pi) - \pi$$

(6.117ciii) $$\sum_{n=1}^{\infty}(-1)^n \frac{si(n\pi)}{n} = \frac{\pi}{2}\log 2 - \frac{\pi}{2}$$

(6.117civ) $$\sum_{n=1}^{\infty}(-1)^n \frac{si(nx)}{n} = \frac{\pi}{2}\log 2 - \frac{x}{2}$$

(6.92b) $$\sum_{n=1}^{\infty}\frac{\varsigma(2n)}{(2n+1)^2} = \frac{1}{2} - \frac{1}{4\pi}\sum_{n=1}^{\infty}\frac{Si(2n\pi)}{n^2} \quad (?)$$

(6.117j) $$\sum_{n=1}^{\infty}\frac{Si(2n\pi)}{n^2} = \frac{1}{2}\pi^2 \log(2\pi) + 2\pi^2 \log A + \frac{5\pi^2}{36}$$

(6.117jii) $$\sum_{n=1}^{\infty}\frac{(-1)^n Si(n\pi)}{n^2} = -\frac{7\varsigma(3)}{4} + \pi^2 \log A - \frac{1}{4}\pi^2$$

(6.117) $$\sum_{n=1}^{\infty}\frac{Ci(2n\pi)}{n^2} = 2\pi^2\left(\log A - \frac{1}{4}\right) = -2\pi^2\left[\frac{1}{6} + \varsigma'(-1)\right]$$

(6.117b) $$\sum_{n=1}^{\infty}(-1)^n \frac{Ci(n\pi)}{n^2} = \frac{1}{6}\pi^2 \log 2 + \frac{\pi^2}{24} + \pi^2 \varsigma'(-1)$$

(6.130v) $$\sum_{n=1}^{\infty}\frac{Ci(n\pi)}{n^2} = -\frac{1}{6}\pi^2 \log 2 - \frac{17}{24}\pi^2 - 2\pi^2 \varsigma'(-1)$$

(6.117fi) $$\int_0^x \log \Gamma(a)\, da =$$

$$= \frac{1}{2}x\log(2\pi) + \frac{1}{4}x\bigl[2 - x + 2(x-1)\log x\bigr] - \frac{1}{2}x^2 + \frac{1}{4\pi}\sum_{n=1}^{\infty}\frac{\sin(2n\pi x)}{n^2} + \frac{1}{12}\log x$$

$$- \frac{1}{2\pi^2}\sum_{n=1}^{\infty}\frac{1}{n^2}\bigl[\cos(2n\pi x)Ci(2n\pi x) + \sin(2n\pi x)Si(2n\pi x)\bigr] + \frac{1}{12} - \varsigma'(-1)$$

(6.117i) $$\frac{1}{2}x\log(2\pi) + \frac{1}{4}x\bigl(2 - x + 2(x-1)\log x\bigr) - \frac{1}{2}x^2 + \frac{1}{4\pi}x\sum_{n=1}^{\infty}\frac{\sin(2n\pi x)}{n^2}$$

$$+ \frac{1}{12}x(x-1)(2x-1)\log x - \frac{5}{36}x^3 + \frac{1}{8}x^2 + \frac{1}{12}x$$

$$+ \frac{1}{12}B_3(x) - \frac{1}{2}\varsigma'(-2, x) + \frac{1}{2}\varsigma'(-2)$$

$$+ \frac{x}{2\pi^2}\sum_{n=1}^{\infty}\frac{1}{n^2}[\sin(2n\pi x)Ci(2n\pi x) - \cos(2n\pi x)Si(2n\pi x)]$$

$$= \frac{1}{2}\left(\frac{1}{4} - 2\log A\right)x + \frac{1}{2}\left(\frac{1}{2}\log(2\pi) - \frac{1}{4}\right)x^2$$

$$+ \frac{1}{2}x^2 \log \Gamma(x) - \frac{1}{2}\log G(1+x) - \log \Gamma_3(1+x)$$

(6.117k) $\quad -\log G(1+x) + x\log \Gamma(x) = \frac{1}{4}x\left[-x + 2(x-1)\log x\right] + \frac{Cl_2(2\pi x)}{4\pi}$

$$-\frac{1}{2\pi^2}\sum_{n=1}^{\infty}\frac{1}{n^2}\left[\cos(2n\pi x)Ci(2n\pi x) + \sin(2n\pi x)Si(2n\pi x)\right] + \frac{1}{12} - \varsigma'(-1)$$

(6.117e) $\qquad \sum_{n=1}^{\infty}\frac{Si(2n\pi)}{n^3} = \frac{1}{18}\pi^3$

(6.91) $\qquad \sum_{n=1}^{\infty}\frac{Si(n\pi)}{n^3} = \frac{5}{72}\pi^3$

(6.94ki) $\qquad \sum_{n=1}^{\infty}(-1)^n \frac{Si(n\pi)}{n^3} = -\frac{1}{18}\pi^3$

(6.94k) $\qquad \sum_{n=0}^{\infty}\frac{Si[2(2n+1)\pi]}{(2n+1)^3} = -\frac{1}{8}\pi^3$

$$\frac{11}{36}t^3 - \frac{1}{24}t - \frac{1}{4}t^2 - \left(t - \frac{1}{2}\right)\log G(1+t) + \log \Gamma_3(1+t) + \frac{1}{2}t^2 \log \Gamma(t)$$

(6.117m)

$$= \frac{1}{12}t(t-1)(2t-1)\log t - \frac{1}{8\pi^2}\sum_{n=1}^{\infty}\frac{\cos(2n\pi t) - 1}{n^3}$$

$$-\frac{1}{4\pi^3}\sum_{n=1}^{\infty}\frac{1}{n^3}\left[\sin(2n\pi t)Ci(2n\pi t) - \cos(2n\pi t)Si(2n\pi t)\right]$$

(6.117d)

$$\frac{1}{12}x(x-1)(2x-1)\log x - \frac{5}{36}x^3 + \frac{1}{8}x^2 - \frac{1}{8\pi^2}\sum_{n=1}^{\infty}\frac{\cos(2n\pi x) - 1}{n^3}$$

$$-\frac{1}{4\pi^3}\sum_{n=1}^{\infty}\frac{1}{n^3}[\sin(2n\pi x)Ci(2n\pi x) - \cos(2n\pi x)Si(2n\pi x)]$$

$$= -\frac{1}{12}B_3(x) + \frac{1}{2}\varsigma'(-2,x) - \frac{1}{2}\varsigma'(-2)$$

(6.94o) $$\frac{11}{576}\pi^4 + \sum_{n=1}^{\infty}\frac{Ci(n\pi)}{n^4} = \varsigma(4)[\gamma + \log \pi] - \varsigma'(4)$$

**Example 16:**

A more elementary example will now be considered. Letting $p(x) = e^x - 1$ in (6.5) we get

(6.131) $$\frac{1}{2}\int_0^t (e^x - 1)\,dx = \sum_{n=0}^{\infty}\int_0^t (e^x - 1)\cos \alpha n x\,dx$$

We have

$$\int_0^t (e^x - 1)\cos nx\,dx = \frac{1}{1+n^2}(e^t \cos nt - 1) + \frac{[n^2(e^t - 1) - 1]}{n(1+n^2)}\sin nt$$

and in particular we get

$$\int_0^\pi (e^x - 1)\cos nx\,dx = \frac{1}{1+n^2}\left[e^\pi(-1)^n - 1\right]$$

Therefore we have

(6.132) $$e^\pi - \pi - 1 = 2\sum_{n=1}^{\infty}\frac{1}{1+n^2} - 2e^\pi \sum_{n=1}^{\infty}\frac{(-1)^n}{1+n^2}$$

The above analysis may be easily verified by using the following identities which were originally derived by Euler [126aa, p.96] (and see Knopp's book [90, p.207]). An elementary proof was also given by Efthimiou [58a].

(6.133) $$\sum_{n=1}^{\infty}\frac{1}{n^2 - \alpha^2} = \frac{1}{2\alpha^2} - \frac{\pi}{2\alpha \tan \pi\alpha}$$

(6.134) $$\frac{\alpha \pi}{\sin \pi\alpha} = 1 + 2\alpha^2 \sum_{n=1}^{\infty}\frac{(-1)^{n+1}}{n^2 - \alpha^2}$$

Letting $\alpha = \sqrt{-1}$ we obtain

(6.135) $$\sum_{n=1}^{\infty}\frac{1}{n^2 + 1} = -\frac{1}{2} - \frac{\pi}{2\tanh \pi}$$

(6.136) $$\frac{\pi}{\sinh \pi} = 1 - 2\sum_{n=1}^{\infty} \frac{(-1)^{n+1}}{n^2 + 1}$$

and (6.132) is confirmed by a simple substitution of (6.135) and (6.136).

After multiplying (6.133) by $2\alpha$ and integrating we obtain

(6.137) $$\sum_{n=1}^{\infty} \int_a^b \frac{2\alpha}{n^2 - \alpha^2} d\alpha = \sum_{n=1}^{\infty} \int_a^b \frac{2\alpha/n^2}{1 - \alpha^2/n^2} d\alpha = \int_a^b \left[ \frac{1}{\alpha} - \frac{\pi}{\tan \pi \alpha} \right] d\alpha$$

This implies that

$$\sum_{n=1}^{\infty} \log \left[ \frac{1 - a^2/n^2}{1 - b^2/n^2} \right] = \log \frac{b}{a} - \log \frac{\sin \pi b}{\sin \pi a}$$

and with a little rearrangement we get

$$\log \frac{a \sin \pi b}{b \sin \pi a} \prod_{n=1}^{\infty} \left[ \frac{1 - a^2/n^2}{1 - b^2/n^2} \right] = 0 = \log 1$$

Hence we have

$$\frac{\sin \pi a}{a} \frac{1}{\prod_{n=1}^{\infty} \left[ 1 - a^2/n^2 \right]} = \frac{\sin \pi b}{b} \frac{1}{\prod_{n=1}^{\infty} \left[ 1 - b^2/n^2 \right]}$$

Therefore, in the limit as $b \to 0$ we have

$$\frac{\sin \pi a}{a} \frac{1}{\prod_{n=1}^{\infty} \left[ 1 - a^2/n^2 \right]} = \pi$$

and hence we have derived Euler's formula (1.6d)

$$\sin x = x \prod_{n=1}^{\infty} \left[ 1 - x^2/n^2\pi^2 \right]$$

We have (with the following steps being justified by absolute convergence for $|x| < 1$)

$$\log \sin \pi x = \log \pi x + \sum_{n=1}^{\infty} \log \left( 1 - \frac{x^2}{n^2} \right)$$

$$= \log \pi x - \sum_{n=1}^{\infty} \sum_{k=1}^{\infty} \frac{x^{2k}}{kn^{2k}}$$

$$= \log \pi x - \sum_{k=1}^{\infty} \frac{x^{2k}}{k} \sum_{n=1}^{\infty} \frac{1}{n^{2k}}$$

(6.138) $$\log \sin \pi x = \log \pi x - \sum_{k=1}^{\infty} \frac{\varsigma(2k)}{k} x^{2k}$$

Upon integrating the above equation from $x = 0$ to $x = 1/2$ we obtain the series representation (6.92). In passing, it may be noted that letting $x = 1/2$ in Euler's formula (1.6d) results in the Wallis identity.

The well-known Maclaurin expansion for the log gamma function was referred to in (6.108)

$$\log \Gamma(1+x) = -\gamma x + \sum_{n=2}^{\infty} (-1)^n \frac{\varsigma(n)}{n} x^n$$

Letting $x \to -x$ we get for $x < 1$

$$\log \Gamma(1-x) = \gamma x + \sum_{n=2}^{\infty} \frac{\varsigma(n)}{n} x^n$$

and hence we have

$$\log \Gamma(1+x) + \log \Gamma(1-x) = \sum_{n=2}^{\infty} \left[1 + (-1)^n\right] \frac{\varsigma(n)}{n} x^n$$

Thus we obtain

$$\log x + \log \Gamma(x)\Gamma(1-x) = \sum_{n=2}^{\infty} \left[1 + (-1)^n\right] \frac{\varsigma(n)}{n} x^n$$

and using Euler's reflection formula we again get

$$\log \sin \pi x = \log \pi x - \sum_{k=1}^{\infty} \frac{\varsigma(2k)}{k} x^{2k}$$

Differentiating this we obtain (6.48)

(6.139) $$\pi \cot \pi x = \frac{1}{x} - 2 \sum_{k=1}^{\infty} \varsigma(2k) x^{2k-1}$$

Combining (6.139) and (C.42a) enables us to determine $\varsigma(2k)$ in terms of the Bernoulli numbers. Alternatively, we may write

(6.139a) $$\pi \cot \pi x = i\pi x + \frac{2i\pi x}{e^{2i\pi x}-1} = 1 + \sum_{k=2}^{\infty} B_k \frac{(2i\pi x)^k}{k!}$$

and then equate coefficients with (6.139).

**Example 17:**

Letting $p(x) = x \sin \alpha x$ in (6.5) we have (for $\alpha \notin N$)

$$-\frac{1}{2}\int_0^\pi x \sin \alpha x \, dx = \sum_{n=1}^{\infty} \int_0^\pi x \sin \alpha x \cos nx \, dx$$

We have

$$\int x \sin \alpha x \sin nx \, dx = \frac{1}{2}\left[\frac{\cos(\alpha-n)x}{(\alpha-n)^2} - \frac{\cos(\alpha+n)x}{(\alpha+n)^2} + \frac{x\sin(\alpha-n)x}{\alpha-n} - \frac{x\sin(\alpha+n)x}{\alpha+n}\right]$$

and therefore

$$\int_0^\pi x \sin \alpha x \sin nx \, dx =$$

$$\frac{1}{2}\left[\frac{\cos(\alpha-n)\pi}{(\alpha-n)^2} - \frac{\cos(\alpha+n)\pi}{(\alpha+n)^2} + \frac{\pi\sin(\alpha-n)\pi}{\alpha-n} - \frac{\pi\sin(\alpha+n)\pi}{\alpha+n} - \frac{1}{(\alpha-n)^2} + \frac{1}{(\alpha+n)^2}\right]$$

$$= \frac{1}{2}\left[\frac{(-1)^n \cos \alpha\pi}{(\alpha-n)^2} - \frac{(-1)^n \cos \alpha\pi}{(\alpha+n)^2} + \frac{\pi(-1)^n 2n \sin \alpha\pi}{\alpha^2-n^2} - \frac{1}{(\alpha-n)^2} + \frac{1}{(\alpha+n)^2}\right]$$

We then have

(6.140) $$\frac{\pi \cos \alpha\pi}{\alpha} - \frac{\sin \alpha\pi}{\alpha^2} =$$

$$\cos \alpha\pi \sum_{n=1}^{\infty}\left[\frac{(-1)^n}{(\alpha-n)^2} - \frac{(-1)^n}{(\alpha+n)^2}\right] + 2\pi \sin \alpha\pi \sum_{n=1}^{\infty}\left[\frac{(-1)^n n}{\alpha^2-n^2}\right] + \sum_{n=1}^{\infty}\left[\frac{1}{(\alpha+n)^2} - \frac{1}{(\alpha-n)^2}\right]$$

## 7. SOME APPLICATIONS OF THE RIEMANN-LEBESGUE LEMMA

The Riemann-Lebesgue lemma is primarily used as a tool in the rigorous proof of the theory of Fourier series. It does however have further uses as illustrated by the derivation of various integral identities in this series of papers, in particular (2.23), (2.24) and the other identities recorded in Section 6. Some additional applications are considered below: these applications, whilst unlikely to be new, do not appear to have received much exposure in the mathematical literature.

**Example 1:**

Let us consider the familiar trigonometric identities which we used in Section 6.

(7.1) $$\frac{1}{2} = \sum_{n=0}^{N} \cos\alpha nt - \frac{\sin\alpha(N+1/2)t}{2\sin(\alpha t/2)}$$

(7.2) $$\frac{1}{2}\cot(\alpha t/2) = \sum_{n=1}^{N} \sin\alpha nt + \frac{\cos\alpha(N+1/2)t}{2\sin(\alpha t/2)}$$

In (7.1) let $\alpha = 1$ and integrate over the interval $J = [a, x]$ where $0 < a < x < 2\pi$ to obtain

(7.3) $$\frac{1}{2}\int_a^x dt = \sum_{n=0}^{N}\int_a^x \cos nt\, dt - \int_a^x \frac{\sin(N+1/2)t}{2\sin(t/2)} dt$$

This gives

(7.4) $$\frac{1}{2}(a-x) = \sum_{n=1}^{N} \frac{\sin nx - \sin na}{n} - \int_a^x \frac{\sin(N+1/2)t}{2\sin(t/2)} dt$$

The denominator of the integral in (7.4) is bounded in the specified interval $J$ and hence, in accordance with the Riemann-Lebesgue lemma ("RLL"), we have

$$\lim_{N\to\infty} \int_a^x \frac{\sin(N+1/2)t}{2\sin(t/2)} dt = 0$$

and therefore with $a = \pi$ we obtain the familiar Fourier series

(7.5) $$\frac{1}{2}(\pi - x) = \sum_{n=1}^{\infty} \frac{\sin nx}{n} \quad (0 < x < 2\pi)$$

One of the attractions of this exposition is that it clearly shows why the Fourier series does not represent the function at the specified end points (because the convergence conditions for the RLL would not then be met). In addition, we have an exact relationship for the partial sum and hence this method could be employed for numerical approximations if required.

It should also be noted that it has not been necessary to specify from the outset that the function was periodic outside the interval [0,2π]: the author has never been able to truly fathom why such periodicity should be a requirement in the Dirichlet conditions used in the rigorous treatment of Fourier series (see for example Apostol [13, p.319], Bressoud [34, p.223], Titchmarsh [128, p.401] and Whittaker and Watson [135, p.164]). Why should the nature of a real function outside of its defined finite domain have any effect on the pointwise convergence of the corresponding Fourier series? I would welcome further elucidation on this point.

**Example 2:**

Let us now consider (7.2) with $\alpha = 1$ and integrate over the interval $[a, x]$ where $0 < a < x < 2\pi$ to obtain

(7.6) $$\frac{1}{2}\int_a^x \cot(t/2)dt = \sum_{n=1}^N \int_a^x \sin nt\, dt + \int_a^x \frac{\cos(N+1/2)t}{2\sin(t/2)} dt$$

In the limit as $N \to \infty$ we obtain another Fourier series

(7.7) $$\log\sin(x/2) - \log\sin(a/2) = -\sum_{n=1}^\infty \frac{\cos nx - \cos na}{n}$$

and, upon letting $a = \pi/2$ this becomes the familiar trigonometric series shown in Carslaw's book [41, p.241]

(7.8) $$\log[2\sin(x/2)] = -\sum_{n=1}^\infty \frac{\cos nx}{n} \quad (0 < x < 2\pi)$$

The derivation of this particular Fourier series is certainly simpler than the standard method (as shown in Carslaw [41, p.131] for example) and it is clear that this method could be used to provide simple derivations of a whole variety of Fourier series.

Integrating (7.8) results in

(7.8a) $$\int_0^u \log[2\sin(x/2)]dx = -\sum_{n=1}^\infty \frac{\sin nu}{n^2} \quad (0 \le u \le 2\pi)$$

which is shown in [130, p.148].

**Example 3:**

Now let $\alpha = 2$ in (7.2) to obtain

(7.9) $$\cot t = 2\sum_{n=1}^N \sin 2nt + \frac{\cos(2N+1)t}{\sin t}$$

(7.10) $$= 2\sum_{n=1}^N \sin 2nt + \frac{\cos 2Nt \cos t - \sin 2Nt \sin t}{\sin t}$$

Rearranging this gives

(7.11) $$\cos 2Nt \cot t = \cot t - 2\sum_{n=1}^{N} \sin 2nt + \sin 2Nt$$

Integrating (7.11) we obtain

(7.12) $$\int \cos 2Nt \cot t \, dt = \log \sin t + \sum_{n=1}^{N} \frac{\cos 2nt}{n} - \frac{\cos 2Nt}{2N} + c$$

With integration by parts we have

(7.13) $$\int_{0}^{\pi/2} t^2 \cos 2Nt \cot t \, dt = t^2 \left\{ \log \sin t + \sum_{n=1}^{N} \frac{\cos 2nt}{n} - \frac{\cos 2Nt}{2N} \right\} \Bigg|_{0}^{\pi/2}$$

$$- 2\int_{0}^{\pi/2} t \left\{ \log \sin t + \sum_{n=1}^{N} \frac{\cos 2nt}{n} - \frac{\cos 2Nt}{2N} \right\} dt$$

(7.14)
$$= \frac{\pi^2}{4} \sum_{n=1}^{N} \frac{(-1)^n}{n} - \frac{\pi^2}{8} \frac{(-1)^N}{N} - 2\int_{0}^{\pi/2} t \log \sin t \, dt - \sum_{n=1}^{N} \frac{(-1)^n}{2n^3} + \sum_{n=1}^{N} \frac{1}{2n^3} + \frac{(-1)^N}{4N^3} - \frac{1}{4N^3}$$

As $N \to \infty$ the left hand side of (7.13) $\to 0$ (using RLL), and the right hand side of (7.14) becomes

(7.15) $$= \frac{\pi^2}{4} \log 2 - 2\int_{0}^{\pi/2} t \log \sin t \, dt - \sum_{n=1}^{\infty} \frac{(-1)^n}{2n^3} + \sum_{n=1}^{\infty} \frac{1}{2n^3} = 0$$

Therefore we again end up with the familiar Euler integral identity (1.11)

(7.16) $$\int_{0}^{\pi/2} t \log \sin t \, dt = \frac{7}{16}\varsigma(3) - \frac{\pi^2}{8} \log 2$$

In fact, by the same method we can easily show that

$$\int_{0}^{x} t^2 \cos 2Nt \cot t \, dt = x^2 \left( \sum_{n=1}^{N} \frac{\cos 2nx}{n} - \frac{\cos 2Nx}{2N} + \log \sin x \right) - 2\left( \sum_{n=1}^{N} \frac{\cos 2nx}{4n^3} - \frac{\cos 2Nx}{8N^3} \right)$$

$$- 2\left( \sum_{n=1}^{N} \frac{x \sin 2nx}{2n^2} - \frac{x \sin 2Nx}{4N^2} \right) - 2\int_{0}^{x} t \log \sin t \, dt + \frac{1}{2}\sum_{n=1}^{N} \frac{1}{n^3} - \frac{1}{4N^3}$$

Therefore, as $N \to \infty$ we obtain

(7.17)
$$2\int_0^x t \log \sin t \, dt = \sum_{n=1}^{\infty} \left( \frac{x^2 \cos 2nx}{n} - \frac{\cos 2nx}{2n^3} - \frac{x \sin 2nx}{n^2} \right) + x^2 \log \sin x + \frac{1}{2}\varsigma(3)$$

which is valid in the interval $(-\pi, \pi)$. Prima facie, we cannot differentiate the right hand side of (7.17) because we would end up with a divergent series $\sum_{n=1}^{\infty} x^2 \sin 2nx$ for $|x| \geq 1$: however it can be seen from (6.20) that

$$\int_0^x t^2 \cot t \, dt = \frac{1}{2}\sum_{n=1}^{\infty} \frac{\cos 2nx}{n^3} - x^2 \sum_{n=1}^{\infty} \frac{\cos 2nx}{n} + x \sum_{n=1}^{\infty} \frac{\sin 2nx}{n} - \frac{1}{2}\varsigma(3)$$

### 8. SOME MISCELLANEOUS RESULTS

**Example (i).** It is shown in Knopp's book [90, p.376] that

(8.1) $$x = \sum_{n=0}^{\infty} c_n \frac{\sin^{2n+1} x}{2n+1} \quad , 0 \leq x \leq \pi/2$$

where

(8.2) $$c_n = (-1)^n \binom{-1/2}{n} = \frac{1}{2^{2n}} \binom{2n}{n} = \frac{1.3...(2n-1)}{2.4...(2n)} = \frac{(2n-1)!!}{(2n)!!}$$

and the central binomial sum $\binom{2n}{n}$ is equal to $\frac{(2n)!}{(n!)^2}$. A further proof of (8.1) is given below. By the binomial theorem we have

$$\frac{1}{\sqrt{1-x^2}} = \sum_{n=0}^{\infty} c_n x^{2n} \quad ,(|x| < 1)$$

and upon integration we obtain

$$\sin^{-1} x = \sum_{n=0}^{\infty} \frac{c_n}{2n+1} x^{2n+1}$$

Then, letting $t = \sin^{-1} x$, we get (8.1).

Using (8.1) we have

$$\int_0^{\pi/2} x\,dx = \frac{\pi^2}{8} = \sum_{n=0}^{\infty} \frac{c_n}{2n+1} \int_0^{\pi/2} \sin^{2n+1} x\,dx$$

and with Wiener's formula (8.11h) this becomes

$$\frac{\pi^2}{8} = \sum_{n=0}^{\infty} \frac{1}{2n+1} \frac{(2n-1)!!}{(2n)!!} \frac{(2n)!!}{(2n+1)!!} = \sum_{n=0}^{\infty} \frac{1}{(2n+1)^2}$$

Multiplying (8.1) by $\cot x$ and integrating over the interval $[0, \pi/2]$ we have

(8.3) $$\int_0^{\pi/2} x \cot x\,dx = \sum_{n=0}^{\infty} c_n \int_0^{\pi/2} \frac{\sin^{2n} x \cos x\,dx}{2n+1}$$

(8.3a) $$= \sum_{n=0}^{\infty} c_n \frac{\sin^{2n+1} x}{(2n+1)^2} \bigg|_0^{\pi/2} = \sum_{n=0}^{\infty} \frac{1}{2^{2n}(2n+1)^2} \binom{2n}{n}$$

We have already evaluated the integral in (3.1), and we therefore have

(8.4) $$\frac{\pi}{2} \log 2 = \sum_{n=0}^{\infty} \frac{1}{2^{2n}(2n+1)^2} \binom{2n}{n}$$

Ramanujan's second paper [76, p.16] was published in the Journal of the Indian Mathematical Society in 1912 and had the enigmatic title "On Question 330 of Professor Sanjana": inter alia, it contained a result similar to (8.4) but curiously the factor of $2^{2n}$ does not appear in the denominator.

More generally we have

$$\int_0^t x \cot x\,dx = \sum_{n=0}^{\infty} c_n \int_0^t \frac{\sin^{2n} x \cos x\,dx}{2n+1}$$

(8.4a) $$\int_0^t x \cot x\,dx = \sum_{n=0}^{\infty} c_n \frac{\sin^{2n+1} t}{(2n+1)^2}$$

Wiener [138aii] gives an elementary and elegant method for expansions of powers of $\sin t$ and $\cos t$. Letting $z = \cos t + i \sin t$ we have

$$\sin^{2n+1} t = \left[\frac{1}{2i}\left(z - \frac{1}{z}\right)\right]^{2n+1}$$

we have

$$\sin^{2n+1} t = i \frac{(-1)^n}{2^{2n+1}} \sum_{k=0}^{2n+1} (-1)^k \binom{2n+1}{k} z^{2n+1-2k}$$

Since $z^{2n+1-2k} = \cos(2n+1-2k)t + i\sin(2n+1-2k)t$ we have

$$\sin^{2n+1} t = i\frac{(-1)^n}{2^{2n+1}}\sum_{k=0}^{2n+1}(-1)^k\binom{2n+1}{k}\cos(2n+1-2k)t$$

$$-\frac{(-1)^n}{2^{2n+1}}\sum_{k=0}^{2n+1}(-1)^k\binom{2n+1}{k}\sin(2n+1-2k)t$$

Equating real and imaginary parts, we immediately see that

$$\frac{(-1)^n}{2^{2n+1}}\sum_{k=0}^{2n+1}(-1)^k\binom{2n+1}{k}\cos(2n+1-2k)t = 0$$

$$\frac{(-1)^{n+1}}{2^{2n+1}}\sum_{k=0}^{2n+1}(-1)^k\binom{2n+1}{k}\sin(2n+1-2k)t = \sin^{2n+1} t$$

We have

$$\sum_{k=0}^{2n+1}(-1)^k\binom{2n+1}{k}\sin(2n+1-2k)t =$$

$$\sum_{k=0}^{n}(-1)^k\binom{2n+1}{k}\sin(2n+1-2k)t + \sum_{k=n+1}^{2n+1}(-1)^k\binom{2n+1}{k}\sin(2n+1-2k)t$$

Since $\binom{i}{j} = \binom{i}{i-j}$ we get

$$\sum_{k=n+1}^{2n+1}(-1)^k\binom{2n+1}{k}\sin(2n+1-2k)t = \sum_{k=n+1}^{2n+1}(-1)^k\binom{2n+1}{2n+1-k}\sin(2n+1-2k)t$$

Designating $j = 2n+1-k$ this becomes

$$= \sum_{j=n}^{0}(-1)^{2n+1-j}\binom{2n+1}{j}(-1)\sin(2n+1-2k)t$$

$$= \sum_{j=0}^{n}(-1)^j\binom{2n+1}{j}\sin(2n+1-2k)t$$

Therefore we obtain

$$\sum_{k=0}^{2n+1}(-1)^k\binom{2n+1}{k}\sin(2n+1-2k)t = 2\sum_{k=0}^{n}(-1)^k\binom{2n+1}{k}\sin(2n+1-2k)t$$

and thus, in accordance with G&R [74, p.30], we have

$$\sin^{2n+1}t = \frac{(-1)^n}{2^{2n}}\sum_{k=0}^{n}(-1)^k\binom{2n+1}{k}\sin(2n+1-2k)t$$

Using the familiar trigonometric identity

$$\sin(2n+1-2k)t = \sin 2nt\cos(1-2k)t + \cos 2nt\sin(1-2k)t$$

we have
$$\sin^{2n+1}t =$$

$$\frac{(-1)^n}{2^{2n}}\sin 2nt\sum_{k=0}^{n}(-1)^k\binom{2n+1}{k}\cos(1-2k)t + \frac{(-1)^n}{2^{2n}}\cos 2nt\sum_{k=0}^{n}(-1)^k\binom{2n+1}{k}\sin(1-2k)t$$

This then gives us

$$\sum_{n=0}^{\infty}c_n\frac{\sin^{2n+1}t}{(2n+1)^2} =$$

$$\sum_{n=0}^{\infty}c_n\frac{(-1)^n}{(2n+1)^2 2^{2n}}\sin 2nt\sum_{k=0}^{n}(-1)^k\binom{2n+1}{k}\cos(2k-1)t -$$

$$\sum_{n=0}^{\infty}c_n\frac{(-1)^n}{(2n+1)^2 2^{2n}}\cos 2nt\sum_{k=0}^{n}(-1)^k\binom{2n+1}{k}\sin(2k-1)t$$

Using (6.5a) we have for $t \in [0,\pi)$

$$\int_0^t x\cot x\,dx = 2\sum_{n=1}^{\infty}\int_0^t x\sin 2nx\,dx$$

$$= \frac{1}{2}\sum_{n=1}^{\infty}\frac{\sin 2nt}{n^2} - t\sum_{n=1}^{\infty}\frac{\cos 2nt}{n}$$

Therefore we have

$$\sum_{n=0}^{\infty}c_n\frac{(-1)^n}{(2n+1)^2 2^{2n}}\sin 2nt\sum_{k=0}^{n}(-1)^k\binom{2n+1}{k}\cos(2k-1)t -$$

$$\sum_{n=0}^{\infty}c_n\frac{(-1)^n}{(2n+1)^2 2^{2n}}\cos 2nt\sum_{k=0}^{n}(-1)^k\binom{2n+1}{k}\sin(2k-1)t$$

$$= \frac{1}{2}\sum_{n=1}^{\infty}\frac{\sin 2nt}{n^2} - t\sum_{n=1}^{\infty}\frac{\cos 2nt}{n}$$

I initially thought that we could now employ Cantor's theorem (1870) on the uniqueness of Fourier series to equate the coefficients of $\sin 2nt$ in both series representations (and also for $\cos 2nt$) but this does not appear to be appropriate.

With $t = \pi/2$ we get for the left-hand side

$$\sum_{n=0}^{\infty} c_n \frac{1}{(2n+1)^2 2^{2n}} \sum_{k=0}^{n}\binom{2n+1}{k}$$

and we see that

$$\sum_{k=0}^{2n+1}\binom{2n+1}{k} = \sum_{k=0}^{n}\binom{2n+1}{k} + \sum_{k=n+1}^{2n+1}\binom{2n+1}{k} = \sum_{k=0}^{n}\binom{2n+1}{k} + \sum_{k=n+1}^{2n+1}\binom{2n+1}{2n+1-k}$$

$$= 2\sum_{k=0}^{n}\binom{2n+1}{k}$$

Using the fact that $(1+1)^{2n+1} = 2^{2n+1} = \sum_{k=0}^{2n+1}\binom{2n+1}{k}$, we then see that

$$\sum_{k=0}^{n}\binom{2n+1}{k} = 2^{2n}$$

and the left-hand side then becomes $\sum_{n=0}^{\infty} c_n \frac{1}{(2n+1)^2}$. The right-hand side is equal to $\frac{\pi}{2}\log 2$ and hence we have regained equation (8.4).

From (8.4a) we have

$$\int_{0}^{t} x \cot x\, dx = \sum_{n=0}^{\infty} c_n \frac{\sin^{2n+1} t}{(2n+1)^2}$$

Using mathematical induction it is relatively easy to prove that

$$\sin^{2n+1} t = \sum_{k=0}^{n} A_k \sin(2k+1)t$$

and hence we have

$$\int_0^t x \cot x \, dx = \sum_{n=0}^{\infty} \frac{c_n}{(2n+1)^2} \sum_{k=0}^{n} A_k \sin(2k+1)t$$

With $t = \pi/2$ we see that

$$1 = \sum_{k=0}^{n} A_k \sin \frac{(2k+1)\pi}{2} = \sum_{k=0}^{n} (-1)^{k+1} A_k$$

We have

$$\int_0^u dt \int_0^t x \cot x \, dx = \sum_{n=0}^{\infty} \frac{c_n}{(2n+1)^2} \int_0^u \sin^{2n+1} t \, dt$$

and using (8.11g) this becomes

$$= \sum_{n=0}^{\infty} \frac{c_n}{(2n+1)^2} \sum_{k=0}^{n} \binom{n}{k} (-1)^{k+1} \frac{\left[\cos^{2k+1} u - 1\right]}{2k+1}$$

As mentioned above, using (6.5a) we have for $t \in [0, \pi)$

$$\int_0^t x \cot x \, dx = \frac{1}{2} \sum_{n=1}^{\infty} \frac{\sin 2nt}{n^2} - t \sum_{n=1}^{\infty} \frac{\cos 2nt}{n}$$

and therefore we obtain

$$\int_0^u dt \int_0^t x \cot x \, dx = \frac{1}{4} \sum_{n=1}^{\infty} \frac{1}{n^3} - \frac{1}{4} \sum_{n=1}^{\infty} \frac{\cos 2nu}{n^3} - \frac{1}{2} t \sum_{n=1}^{\infty} \frac{\sin 2nu}{n^2} - \frac{1}{4} \sum_{n=1}^{\infty} \frac{\cos 2nu}{n^2} + \frac{1}{4} \sum_{n=1}^{\infty} \frac{1}{n^2}$$

We then have

(8.4b) $$\frac{1}{4}\varsigma(2) + \frac{1}{4}\varsigma(3) - \frac{1}{4} \sum_{n=1}^{\infty} \frac{\cos 2nu}{n^3} - \frac{1}{2} t \sum_{n=1}^{\infty} \frac{\sin 2nu}{n^2} - \frac{1}{4} \sum_{n=1}^{\infty} \frac{\cos 2nu}{n^2} =$$

$$\sum_{n=0}^{\infty} \frac{c_n}{(2n+1)^2} \sum_{k=0}^{n} \binom{n}{k} (-1)^{k+1} \frac{\left[\cos^{2k+1} u - 1\right]}{2k+1}$$

**Example (ii).** Let us now consider the integral

(8.5) $$J = \int_0^{\pi/2} x^2 \cot x \, dx = x^2 \log \sin x \Big|_0^{\pi/2} - 2 \int_0^{\pi/2} x \log \sin x \, dx$$

$$= -2 \int_0^{\pi/2} x \log \sin x \, dx$$

Multiplying (8.1) by $x \cot x$ and integrating, $J$ may be expressed as

(8.6) $$J = \sum_{n=0}^{\infty} c_n \int_0^{\pi/2} \frac{x \sin^{2n} x \cos x \, dx}{2n+1}$$

Integration by parts gives us

(8.7) $$\int_0^{\pi/2} x \sin^{2n} x \cos x \, dx = \left. \frac{x \sin^{2n+1} x}{2n+1} \right|_0^{\pi/2} - \int_0^{\pi/2} \frac{\sin^{2n+1} x \, dx}{2n+1}$$

The latter integral may be evaluated using the Wallis integral (see [25, p.113] and [137])

(8.8) $$\int_0^{\pi/2} \sin^{2n+1} x \, dx = \frac{1}{c_n (2n+1)}$$

and we may therefore write

(8.9) $$J = \sum_{n=0}^{\infty} c_n \left\{ \frac{\pi}{2(2n+1)^2} - \frac{1}{c_n (2n+1)^3} \right\}$$

$$= \frac{\pi}{2} \sum_{n=0}^{\infty} c_n \frac{1}{(2n+1)^2} - \sum_{n=0}^{\infty} \frac{1}{(2n+1)^3}$$

(8.10) $$= \frac{\pi}{4} \log 2 - \sum_{n=0}^{\infty} \frac{1}{(2n+1)^3}$$

where, in the last part, we used equation (8.4). Combining (8.5) and (8.10) we deduce

(8.11) $$2 \int_0^{\pi/2} x \log \sin x \, dx = \sum_{n=0}^{\infty} \frac{1}{(2n+1)^3} - \frac{\pi^2}{4} \log 2$$

and this is simply a mildly disguised form of Euler's 1772 integral (1.11) because $\frac{7}{8} \varsigma(3) = \sum_{n=0}^{\infty} \frac{1}{(2n+1)^3}$, as can be seen from (1.1) of Volume I.

It should be noted that if we integrate (8.1) over the interval $[0, \pi/2]$ and use (8.8) then we obtain a simple proof that

$$\frac{\pi^2}{8} = \sum_{n=0}^{\infty} \frac{1}{(2n+1)^2}$$

Letting $x = \sin^{-1} t$ in (8.1) is equivalent to one of Euler's more rigorous derivations of the Basle identity (see Kimble's paper "Euler's Other Proof" [84]).

Wiener et al. have given a useful approach to evaluating trigonometric integrals in [138aa]: their paper includes the following

$$\int_0^t \cos^{2n+1} x \, dx = \int_0^t \cos^{2n} x \cos x \, dx = \int_0^t (1-\sin^2 x)^n x \cos x \, dx$$

$$= \int_0^t \sum_{k=0}^n \binom{n}{k} (-1)^k \sin^{2k} x \cos x \, dx$$

We therefore have

(8.11a) $$\int_0^t \cos^{2n+1} x \, dx = \sum_{k=0}^n \binom{n}{k} (-1)^k \frac{\sin^{2k+1} t}{2k+1}$$

With $t = \pi/2$ we get

(8.11b) $$\int_0^{\pi/2} \cos^{2n+1} x \, dx = \sum_{k=0}^n \binom{n}{k} \frac{(-1)^k}{2k+1}$$

as compared with the frequently quoted form of the Wallis integral formula

(8.11c) $$\int_0^{\pi/2} \cos^{2n+1} x \, dx = \frac{(2n)!!}{(2n+1)!!}$$

We therefore have the following identity which is reported in [90, p.270]

(8.11d) $$\sum_{k=0}^n \binom{n}{k} (-1)^k \frac{1}{2k+1} = \frac{(2n)!!}{(2n+1)!!} = \frac{[2^n n!]^2}{(2n+1)!}$$

Completing the summation we get

$$\sum_{n=0}^\infty \frac{1}{2^n} \int_0^{\pi/2} \cos^{2n+1} x \, dx = \sum_{n=0}^\infty \frac{1}{2^n} \sum_{k=0}^n \binom{n}{k} \frac{(-1)^k}{2k+1}$$

Interchanging the order of summation and integration we obtain

$$\int_0^{\pi/2} \sum_{n=0}^\infty \frac{\cos^{2n+1} x}{2^n} \, dx = \int_0^{\pi/2} \frac{\cos x}{\left[1 - \frac{\cos^2 x}{2}\right]} \, dx = 2\int_0^{\pi/2} \frac{\cos x}{[1+\sin^2 x]} \, dx$$

and the obvious substitution $u = \sin x$ provides us with

$$= 2\int_0^1 \frac{1}{\left[1+u^2\right]} du = \tan^{-1} 1 = \frac{\pi}{2}$$

Therefore we deduce that

(8.11e) $$\sum_{n=0}^{\infty} \frac{1}{2^n} \sum_{k=0}^{n} \binom{n}{k} \frac{(-1)^k}{2k+1} = \frac{\pi}{2}$$

and

(8.11f) $$\sum_{n=0}^{\infty} \frac{1}{2^n} \frac{(2n)!!}{(2n+1)!!} = \frac{\pi}{2}$$

and this may be a consequence of Euler's transformation of series technique.

Similarly we also have

(8.11g) $$\int \sin^{2n+1} x \, dx = \sum_{k=0}^{n} \binom{n}{k} (-1)^{k+1} \frac{\cos^{2k+1} t}{2k+1}$$

and in particular we get

(8.11h) $$\int_0^{\pi/2} \sin^{2n+1} x \, dx = \sum_{k=0}^{n} \binom{n}{k} \frac{(-1)^k}{2k+1}$$

Completing the summation we get for $\sin x / 2 < 1$ we have

$$\sum_{n=1}^{\infty} \frac{1}{2^n} \int_0^{\pi/2} \sin^{2n+1} x \, dx = \int_0^{\pi/2} \frac{\sin^3 x}{1+\cos^2 x} dx = \sum_{n=1}^{\infty} \frac{1}{2^n} \sum_{k=0}^{n} \binom{n}{k} \frac{(-1)^k}{2k+1}$$

With the substitution $u = \cos x$ we get

$$\int_0^{\pi/2} \frac{\sin^3 x}{1+\cos^2 x} dx = \int_1^0 \left(1 - \frac{2}{1+u^2}\right) du = \frac{\pi}{2}$$

and hence we obtain (8.11e) again

(8.11i) $$\frac{\pi}{2} = \sum_{n=1}^{\infty} \frac{1}{2^n} \sum_{k=0}^{n} \binom{n}{k} \frac{(-1)^k}{2k+1}$$

The authors, Wiener et al [138ai], also show with the substitution $u = \tan x$ that

$$\int \tan^{2n+1} x \, dx = \int \frac{u^{2n+1}}{1+u^2} du = \int \left( u^{2n-1} - u^{2n-3} + u^{2n-5} - \ldots + (-1)^n \frac{u}{1+u^2} \right)$$

and hence

(8.11j) $$\int_0^t \tan^{2n+1} x \, dx = (-1)^n \log \sec t + (-1)^n \sum_{k=1}^n (-1)^k \frac{\tan^{2k} t}{2k}$$

With $t = \pi/4$ we obtain

$$(-1)^n \int_0^{\pi/4} \tan^{2n+1} x \, dx = \frac{1}{2} \log 2 + \frac{1}{2} \sum_{k=1}^n \frac{(-1)^k}{k}$$

and as $n \to \infty$ we see that $\log 2 = -\sum_{k=1}^{\infty} \frac{(-1)^k}{k}$.

Completing the summation we get for $\tan x / 2 < 1$

$$\sum_{n=1}^{\infty} \frac{1}{2^n} \int_0^t \tan^{2n+1} x \, dx = \log \sec t \sum_{n=1}^{\infty} \frac{(-1)^n}{2^n} + \sum_{n=1}^{\infty} \frac{(-1)^n}{2^n} \sum_{k=1}^n (-1)^k \frac{\tan^{2k} t}{2k}$$

The geometric series converges to a surprisingly simple form and we get

$$\sum_{n=1}^{\infty} \frac{1}{2^n} \int_0^t \tan^{2n+1} x \, dx = -2 \int_0^t \cot x \, dx = -2 \log \cos t$$

We then find an interesting series expansion for $\log \cos t$ in the range $[0, \pi/2)$

(8.11k) $$\log \cos t = \frac{3}{7} \sum_{n=1}^{\infty} \frac{(-1)^n}{2^n} \sum_{k=1}^n (-1)^k \frac{\tan^{2k} t}{2k}$$

With $t = \pi/4$ this becomes

(8.11l) $$\log 2 = \frac{6}{7} \sum_{n=1}^{\infty} \frac{(-1)^n}{2^{n+1}} \sum_{k=1}^n \frac{(-1)^{k+1}}{k}$$

and using (4.4.7) this equates to

(8.11m) $$\log 2 = \frac{6}{7} \sum_{n=1}^{\infty} (-1)^n \frac{H_n}{2^{n+1}}$$

**Example (iii).** We have already seen from (6.28) that

$$(8.12) \qquad 2G = \int_0^{\pi/2} \frac{x}{\sin x} dx$$

We shall now evaluate the integral in a different way: using (8.1) we have

$$(8.13) \qquad \int_0^{\pi/2} \frac{x}{\sin x} dx = \sum_{n=0}^{\infty} \frac{c_n}{2n+1} \int_0^{\pi/2} \sin^{2n} x \, dx$$

From [137] we have the Wallis formula

$$(8.14) \qquad \int_0^{\pi/2} \sin^{2n} x \, dx = \frac{\pi}{2} \frac{(2n-1)!!}{(2n)!!} = \frac{\pi}{2} c_n$$

and we therefore have

$$(8.15) \qquad 2G = \frac{\pi}{2} \sum_{n=0}^{\infty} \frac{c_n^2}{2n+1}$$

or

$$(8.15a) \qquad G = \frac{\pi}{4} \sum_{n=0}^{\infty} \frac{1}{(2n+1)2^{4n}} \binom{2n}{n}^2 = \frac{\pi}{4} \sum_{n=0}^{\infty} \frac{[(2n)!]^2}{(2n+1)[(n)!]^4 \, 2^{4n}}$$

Is this a new result? Unfortunately, the answer is no because the identity has actually been proved several times before: see, for example, a paper by Adamchik in 2001 [6]. The formula does have some history: on 16 January 1913 a letter, accompanied by twelve pages of dense mathematical formulae, was sent by a 26 year old Indian clerk to a renowned English mathematician [76]. As is well known, the writer was Srinivasa Ramanujan (1887-1920) and the Cambridge mathematician was G.H. Hardy (1877-1947) .On one of these pages Ramanujan stated a formula similar to (8.15a) and, as usual, this was stated without giving any proof. During the course of long discussions between Hardy and his collaborator, J.E. Littlewood, one of these pages was mislaid and never recovered: fortunately, it was not the page referred to above. Ramanujan's proof of (8.15a) can now be seen in [21, p.39] and [110].

It is perhaps not so well known that, in addition to writing to Hardy, Ramanujan also wrote to two other prominent mathematicians (Hardy, who was an ardent cricket fan, never did reveal their identities during his lifetime, or at least, never in writing). In the biography, The man who knew infinity: A Life of the Genius Ramanujan, Kanigel [83, p.170] reports that the other mathematicians were H.F.Baker and E.W. Hobson, both of whom were Senior Wranglers at Cambridge (Hobson held the Sadleirian chair of pure mathematics and, upon his retirement, was succeeded in this position by Hardy). They failed however to recognise Ramanujan's genius: it would be interesting to discover, even at this late juncture, if any of the papers which they received from him have survived the passage of time. Interestingly, both Baker and Hobson were signatories to the document recommending Ramanujan's election to the Royal Society in 1918. Further information regarding Ramanujan, and facsimile copies of his Notebooks, are contained in a website maintained by K. Srinivasa Rao [111].

Let us now continue on in the same vein: since $\int \frac{dx}{\sin^2 x} = -\cot x$, using integration by parts we get

$$\int \frac{x^2 \, dx}{\sin^2 x} = -x^2 \cot x + 2\int x \cot x \, dx$$

and therefore $\int_0^{\pi/2} \frac{x^2 \, dx}{\sin^2 x} = 2\int_0^{\pi/2} x \cot x \, dx$

Using (8.1) we obtain

$$\int_0^{\pi/2} \frac{x^2}{\sin^2 x} \, dx = \sum_{n=0}^{\infty} \frac{c_n}{2n+1} \int_0^{\pi/2} x^2 \sin^{2n-1} x \, dx$$

and the latter integral may be evaluated by using the decomposition formula in [74, p.30]

(8.15b) $\quad \sin^{2n-1} x = \frac{1}{2^{2n-2}} \left[ \sum_{k=0}^{n-1} (-1)^{n+k-1} \binom{2n-1}{k} \sin(2n-2k-1)x \right]$

In [74, p.53] it is reported that, using the Maclaurin series, we have for $|x| < \pi/2$

(8.15c) $\quad \log \cos x = \log \sqrt{1-\sin^2 x} = \frac{1}{2}\log(1-\sin^2 x) = -\frac{1}{2}\sum_{n=1}^{\infty} \frac{\sin^{2n} x}{n}$

and we then obtain

$$\int_0^{\pi/4} \log \cos x \, dx = -\frac{1}{2}\sum_{n=1}^{\infty} \frac{1}{n} \int_0^{\pi/4} \sin^{2n} x \, dx$$

From (6.107b) we have

$$\int_0^{\pi/4} \log \cos x \, dx = \frac{G}{2} - \frac{\pi}{4}\log 2$$

We have [74, p.30]

$$(8.15d) \qquad \sin^{2n} x = \frac{1}{2^{2n}} \left[ \sum_{k=0}^{n-1} (-1)^{n-k} 2\binom{2n}{k} \cos 2(n-k)x + \binom{2n}{n} \right]$$

Therefore we get

$$\int_0^{\pi/4} \sin^{2n} x \, dx = \frac{1}{2^{2n}} \left[ \sum_{k=0}^{n-1} (-1)^{n-k} \binom{2n}{k} \frac{\sin[(n-k)\pi/2]}{n-k} + \frac{\pi}{4}\binom{2n}{n} \right]$$

and hence we get

$$\frac{\pi}{2}\log 2 - G = \sum_{n=1}^{\infty} \frac{(-1)^n}{n2^{2n}} \left[ \sum_{k=0}^{n-1} (-1)^k \binom{2n}{k} \frac{\sin[(n-k)\pi/2]}{n-k} \right] + \frac{\pi}{4} \sum_{n=1}^{\infty} \frac{1}{n2^{2n}} \binom{2n}{n}$$

It is interesting in this connection to note that Ross [114a] has proved that

$$(8.15e) \qquad \sum_{n=1}^{\infty} \frac{1.3.5....(2n-1)}{n2^{2n} n!} = 2\log 2$$

but Lehmer's observation [97] prima facie appeared to be more germane: he easily showed that

$$(8.15f) \qquad \sum_{n=1}^{\infty} \frac{1}{n}\binom{2n}{n} x^n = 2\log\left[\frac{1-\sqrt{1-4x}}{2x}\right]$$

and thus

$$\sum_{n=1}^{\infty} \frac{1}{n}\binom{2n}{n} y^{2n} = 2\log\left[\frac{1-\sqrt{1-4y^2}}{2y^2}\right]$$

With $y = 1/2$ we get

$$(8.15g) \qquad \sum_{n=1}^{\infty} \frac{1}{n2^{2n}}\binom{2n}{n} = 2\log 2$$

Since $1.3.5....(2n-1) = \frac{(2n)!}{2.4.6...2n} = \frac{(2n)!}{2^n n!}$ it is clear that the Ross and Lehmer formulae are in fact equivalent. Raina and Ladda [109a] have pointed out how (8.15e) may be derived from (4.3.32).

Hence we deduce that

(8.15h) $$G = \sum_{n=1}^{\infty} \frac{(-1)^{n+1}}{n2^{2n}} \left[ \sum_{k=0}^{n-1} (-1)^k \binom{2n}{k} \frac{\sin[(n-k)\pi/2]}{n-k} \right]$$

In an attempt to simplify the summand, I considered the following analysis

$$\sum_{k=0}^{n-1} (-1)^k \binom{2n}{k} \frac{\sin[(n-k)\pi/2]}{n-k} = 2\sum_{k=0}^{n-1} \int_0^{\pi/4} (-1)^k \binom{2n}{k} \cos 2[(n-k)x] dx$$

Making reference to (8.15d) we have

$$2^{2n} \sin^{2n} x - \binom{2n}{n} = 2(-1)^n \left[ \sum_{k=0}^{n-1} (-1)^k \binom{2n}{k} \cos 2(n-k)x \right]$$

but proceeding further in this direction will unfortunately prove to be circular in nature.

Another incomplete alternative approach is noted below

$$\sum_{k=0}^{n-1} (-1)^k \binom{2n}{k} \frac{\sin[(n-k)\pi/2]}{n-k} = \sum_{k=0}^{n-1} \int_0^{\pi/2} (-1)^k \binom{2n}{k} \cos[(n-k)x] dx$$

$$= \text{Re} \left\{ \sum_{k=0}^{n-1} \int_0^{\pi/2} (-1)^k \binom{2n}{k} e^{[i(n-k)x]} dx \right\}$$

$$= \text{Re} \left\{ \sum_{k=0}^{n-1} \int_0^{\pi/2} e^{inx} (-1)^k \binom{2n}{k} [e^{-ix}]^k dx \right\}$$

$$= \text{Re} \left\{ \sum_{k=0}^{n-1} \int_0^{\pi/2} e^{inx} \binom{2n}{k} [-e^{-ix}]^k dx \right\}$$

$$= \text{Re} \left\{ \int_0^{\pi/2} e^{inx} \sum_{k=0}^{n-1} \binom{2n}{k} [-e^{-ix}]^k dx \right\}$$

If we multiply (8.15c) by $\sin x$ and integrate we obtain

$$\int_0^t \sin x \log \cos x \, dx = -\frac{1}{2} \sum_{n=1}^{\infty} \frac{1}{n} \int_0^t \sin^{2n+1} x \, dx$$

and, using the Wiener formula (8.11g), this becomes

$$\int_0^t \sin x \log \cos x\, dx = -\frac{1}{2}\sum_{n=1}^{\infty}\frac{1}{n}\left[\sum_{k=0}^{n}\binom{n}{k}(-1)^{k+1}\frac{\cos^{2k+1} t}{2k+1} - \sum_{k=0}^{n}\binom{n}{k}(-1)^{k+1}\frac{1}{2k+1}\right]$$

Integration by parts also gives us

$$\int_0^t \sin x \log \cos x\, dx = -\cos x \log \cos x + \cos x \Big|_0^t$$

$$= -\cos t \log \cos t + \cos t - 1$$

and we therefore have

(8.15i) $\quad \dfrac{1}{2}\sum_{n=1}^{\infty}\dfrac{1}{n}\sum_{k=0}^{n}\binom{n}{k}(-1)^k \dfrac{\left[\cos^{2k+1} t - 1\right]}{2k+1} = \cos t(1 - \log \cos t) - 1$

With $t = \pi/2$ we obtain

$$\frac{1}{2}\sum_{n=1}^{\infty}\frac{1}{n}\sum_{k=0}^{n}\binom{n}{k}\frac{(-1)^k}{2k+1} = 1$$

From an elementary exercise in [108, p.229] we have

$$\log \tan\left(\frac{\pi}{4} + \frac{x}{2}\right) = \sum_{n=0}^{\infty} \frac{\sin^{2n+1} x}{2n+1}$$

and therefore upon integration we obtain

(8.15j) $\quad \displaystyle\int_0^t \log \tan\left(\frac{\pi}{4} + \frac{x}{2}\right) dx = \sum_{n=0}^{\infty}\frac{1}{2n+1}\int_0^t \sin^{2n+1} x\, dx$

$$= \sum_{n=0}^{\infty}\frac{1}{2n+1}\sum_{k=0}^{n}\binom{n}{k}(-1)^{k+1}\frac{\left[\cos^{2k+1} t - 1\right]}{2k+1}$$

**Example (iv).** Using integration by parts and (6.107d) we have

(8.16) $\quad \displaystyle\int_0^{\pi/4} x \cot x\, dx = x \log \sin x \Big|_0^{\pi/4} - \int_0^{\pi/4} \log \sin x\, dx$

$$= \frac{\pi}{8}\log 2 + \frac{G}{2}$$

$$= \sum_{n=0}^{\infty} c_n \int_0^{\pi/4} \frac{\sin^{2n} x \cos x \, dx}{2n+1} = \frac{1}{\sqrt{2}} \sum_{n=0}^{\infty} \frac{1}{2^{3n}(2n+1)^2} \binom{2n}{n}$$

This then gives

(8.17) $$\frac{\pi}{8} \log 2 + \frac{G}{2} = \frac{1}{\sqrt{2}} \sum_{n=0}^{\infty} \frac{1}{2^{3n}(2n+1)^2} \binom{2n}{n}$$

D.M. Bradley [33] has provided a very detailed list (and proofs) of more than 70 integrals and series involving Catalan's constant on his website and the above formula appears there as item (68). This identity was also due to Ramanujan.

**Example (v).** The inspiration for the following example came from D.H. Lehmer's paper, "Interesting series involving the central binomial coefficient" [97]. By the binomial theorem we have

(8.18) $$\sum_{n=0}^{\infty} \frac{1}{2^{2n}} \binom{2n}{n} x^n = \frac{1}{\sqrt{1-x}} \quad ,(|x|<1)$$

Replacing $x$ by $x^2$ in (8.18) and integrating over the range $[0, t]$, we obtain (after dividing the result by $t$)

(8.19) $$\sum_{n=0}^{\infty} \frac{1}{2^{2n}(2n+1)} \binom{2n}{n} t^{2n} = \frac{\sin^{-1} t}{t} \quad ,(|t| \leq 1)$$

We now multiply (8.19) by $\sin^{-1} t$, and integrate over the range $[0,1]$, to give

(8.20) $$\sum_{n=0}^{\infty} \int_0^1 \frac{1}{2^{2n}(2n+1)} \binom{2n}{n} t^{2n} \sin^{-1} t \, dt = \int_0^1 \frac{(\sin^{-1} t)^2}{t} dt$$

(and using the substitution $t = \sin x$, this is equivalent to the formula (8.3)). The series in (8.20) converges absolutely and uniformly on the unit closed disc [96] and hence we may interchange the order of integration and summation. With the substitution $t = \sin x$ we have

(8.21) $$I_n = \int_0^1 t^{2n} \sin^{-1} t \, dt = \int_0^{\pi/2} x \sin^{2n} x \cos x \, dx$$

$$= \frac{x \sin^{2n+1} x}{2n+1} \Big|_0^{\pi/2} - \frac{1}{2n+1} \int_0^{\pi/2} \sin^{2n+1} x \, dx$$

We therefore have

$$(8.22) \qquad \int_0^1 t^{2n} \sin^{-1} t \, dt = \frac{1}{2n+1}\left[\frac{\pi}{2} - \frac{2^n n!}{(2n+1)!!}\right]$$

where we have again used the Wallis integral (8.8) ($I_n$ is also contained in G&R [74, p.597]).

The left-hand side of (8.20) may now be written as

$$= \sum_{n=0}^{\infty} \frac{1}{2^{2n}(2n+1)^2} \binom{2n}{n} \left[\frac{\pi}{2} - \frac{2^n n!}{(2n+1)!!}\right]$$

$$= \frac{\pi}{2} \sum_{n=0}^{\infty} \frac{1}{2^{2n}(2n+1)^2} \binom{2n}{n} - \sum_{n=0}^{\infty} \frac{(2n)!}{2^n n!(2n+1)^2 (2n+1)!!}$$

$$(8.23) \qquad = \frac{\pi^2}{4} \log 2 - \sum_{n=0}^{\infty} \frac{1}{(2n+1)^3} = \frac{\pi^2}{4} \log 2 - \frac{7}{8}\varsigma(3)$$

where we have used (8.4) and the fact that $(2n+1)!! = \frac{(2n+1)!}{2^n n!}$.

Now, using the substitution $t = \sin\theta$ in the integral in (8.20) we get

$$(8.24) \qquad \int_0^1 \frac{(\sin^{-1} t)^2}{t} dt = \int_0^{\pi/2} \theta^2 \cot\theta \, d\theta = -2\int_0^{\pi/2} \theta \log \sin\theta \, d\theta$$

Comparing (8.23) and (8.24) we obtain yet another proof of the 1772 Euler integral shown in (1.11).

$$\int_0^{\pi/2} \theta \log \sin\theta \, d\theta = \frac{7}{16}\varsigma(3) - \frac{\pi^2}{8} \log 2$$

Equation (8.21) was the starting point for a short paper written by Boo Rim Choe [44] in 1987 where he gave an elementary proof that $\sum_{n=1}^{\infty} \frac{1}{n^2} = \frac{\pi^2}{6}$. A similar method was employed by J.A. Ewell [65] a few years later when he found a "new" fast converging series representation for $\varsigma(3)$.

$$(8.25) \qquad \varsigma(3) = \frac{\pi^2}{7}\left\{1 - 4\sum_{n=1}^{\infty} \frac{\varsigma(2n)}{(2n+1)(2n+2)2^{2n}}\right\}$$

or equivalently

$$\varsigma(3) = -\frac{4\pi^2}{7}\left\{\sum_{n=0}^{\infty}\frac{\varsigma(2n)}{(2n+1)(2n+2)2^{2n}}\right\}$$

since $\varsigma(0) = -1/2$ (see equation (3.11a)). However, as reported by Srivastava in [125] it was subsequently ascertained that (8.25) was not new: it had in fact been published by Euler in 1772 in his paper Exercitationes Analyticae [126, p.289] and was also rediscovered by Ramaswami in 1934. I also rediscovered this identity in 2004 (see equation (6.71)).

**Example (vi).** Using Bürmann's theorem, it is an exercise in Whittaker & Watson [135, p.130] to prove

(8.26) $$x^2 = \sin^2 x + \left(\frac{2}{3}\right)\frac{1}{2}\sin^4 x + \left(\frac{2.4}{3.5}\right)\frac{1}{3}\sin^6 x + ...$$

(8.27) $$= \sum_{n=1}^{\infty} A_n \sin^{2n} x$$

where $$A_n = \frac{(2n-2)!!}{(2n-1)!!}\frac{1}{n}$$

Multiplying (8.26) by $\cot x$ and integrating, we obtain

(8.28) $$\int_0^{\pi/2} x^2 \cot x\, dx = \sum_{n=1}^{\infty} A_n \int_0^{\pi/2} \sin^{2n-1} x \cos x\, dx$$

(8.29) $$= \frac{1}{2}\sum_{n=1}^{\infty}\frac{1}{n^2}\frac{(2n-2)!!}{(2n-1)!!}$$

(8.30) $$= \frac{1}{8}\sum_{n=1}^{\infty}\frac{2^{2n}[(n-1)!]^2}{n^2(2n-1)!}$$

where we have used the definitions of the double factorials

(8.31) $$(2n)!! = 2.4...(2n) = 2^n n! \qquad (2n+1)!! = 1.3.5...(2n+1) = \frac{(2n+1)!}{2^n n!}$$

Therefore, referring to (8.24) we have shown that

(8.31a) $$\varsigma(3) - 2\pi^2 \log 2 = \sum_{n=1}^{\infty}\frac{2^{2n}[(n-1)!]^2}{n^2(2n-1)!}$$

Let us now divide (8.27) by $\sin x$ and integrate to obtain

$$\int_0^{\pi/2} \frac{x^2}{\sin x}\,dx = \sum_{n=1}^{\infty} A_n \int_0^{\pi/2} \sin^{2n-1} x\,dx$$

In (6.29) we have previously shown that

$$\int_0^{\pi/2} \frac{x^2}{\sin x}\,dx = 2\pi G - \frac{7}{2}\varsigma(3)$$

and the Wiener integral (8.11h) gives us [25, p.113]

$$\int_0^{\pi/2} \sin^{2n-1} x\,dx = \sum_{k=0}^{n-1} \binom{n-1}{k} \frac{(-1)^k}{2k+1} = \frac{(2n-2)!!}{(2n-1)!!}$$

Hence we obtain

$$2\pi G - \frac{7}{2}\varsigma(3) = \sum_{n=1}^{\infty} \frac{1}{n}\left[\frac{(2n-2)!!}{(2n-1)!!}\right]^2$$

or alternatively

(8.31b) $$2\pi G - \frac{7}{2}\varsigma(3) = \frac{1}{16}\sum_{n=1}^{\infty} \frac{2^{4n}\left[(n-1)!\right]^4}{n\left[(2n-1)!\right]^2}$$

**Example (vii).** Ghusayni [70] used the following well-known identity (see [26] and [102]): this identity was known to Euler (see [91]).

(8.32) $$2(\sin^{-1} x)^2 = \sum_{n=1}^{\infty} \frac{(2x)^{2n}}{n^2}\binom{2n}{n}^{-1}, \quad -\frac{1}{2} \leq x \leq \frac{1}{2}$$

and, dividing by $x$ and integrating, we obtain

(8.33) $$2\int_0^{1/2} \frac{(\sin^{-1} x)^2}{x}\,dx = \frac{1}{2}\sum_{n=1}^{\infty} \frac{1}{n^3}\binom{2n}{n}^{-1}$$

Using the substitution $x = \sin t$, and integration by parts, we obtain (as in (8.24))

(8.34) $$2\int_0^{1/2} \frac{(\sin^{-1} x)^2}{x}\,dx = 2\int_0^{\pi/6} t^2 \cot t\,dt$$

As in (1.12) using integration by parts, it is easily seen that

(8.34a)
$$\int_0^{\pi/6} t^2 \cot t\, dt = t^2 \log \sin t \Big|_0^{\pi/6} - 2\int_0^{\pi/6} t \log \sin t\, dt = \frac{\pi^2}{36}\log \sin(\pi/6) - 2\int_0^{\pi/6} t \log \sin t\, dt$$

$$= -\frac{\pi^2}{36}\log 2 - 2\int_0^{\pi/6} t \log \sin t\, dt$$

Then using the basic identity (6.5a)

$$\frac{1}{2}\int_a^b p(x)\cot x\, dx = \sum_{n=0}^{\infty}\int_a^b p(x)\sin 2nx\, dx$$

we obtain

(8.35) $$2\int_0^{1/2}\frac{(\sin^{-1} x)^2}{x}dx = 4\sum_{n=0}^{\infty}\int_0^{\pi/6} x^2 \sin 2nx\, dx$$

We have

$$4\int_0^{\pi/6} x^2 \sin 2nx\, dx = \left(\frac{1}{n^3} - \frac{2x^2}{n}\right)\cos 2nx + \frac{2x\sin 2nx}{n^2}\Big|_0^{\pi/6}$$

$$= \frac{1}{n^3}\cos(n\pi/3) - \frac{1}{n^3} + \frac{\pi}{3n^2}\sin(n\pi/3)$$

Accordingly, we get

$$4\sum_{n=1}^{\infty}\int_0^{\pi/6} x^2 \sin 2nx\, dx = \sum_{n=1}^{\infty}\left[\frac{1}{n^3}\cos(n\pi/3) - \frac{1}{n^3} + \frac{\pi}{3n^2}\sin(n\pi/3)\right]$$

$$= \sum_{n=1}^{\infty}\frac{\cos(n\pi/3)}{n^3} + \frac{\pi}{3}\sum_{n=1}^{\infty}\frac{\sin(n\pi/3)}{n^2} - \varsigma(3)$$

As reported by Lewin [100] and Srivastava and Tsumura ([125a] and [126, p.293]), we have

(8.36a) $$\sum_{n=1}^{\infty}\frac{\cos(n\pi/3)}{n^s} = \frac{1}{2}(6^{1-s} - 3^{1-s} - 2^{1-s} + 1)\varsigma(s)$$

(8.36b) $$\sum_{n=1}^{\infty}\frac{\sin(n\pi/3)}{n^s} = \sqrt{3}\left\{\frac{3^{-s}-1}{2}\varsigma(s) + 6^{-s}\left[\varsigma\left(s,\frac{1}{6}\right) + \varsigma\left(s,\frac{1}{3}\right)\right]\right\}$$

Hence we have with $s = 3$ and $s = 2$ respectively

(8.36c) $$\sum_{n=1}^{\infty}\frac{\cos(n\pi/3)}{n^3} = \frac{1}{3}\varsigma(3)$$

(8.36d) $$\sum_{n=1}^{\infty}\frac{\sin(n\pi/3)}{n^2} = \sqrt{3}\left\{\frac{3^{-2}-1}{2}\varsigma(2) + 6^{-2}\left[\varsigma\left(2,\frac{1}{6}\right) + \varsigma\left(2,\frac{1}{3}\right)\right]\right\}$$

Using (8.36d) we then obtain Ghusayni's result [70]

(8.37) $$\varsigma(3) = \frac{\pi}{2}\sum_{n=1}^{\infty}\frac{\sin(n\pi/3)}{n^2} - \frac{3}{4}\sum_{n=1}^{\infty}\frac{1}{n^3}\binom{2n}{n}^{-1}$$

(8.37a) $$\varsigma(3) = \frac{\pi}{2}\sqrt{3}\left\{\frac{3^{-2}-1}{2}\varsigma(2) + 6^{-2}\left[\varsigma\left(2,\frac{1}{6}\right) + \varsigma\left(2,\frac{1}{3}\right)\right]\right\} - \frac{3}{4}\sum_{n=1}^{\infty}\frac{1}{n^3}\binom{2n}{n}^{-1}$$

From the above we also have

(8.38) $$\sum_{n=1}^{\infty}\frac{1}{n^3}\binom{2n}{n}^{-1} = -2\int_0^{\pi/3} t\log[2\sin(t/2)]\,dt = -8\int_0^{\pi/6} x\log[2\sin x]\,dx$$

and I first came across this result in van der Poorten's 1979 paper "Some wonderful formulae…an introduction to Polylogarithms" [131b]. Using (7.17) we can also evaluate $\int_0^{\pi/6} x\log[2\sin x]\,dx$.

Using (8.32) and letting $t = \sin^{-1} x$ we obtain

$$t^2 = \frac{1}{2}\sum_{n=1}^{\infty}\frac{2^{2n}\sin^{2n} t}{n^2}\binom{2n}{n}^{-1}$$

and, multiplying by $\cot t$ and integrating, we have

$$\int_0^{\pi/6} t^2 \cot t\,dt = \frac{1}{2}\sum_{n=1}^{\infty}\frac{2^{2n}}{n^2}\binom{2n}{n}^{-1}\int_0^{\pi/6}\sin^{2n-1} t\cos t\,dt$$

$$= \frac{1}{4}\sum_{n=1}^{\infty}\frac{1}{n^3}\binom{2n}{n}^{-1}$$

$$= -\frac{\pi^2}{36}\log 2 - 2\int_0^{\pi/6} t\log\sin t\,dt$$

where, in the final part, we have used (8.34a). This therefore proves (8.38).

In a follow-up paper in 2000, Ghusayni [70], having noted an earlier paper [76a], reported that

(8.39) $$\sum_{n=1}^{\infty}\frac{\sin(n\pi/3)}{n^2}=\frac{\sqrt{3}}{2}\left\{\frac{1}{1^2}+\frac{1}{2^2}-\frac{1}{4^2}-\frac{1}{5^2}+\frac{1}{7^2}+\frac{1}{8^2}-\cdots\right\}$$

$$=\frac{\sqrt{3}}{2}\left\{-\frac{2}{9}\pi^2+\frac{1}{3}\psi^{(1)}\left(\frac{1}{3}\right)\right\}$$

(8.40) $$\frac{1}{1^2}+\frac{1}{2^2}-\frac{1}{4^2}-\frac{1}{5^2}+\frac{1}{7^2}+\frac{1}{8^2}-\cdots=-\frac{2}{27}\pi^2-2\int_0^1\frac{\log x}{1+x^3}dx$$

and Ghusayni concluded that

(8.41) $$\varsigma(3)=-\frac{\sqrt{3}}{18}\pi^3+\frac{3\sqrt{3}}{4}\pi\sum_{n=1}^{\infty}\frac{1}{(3n-2)^2}-\frac{3}{4}\sum_{n=1}^{\infty}\frac{1}{n^3}\binom{2n}{n}^{-1}$$

I understand from Moen's paper [103a] that Zucker [142a] has evaluated $\sum_{n=1}^{\infty}\frac{1}{n^3}\binom{2n}{n}^{-1}$ in 1985 (but I have not been able to obtain a copy of that paper). Reference should also be made to the recent paper "Certain series related to the triple sine function" by Koyama and Kurokawa [93]. The Mathworld website for the central binomial coefficient reports the formula

(8.42) $$\sum_{n=1}^{\infty}\frac{1}{n^3}\binom{2n}{n}^{-1}=\frac{1}{18}\pi\sqrt{3}\left[\psi^{(1)}\left(\frac{1}{3}\right)-\psi^{(1)}\left(\frac{2}{3}\right)\right]-\frac{4}{3}\varsigma(3)$$

We have the reflection formula

(8.42a) $$\psi^{(k)}(1-x)+(-1)^{k+1}\psi^{(k)}(x)=(-1)^k\pi\frac{d^k}{dx^k}\cot\pi x$$

and therefore we see that

(8.42b) $$\psi^{(1)}\left(\frac{1}{3}\right)+\psi^{(1)}\left(\frac{2}{3}\right)=\pi^2/\sin^2\left(\frac{\pi}{3}\right)=2\pi^2$$

We then obtain

(8.42c) $$\sum_{n=1}^{\infty}\frac{1}{n^3}\binom{2n}{n}^{-1}=\frac{1}{9}\pi\sqrt{3}\,\psi^{(1)}\left(\frac{1}{3}\right)-\frac{4}{3}\varsigma(3)+\frac{1}{9}\pi^3\sqrt{3}$$

We have

$$\psi^{(k)}\left(\frac{p}{q}\right) = (-1)^{k+1} k! \varsigma\left(k+1, \frac{p}{q}\right)$$

and therefore we see that

$$\psi'\left(\frac{1}{3}\right) = \varsigma\left(2, \frac{1}{3}\right) = 9 \sum_{n=0}^{\infty} \frac{1}{(3n+1)^2}$$

or equivalently

$$\psi'\left(\frac{1}{3}\right) = 9 \sum_{n=1}^{\infty} \frac{1}{(3n-2)^2}$$

With regard to (8.40)

$$\frac{1}{1^2} + \frac{1}{2^2} - \frac{1}{4^2} - \frac{1}{5^2} + \frac{1}{7^2} + \frac{1}{8^2} - \ldots = -\frac{2}{27}\pi^2 - 2\int_0^1 \frac{\log x}{1+x^3} dx$$

Knopp [90, p.236] has recorded that

$$\int \frac{1}{1+x^3} dx = \frac{1}{3}\log(1+x) - \frac{1}{6}\log(1-x+x^2) + \frac{1}{\sqrt{3}} \tan^{-1}\left[\frac{2x-1}{\sqrt{3}}\right]$$

$$= \frac{1}{6}\log(1+x) - \frac{1}{6}\log(1+x^3) + \frac{1}{\sqrt{3}} \tan^{-1}\left[\frac{2x-1}{\sqrt{3}}\right]$$

Therefore, using integration by parts we have

$$\int_a^1 \frac{\log x}{1+x^3} dx = -\left[\frac{1}{6}\log(1+a) - \frac{1}{6}\log(1+a^3) + \frac{1}{\sqrt{3}} \tan^{-1}\left[\frac{2a-1}{\sqrt{3}}\right]\right] \log a$$

$$- \int_a^1 \left[\frac{1}{6}\log(1+x) - \frac{1}{6}\log(1+x^3) + \frac{1}{\sqrt{3}} \tan^{-1}\left[\frac{2x-1}{\sqrt{3}}\right]\right] \frac{dx}{x}$$

and we have designated the lower limit as $a$ so as to hopefully remove a singularity in the integrated part at $x = 0$.

We have

$$\int \frac{\log(1+x)}{x} dx = -Li_2(-x)$$

and
$$\int \frac{\log(1+x^3)}{x} dx = -Li_2(-x^2)$$

Therefore we get

$$\int_0^1 \left[ \frac{1}{6}\log(1+x) - \frac{1}{6}\log(1+x^3) \right] \frac{dx}{x} = 0$$

The final integral is

$$\frac{1}{\sqrt{3}} \int_a^1 \tan^{-1}\left[\frac{2x-1}{\sqrt{3}}\right] \frac{dx}{x} = A \int_b^A \frac{\tan^{-1} t}{t+A} dt$$

where $A = 1/\sqrt{3}$ and $b = A(2a-1)$.

The Wolfram Integrator gives the following result

$$2\int \frac{\tan^{-1} x}{x+A} dx = 2\tan^{-1} x \left[ \frac{1}{2}\log(1+x^2) + \log\left[\sin\left(\tan^{-1} A\right) + \tan^{-1} x\right]\right]$$

$$-\frac{i}{4}\left[\pi - 2\tan^{-1} x\right]^2 - i\left[\tan^{-1} A + \tan^{-1} x\right]^2$$

$$+\left[\pi - 2\tan^{-1} x\right] \log\left[1 - \exp i\left[\pi - 2\tan^{-1} x\right]\right]$$

$$+2\left[\tan^{-1} A + \tan^{-1} x\right] \log\left[1 - \exp 2i\left[\tan^{-1} A + \tan^{-1} x\right]\right]$$

$$-\left[\pi - 2\tan^{-1} x\right] \log\left[2\sin\frac{1}{2}\left[\pi - 2\tan^{-1} x\right]\right]$$

$$-2\left[\tan^{-1} A + \tan^{-1} x\right] \log\left[2\sin\left[\tan^{-1} A + \tan^{-1} x\right]\right]$$

$$-iLi_2\left[\exp i\left[\pi - 2\tan^{-1} x\right]\right] + iLi_2\left[\exp 2i\left[\tan^{-1} A + \tan^{-1} x\right]\right]$$

We have the specific integral with $a = 0$ (whereupon $b = -A$)

$$2\int_{-A}^A \frac{\tan^{-1} x}{x+A} dx = 2\tan^{-1} A \left[\frac{1}{2}\log(1+A^2) + \log\left[\sin\left(\tan^{-1} A\right) + \tan^{-1} A\right]\right]$$

$$+2\tan^{-1} A \left[\frac{1}{2}\log(1+A^2) + \log\left[-\sin\left(\tan^{-1} A\right) - \tan^{-1} A\right]\right]$$

$$-\frac{i}{4}\left[\pi - 2\tan^{-1} A\right]^2 + \frac{i}{4}\left[\pi + 2\tan^{-1} A\right]^2 - i\left[2\tan^{-1} A\right]^2$$

$$+\left[\pi - 2\tan^{-1} A\right]\log\left[1 - \exp i\left[\pi - 2\tan^{-1} A\right]\right]$$

$$-\left[\pi + 2\tan^{-1} A\right]\log\left[1 - \exp i\left[\pi + 2\tan^{-1} A\right]\right]$$

$$+4\tan^{-1} A \log\left[1 - \exp 4i\tan^{-1} A\right]$$

$$-\left[\pi - 2\tan^{-1} A\right]\log\left[2\sin\frac{1}{2}\left[\pi - 2\tan^{-1} A\right]\right]$$

$$+\left[\pi + 2\tan^{-1} A\right]\log\left[2\sin\frac{1}{2}\left[\pi + 2\tan^{-1} A\right]\right]$$

$$-4\tan^{-1} A \log\left[2\sin\left[2\tan^{-1} A\right]\right]$$

$$-iLi_2\left[\exp i\left[\pi - 2\tan^{-1} A\right]\right] + iLi_2\left[\exp i\left[\pi + 2\tan^{-1} A\right]\right]$$

$$+iLi_2\left[\exp 4i\tan^{-1} A\right] - iLi_2(1)$$

We have

$$-iLi_2\left[\exp i\left[\pi - 2\tan^{-1} A\right]\right] + iLi_2\left[\exp i\left[\pi + 2\tan^{-1} A\right]\right]$$

$$= -iLi_2\left[-\exp i\left[-2\tan^{-1} A\right]\right] + iLi_2\left[-\exp i\left[2\tan^{-1} A\right]\right]$$

$$= 2i\sum_{n=1}^{\infty}\frac{(-1)^n \sin\left[2\tan^{-1} A\right]}{n^2} = 2i\sin\left[2\tan^{-1} A\right]\sum_{n=1}^{\infty}\frac{(-1)^n}{n^2}$$

and this is purely imaginary.

Unfortunately, I do not see any simplified expression emerging from this analysis. Using the fact that

$$\sin\left(\tan^{-1} A\right) = \frac{A}{\sqrt{1+A^2}}$$

may be of some assistance. Further work is required here.

□

In 1978, in his celebrated proof of the irrationality of $\varsigma(3)$, Apéry used the following formula (see [11] and [132])

(8.43) $$\sum_{k=1}^{\infty} \frac{a_1 a_2 \ldots a_{k-1}}{(x+a_1)(x+a_2)\ldots(x+a_k)} = \frac{1}{x}$$

to show that

(8.44) $$\varsigma(3) = \frac{5}{2} \sum_{n=1}^{\infty} \frac{(-1)^{n+1}}{n^3} \binom{2n}{n}^{-1}$$

As pointed out by Srivastava in [125], it was not actually Apéry who discovered (8.44): it was actually Markov (the inventor of Markov chains) who found it in 1890 in connection with his series convergence acceleration technique. Further historical information on Markov's technique and its relationship with the much more modern WZ method is given in the paper by Kondratieva and Sadov [91b].

We also know from Tolstov's book on Fourier series [130, p.148] that

$$\sum_{n=1}^{\infty} \frac{\cos nx}{n^3} = \varsigma(3) + \int_0^x dt \int_0^t \log 2\sin(u/2) du \quad \text{for } x \in [0, 2\pi]$$

and

$$\sum_{n=1}^{\infty} \frac{\sin nx}{n^2} = -\int_0^x \log 2\sin(u/2) du \quad \text{for } x \in [0, 2\pi]$$

We now multiply (8.19) by $\sin^{-1} t$, and integrate over the range [0,1/2], to give

(8.45) $$\sum_{n=0}^{\infty} \int_0^{1/2} \frac{1}{2^{2n}(2n+1)} \binom{2n}{n} t^{2n} \sin^{-1} t \, dt = \int_0^{1/2} \frac{(\sin^{-1} t)^2}{t} dt = \frac{1}{4} \sum_{n=1}^{\infty} \frac{1}{n^3} \binom{2n}{n}^{-1}$$

With the substitution $t = \sin x$ we have

$$\int_0^{1/2} t^{2n} \sin^{-1} t \, dt = \int_0^{\pi/6} x \sin^{2n} x \cos x \, dx$$

$$= \frac{x \sin^{2n+1} x}{2n+1} \bigg|_0^{\pi/6} - \frac{1}{2n+1} \int_0^{\pi/6} \sin^{2n+1} x \, dx$$

We have from G&R [74, p.149]

$$\int \sin^{2n+1} x\, dx = \frac{1}{2^{2n}} (-1)^{n+1} \sum_{k=0}^{n} (-1)^k \binom{2n+1}{k} \frac{\cos(2n+1-2k)x}{(2n+1-2k)}$$

The Wiener formula is more elegant here:

$$\int_0^{\pi/6} \sin^{2n+1} x\, dx = \sum_{k=0}^{n} \binom{n}{k} (-1)^{k+1} \frac{\cos^{2k+1} t}{2k+1} \Bigg|_0^{\pi/6}$$

but I don't really see anything of great beauty emerging from this particular proliferation of formulae.

In [90, p.266], Knopp proved the following series expansion and used it to determine the value of $\varsigma(2)$

(8.46) $$(\sin^{-1} x)^2 = \frac{1}{2} \sum_{n=1}^{\infty} \frac{[(n-1)!]^2 (2x)^{2n}}{(2n)!}, \quad (|x| \leq 1)$$

This formula is the same series expansion employed by Ghusayni in (8.32) above. We obtain the following by integration

(8.47) $$\int_0^{1/2} \frac{(\sin^{-1} x)^2 dx}{x} = \frac{1}{2} \sum_{n=1}^{\infty} \int_0^{1/2} \frac{[(n-1)!]^2 (2x)^{2n-1} dx}{(2n)!}$$

(8.47a) $$= \frac{1}{4} \sum_{n=1}^{\infty} \frac{[(n-1)!]^2}{n(2n)!}$$

The formula (8.47a) was actually printed in [90, p.266] with the factorial sign in the numerator being inadvertently omitted and this omission momentarily led me to believe that there may be another remarkable mathematical relationship between $e$ and $\pi$!

**Example (viii).** David Bierens de Haan (1822-1895) provided several proofs of the following integral in his 1862 book "Exposé de la Théorie, Propriétés, des formules de transformation, et des méthodes d'évaluation des intégrales définies" [55a], a copy of which is available on the internet courtesy of the Michigan Historical Mathematics Collection.

(8.48) $$\int_0^{\pi/2} \cos^{p-1} x \cos ax\, dx = \frac{\pi}{2^p} \frac{\Gamma(p)}{\Gamma\left(\frac{p+a+1}{2}\right) \Gamma\left(\frac{p-a+1}{2}\right)}$$

This formula was used by the Borweins [27] who also provided a further proof of it using contour integration.

Differentiating (8.48) with respect to $a$ we immediately obtain

$$(8.48a) \quad \int_0^{\pi/2} x\cos^{p-1} x \sin ax\, dx = \frac{\pi}{2^{p+1}} \Gamma(p) \frac{\psi\left(\frac{p+a+1}{2}\right) - \psi\left(\frac{p-a+1}{2}\right)}{\Gamma\left(\frac{p+a+1}{2}\right) \Gamma\left(\frac{p-a+1}{2}\right)}$$

where $p > 0$ and $-(p+1) < a < p+1$ and $\Gamma(x)$ and $\psi(x)$ are the gamma and digamma functions.

At a very early stage of this series of papers, I made reference to the above integral (which appears in Gradshteyn and Ryzhik [74, 3.832 1]) and used it to corroborate one of my results in (3.8a). Coffey has also used the same integral, but to much more advantage!

Part of the following analysis is based on Coffey's recent paper "On some log-cosine integrals related to $\varsigma(3), \varsigma(4)$ and $\varsigma(6)$" [45b].

For convenience, letting $p_1 = \left(\frac{p+a+1}{2}\right)$ and $p_2 = \left(\frac{p-a+1}{2}\right)$ we obtain by differentiating (8.48a)

$$(8.49) \quad \frac{\partial}{\partial a} I(a,p) = \frac{\pi}{2^{p+2}} \frac{\Gamma(p)}{\Gamma(p_1)\Gamma(p_2)} \left\{ -[\psi(p_1) - \psi(p_2)]^2 + \psi'(p_1) + \psi'(p_2) \right\}$$

$$= \int_0^{\pi/2} x^2 \cos^{p-1} x \cos ax\, dx$$

Putting $a = 0$ in (8.49) we get

$$(8.50) \quad J_p = \frac{\partial}{\partial a} I(a,p) \bigg|_{a=0} = \int_0^{\pi/2} x^2 \cos^{p-1} x\, dx = \frac{\pi}{2^{p+1}} \frac{\Gamma(p)}{\Gamma^2[(p+1)/2]} \psi'[(p+1)/2]$$

Then, differentiating (8.50) with respect to the parameter $p$, we have

$$\frac{d}{dp} J_p = \int_0^{\pi/2} x^2 \cos^{p-1} x \log \cos x\, dx$$

$$= \frac{\pi}{2^{p+1}} \frac{\Gamma(p)}{\Gamma^2[(p+1)/2]} \psi'[(p+1)/2] \left\{ -\log 2 + \psi(p) - \psi[(p+1)/2] + \frac{\psi''[(p+1)/2]}{2\psi'[(p+1)/2]} \right\}$$

Evaluating this at $p = 1$, and using (E.22j), we obtain

$$\int_0^{\pi/2} x^2 \log \cos x \, dx = \frac{\pi^3}{24}\left[-\log 2 + \frac{\psi''(1)}{\psi'(1)}\right]$$

and hence we have

(8.51) $$\int_0^{\pi/2} x^2 \log \cos x \, dx = -\frac{\pi^3}{24}\log 2 - \frac{\pi}{2}\varsigma(3)$$

Similarly, Coffey obtains

$$\frac{d^2}{dp^2}J_p = \int_0^{\pi/2} x^2 \cos^{p-1} x \log^2 \cos x \, dx$$

$$= \frac{\pi}{2^{p+1}}\frac{\Gamma(p)}{\Gamma^2[(p+1)/2]}\psi'[(p+1)/2]\left\{\begin{array}{l}\left[-\log 2 + \psi(p) - \psi[(p+1)/2] + \frac{\psi''[(p+1)/2]}{2\psi'[(p+1)/2]}\right]^2 \\ +\psi'(p) - \frac{1}{2}\psi'[(p+1)/2] + \frac{1}{2}\frac{d}{dp}\frac{\psi''[(p+1)/2]}{\psi'[(p+1)/2]}\end{array}\right\}$$

and in particular we have

(8.52) $$\int_0^{\pi/2} x^2 \log^2 \cos x \, dx = \frac{\pi}{1440}\left[11\pi^4 + 60\pi^2 \log^2 2 + 720\varsigma(3)\log 2\right]$$

Using the substitution $x = \theta/2$ we see that

$$\int_0^{\pi} \theta^2 \log^2[2\cos(\theta/2)]d\theta = 8\int_0^{\pi/2} x^2(\log 2 + \log \cos x)^2 dx$$

$$= 8\log^2 2 \int_0^{\pi/2} x^2 dx + 16\log 2 \int_0^{\pi/2} x^2 \log \cos x \, dx + 8\log^2 2 \int_0^{\pi/2} x^2 \log^2 \cos x \, dx$$

Employing (8.51) and (8.52) we find that

(8.52a) $$\int_0^{\pi} \theta^2 \log^2[2\cos(\theta/2)]d\theta = \frac{11\pi^4}{180} = \frac{11}{4}\varsigma(4)$$

As mentioned by Coffey, this integral was first evaluated by D. Borwein and J. M. Borwein in their 1995 paper, "On an Intriguing Integral and Some Series Related to $\varsigma(4)$" [27] and used by them to prove (4.2.42).

$$\sum_{n=1}^{\infty}\frac{\left(H_n^{(1)}\right)^2}{n^2}=\frac{17\pi^4}{360}=\frac{17}{4}\varsigma(4)$$

In (3.46a) we showed that

(8.53) $$\sum_{n=1}^{\infty}\frac{\left(H_n^{(1)}\right)^2}{n}x^n=-\frac{1}{3}\log^3(1-x)+Li_3(x)-Li_2(x)\log(1-x)$$

where we note that the series obviously does not converge at $x=1$. Using the same method employed by the Borweins [27] and Coffey [45b], we try to obtain a Fourier-type series from (8.53) by making the substitution $x=e^{it}$. We have

$$\log(1-e^{it})=\log e^{it/2}(e^{-it/2}-e^{it/2})=\log\left[-2i\sin(t/2)e^{it/2}\right]$$

$$=\log\left[2\sin(t/2)e^{it/2}\right]+i\frac{t}{2}-\log i$$

$$=\log\left[2\sin(t/2)\right]+\frac{i}{2}(t-\pi)$$

Therefore we have

$$\log^2(1-e^{it})=\log^2\left[2\sin(t/2)\right]-\frac{1}{4}(t-\pi)^2+i(t-\pi)\log\left[2\sin(t/2)\right]$$

which was used in [27] and

$$\log^3(1-e^{it})=\log^3\left[2\sin(t/2)\right]-\frac{3}{4}(t-\pi)^2\log\left[2\sin(t/2)\right]$$

$$+i\frac{3}{2}(t-\pi)\log^2\left[2\sin(t/2)\right]-i\frac{1}{8}(t-\pi)^3$$

Then, by substituting the latter formula in (8.53) we get

(8.53a)

$$\sum_{n=1}^{\infty}\frac{\left(H_n^{(1)}\right)^2}{n}\cos nt=-\frac{1}{3}\left\{\log^3\left[2\sin(t/2)\right]-\frac{3}{4}(t-\pi)^2\log\left[2\sin(t/2)\right]\right\}$$

$$+\sum_{n=1}^{\infty}\frac{\cos nt}{n^3}-\log\left[2\sin(t/2)\right]\sum_{n=1}^{\infty}\frac{\cos nt}{n^2}+\frac{1}{2}(t-\pi)\sum_{n=1}^{\infty}\frac{\sin nt}{n^2}$$

(8.53b)

$$\sum_{n=1}^{\infty}\frac{\left(H_n^{(1)}\right)^2}{n}\sin nt = -\frac{1}{3}\left\{\frac{3}{2}(t-\pi)\log^2[2\sin(t/2)]-\frac{1}{8}(t-\pi)^2\right\}+\sum_{n=1}^{\infty}\frac{\sin nt}{n^3}$$

$$-\frac{1}{2}(t-\pi)\sum_{n=1}^{\infty}\frac{\cos nt}{n^2}-\log[2\sin(t/2)]\sum_{n=1}^{\infty}\frac{\sin nt}{n^2}$$

Unfortunately, we do not quite obtain Fourier series and hence we are unable to apply Parseval's theorem to the above identities.

It may also be useful to integrate (8.53a).

Differentiating (8.48) with respect to $p$ we immediately obtain

(8.54) $\displaystyle\int_0^{\pi/2}\cos^{p-1}x\cos ax\log\cos x\,dx =$

$$-\frac{\pi\log 2}{2^p}\frac{\Gamma(p)}{\Gamma\left(\frac{p+a+1}{2}\right)\Gamma\left(\frac{p-a+1}{2}\right)}+\frac{\pi}{2^p}\frac{\Gamma'(p)}{\Gamma\left(\frac{p+a+1}{2}\right)\Gamma\left(\frac{p-a+1}{2}\right)}$$

$$-\frac{\pi}{2^{p+1}}\Gamma(p)\frac{\Gamma\left(\frac{p+a+1}{2}\right)\Gamma'\left(\frac{p-a+1}{2}\right)+\Gamma'\left(\frac{p+a+1}{2}\right)\Gamma\left(\frac{p-a+1}{2}\right)}{\Gamma\left(\frac{p+a+1}{2}\right)\Gamma\left(\frac{p-a+1}{2}\right)}$$

$$=-\frac{\pi\log 2}{2^p}\frac{\Gamma(p)}{\Gamma\left(\frac{p+a+1}{2}\right)\Gamma\left(\frac{p-a+1}{2}\right)}+\frac{\pi}{2^p}\frac{\Gamma'(p)}{\Gamma\left(\frac{p+a+1}{2}\right)\Gamma\left(\frac{p-a+1}{2}\right)}$$

$$-\frac{\pi}{2^{p+1}}\Gamma(p)\left[\psi\left(\frac{p-a+1}{2}\right)+\psi\left(\frac{p+a+1}{2}\right)\right]$$

With $p=2$ and $a=1$ we obtain

$$\int_0^{\pi/2}\cos^2 x\log\cos x\,dx=\frac{\pi}{4}\left[-\log 2+\Gamma'(2)-\frac{1}{2}\psi(1)-\frac{1}{2}\psi(2)\right]$$

We have from (4.3.16) $\psi(n+1)=H_n-\gamma$ and hence

$\psi(1)=-\gamma$ and $\psi(2)=1-\gamma$

and from (E.20) we have

$$\psi(m) = \frac{\Gamma'(m)}{\Gamma(m)} = -\gamma - \sum_{k=0}^{\infty}\left(\frac{1}{m+k} - \frac{1}{k+1}\right)$$

Hence

$$\Gamma'(2) = \psi(2) = 1 - \gamma$$

Therefore we get

(8.55) $$\int_{0}^{\pi/2} \cos^2 x \log \cos x \, dx = \frac{\pi}{4}\left[\frac{1}{2} - \log 2\right]$$

which is in agreement with the more general formula given in G&R [74, p.581].

We have from (2.25) and (2.26) for suitably behaved functions

(8.56a) $$\sum_{n=1}^{\infty}\int_{a}^{b} p(x)\cos^n x \, \cos nx \, dx = 0$$

(8.56b) $$\sum_{n=1}^{\infty}\int_{a}^{b} p(x)\cos^n x \, \sin nx \, dx = \int_{a}^{b} p(x)\cot x \, dx$$

The existence of the above identities is formally explained in a neat way by the following analysis. From (2.9) we have

(8.56c) $$2^n \cos^n x(\cos nx + i\sin nx) = \left(1 + e^{2ix}\right)^n$$

which gives us the infinite geometric series

$$\sum_{n=1}^{\infty} \cos^n x(\cos nx + i\sin nx) = \sum_{n=1}^{\infty}\left(\frac{1+e^{2ix}}{2}\right)^n$$

provided $\left|\dfrac{1+e^{2ix}}{2}\right| < 1$. We obtain

$$\sum_{n=1}^{\infty}\left(\frac{1+e^{2ix}}{2}\right)^n = \frac{1+e^{2ix}}{1-e^{2ix}} = i\cot x$$

Hence equating real and imaginary parts we get (ignoring issues of convergence!)

$$\sum_{n=1}^{\infty} \cos^n x \cos nx = 0$$

$$\sum_{n=1}^{\infty} \cos^n x \sin nx = \cot x$$

Much later, I discovered that the following formulae are contained in Ramanujan's Notebook [21, Part I, p.246] for $x \leq \pi/2$

$$2^n \cos^n x \cos(a+n)x = \sum_{k=1}^{\infty} \binom{n}{k} \cos(a+2k)x$$

$$2^n \cos^n x \sin(a+n)x = \sum_{k=1}^{\infty} \binom{n}{k} \sin(a+2k)x$$

and it may be noted that the second formula may be obtained by differentiating the first one with respect to $a$.

From (8.49) we have

$$\frac{\partial}{\partial a} I(a,p) = \frac{\pi}{2^{p+2}} \frac{\Gamma(p)}{\Gamma(p_1)\Gamma(p_2)} \left\{ -[\psi(p_1) - \psi(p_2)]^2 + \psi'(p_1) + \psi'(p_2) \right\}$$

$$= \int_0^{\pi/2} x^2 \cos^{p-1} x \cos ax \, dx$$

Therefore, with $a = n$ and $p = n+1$ we get

(8.57) $$\int_0^{\pi/2} x^2 \cos^n x \cos nx \, dx = \frac{\pi}{2^{n+3}} \left\{ -[\psi(n+1) - \psi(1)]^2 + \psi'(n+1) + \psi'(1) \right\}$$

We have from (4.3.16) and (E.16a)

(8.57a) $$\psi(n+1) = H_n^{(1)} - \gamma$$

$$\psi(n+1) - \psi(1) = H_n^{(1)}$$

$$\psi'(m) = \sum_{n=0}^{\infty} \frac{1}{(m+n)^2} = \varsigma(2) - H_{m-1}^{(2)}$$

$$\psi'(1) = \varsigma(2)$$

$$\psi'(n+1) - \psi'(1) = -H_n^{(2)}$$

$$\psi''(x) = -2\sum_{k=0}^{\infty} \frac{1}{(x+k)^3}$$

$$\psi''(n+1) = -2\left\{\frac{1}{(n+1)^3} + ...\right\} = -2\left\{\varsigma(3) - H_n^{(3)}\right\}$$

$$\psi''(1) = -2\varsigma(3)$$

$$\psi''(n+1) - \psi''(1) = 2H_n^{(3)}$$

and we therefore get

(8.58) $$\int_0^{\pi/2} x^2 \cos^n x \cos nx \, dx = \frac{\pi}{2^{n+3}}\left\{-\left[H_n^{(1)}\right]^2 + 2\varsigma(2) - H_n^{(2)}\right\}$$

Completing the summation we have

$$\sum_{n=1}^{\infty} \int_0^{\pi/2} x^2 \cos^n x \cos nx \, dx = \sum_{n=1}^{\infty} \frac{\pi}{2^{n+3}}\left\{-\left[H_n^{(1)}\right]^2 + 2\varsigma(2) - H_n^{(2)}\right\}$$

$$= -\frac{\pi}{8}\left[\sum_{n=1}^{\infty} \frac{\left[H_n^{(1)}\right]^2}{2^n} + \sum_{n=1}^{\infty} \frac{H_n^{(2)}}{2^n} - 2\varsigma(2)\right]$$

We recall from (3.33) that

$$\sum_{n=1}^{\infty} H_n^{(r)} x^n = \frac{Li_r(x)}{1-x} \quad , x \in [0,1)$$

which gives us

$$\sum_{n=1}^{\infty} \frac{H_n^{(r)}}{2^n} = 2Li_r(1/2)$$

We also have from (3.35)

$$\frac{\log^2(1-x) + Li_2(x)}{1-x} = \sum_{n=1}^{\infty} \left(H_n^{(1)}\right)^2 x^n \quad , x \in [0,1)$$

which results in

$$2\log^2 2 + 2Li_2(1/2) = \sum_{n=1}^{\infty} \frac{\left(H_n^{(1)}\right)^2}{2^n}$$

Therefore we have

(8.58a) $$\sum_{n=1}^{\infty} \frac{\left(H_n^{(1)}\right)^2}{2^n} + \sum_{n=1}^{\infty} \frac{H_n^{(2)}}{2^n} - 2\varsigma(2) = 2\log^2 2 + 4Li_2(1/2) - 2\varsigma(2) = 0$$

where we have employed the dilogarithm identity (3.43a).

A further differentiation results in

$$\frac{\partial^2}{\partial a^2} I(a,p) = \int_0^{\pi/2} x^3 \cos^{p-1} x \sin ax \, dx$$

We also have

$$\frac{2^{p+2}}{\pi} \frac{\Gamma(p_1)\Gamma(p_2)}{\Gamma(p)} \frac{\partial^2}{\partial a^2} I(a,p) =$$

$$\left\{ -2[\psi(p_1) - \psi(p_2)] \left[ \frac{1}{2}\psi'(p_1) + \frac{1}{2}\psi'(p_2) \right] + \frac{1}{2}\psi''(p_1) - \frac{1}{2}\psi''(p_2) \right\}$$

$$+ \frac{1}{\Gamma(p_1)\Gamma(p_2)} \left\{ -[\psi(p_1) - \psi(p_2)]^2 + \psi'(p_1) + \psi'(p_2) \right\} \frac{1}{2} [\Gamma'(p_1)\Gamma(p_2) - \Gamma(p_1)\Gamma'(p_2)]$$

and this may be written as

$$\frac{2^{p+2}}{\pi} \frac{\Gamma(p_1)\Gamma(p_2)}{\Gamma(p)} \frac{\partial^2}{\partial a^2} I(a,p) =$$

$$= \frac{\pi}{2^{p+3}} \frac{\Gamma(p)}{\Gamma(p_1)\Gamma(p_2)} \left\{ -2[\psi(p_1) - \psi(p_2)][\psi'(p_1) + \psi'(p_2)] + \psi''(p_1) - \psi''(p_2) \right\}$$

$$+ \frac{\pi}{2^{p+3}} \frac{\Gamma(p)}{\Gamma(p_1)\Gamma(p_2)} \left\{ -[\psi(p_1) - \psi(p_2)]^2 + \psi'(p_1) + \psi'(p_2) \right\} [\psi(p_1) - \psi(p_2)]$$

As before, with $a = n$ and $p = n+1$ we get

$$\frac{\partial^2}{\partial a^2} I(n, n+1) =$$

$$= \frac{\pi}{2^{n+4}} \left\{ -2[\psi(n+1) - \psi(1)][\psi'(n+1) + \psi'(1)] + \psi''(n+1) - \psi''(1) \right\}$$

$$+\frac{\pi}{2^{n+4}}\left\{-\left[\psi(n+1)-\psi(1)\right]^2+\psi'(n+1)+\psi'(1)\right\}\left[\psi'(n+1)-\psi'(1)\right]$$

and this simplifies to

$$=\frac{\pi}{2^{n+4}}\left\{-2H_n^{(1)}\left[2\varsigma(2)-H_n^{(2)}\right]+2H_n^{(3)}\right\}$$

$$-\frac{\pi}{2^{n+4}}\left\{-\left[H_n^{(1)}\right]^2+2\varsigma(2)-H_n^{(2)}\right\}H_n^{(2)}$$

We then have

$$\frac{\partial^2}{\partial a^2}I(n,n+1)=$$

$$\frac{\pi}{2^{n+3}}\left\{-H_n^{(1)}\left[2\varsigma(2)-H_n^{(2)}\right]+H_n^{(3)}\right\}-\frac{\pi}{2^{n+4}}\left\{-\left[H_n^{(1)}\right]^2+2\varsigma(2)-H_n^{(2)}\right\}H_n^{(2)}$$

At first glance, the "dimensions" of the above equation in terms of $H_n^{(r)}$ appear to be incorrect, but sanity is restored by noting from (8.58a) that, after making the summation, the sum in the second set of braces is identically zero.

We have from (8.56b)

$$\sum_{n=1}^{\infty}\int_0^{\pi/2} x^3 \cos^n x \sin nx\, dx = \int_0^{\pi/2} x^3 \cot x\, dx$$

and completing the summation we have

$$\sum_{n=1}^{\infty}\int_0^{\pi/2} x^3 \cos^n x \cos nx\, dx = -\frac{\pi}{4}\varsigma(2)\sum_{n=1}^{\infty}\frac{H_n^{(1)}}{2^n}+\frac{\pi}{8}\sum_{n=1}^{\infty}\frac{H_n^{(1)}H_n^{(2)}}{2^n}+\frac{\pi}{8}\sum_{n=1}^{\infty}\frac{H_n^{(3)}}{2^n}$$

We already know from (3.8c) that $\sum_{n=1}^{\infty}\frac{H_n^{(1)}}{2^n}=2\log 2$ and, as noted above, we also have $\sum_{n=1}^{\infty}\frac{H_n^{(3)}}{2^n}=2Li_3(1/2)$: we are therefore positioned to evaluate $\sum_{n=1}^{\infty}\frac{H_n^{(1)}H_n^{(2)}}{2^n}$. The integral $\int_0^{\pi/2} x^3 \cot x\, dx$ may be easily ascertained by using (6.5a) but, coincidentally, I note that Amigó [7a] has already determined this integral by an entirely different method. We have

(8.59) $$\int_0^{\pi/2} x^3 \cot x \, dx = \frac{1}{8}\pi^3 \log 2 - \frac{9}{10}\pi \varsigma(3)$$

and therefore we get

(8.60) $$\sum_{n=1}^{\infty} \frac{H_n^{(1)} H_n^{(2)}}{2^n} = \pi^2 \log 2 - \frac{72}{10}\varsigma(3) + 4\varsigma(2)\log 2 - 2Li_3(1/2)$$

$$= \frac{11}{6}\pi^2 \log 2 - \frac{179}{20}\varsigma(3) - \frac{1}{3}\log^3 2$$

(where we have used Landen's formula (3.43b)).

We also have using (3.23)

$$\sum_{n=1}^{\infty} \frac{H_n^{(1)} H_n^{(2)}}{2^n} = \sum_{n=1}^{\infty} \frac{H_n^{(1)}}{2^n} \sum_{k=1}^{n} \frac{1}{k^2} = \sum_{n=1}^{\infty} \frac{1}{n^2} \sum_{k=n}^{\infty} \frac{H_k^{(1)}}{2^k}$$

$$= \sum_{n=1}^{\infty} \frac{1}{n^2} \left[ \sum_{k=1}^{\infty} \frac{H_k^{(1)}}{2^k} - \sum_{k=1}^{n-1} \frac{H_k^{(1)}}{2^k} \right]$$

$$= \varsigma(2) \log 2 - \sum_{n=1}^{\infty} \frac{1}{n^2} \sum_{k=1}^{n-1} \frac{H_k^{(1)}}{2^k}$$

□

Being rather brave, we venture to take the third derivative and obtain

$$\frac{\partial^3}{\partial a^3} I(a, p) = \int_0^{\pi/2} x^4 \cos^{p-1} x \cos ax \, dx$$

and this is also equal to

$$\frac{2^{p+4}}{\pi} \frac{\Gamma(p_1)\Gamma(p_2)}{\Gamma(p)} \int_0^{\pi/2} x^4 \cos^{p-1} x \cos ax \, dx =$$

$$\begin{cases} -2[\psi(p_1)-\psi(p_2)][\psi''(p_1)-\psi''(p_2)] \\ -2[\psi'(p_1)+\psi'(p_2)][\psi'(p_1)+\psi'(p_2)] + \psi'''(p_1)+\psi'''(p_2) \end{cases}$$

$$+ \{-2[\psi(p_1)-\psi(p_2)][\psi'(p_1)+\psi'(p_2)] + \psi''(p_1)-\psi''(p_2)\}[\psi'(p_1)-\psi'(p_2)]$$

$$+\left\{-[\psi(p_1)-\psi(p_2)]^2+\psi'(p_1)+\psi'(p_2)\right\}[\psi''(p_1)+\psi''(p_2)]$$

$$+\left\{-2[\psi(p_1)-\psi(p_2)][\psi'(p_1)+\psi'(p_2)]+\psi''(p_1)-\psi''(p_2)\right\}[\psi'(p_1)-\psi'(p_2)]$$

$$+\left\{-[\psi(p_1)-\psi(p_2)]^2+\psi'(p_1)+\psi'(p_2)\right\}[\psi'(p_1)-\psi'(p_2)]^2$$

Substituting $a=n$ and $p=n+1$ (and noting that the second line of the above equation is equal to the forth line) we get

$$=\frac{\pi}{2^{n+5}}\left\{-2[\psi(n+1)-\psi(1)][\psi''(n+1)-\psi''(1)]-2[\psi'(n+1)+\psi'(1)]^2+\psi'''(n+1)+\psi'''(1)\right\}$$

$$+\frac{2\pi}{2^{n+5}}\left\{-2[\psi(n+1)-\psi(1)][\psi'(n+1)+\psi'(1)]+\psi''(n+1)-\psi''(1)\right\}[\psi'(n+1)-\psi'(1)]$$

$$+\frac{\pi}{2^{n+5}}\left\{-[\psi(n+1)-\psi(1)]^2+\psi'(n+1)+\psi'(1)\right\}[\psi''(n+1)+\psi''(1)]$$

$$+\frac{\pi}{2^{n+5}}\left\{-[\psi(n+1)-\psi(1)]^2+\psi'(n+1)+\psi'(1)\right\}[\psi'(n+1)-\psi'(1)]^2$$

This then becomes

$$\frac{2^{n+5}}{\pi}\frac{\partial^3}{\partial a^3}I(a,p)=$$

$$-2H_n^{(1)}H_n^{(3)}-2\left[2\varsigma(2)-H_n^{(2)}\right]^2+12\varsigma(4)-6H_n^{(4)}$$

$$-2\left\{-2H_n^{(1)}\left[2\varsigma(2)-H_n^{(2)}\right]+2H_n^{(3)}\right\}H_n^{(2)}$$

$$+\left\{-\left[H_n^{(1)}\right]^2+2\varsigma(2)-H_n^{(2)}\right\}\left\{-4\varsigma(3)+2H_n^{(3)}\right\}$$

$$+\left\{-\left[H_n^{(1)}\right]^2+2\varsigma(2)-H_n^{(2)}\right\}\left[H_n^{(2)}\right]^2$$

Reference to (8.58a) shows that, after making the relevant summation, the series comprising the third and forth lines of the above equation will respectively become equal to zero and hence we get

$$\frac{2^{n+5}}{\pi}\frac{\partial^3}{\partial a^3}I(a,p)=$$

$$-2H_n^{(1)}H_n^{(3)} - 2\left[2\varsigma(2) - H_n^{(2)}\right]^2 + 12\varsigma(4) - 6H_n^{(4)}$$

$$-2\left\{-2H_n^{(1)}\left[2\varsigma(2) - H_n^{(2)}\right] + 2H_n^{(3)}\right\}H_n^{(2)} \text{ (plus other bits which we may}$$

disregard).

The first line of the above equation is of dimension 4 in terms of $H_n^{(r)}$ whereas curiously the second line has a dimension of 5.

Completing the summation we have

$$\sum_{n=1}^{\infty}\int_0^{\pi/2} x^4 \cos^n x \cos nx\, dx = \sum_{n=1}^{\infty} \frac{\pi}{2^{n+5}}\left\{-2H_n^{(1)}H_n^{(3)} - 2\left[2\varsigma(2) - H_n^{(2)}\right]^2 + 12\varsigma(4) - 6H_n^{(4)}\right\}$$

$$+\sum_{n=1}^{\infty} \frac{\pi}{2^{n+5}}\left\{-2\left\{-2H_n^{(1)}\left[2\varsigma(2) - H_n^{(2)}\right] + 2H_n^{(3)}\right\}H_n^{(2)}\right\}$$

Purely based on "dimensional" grounds, my conjecture is that

(8.61) $$\sum_{n=1}^{\infty} \frac{1}{2^n}\left\{-2H_n^{(1)}H_n^{(3)} - 2\left[2\varsigma(2) - H_n^{(2)}\right]^2 + 12\varsigma(4) - 6H_n^{(4)}\right\} = 0$$

and

(8.62) $$\sum_{n=1}^{\infty} \frac{1}{2^n}\left\{-2\left\{-2H_n^{(1)}\left[2\varsigma(2) - H_n^{(2)}\right] + 2H_n^{(3)}\right\}H_n^{(2)}\right\} = 0$$

It is clear that higher derivatives will produce additional identities and probably an algebraic nightmare!

□

The application of the method of Section 6 to integrands containing the term

$$\frac{1}{1-ue^{ix}} = \frac{1-ue^{-ix}}{1-2u\cos x + u^2}$$

will be considered in a (much shorter) follow-up paper.

**REFERENCES**

[1]  M. Abramowitz and I.A. Stegun (Eds.), Handbook of Mathematical Functions with Formulas, Graphs and Mathematical Tables. Dover, New York, 1970.
 http://www.math.sfu.ca/~cbm/aands/

[2]  V.S.Adamchik, On Stirling Numbers and Euler Sums.
J.Comput.Appl.Math.79, 119-130, 1997.
http://www-2.cs.cmu.edu/~adamchik/articles/stirling.htm

[2a] V.S.Adamchik, A Class of Logarithmic Integrals. Proceedings of the 1997 International Symposium on Symbolic and Algebraic Computation.
ACM, Academic Press, 1-8, 2001.
http://www-2.cs.cmu.edu/~adamchik/articles/logs.htm

[3]  V.S.Adamchik and H.M. Srivastava, Some Series of the Zeta and Related Functions. Analysis 18, 131-144, 1998.
http://www-2.cs.cmu.edu/~adamchik/articles/sums.htm

[4]  V.S. Adamchik, Polygamma Functions of Negative Order. J.Comp.and Applied Math.100, 191-199, 1998. Polygamma Functions of Negative Order

[5]  V.S. Adamchik, Some definite Integrals Associated with the Riemann Zeta Function. Journal for Analysis and its Applications (ZAA), 19, 831-846, 2000.
 http://www-2.cs.cmu.edu/~adamchik/articles/zaa.htm

[5a] V.S. Adamchik, On the Barnes Function. Proceedings of the 2001 International Symposium on Symbolic and Algebraic Computation. (July 22-25, 2001, London, Canada), Academic Press, 15-20, 2001
http://www-2.cs.cmu.edu/~adamchik/articles/issac01/issac01.pdf

[5b] V.S. Adamchik, Symbolic and numeric computations of the Barnes function. Computer Physics Communications, 157 (2004) 181-190.

[6]  V.S.Adamchik, Certain Series Associated with Catalan's Constant. Journal for Analysis and its Applications (ZAA), 21, 3(2002), 817-826.
http://www-2.cs.cmu.edu/~adamchik/articles/csum.html

[6a] V.S.Adamchik, Contributions to the Theory of the Barnes Function. Computer Physics Communications, 2003.
http://www-2.cs.cmu.edu/~adamchik/articles/barnes1.pdf

[6b] V.S.Adamchik, Symbolic and numeric computations of the Barnes function. Computer Physics Comms., 157 (2004) 181-190.
Symbolic and numeric computations of the Barnes function

[6c] V.S.Adamchik, The multiple gamma function and its application to computation of series. The Ramanujan Journal, 9, 271-288, 2005.

[6x] O.R. Ainsworth and L.W. Howell, The generalized Euler-Mascheroni constants.

NASA Centre for AeroSpace Information (CASI)
NASA-TP-2264 ;NAS 1.60.2264, 1984. View PDF File

[6y] O.R. Ainsworth and L.W. Howell, An integral representation of the generalized Euler-Mascheroni constants.
NASA Centre for AeroSpace Information (CASI)
NASA-TP-2456 ;NAS 1.60.2456, 1985. View PDF File

[6ai] S. Akiyama and Y. Tanigawa, Multiple zeta values at non-positive integers.
The Ramanujan Journal, Vol.5, No.4, 327-351, 2001.

[6aii] U. Alfred, The Amateur Mathematician.Math.Mag, 34, 311-315, 1961.

[6aiii] J.-P. Allouche, J. Shallit and J. Sondow, Summation of Series Defined by Counting Blocks of Digits. math.NT/0512399 [abs, ps, pdf, other] 2005.
J. Number Theory 123 (2007) 133-143

[7] J.M. Amigó, Relations among sums of reciprocal powers.
Israel Journal of Math.124, 177-184, 2001.

[7a] J.M. Amigó, Relations among sums of reciprocal powers II.
http://www.crm.es/Publications/01/483.pdf

[8] P. Amore, Convergence acceleration of series through a variational approach.
Math-ph/0408036 [abs, ps, pdf, other]

[8a] G.E. Andrews, R. Askey and R. Roy, Special Functions.
Cambridge University Press, Cambridge, 1999.

[9] J. Anglesio, A fairly general family of integrals.
Amer.Math.Monthly, 104, 665-666, 1997.

[10] F. Apéry, Roger Apéry, 1916-1999: A Radical mathematician.
The Mathematical Intelligencer, 18, No.2, 54-61, 1996.
*Roger Apéry, 1916-1999 : A Radical Mathematician*

[11] R. Apéry, Irrationalité de $\varsigma(2)$ et $\varsigma(3)$ in Journées Arithmétiques de Luminy (Colloq. Internat. CRNS, Centre Univ. Luminy, 1978).
Astérisque, 61, Soc. Math. France, Paris11-13, 1979.

[12] T.M. Apostol, Another Elementary Proof of Euler's Formula for $\varsigma(2n)$.
Amer.Math.Monthly, 80,425-431, 1973.

[13] T.M. Apostol, Mathematical Analysis, Second Ed., Addison-Wesley Publishing Company, Menlo Park (California), London and Don Mills (Ontario), 1974.

[14] T.M. Apostol, Introduction to Analytic Number Theory.
Springer-Verlag, New York, Heidelberg and Berlin, 1976.

[14aa] T.M. Apostol, Formulas for Higher Derivatives of the Riemann Zeta Function.

Math. of Comp., 169, 223-232, 1985.

[14a] T.M. Apostol, An Elementary View of Euler's Summation Formula.
Amer.Math.Monthly, 106, 409-418, 1999.

[14b] T.M. Apostol, Remark on the Hurwitz zeta function.
Proc.Amer.Math.Soc., 5, 690-693, 1951.

[15] R. Ayoub, Euler and the Zeta Function.
Amer.Math.Monthly, 81, 1067-1086, 1974.

[16] D.H. Bailey, J.M. Borwein and R.Girgensohn, Experimental Evaluation of Euler Sums.
[Experimental Evaluation of Euler Sums - Bailey, Borwein, Girgensohn (1994)](#)

[16a] D.H. Bailey, J.M. Borwein, V. Kapoor and E. Weisstein, Ten problems in experimental mathematics. Amer.Math.Monthly, 481-509, 2006.

[17] D.H. Bailey, P.B. Borwein and S.Plouffe, On the rapid computation of various polylogarithmic constants. Mathematics of Computation, 66(218), 903-913, 1997.
[On the Rapid Computation of Various Polylogarithmic Constants](#)

[17aa] E.W. Barnes, The theory of the G-function. Quart. J. Math.31, 264-314, 1899.

[17a] R.G. Bartle, The Elements of Real Analysis. 2$^{nd}$ Ed. John Wiley & Sons Inc, New York, 1976.

[18] E.T. Bell, Men of Mathematics. Simon & Schuster, New York, 1937.

[19] B.C. Berndt, Elementary Evaluation of $\varsigma(2n)$. Math.Mag.48, 148-153, 1975.

[20] B.C. Berndt, The Gamma Function and the Hurwitz Zeta Function.
Amer.Math. Monthly, 92,126-130, 1985.

[21] B.C. Berndt, Ramanujan's Notebooks. Parts I-III, Springer-Verlag, 1985-1991.

[22] J.L.F. Bertrand, Traité de Calcul Différentiel et de Calcul Intégral (1864).
http://math-sahel.ujf-grenoble.fr/LiNuM/TM/Gallica/S099558.html
http://math-sahel.ujf-grenoble.fr/LiNuM/TM/Gallica/S099559.html

[23] F. Beukers, A note on the irrationality of $\varsigma(2)$ and $\varsigma(3)$.
Bull. London Math.Soc.11, 268-272, 1979.

[23aa] L. Berggren, J. Borwein and P. Borwein, Pi: A Source Book.
Springer-Verlag, New York, 1997.

[23a] M.G. Beumer, Some special integrals. Amer.Math.Monthly, 68, 645-647, 1961.

[23aa] J. Billingham and A.C. King, Uniform asymptotic expansions for the Barnes

double gamma function. Proc. R. Soc. Lond. A (1997) 454, 1817-1829.

[24] J. Blümlein, Algebraic Relations between Harmonic Sums and Associated Quantities. Comput.Phys.Commun. 159, 19-54, 2004.
Hep-ph/0311046 Abstract and Postscript and PDF]

[24aa] J. Blümlein, Analytic Continuation of Mellin Transforms up to two-loop Order. Comput.Phys.Commun. 133 (2000) 76-104.
hep-ph/0003100 [abs, ps, pdf, other]

[24a] H.P. Boas and E. Friedman, A simplification in certain contour integrals. Amer.Math.Monthly, 84, 467-468, 1977.

[24a] J. Bohman and C.-E. Fröberg, The Stieltjes Function-Definition and Properties. Math. of Computation, 51, 281-289, 1988.

[24b] E. Bombieri and J.C. Lagarias, Complements to Li's criterion for the Riemann hypothesis. J. Number Theory 77, 274-287 (1999).

[25] G. Boros and V.H. Moll, Irresistible Integrals: Symbolics, Analysis and Experiments in the Evaluation of Integrals. Cambridge University Press, 2004.

[25a] G. Boros, O. Espinosa and V.H. Moll, On some families of integrals solvable in terms of polygamma and negapolygamma functions.2002.
math.CA/0305131 [abs, ps, pdf, other]

[26] J.M. Borwein and P. Borwein, Pi and the AGM.Wiley-Interscience, New York, 1987.

[27] D. Borwein and J.M. Borwein, On an Intriguing Integral and Some Series Related to $\varsigma(4)$ . Proc. Amer. Math. Soc. 123, 1191-1198, 1995.
http://www.math.uwo.ca/~dborwein/cv/zeta4.pdf

[28] D. Borwein, J.M. Borwein and R. Girgensohn, Explicit Evaluations of Euler Sums. Proc. Edinburgh Math. Soc. (2) 38, 277-294, 1995.
Explicit evaluation of Euler sums - Borwein, Borwein, Girgensohn (1994)

[28a] J.M. Borwein and R. Girgensohn, Evaluation of triple Euler sums. Electron. J. Combin., 3:1-27, 1996.
Evaluation Of Triple Euler Sums - Jonathan Borwein (1995)

[29] J.M. Borwein, D.M. Bradley and R.E. Crandall, Computational Strategies for the Riemann Zeta Function. J. Comput. Appl. Math. 123, 247-296, 2000.
http://eprints.cecm.sfu.ca/archive/00000211/

[30] J.M. Borwein, D.M. Bradley, D.J. Broadhurst and P. Lisoněk, Special Values of Multiple Polylogarithms.
http://arxiv.org/abs/math/9910045

[30a] J.M. Borwein and D.M. Bradley, Thirty-two Goldbach Variations.

math.NT/0502034 [abs, ps, pdf, other] (to appear in International Journal of Number Theory), 2005.

[30b] J.M. Borwein, I.J. Zucker and J. Boersma. Evaluation of character Euler sums. http://eprints.cecm.sfu.ca/archive/00000255/ 2004.

[30c] M.T. Boudjelkha, A proof that extends Hurwitz formula into the critical strip. Applied Mathematics Letters, 14 (2001) 309-403.

[31] K.N. Boyadzhiev, Consecutive evaluation of Euler sums, Internat. J. Math. Sci., 29:9 (2002), 555-561

[32] K.N. Boyadzhiev Evaluation of Euler-Zagier sums, Internat. J. Math. Math. Sci., 27:7 (2001) 407-412

[32a] P. Bracken; C. Wenchang and D.C.L. Veliana, Summing Inverted Binomial Coefficients. Math.Mag., 77, 398-399, 2004.

[33] D.M. Bradley, Representations of Catalan's constant, (an unpublished catalogue of formulae for the alternating sum of the reciprocals of the odd positive squares), 1998.
http://germain.umemat.maine.edu/faculty/bradley/papers/pub.html

[33] D.M. Bradley, A class of series acceleration formulae for Catalan's constant. The Ramanujan Journal, Vol. 3, Issue 2, 159-173, 1999.
http://germain.umemat.maine.edu/faculty/bradley/papers/rj.pdf

[33b] M. Brede, Eine reihenentwicklung der vervollständigten und ergänzten Riemannschen zetafunktion und verwandtes
http://www.mathematik.uni-kassel.de/~koepf/Diplome/brede.pdf

[34] D. Bressoud, A Radical Approach to Real Analysis. The Mathematical Association of America, 1994.

[35] W.E. Briggs, S. Chowla, A.J. Kempner and W.E. Mientka, On some infinite series. Scripta Math, 21, 28-30, 1955.

[35a] W.E. Briggs and S. Chowla, The power series coefficients of $\varsigma(s)$. Amer. Math. Monthly, 62, 323-325, 1955.

[35b] W.E. Briggs, Some constants associated with the Riemann zeta-function. (1955-1956), Michigan Math. J. 3, 117-121.

[36] D.J. Broadhurst, Polylogarithmic ladders, hypergeometric series and the ten millionth digits of $\varsigma(3)$ and $\varsigma(5)$.1998.
math.CA/9803067 [abs, ps, pdf, other]

[36a] K.A. Broughan, Vanishing of the integral of the Hurwitz zeta function. Bull. Austral. Math. Soc. 65 (2002) 121-127.
Vanishing of the integral of the Hurwitz zeta function

[36b] T.J.I'a Bromwich, Introduction to the theory of infinite series. $2^{nd}$ edition Macmillan & Co Ltd, 1965.

[37] R.G. Buschman. Math. Mag.Vol.32, p107-108, 1958.

[38] P.L. Butzer, C. Markett and M. Schmidt, Stirling numbers, central factorial numbers and representations of the Riemann zeta function. Resultate Math.19, 257-274, 1991.

[38a] B. Candelpergher, M.A. Coppo and E. Delabaere, La Sommation de Ramanujan. L'Enseignement Mathématique, 43, 93-132, 1997.
[PS] La sommation de Ramanujan

[39] L. Carlitz, Eulerian Numbers and Polynomials. Math.Mag.32, 247-260, 1959.

[40] A.G. Caris, Amer.Math.Monthly, 21,336-340, 1914.

[41] H.S. Carslaw, Introduction to the theory of Fourier Series and Integrals. Third Ed. Dover Publications Inc, 1930.

[42] D. Castellanos, The Ubiquitous $\pi$. Math.Mag.61, 67-98, 1988.

[43] P. Cerone, M.A. Chaudhry, G. Korvin and A. Qadir, New Inequalities involving the Zeta Function. Journal of Inequalities in Pure and Applied Mathematics. Vol.5, No.2, Article 43, 2004.
http://jipam.vu.edu.au/images/130_03_JIPAM/130_03.pdf

[43a] Chao-Ping Chen and Feng Qi, The best bounds of the harmonic sequence. math.CA/0306233 [abs, ps, pdf, other],2003.

[43b] Hwang Chien-Lih, Relations between Euler's constant, Riemann's zeta function and Bernoulli numbers. Math. Gaz., 89, 57-59, 2005.

[43c] H. Chen, A power series and its applications.
Int. J. of Math. Educ. in Science and Technology, 37:3,362-368 (2005).

[43d] H. Chen and P. Khalili, On a class of logarithmic integrals.
Int. J. of Math. Educ. in Science and Technology, 37:1,119-125 (2006).

[43d] Y.J. Cho, M. Jung, J. Choi and H.M. Srivastava, Closed-form evaluations of definite integrals and associated infinite series involving the Riemann zeta function.Int. J. Comput. Math., 83,Nos. 5-6,461-472, 2006.

[44] Boo Rim Choe, An Elementary Proof of $\sum_{n=1}^{\infty} 1/n^2 = \pi^2/6$.
Amer.Math.Monthly, 94,662-663, 1987.

[45] J. Choi, H.M. Srivastava and V.S. Adamchik, Multiple Gamma and Related Functions. The Ramanujan Journal, 2003.

[45aa] J. Choi and H.M. Srivastava, Certain classes of series involving the Zeta function. J.Math.Anal.Appl., 231, 91-117,1999.

[45ab] J. Choi and H.M. Srivastava, Certain classes of series associated with the Zeta function and multiple gamma functions.
J. Comput. Appl. Math., 118, 87-109, 2000.

[45ac] J. Choi, Y.J. Cho and H.M. Srivastava, Series involving the Zeta function and multiple Gamma functions. Appl.Math.Comput.159 (2004)509-537.

[45aci] J. Choi and H.M. Srivastava, Explicit evaluation of Euler and related sums. The Ramanujan Journal, 10, 51-70, 2005.

[45ad] B.K. Choudhury, The Riemann zeta function and its derivatives.
Proc. Roy. Soc. Lond. A (1995) 450, 477-499.

[45ae] V.O. Choulakian; K.F. Anderson, Series of sine integrals.
Amer.Math.Monthly, 105, 474-475, 1998.

[45b] M.W. Coffey, On some log-cosine integrals related to $\varsigma(3), \varsigma(4)$ and $\varsigma(6)$.
J. Comput. Appl. Math., 159, 205-215, 2003.

[45c] M.W. Coffey, New results on the Stieltjes constants: Asymptotic and exact evaluation. J. Math. Anal. Appl., 317 (2006)603-612.
math-ph/0506061 [abs, ps, pdf, other]

[45d] M.W. Coffey, On one-dimensional digamma and polygamma series related to the evaluation of Feynman diagrams.
J. Comput.Appl. Math, 183, 84-100, 2005.
math-ph/0505051 [abs, ps, pdf, other]

[45e] M.W. Coffey, New summation relations for the Stieltjes constants
Proc. R. Soc. A ,462, 2563-2573, 2006.

[45f] M.W. Coffey, Toward verification of the Riemann Hypothesis: Application of the Li criterion. math-ph/0505052 [abs, ps, pdf, other],2005. Math. Phys. Analysis and Geometry, 8, 211-255, 2005.

[45g] M.W. Coffey, Polygamma theory, the Li/Keiper constants and validity of the Riemann Hypothesis. math-ph/0507042 [abs, ps, pdf, other],2005.

[45h] M.W. Coffey, A set of identities for alternating binomial sums arising in computing applications. math-ph/0608049 [abs, ps, pdf, other],2006.

[45i] M.W. Coffey, The Stieltjes constants, their relation to the $\eta_j$ coefficients, and representation of the Hurwitz zeta function.
arXiv:math-ph/0612086 [ps, pdf, other], 2007.

[45j] M.W. Coffey, Series of zeta values, the Stieltjes constants, and a sum $S_\gamma(n)$.
arXiv:0706.0345 [ps, pdf, other], 2007.

[46] S.W. Coffman, B. Shawer, H. Kappus, B.C. Berndt, An Infinite Series with Harmonic Numbers. Math. Mag., 60, 118-119, 1987.

[46a] G. Cognola, E. Elizalde and K. Kirsten, Casimir Energies for Spherically Symmetric Cavities. J.Phys. A34 (2001) 7311-7327
hep-th/9906228 [abs, ps, pdf, other]

[46aa] C.B. Collins, The role of Bell polynomials in integration.
J. Comput. Appl. Math. 131 (2001) 195-211.

[46b] M.A. Coppo, Nouvelles expressions des constantes de Stieltjes.
Expositiones Mathematicae 17, No. 4, 349-358 (1999).

[47] F. Lee Cook, A Simple Explicit Formula for the Bernoulli Numbers.
The Two-Year College Mathematics Journal 13, 273-274, 1982.

[48] R. E. Crandall and J. P. Buhler, On the evaluation of Euler sums.
Experimental Math., 3 (1994), no. 4, 275–285
Full text (Postscript)

[48a] D. Cvijović and J. Klinowski, Closed-form summation of some trigonometric series. Math. Comput., 64, 205-210, 1995.

[49] D. Cvijović and J. Klinowski, New rapidly convergent series representations for $\varsigma(2n+1)$. Proc. Amer. Math. Soc. 125, 1263-1271, 1997.
http://www.ams.org/proc/1997-125-05/S0002-9939-97-03795-7/home.html

[49a] D. Cvijović, The Haruki-Rassias and related integral representations of the Bernoulli and Euler polynomials. J. Math. Anal. Appl. (to appear) 2007.

[49b] D. Cvijović, New integral representations of the polylogarithm function.
Proc. R. Soc. A (2007), 463, 897-905.

[50] D. Cvijović and J. Klinowski, Integral Representations of the Riemann Zeta Function for Odd-Integer Arguments. J.Comput.Appl.Math.142, 435-439, 2002.

[50a] O.T. Dasbach, Torus Knot complements: A natural series for the natural logarithm. math.GT/0611027 [abs, ps, pdf, other].

[51] M. Dalai, Recurrence Relations for the Lerch $\Phi$ Function and Applications.
math.NT/0411087 [abs, ps, pdf, other] 2004.

[51a] A. Das and G.V. Dunne, Large-order Perturbation Theory and de Sitter/Anti de Sitter Effective Actions, Phys.Rev. D74 (2006) 044029
hep-th/0607168 [abs, ps, pdf, other]

[51b] A.I. Davydychev and M. Yu. Kalmykov, New results for the epsilon-expansion of certain one-, two- and three-loop Feynman diagrams Nucl.Phys. B605 (2001) 266-318
arXiv:hep-th/0012189 [ps, pdf, other]

[52] R. Dedekind, Üeber ein Eulerische Integral. J. Reine Ang. Math., Vol.45, 1853.
http://www.digizeitschriften.de/no_cache/home/open-access/nach-zeitschriftentiteln/

[53] J. B. Dence, Development of Euler Numbers. Missouri Journal of Mathematical Sciences, 9, 1-8, 1997.   A Development of Euler Numbers

[53a] A. Devoto and D.W. Duke, Table of integrals and formulae for Feynman diagram calculations. Florida State University, FSU-HEP-831003, 1983.
http://www.csit.fsu.edu/~dduke/integrals.htm

[54] K. Dilcher, Some $q$-Series Identities Related to Divisor Functions. Discrete Math. 145, 83-93, 1995.

[54a] K. Dilcher, Generalized Euler constants for arithmetical progressions Math. of Comp.,Vol.59,No.199,259-282,1992.

[55] P.J. de Doelder, On some series containing $\psi(x)-\psi(y)$ and $(\psi(x)-\psi(y))^2$ for certain values of $x$ and $y$. J. Comput. Appl. Math. 37, 125-141, 1991.

[55a] D. Bierens de Haan, Exposé de la Théorie, Propriétés, des formules de transformation, et des méthodes d'évaluation des intégrales définies, C.G. Van der Post, Amsterdam, 1862. Available on the internet at the University of Michigan Historical Mathematics Collection.
http://www.hti.umich.edu/u/umhistmath/

[55b] B. Doyon, J. Lepowsky and A. Milas, Twisted vertex operators and Bernoulli polynomials. math.QA/0311151 [abs, ps, pdf, other], 2005.

[56] W. Dunham, Euler, The Master of Us All. Mathematical Association of America, Washington, DC, 1999.

[56a] J. Duoandikoetxea, A sequence of polynomials related to the evaluation of the Riemann zeta function. Math. Mag, 80, No. 1, 38-45, 2007.

[57] H.M. Edwards, Riemann's Zeta Function. Academic Press, New York and London, 1974.

[58] H.M. Edwards, Fermat's Last Theorem: A Genetic Introduction to Algebraic Number Theory. Springer-Verlag, 1977.

[58a] C.J. Efthimiou, Finding exact values for infinite series.
Math. Mag. 72, 45-51, 1999. arXiv:math-ph/9804010 [ps, pdf, other]

[58aa] C.J. Efthimiou, Trigonometric Series via Laplace Transforms.

      arXiv:0707.3590 [ps, pdf, other] 2007.

[58b] A. Erdélyi, W. Magnus, F. Oberhettinger and F.G. Tricomi.
      Higher Transcendental Functions, Volume I, McGraw-Hill Book Company,
      Inc, 1953.

[58c] E. Elizalde, Derivative of the generalised Riemann zeta function $\varsigma(z,q)$ at
      $z = -1$. J. Phys. A Math. Gen. (1985) 1637-1640

[58ci] E. Elizalde and A. Romeo, An integral involving the generalized zeta function.
      Internat. J. Maths. & Maths. Sci. Vol.13, No.3, (1990) 453-460.
      http://www.hindawi.com/GetArticle.aspx?doi=10.1155/S0161171290000679&e=CTA

[58d] A. Erdélyi et al, Tables of Integral Transforms. McGraw-Hill Book Company,
      New York, 1964.

[59] O. Espinosa and V.H. Moll, On some integrals involving the Hurwitz zeta
      function: Part I. The Ramanujan Journal, 6,150-188, 2002.
      http://arxiv.org/abs/math.CA/0012078

[60] O. Espinosa and V. H. Moll. On some integrals involving the Hurwitz zeta
      function: Part 2. The Ramanujan Journal, 6,449-468, 2002. ps    pdf

[61] L. Euler, Demonstratio insignis theorematis numerici circa unicias potestatum
      binomialium .Nova Acta Acad. Sci. Petropol.15 (1799/1802), 33-43; reprinted in
      Opera Omnia, Ser. I,Vol. 16(2), B.G. Teubner,Leipzig,1935, pp.104-116.

[62] L. Euler, Institutiones Calculi Differentialis, Petrograd, 1755, pp. 487-491.

[63] The Euler Archive. Website http://www.eulerarchive.org/

[64] Russell Euler, Evaluating a Family of Integrals. Missouri Journal of
      Mathematical   Sciences 9, 1-4, 1997. Evaluating a Family of Integrals

 [65] J.A. Ewell, A New Series Representation for $\varsigma(3)$.
      Math.Monthly, 97, 219-220, 1990.

[65a] O.Furdui, College Math. Journal, 38, No.1, 61, 2007

[66] E. Fischer, Intermediate Real Analysis. Springer-Verlag, New York, 1983.

[67] P. Flajolet, X. Gourdon and P. Dumas, Mellin Transforms and Asymptotics:
      Harmonic sums. *Theoretical Computer Science,* vol. 144 (1-2), pp. 3-58, 1995.

 [68] P. Flajolet and R. Sedgewick, Mellin Transforms and Asymptotics: Finite
Differences and Rice's Integrals.Theor.Comput.Sci.144, 101-124, 1995.
Mellin Transforms and Asymptotics : Finite Differences and Rice's Integrals   (117kb),

[69] P. Flajolet and B. Salvy, Euler Sums and Contour Integral Representations
(115kb),   (INRIA, RR2917), June 1996. The final version appears in *Journal of*

*Experimental Mathematics*, volume **7**(1), 1998, pp. 15-35, where it is available electronically, by courtesy of the publisher.

[69aa] J.Fleisher, A.V. Kotikov and O.L. Veretin, Analytic two-loop results for self energy- and vertex-type diagrams with one non-zero mass**.** hep-ph/9808242 [abs, ps, pdf, other] Nucl.Phys. B547 (1999) 343-374.

[69a] P. Freitas, Integrals of polylogarithmic functions, recurrence relations and associated Euler sums. Math.CA/0406401 [abs, ps, pdf, other] 2004.

[69aa] P.G.O. Freund and A.V. Zabrodin, A Hierarchical Array of Integrable Models. J.Math.Phys. 34 (1993) 5832-5842. hep-th/9208033 [abs, ps, pdf, other]

[69b] R. Gastmans and W. Troost, On the evaluation of polylogarithmic integrals. Simon Stevin, 55, 205-219, 1981.

[69c] C. Georghiou and A.N. Philippou, Harmonic sums and the Zeta function. Fibonacci Quart., 21, 29-36, 1983.

[70] Badih Ghusayni Some Representations of zeta(3).Missouri Journal of Mathematical Sciences 10, 169-175, 1998.

[70] Badih Ghusayni. Euler-type formula using Maple. Palma Research Journal, 7, 175-180, 2001.
http://www.ndu.edu.lb/academics/palma/20010701/vol7is1a17.doc

[70aa] J. Ginsburg, Note on Stirling's Numbers. Amer.Math.Monthly, 35, 77-80, 1928.

[70ab] M.L. Glasser, Evaluation of some integrals involving the $\psi$ - function. Math. of Comp., Vol.20, No.94, 332-333, 1966.

[70a] M.A. Glicksman, Euler's Constant. Amer.Math.Monthly, 50, 575, 1943.

[71] R.W. Gosper, $\int_{n/4}^{m/6} \log \Gamma(z) dz$. In Special Functions, q-series and related topics. Amer.Math.Soc.Vol. 14.

[71a] T.H. Gronwall, The gamma function in integral calculus. Annals of Math., 20, 35-124, 1918.

[72] H.W. Gould, Some Relations involving the Finite Harmonic Series. Math.Mag., 34,317-321, 1961.

[73] H.W. Gould, Combinatorial Identities.Rev.Ed.University of West Virginia, U.S.A., 1972.

[73a] H.W. Gould, Explicit formulas of Bernoulli Numbers.

        Amer.Math.Monthly, 79, 44-51, 1972.

[73b] H.W. Gould, Stirling Number Representation Problems.
     Proc. Amer. Math. Soc., 11, 447-451, 1960.

[74] I.S. Gradshteyn and I.M. Ryzhik, Tables of Integrals, Series and Products.
    Sixth Ed., Academic Press, 2000.
    Errata for Sixth Edition http://www.mathtable.com/errata/gr6_errata.pdf

[75] R.L. Graham, D.E. Knuth and O. Patashnik, Concrete Mathematics. Second Ed.
    Addison-Wesley Publishing Company, Reading, Massachusetts, 1994.

[75a] R. Greenberg, D.C.B. Marsh and A.E. Danese, A Zeta-function Summation.
     Amer.Math.Monthly, 74, 80-81, 1967.

[75aa] J. Guillera and J. Sondow, Double integrals and infinite products for some
     classical constants via analytic continuations of Lerch's transcendent.2005.
     math.NT/0506319 [abs, ps, pdf, other]

[76] G.H. Hardy et al., Collected Papers of Srinivasa Ramanujan.Cambridge
    University Press, Cambridge, 1927; reprinted by Chelsea, 1962; reprinted by
    American Mathematical Society, 2000.
    http://www.imsc.res.in/~rao/ramanujan/CamUnivCpapers/collectedright1.htm

[76aa] G.H. Hardy, Divergent Series. Chelsea Publishing Company, New York, 1991.

[76a] F. Haring; G.T.Nelson; G.Bach. $\varsigma(n)$, $\psi^{(n)}$ and an Infinite Series.
     Amer.Math.Monthly, 81, 180-181, 1974.

[76b] M. Hashimoto, S. Kanemitsu, T. Tanigawa, M. Yoshimoto and W.-P.Zhang,
     On some slowly convergent series involving the Hurwitz zeta function.2002.
     http://www.imsc.res.in/~icsf2002/papers/tanigawa.pdf

[76c] F.E. Harris, Spherical Bessel expansions of sine, cosine and exponential
     integrals. Appl. Numerical Math, 34 (2000) 95-98.

[77] H. Hasse, Ein Summierungsverfahren für Die Riemannsche $\varsigma$ - Reithe.
    Math.Z.32, 458-464, 1930.
    http://dz-srv1.sub.uni-goettingen.de/sub/digbib/loader?ht=VIEW&did=D23956&p=462

[78] J. Havil, Gamma: Exploring Euler's Constant. Princeton University Press,
    Princeton, NJ, 2003.

[79] Hernández, V. Solution IV of Problem 10490: A Reciprocal Summation
     Identity. *Amer. Math. Monthly* 106, 589-590, 1999.

[79a] M.Hirschhorn, A new formula for Pi. Austral.Math.Soc.Gazette, 25, 82-83,
     1998.

[80] M.E. Hoffman, Quasi-symmetric functions and mod *p* multiple harmonic sums,

http://arxiv.org/PS_cache/math/pdf/0401/0401319.pdf ,2004.

[80a] K. Ireland and M. Rosen, A Classical Introduction to Modern Number Theory. Second edition, Springer-Verlag New York Inc, 1990.

[81] A. Ivić, The Riemann Zeta- Function: Theory and Applications. Dover Publications Inc, 2003.

[82] W.P. Johnson, The Curious History of Faà di Bruno's Formula. Amer.Math.Monthly 109,217-234, 2002.

[82aa] M. Kamela and C.P. Burgess, Massive-Scalar Effective Actions on Anti-de Sitter Spacetime.
Can.J.Phys. 77 (1999) 85-99. hep-th/9808107 [abs, ps, pdf, other]

[82a] M. Kaneko, The Akiyama-Tanigawa algorithm for Bernoulli numbers. Journal of Integer Sequences, Vol. 3, Article 00.2.9, 2000.
http://www.cs.uwaterloo.ca/journals/JIS/VOL3/KANEKO/AT-kaneko.pdf

[82b] S. Kanemitsu, M. Katsurada and M. Yoshimoto, On the Hurwitz-Lerch zeta function.   Aequationes Math. 59 (2000) 1-19.

[83] R. Kanigel, The Man Who Knew Infinity: A Life of the Genius Ramanujan.Charles Scribners' Sons, New York, 1991.

[83a] J.B. Keiper, power series expansions of Riemann's $\xi$ function.
Math. Comp.58, 765-773 (1992).

[84] G. Kimble, Euler's Other Proof. Math. Mag., 60,282, 1977.

[84a] K. Kimoto and M. Wakayama, Apéry-like numbers arising from special values of spectral zeta functions for non-commutative harmonic oscillators. 2006. math.NT/0603700 [abs, ps, pdf, other]

[85] A.N. Kirillov, Dilogarithm Identities. *Progress of Theor. Phys. Suppl.* 118, 61-142, 1995. http://arxiv.org/abs/hep-th/9408113

[86] P. Kirschenhofer, A Note on Alternating Sums. The Electronic Journal of Combinatorics 3 (2), #R7, 1996. R7: Peter Kirschenhofer

[86a] M.S. Klamkin; H.F. Sandham; M.R. Spiegel.
 Amer. Math. Monthly, 62, 588-590, 1955.

[87] M. Kline, Mathematical Thought from Ancient to Modern Times.Vol.2, Oxford University Press, 1972.

[88] P. Knopf, The Operator $\left(x\dfrac{d}{dx}\right)^n$ and its Application to Series.

Math.Mag.76, 364-371, 2003.

[89] K.Knopp, Theory of Functions. Dover, New York, 1996.

[90] K. Knopp, Theory and Application of Infinite Series. Second English Ed.Dover Publications Inc, New York, 1990.

[90a] D.E. Knuth, Euler's constant to 1271 places. Math. of Computation, 16, 275-281, 1962.

[90b] D.E. Knuth, The Art of Computer Programming, Vol. 1, Addison Wesley, 1977.

[91] M. Koecher: Lettters, Math.Intell.2, 62-64,1980.

[91a] K.S. Kölbig, Nielsen's generalised polylogarithms. SIAM J. Math.Anal.Vol.17, No.5, 1232-1258, 1986.

[91aa] K.S. Kölbig, Some infinite integrals with powers of logarithms and the complete Bell polynomials. J. Comput.Appl. Math, 69 (1996), 39-47.

[91b] Kondratieva and Sadov, Markov's Transformation of series and the WZ method math.CA/0405592 [abs, ps, pdf, other], 2004.

[92] R.A. Kortram, Simple proofs for

$$\sum_{k=1}^{\infty}\frac{1}{k^2} = \frac{\pi^2}{6} \text{ and } \sin x = x\prod_{k=1}^{\infty}\left(1-\frac{x^2}{k^2\pi^2}\right).$$ Math. Mag.69, 122-125, 1996.

[93] S. Koyama and N. Kurokawa, Certain Series Related to the Triple Sine Function. http://www.math.keio.ac.jp/local/koyama/papers/English/series.pdf

[93a] S. Koyama and N. Kurokawa, Kummer's formula for the multiple gamma functions. Presented at the conference on Zetas and Trace Formulas in Okinawa, November, 2002.
www.math.titech.ac.jp/~tosho/Preprints/pdf/128.pdf

[94] E.E. Kummer, Beitrag zur Theorie der Function $\Gamma(x) = \int_0^{\infty} e^{-v}v^{x-1}dv$.

J. Reine Angew.Math., 35, 1-4, 1847.
http://www.digizeitschriften.de/index.php?id=132&L=2

[94a] J. Landen, A new method of computing sums of certain series. Phil.Trans.R.Soc.Lond., 51, 553-565, 1760.

[94aa] H. Langman; J.F. Locke; C.W. Trigg. Amer. Math. Monthly, 43, 196-197, 1936.

[94b] J. Landen, Mathematical Memoirs, 1, 1780.

[95] P.J. Larcombe, E.J. Fennessey and W.A. Koepf, Integral proofs of Two Alternating Sign Binomial Coefficient Identities.

http://citeseer.ist.psu.edu/598454.html

[96] Kee-Wai Lau, Some Definite Integrals and Infinite Series. Amer.Math.Monthly 99, 267-271, 1992.

[97] D.H. Lehmer, Interesting series involving the central binomial coefficient. Amer.Math.Monthly 92,449-457, 1985.

[98] M.E. Levenson, J.F. Locke and H. Tate, Amer.Math.Monthly, 45, 56-58, 1938.

[99] M.E. Levenson, A recursion formula for $\int_0^\infty e^{-t} \log^{n+1} t\, dt$. Amer.Math.Monthly, 65, 695-696, 1958.

[100] L. Lewin, Polylogarithms and Associated Functions. Elsevier (North-Holland), New York, London and Amsterdam, 1981.

[101] L. Lewin (Editor), Structural Properties of Polylogarithms. (Mathematical Surveys and Monographs, Vol.37), American Mathematical Society, Providence, Rhode Island, 1991.

[101i] X.-J. Li, The positivity of a sequence of numbers and the Riemann Hypothesis. J. Number Th., 65, 325-333, 1997

[101aa] G.J. Lodge; R. Breusch. Riemann Zeta Function. Amer.Math.Monthly, 71, 446-447, 1964.

[101ab] J.L. Lopez, Several series containing gamma and polygamma functions. J. Comput. Appl. Math, 90, (1998), 15-23.

[101a] M. Lutzky, Evaluation of some integrals by contour integration. Amer.Math.Monthly, 77, 1080-1082, 1970.

[101aa] T. Mansour, Gamma function, Beta function and combinatorial identities. math.CO/0104026 [abs, ps, pdf, other], 2001.

[101b] L.C. Maximon, The dilogarithm function for complex argument. Proceedings: Mathematical, Physical and Engineering Sciences, The Royal Society, 459 (2339), 2807-2819, 2003.
http://www.pubs.royalsoc.ac.uk/proc_phys_homepage.shtml
http://www.math.uio.no/~didier/dilog.pdf

[101c] K. Maślanka, Effective method of computing Li's coefficients and their properties. math.NT/0402168 [abs, ps, pdf, other]

[101d] K. Maślanka, An explicit formula relating Stieltjes constants and Li's numbers. math.NT/0406312 [abs, ps, pdf, other]

[102] Z.R. Melzak, Companion to Concrete Mathematics.Wiley-Interscience, New

York, 1973.

[102a] M. Milgram, On Some Sums of Digamma and Polygamma functions.
math.CA/0406338 [abs, pdf]

[103] J. Miller and V.S. Adamchik, Derivatives of the Hurwitz Zeta Function for Rational Arguments. *J. Comp. and Applied Math.,* 100(1998), 201--206.
Derivatives of the Hurwitz Zeta Function for Rational Arguments

[103a] C. Moen, Infinite series with binomial coefficients.
Math. Mag., 64, 53-55, 1991.

[103ab] C.Musès, Some new considerations on the Bernoulli numbers, the factorial function and Riemann's zeta function.
Appl. Math. and Comput.113 (2000) 1-21.

[103ac] H. Muzaffar, Some interesting series arising from the power series expansion of $\left(\sin^{-1} x\right)^{q}$. Int. J. of Math. and Math. Sci. 2005:14(2005) 2329-2336.

[103ai] T.S. Nanjundiah, Van der Pol's expressions for the gamma function.
Proc.Amer.Math.Soc., 9, 305-307, 1958.

[104] C. Nash and D. O'Connor, Determinants of Laplacians, the Ray-Singer torsion on lens spaces and the Riemann zeta function.J.Math.Phys.36, 1462-1505,1995.
http://arxiv.org/pdf/hep-th/9212022

[104a] N. Nielsen,Theorie des Integrallogarithmus und verwanter tranzendenten 1906.
http://www.math.uni-bielefeld.de/~rehmann/DML/dml_links_author_H.html

[104b] N. Nielsen, Die Gammafunktion. Chelsea Publishing Company, Bronx and New York, 1965.

[105] N.E. Nörlund, Vorlesungen über Differenzenrechnung.Chelsea, 1954.
http://dz-srv1.sub.uni-goettingen.de/cache/browse/AuthorMathematicaMonograph,WorkContainedN1.html

[105(i)] N.E. Nörlund, Leçons sur les séries d'interpolation.
Paris, Gauthier-Villars, 1926.

[105(ii)] O.M. Ogreid and P. Osland, Some infinite series related to Feynman diagrams. math-ph/0010026 [abs, ps, pdf, other]

[105(iii)] O.M. Ogreid and P. Osland, More Series related to the Euler Series.
hep-th/9904206 [abs, ps, pdf, other]

[105(iv)] D. Oprisa and S. Stieberger, Six Gluon Open Superstring Disk Amplitude, Multiple Hypergeometric Series and Euler-Zagier Sums. 2005.
hep-th/0509042 [abs, ps, pdf, other]

[105(v)] T.J. Osler, An introduction to the zeta function. Preprint 2004.

[An Introduction to the Zeta Function](#)

[105a] A. Panholzer and H. Prodinger, Computer-free evaluation of a double infinite
sum via Euler sums, 2005.
http://math.sun.ac.za/~prodinger/abstract/abs_218.htm

[105aa] R. Pemantle and C. Schneider, When is 0.999...equal to 1?
math.CO/0511574 [abs, ps, pdf, other]

[106] R. Penrose, The Road to Reality: A Complete Guide to the Laws of the
Universe. Jonathan Cape, London, 2004.

[107] S. Plouffe, Plouffe's Inverter. http://pi.lacim.uqam.ca/eng/

[108] R.I. Porter, Further Mathematics. Bell and Sons Ltd, London, 1963.

[108a] G. Póyla and G. Szegö, Problems and Theorems in Analysis, Vol.I
Springer-Verlag, New York 1972.

[109] H. Prodinger, A q-Analogue of a Formula of Hernandez Obtained by Inverting
a Result of Dilcher. *Austral. J. Combin.* 21, 271-274, 2000.

[109i] A.P.Prudnikov, Yu.A Brychkov and O.I Marichev, *Integrals and series,
volume I: elementary functions*. New York, NY: Gordon and Breach,1986.

[109a] R.K. Raina and R.K. Ladda, A new family of functional series relations
involving digamma functions.
Ann. Math. Blaise Pascal, Vol. 3, No. 2, 1996, 189-198.
http://www.numdam.org/item?id=AMBP_1996__3_2_189_0

[110] Srinivasa Ramanujan, Notebooks of Srinivasa Ramanujan, Vol.1, Tata Institute
of Fundamental Research, Bombay, 1957.

[110aa] H. Rademacher, Topics in Analytic Number Theory.
Springer-Verlag, 1973.

[110a] S.K. Lakshamana Rao, On the Sequence for Euler's Constant.
Amer.Math.Monthly, 63, 572-573, 1956.

[110b] K. Roach, Generalization of Adamchik's formulas. 2005.
http://www.planetquantum.com/Notes/Adamchik97Review.pdf

[111] K.Srinivasa Rao, Ramanujan's Notebooks and other material.
http://www.imsc.res.in/~rao/ramanujan/index.html

[111aa] G. E. Raynor ,On Serret's integral formula..
Bull. Amer. Math. Soc. Volume 45, Number 12, Part 1 (1939), 911-917**.**
On Serret's integral formula

[111a] H. Ruben, A Note on the Trigamma Function.
Amer.Math.Monthly, 83, 622-623, 1976.

[112] G.F.B. Riemann, Üeber die Anzahl der Primzahlen unter einer gegebenen Grösse.Monatsber.Königl.Preuss.Akad.Wiss.,Berlin,671-680,1859.
http://www.maths.tcd.ie/pub/HistMath/People/Riemann/Zeta/

[113] T. Rivoal, La function Zeta de Riemann prend une infinité de valeurs irrationnelles aux entiers impairs. Comptes Rendus Acad.Sci.Paris 331,267-270, 2000.
http://arxiv.org/abs/math/0008051

[114] T. Rivoal, Irrationalité d'au moins un des neuf nombres ς(5),ς(7),…,ς(21). Acta Arith. 103:2 (2002), 157-167 (E-print math.NT/0104221).

[114aa] T. Rivoal, Polynômes de type Legendre et approximations de la constant d'Euler. Note (2005).  DVI, PS, PDF
http://www-fourier.ujf-grenoble.fr/

[114a] B. Ross, Serendipity in mathematics or how one is led to discover that
$$\sum_{n=1}^{\infty} \frac{1.3.5....(2n-1)}{n2^{2n} n!} = \frac{1}{2} + \frac{3}{16} + \frac{15}{144} + ... = \log 4$$
Amer.Math.Monthly, 90, 562-566, 1983.

[115] W. Rudin, Principles of Mathematical Analysis. Third Ed.McGraw-Hill Book Company, 1976.

[116] D.C. Russell, Another Eulerian-Type Proof. Math. Mag., 64, 349, 1991.

[116aa] G. Rutledge and R.D. Douglass, Evaluation of $\int_0^1 \frac{\log u}{u} \log^2(1+u)\, du$ and related definite integrals. Amer.Math.Monthly, 41, 29-36, 1934.

[116ab] G. Rutledge and R.D. Douglass, Tables of Definite Integrals.
Amer.Math.Monthly, 45, 525, 1938.

[116a] G. Rzadkowski, A Short Proof of the Explicit Formula for Bernoulli Numbers.
Amer.Math.Monthly, 111, 432-434, 2004.

[117] H.F. Sandham, A Well-Known Integral. Amer.Math.Monthly, 53, 587, 1946.

[118] H.F. Sandham, Advanced Problems 4353, Amer.Math.Monthly, 56, 414, 1949.

[118aa] H.F. Sandham; E. Trost, Amer.Math.Monthly, 58, 705-706, 1951.

[118a] Z.Sasvári, An Elementary Proof of Binet's Formula for the Gamma Function.
Amer.Math.Monthly, 106, 156-158, 1999.

[119] P. Sebah and X. Gourdon, Introduction to the Gamma Function.

[PDF] Introduction to the Gamma Function

[119a] P. Sebah and X. Gourdon, The Riemann zeta function $\varsigma(s)$ :Generalities
http://numbers.computation.free.fr/Constants/Miscellaneous/zetageneralities.pdf

[119b] J. Ser, Sur une expression de la fonction $\varsigma(s)$ de Riemann.
Comptes Rendus, 182, 1075-1077, 1926.
http://gallica.bnf.fr/Catalogue/noticesInd/FRBNF34348108.htm#listeUC

[120] L.-C. Shen, Remarks on some integrals and series involving the Stirling numbers and $\varsigma(n)$. Trans. Amer. Math. Soc. 347, 1391-1399, 1995.

[120a] R. Sitaramachandrarao, A formula of S.Ramanujan.
J.Number Theory 25, 1-19, 1987.

[120ai] W.D. Smith, A "good" problem equivalent to the Riemann hypothesis.

http://math.temple.edu/~wds/homepage/riemann2.pdf

[120aa] A. Snowden, Collection of Mathematical Articles. 2003.
http://www.math.princeton.edu/~asnowden/math-cont/dorfman.pdf

[121] J. Sondow, Analytic Continuation of Riemann's Zeta Function and Values at Negative Integers via Euler's Transformation of Series.Proc.Amer.Math.Soc. 120,421-424, 1994.
http://home.earthlink.net/~jsondow/id5.html

[122] J. Sondow, Zeros of the Alternating Zeta Function on the Line $\text{Re}(s)=1$.
Amer.Math.Monthly, 110, 435-437, 2003.
math.NT/0209393 [abs, ps, pdf]

[123] J. Sondow, An Infinite Product for $e^\gamma$ via Hypergeometric Formulas for Euler's Constant $\gamma$. 2003(preprint) http://arxiv.org/abs/math.CA/0306008

[123a] J. Sondow, A faster product for $\pi$ and a new integral for $\log\frac{\pi}{2}$.
Math.NT/0401406 [abs, pdf] Amer. Math. Monthly 112 (2005) 729-734.

[123aa] J. Sondow, Double Integrals for Euler's Constant and ln(4/Pi) and an Analog of Hadjicostas's Formula. Amer.Math.Monthly, 112, 61-65, 2005.
math.CA/0211148 [abs, pdf]

[123ab] J. Sondow and P. Hadjicostas, The Generalized-Euler-Constant Function $\gamma(z)$ and a Generalization of Somos's Quadratic Recurrence Constant.
math.CA/0610499 [abs, ps, pdf, other], 2006.
J. Math. Anal. Appl. 332 (2007) 292-314.

[123ac] A Speiser, Geometrisches zur Riemannschen zetafunktion.
Math. Ann. 110 (1934).

[123b] W. Spence, An essay on the theory of various orders of logarithmic transcendents. 1809.

[123bi] J. Spieß, Some identities involving harmonic numbers. Math. of Computation, 55, No.192, 839-863, 1990.

[123c] W.G. Spohn; A.S. Adikesavan; H.W. Gould. Amer.Math.Monthly, 75, 204-205,1968.

[124] E.L. Stark, The Series $\sum_{k=1}^{\infty} k^{-s}, s = 2,3,4,...,$ Once More. Math. Mag., 47,197-202, 1974.

[125] H.M. Srivastava, Some Families of Rapidly Convergent Series Representations for the Zeta Function. Taiwanese Journal of Mathematics, Vol.4, No.4, 569-598, 2000.
http://www.math.nthu.edu.tw/~tjm/abstract/0012/tjm0012_3.pdf

[125a] H.M. Srivastava and H. Tsumura, A certain class of rapidly convergent series representations for $\varsigma(2n+1)$ .J. Comput. Appl. Math., 118, 323-325, 2000.

[125aa] H.M. Srivastava, M.L. Glasser and V.S. Adamchik. Some definite integrals associated with the Riemann zeta function.
Z. Anal.Anwendungen, 129, 77-84, 2000.

[126] H.M. Srivastava and J. Choi, Series Associated with the Zeta and Related Functions. Kluwer Academic Publishers, Dordrecht, the Netherlands, 2001.

[126a] G. Stephenson, Mathematical Methods for Science Students.7[th] Ed. Longman Group Limited, London, 1970.

[127] P.M. Stevenson, Phys. Rev. D 23, 2916, 1981.

[127a] The Mactutor History of Mathematics archive.
http://www-history.mcs.st-andrews.ac.uk/Mathematicians/Faa_di_Bruno.html

[128] E.C. Titchmarsh, The Theory of Functions.2[nd] Ed., Oxford University Press, 1932.

[129] E.C. Titchmarsh, The Zeta-Function of Riemann. Oxford University (Clarendon) Press, Oxford, London and New York, 1951; Second Ed. (Revised by D.R. Heath- Brown), 1986.

[130] G.P. Tolstov, Fourier Series. (Translated from the Russian by R.A. Silverman) Dover Publications Inc, New York, 1976.

[130a] D.B. Tyler; P.R. Chernoff; R.B. Nelsen. An old sum reappears. Amer.Math.Monthly, 94, 466-468, 1987.

[131a] A. van der Poorten, Some wonderful formulae…footnotes to Apéry's proof of the irrationality of ς(3). Séminaire Delange-Pisot-Poitou (Théorie des Nombres) 29,1-7, 1978/1979.
http://www.ega-math.narod.ru/Apery1.htm

[131b] A. van der Poorten, Some wonderful formulae…an introduction to Polylogarithms. Queen's Papers in Pure Appl.Math.54, 269-286, 1979.
http://www.ega-math.narod.ru/Apery2.htm#ref10txt

[132] A. van der Poorten, A proof that Euler missed… Apéry's proof of the irrationality of ς(3). Math. Intelligencer 1, 195-203, 1979.

[133] A. van der Poorten, Notes on Fermat's Last Theorem. John Wiley & Sons Inc., New York, 1996.

[133] J.A.M. Vermaseren, Harmonic sums, Mellin transforms and Integrals. Int.J.Mod.Phys. A14 (1999) 2037-2076
http://arXiv.org/abs/hep-ph/9806280

[133a] M.B. Villarino, Ramanujan's approximation to the $n$ th partial sum of the harmonic series. Math.CA/0402354 [abs, ps, pdf, other]

[133b] A. Voros, Special functions, spectral functions and the Selberg zeta function. Comm. Math. Phys.110, 439-465, 1987.

[134] E.W. Weisstein, Dilcher's Formula. From Mathworld-A Wolfram Web Resource.
http://mathworld.wolfram.com/DilchersFormula.html

[135] E.T. Whittaker and G.N. Watson, A Course of Modern Analysis: An Introduction to the General Theory of Infinite Processes and of Analytic Functions; With an Account of the Principal Transcendental Functions. Fourth Ed., Cambridge University Press, Cambridge, London and New York, 1963.

[136] B. Wiener and J. Wiener, DeMoivre's Formula to the Rescue. Missouri Journal of Mathematical Sciences, 13, 1-9, 2001.

[137] J.Wiener, An Analytical Approach to a Trigonometric Integral . Missouri Journal of Mathematical Sciences 2, 75-77, 1990.

[138] J. Wiener Integration of Rational Functions by the Substitution $x = u^{-1}$ Missouri Journal of Mathematical Sciences.

[138a] J. Wiener, Differentiation with respect to a parameter. The College Mathematics Journal, 32, 180-184, 2001.

[138ai] J. Wiener, D.P. Skow and W. Watkins, Integrating powers of trigonometric functions . Missouri Journal of Mathematical Sciences, 3(1992), 55-61.
[PS] Integrating Powers of Trigonometric Functions

[138aii] J. Wiener, Integrals of $\cos^{2n} x$ and $\sin^{2n} x$.
The College Mathematics Journal, 31, 60-61, 2000.

[138b] H.S. Wilf, The asymptotic behaviour of the Stirling numbers of the first kind.
Journal of Combinatorial Theory Series A, 64, 344-349, 1993.
http://www.mathnet.or.kr/papers/Pennsy/Wilf/stirling.pdf

[139] S. Wolfram, The Integrator. http://integrals.wolfram.com/

[139a] Li Yingying, On Euler's Constant-Calculating Sums by Integrals.
Amer. Math. Monthly, 109, 845-850, 2002.

[139b] Wu Yun-Fei, New series involving the zeta function.
IJMMS 28:7 (2001) 403-411
[PDF] New series involving the zeta function

[140] D. Zagier, The Remarkable Dilogarithm. Jour.Math.Phy.Sci, 22,131-145, 1988.

[141] D. Zeilberger, Closed Form (pun intended!).Contemporary Mathematics,
143,579- 608, 1993.
http://www.math.rutgers.edu/~zeilberg/mamarim/mamarimPDF/pun.pdf

[142] D. Zeilberger, Computerized Deconstruction. Advances in Applied
Mathematics, 30, 633-654, 2003.
http://www.math.rutgers.edu/~zeilberg/mamarim/mamarimPDF/derrida.pdf

[142aa] Zhang Nan-Yue and K.S. Williams, Values of the Riemann zeta function and
integrals involving $\log\left(2\sinh\frac{\theta}{2}\right)$ and $\log\left(2\sin\frac{\theta}{2}\right)$.
Pacific J. Math., 168, 271-289, 1995.
http://projecteuclid.org/Dienst/UI/1.0/Summarize/euclid.pjm/1102620561

[142a] I.J. Zucker, On the series $\sum_{k=1}^{\infty}\binom{2k}{k}^{-1} k^{-n}$ and related sums.
J. Number Theory, 20, 92-102, 1985.

[142b] De-Yin Zheng, Further summation formulae related to generalized harmonic
numbers *Journal of Mathematical Analysis and Applications*, *In Press, Corrected Proof, Available online 12 February 2007.*

[143] W. Zudilin, *One of the numbers ζ(5), ζ(7), ζ(9), ζ(11) is irrational*, Uspekhi
Mat. Nauk [Russian Math. Surveys] 56:4 (2001), 149--150 (pdf, gzip ps).
Full details of various papers relating to the (assumed) irrationality of ζ(2n+1)
are contained in Zudilin's website http://wain.mi.ras.ru/zw/

[144] W. Zudilin, An elementary proof of Apéry's theorem.
Math.NT/0202159 [abs, ps, pdf, other] 2002.

[145] A. Zygmund, Trigonometric Sums. Cambridge Mathematical Library, 2002.

Donal F. Connon
Elmhurst
Dundle Road
Matfield
Kent TN12 7HD
dconnon@btopenworld.com